\documentclass[12pt]{article}
\usepackage{amssymb,amsmath}
\usepackage{graphicx}
\usepackage{a4}
\numberwithin{equation}{section}
\begin{document}
 \baselineskip16pt

\title{\bf On Conjectures of Minkowski and Woods for $n=9$}
\author{Leetika Kathuria{\footnote{The author acknowledges the support of CSIR, India. The paper forms a part of her
Ph.D. dissertation accepted by Panjab University, Chandigarh.  }}~ and  Madhu Raka
 \\ \small{\em Centre for Advanced Study in Mathematics}\\
\small{\em Panjab University, Chandigarh-160014, INDIA}\\
\date{}}

\maketitle
{\abstract{$~~$\par Let $\mathbb{R}^n$ be the n-dimensional Euclidean space with  $O$ as the origin. Let $\wedge$ be a
 lattice of determinant $1$ such that  there is a sphere $|X|<R$ which contains no
point of $\wedge$ other than $O$ and has $n$ linearly independent points of $\wedge$ on its boundary. A well known
conjecture in the geometry of numbers asserts that any closed sphere in $\mathbb{R}^n $ of  radius $ \sqrt{n/4}$
contains a point of $\wedge$. This is known to be true for $n\leq 8$. Here we prove a more general conjecture of Woods
for $n=9$ from which this conjecture follows in $\mathbb{R}^9$. Together with a result of C. T. McMullen (2005), the
long standing conjecture of Minkowski  follows for $n=9$.
 }}
\section{Introduction}
Let
$L_{i}=a_{i1}x_{1}+\cdots+ a_{in}x_{n}$, $1\leq i\leq n$ be $n$  real linear forms in $n$ variables $x_{1},\ldots, x_{n}$ and having determinant $\Delta=\mbox{det}(a_{ij})\neq0$. The following conjecture is attributed to H. Minkowski [18]:\\
\emph{For any given real numbers $c_1,\ldots, c_n$, there exists integers $x_1,\ldots,x_n$ such that
\begin{equation}
\mid(L_{1}+c_{1})\cdots(L_{n}+c_{n})\mid \leqslant \frac{1}{2^{n}}\mid \Delta \mid.
\end{equation}
Equality is necessary if and only if after a suitable unimodular transformation the linear forms $L_{i}$ have the form $2c_{i}x_{i}$ for $1\leq i \leq n$.\\}
It is obvious that equality occurs in (1.1) for the cases mentioned in the conjecture.\vspace{2mm}\\
This result is trivial for $n = 1$. For $n = 2$, a proof was first given in 1899 by Minkowski. Several mathematicians such as Mordell, Landau, Perron, Pall, Macbeath, Sawyer, Cassels have obtained a variety of proofs, partly in an effort to find a proof which would generalize to higher dimensions. For a detailed history, see Gruber and Lekkerkerker [9], Bambah et al [1] and Hans-Gill et al [11]. Minkowski's conjecture has so far been proved for $n\leq 8$.
\vspace{2mm}\\
\noindent For $n = 3$,  following three approaches have been tried:\vspace{1mm}\\
I.~~~~~~Remak-Davenport Approach.\\
II.~~~~~Birch and Swinnerton-Dyer Approach.\\
III.~~~~DOTU-matrix Approach by Narzullaev.\vspace{2mm}\\
Only the Remak-Davenport Approach has been successfully extended to give proofs for $4\leq n \leq 8$. It consists of proving the following two conjectures:\vspace{2mm}\\
\textbf{Conjecture I.} \emph{For any lattice $\Lambda$ in $\mathbb{R}^{n}$ there is an ellipsoid
$E : a_{1}x_{1}^{2}+\cdots+ a_{n}x_{n}^{2} < 1$
which contains no point of $\Lambda$ other than $O$ but has $n$ linearly independent points of $\Lambda$ on its boundary}.\vspace{2mm}\\
\textbf{Conjecture II.} \emph{If $\Lambda$ is a lattice of determinent $1$ and there is a sphere $|X|<R$ which contains no point of $\Lambda$ other than O and has $n$ linear independent points of $\Lambda$ on its boundary then $\Lambda$ is a covering lattice for the closed sphere of radius $\sqrt{n}/2$. Equivalently every closed sphere of radius $\sqrt{n}/2$ lying in $\mathbb{R}^{n}$ contains a point of $\Lambda$}.
\par In 2005, McMullen [17] proved Conjecture I for all $n \geq 3$, using a result of Birch and Swinnerton-Dyer [4]. So to prove Minkowski's Conjecture, it is enough to prove Conjecture II, which is known to be true for $n\leq8$. For a detailed history of these conjectures, see Hans-Gill et al [11], [13]. Using stable lattices,  Shapira and  Weiss [20]
have given a different proof of Minkowski's Conjecture for $ n\leq 7$.\vspace{2mm}

\par \textbf{Here we shall prove Conjecture II for $n=9$, thereby proving Minkowski's Conjecture for $n=9$}.\vspace{2mm}

\par Woods [22], [23] formulated a conjecture from which conjecture-II follows immediately. To state Woods' conjecture, we need to introduce some terminology :
\par Let $\mathbb{L}$ be a lattice in $\mathbb{R}^{n}$. By the reduction theory of quadratic forms introduced by Korkine and Zolotareff [16], a cartesian co-ordinate system may be chosen in $\mathbb{R}^{n}$ in such a way that $\mathbb{L}$ has a basis of the form $$(A_1,0,0,\ldots,0),(a_{2,1},A_2,0,\ldots,0),\ldots,(a_{n,1},a_{n,2},\ldots,a_{n,n-1},A_n),$$ where $A_1, A_2,\ldots,A_n$ are all positive and further for each $i=1,2,\ldots,n$ any two points of the lattice in $\mathbb{R}^{n-i+1}$ with basis
$$(A_i,0,0,\ldots,0),(a_{i+1,i},A_{i+1},0,\ldots,0),\ldots,
(a_{n,i},a_{n,i+1},\ldots,a_{n,n-1},A_n)$$
 are at a distance atleast $A_i$ apart. Such a basis of $\mathbb{L}$ is called a reduced basis.\vspace{2mm}\\
\textbf{Conjecture III (Woods)}:\emph{ If $A_{1}A_{2}\cdots A_{n}=1$ and $A_{i}\leqslant A_{1}$ for each $i$ then any closed sphere in $\mathbb{R}^{n}$ of radius  $\sqrt{n}/2$ contains a point of $\mathbb{L}$.}\vspace{2mm}\\
Woods ([21], [22], [23]) proved this conjecture for $4\leqslant n \leqslant 6$. Hans-Gill et al [10] gave a unified proof of Woods' Conjecture for $n \leqslant 6$. Hans-Gill et al ([11], [13]) proved Woods' Conjecture  for $n=7$ and $n=8$ and thus completed the proof of Minkowski's conjecture for $n=7$ and 8. While answering a question of Shapira and Weiss [20], the authors along with Hans-Gill [15]  have given another proof of Woods' Conjecture and hence of Minkowski's Conjecture for $n\leq7$. Hans-Gill et al ([12], [14]) have obtained estimates to Conjectures of Minkowski and Woods for $9\leq n\leq31$. In particular they have proved the weaker result that if hypothesis of Conjecture III holds, then any closed sphere in $\mathbb{R}^9$ of radius $\frac{\sqrt{9.2587472}}{2}$ contains a point of $\mathbb{L}$.\vspace{1mm}\\
Here we shall prove\vspace{2mm}\\
\textbf{Theorem.} Conjecture III is true for $n = 9$.\vspace{1mm}\\
As remarked earlier, this implies Minkowski's Conjecture for $n=9$.
\section{Preliminary Lemmas}
For a unit sphere $S_{n}$ with center $O$ in $\mathbb{R}^{n}$, let $\Delta(S_{n})$ be the \emph{critical determinant} of $S_{n}$, defined as
\vspace{-2mm}$$\Delta(S_{n}) = \inf\{d(\Lambda):\Lambda~\mbox{has no non-zero point in the interior of}~S_{n}\},$$ where $d(\Lambda)$ denotes the determinant of the lattice $\Lambda$.
 \par Let $\mathbb{L}$ be a lattice in $\mathbb{R}^{n}$ reduced in the sense of Korkine and Zolotareff and $A_1, A_2,\ldots,A_n$ be defined as in Section 1. We state below some preliminary lemmas. Lemmas 1 and 2 are due to Woods [21], Lemma 3 is due to Korkine and Zolotareff
[16] and Lemma 4 is due to  Pendavingh and Van Zwam [19].
 In Lemma 5, the cases $n=2$ and $3$ are classical results of Lagrange and
 Gauss; $n=4$ and $5$ are due to Korkine and Zolotareff [16] while $n=6, 7$ and $8$ are due to Blichfeldt [5].\vspace{2mm}\\
\textbf{Lemma 1.} \emph{If $2\Delta(S_{n+1})A_{1}^{n}\geq d(\mathbb{L})$, then any closed sphere of radius
$$R=A_{1}\{ 1-(A_{1}^{n} \Delta (S_{n+1})/ d(\mathbb{L}))^{2}\}^{1/2}$$ in $\mathbb{R}^{n}$ contains a point of $\mathbb{L}$.}\vspace{2mm}\\
\textbf{Lemma 2.} \emph{For a fixed  integer $i$ with $1\leqslant i\leqslant n-1$, denote by $\mathbb{L}_{1}$ the lattice in $\mathbb{R}^{i}$ with the reduced basis $$(A_1,0,0,\ldots,0),(a_{2,1},A_2,0,\ldots,0),\ldots,(a_{i,1},a_{i,2},\ldots,a_{i,i-1},A_i) $$ and denote by $\mathbb{L}_{2}$ the lattice in $\mathbb{R}^{n-i}$ with the reduced basis $$(A_{i+1},0,0,\ldots,0),(a_{i+2,i+1},A_{i+2},0,\ldots,0),\ldots,(a_{n,i+1},a_{n,i+2},\ldots,a_{n,n-1},A_n) $$
If any sphere in $\mathbb{R}^{i}$ of radius $r_{1}$ contains a point of $\mathbb{L}_{1}$ and if any sphere in $\mathbb{R}^{n-i}$ of radius $r_{2}$ contains a point
of $\mathbb{L}_{2}$ then any sphere in $\mathbb{R}^{n}$ of radius $(r_{1}^{2}+r_{2}^{2})^{1/2}$ contains a point of $\mathbb{L}$.}\vspace{2mm}\\
{\noindent \bf Lemma 3.}\emph{ For all relevant $i$, $A_{i+1}^2\geq\frac{3}{4}A_i^2$ and $A_{i+2}^2\geq\frac{2}{3}A_i^2
~$.}\vspace{2mm}\\
{\noindent \bf Lemma 4.}\emph{ For all relevant $i$, $A_{i+4}^2\geq0.46873A_i^2 $~.}\vspace{2mm}\\
{\noindent \bf Lemma 5.}\emph{$~~~\Delta (S_n) = ~1/\sqrt{2}, ~1/2, ~1/2\sqrt{2}, ~\sqrt{3}/8,1/8$ and $1/16  ~for~
 n=~3,~4,~5,$ $6,~7$ and $~8$ ~respectively.}
\section{Plan of the Proof}  We use the notation and approach of Hans-Gill et al [13] but our method of dealing with various inequalities is somewhat different. We assume that Conjecture III is false for $n=9$ and derive a contradiction. Let $\mathbb{L}$ be lattice satisfying the hypothesis of the conjecture for $n=9$ i.e.  $A_{1}A_{2}\cdots A_{9}=1$~\mbox{and}~ $A_{i}\leqslant A_{1}$ for each $i$. Suppose that there exists a closed sphere of radius $\sqrt{9}/2$ in $\mathbb{R}^{9}$ that contains no point of $\mathbb{L}$. Write $A=A_1^2,~ B=A_2^2,~ C=A_3^2,\ldots, I=A_9^2$. So we have $ABCDEFGHI=1$.\vspace{2mm}\\
$~~~$ If
$(\lambda_{1},\lambda_{2},\cdots,\lambda_{s})$ is an ordered partition of $n$, then the conditional inequality arising
from it, by using Lemmas 1 and 2, is also denoted by $(\lambda_{1},\lambda_{2},\cdots,\lambda_{s})$. If the conditions
in an inequality $(\lambda_{1},\lambda_{2},\cdots,\lambda_{s})$ are satisfied then we say that
$(\lambda_{1},\lambda_{2},\cdots,\lambda_{s})$ holds.\vspace{2mm}\\
For example the inequality (1, 1, 1, 1, 1, 1, 1, 2) results in the conditional inequality :
\begin{equation}
{\rm if~~} 2H\geq I~~ {\rm then~~~} A+B+C+D+E+F+G+4H-\frac{2H^2}{I}>9.
\end{equation}
Since $4H-2H^2/I\leq 2I$,  the second inequality in (3.1) gives
\begin{equation}
A+B+C+D+E+F+G+2I>9.
\end{equation}
One may remark here that the condition $2H\geq I$ is necessary only if we want to use inequality (3.1), but it is not necessary if we want to  use the weaker inequality (3.2). This is so because if $2H< I$, using the partition $(1,1)$ in place of $(2)$ for the relevant part, we get the upper bound $H+I$ which is clearly less than $2I$. We shall call inequalities of type (3.2) as \emph{weak} inequalities.\vspace{1mm} \\
$~~~$Sometimes, instead of Lemma 1, we are able to use the fact that Woods Conjecture is true for dimensions less than or equal to $8$. The
use of this is indicated by putting $^{*}$ on the corresponding part of the partition. For example, the inequality $(5^{*},4)$
 is  \begin{equation} {\rm if~~} F^4ABCDE\geq 2~~ {\rm then~~~} 5(ABCDE)^{\frac{1}{5}}+4F-\frac{1}{2}F^{5}ABCDE > 9,\end{equation}
the hypothesis of the conjecture in $5$ variables being satisfied.\vspace{2mm}\\
$~~~~$We observe that the inequalities of the type $(3,1,1,\cdots,1)$, $(2,1,1,\cdots,1)$, $(1,2,1,\cdots,1)$, $(1,2,2,1,\cdots,1)$, $(3,2,1,\cdots,1)$ etc. always hold. \vspace{2mm}\\
$~~~~$We also observe that for positive real numbers
 $X_{1},\cdots, X_{k}$ we have $X_{1}+\cdots +X_{k} \leq(k-1)+X_{1}\cdots X_{k}$ if either all $X_{i}\leq 1$ or all $X_{i} >1$.\vspace{2mm}\\
$~~~~~$Throughout the paper we shall use the following notation:\\
$a=A-1, b=|B-1|,c=|C-1|,d=|D-1|,e=|E-1|,f=|F-1|,g=|G-1|,h=|H-1|,i=|I-1|.$ Also we can assume $A>1$, because if $A\leq 1$, we must
have  $A=B=C=D=E=F=G=H=I=1.$ In this case Woods' Conjecture  can be seen to be true using inequality $(1,1,1,1,1,1,1,1,1)$. Also it is known that $A\leq \gamma_9<2.1326324$(See [6],[7]).\vspace{2mm}\\
$~~~~~$Each of $B,C,\ldots,I$ can either be $>1$ or $\leqslant$ 1. This give rise to $2^8=256$ cases. The case where each of $B,C,\cdots,I$ is $>1$ does not arise as $ABCDEFGHI$ $ =1$. Cases where $H > 1,~I \leq 1$ or where $G>1, ~H\leq 1, ~I \leq 1$   do not arise (see Propositions 1 and 2 of [13]). This settles $64+32$, i.e. 96 cases. Cases in which $B\leq1$ and at most two out of $C,D,E,F,G,H,I$ are greater than 1 also do not arise (see Proposition 3 (i)  of [13]. ) This settles additional $18$ cases. Also the cases in which exactly one of $B,C,D,E,F,G,H,I$ is $\leq 1$ do not arise (see Lemma 5 (i)  of [13]). This settles $7$ more cases. The remaining $134$ cases are listed in Table 1. The inequalities used to get a contradiction for some easy cases are also listed in Table 1.\vspace{2mm}
 \par In Section 4 we have discussed 111 easy cases. Out of the remaining $23$ cases, $18$ cases are somewhat less difficult and have been dealt in Section 5. The 5 most difficult cases are dealt with separately in Section 6 which constitute almost half of the length of the paper.
We would like to remark that in many cases there are alternative ways of proof using different set of inequalities. We have chosen to describe the method which we find convenient.\vspace{2mm}\\
{\noindent \bf Remark 1.} In this paper we need to maximize or minimize frequently functions of several variables. While doing this we shall find it convenient to name the function involved as $\phi(x), \psi(y)$ etc. to indicate that it is being regarded as function of that variable and other variables are kept fixed. When we say that a given function of several variables in $x,y,\cdots$ is an increasing/decreasing function of $x,y,\cdots$, it means that the concerned property holds when function is considered as a function of one variable at a time, all other variables being fixed. Making use of calculus, we reduce the number of variables one by one by replacing it with the values where it can have its optimum value and result in several functions in at most two variables. For functions in one/two variables, we arrive at a contradiction by plotting their 2/3-dimensional graphs using the software Mathematica. \vspace{2mm}\\
{\noindent \bf Remark 2.} Sometimes we have a function $\phi$ symmetric in a number of variables $x_1,x_2,\cdots,x_r$, where $0\leq x_i\leq a$ for $1\leq i\leq r$. We find that second derivative of the function $\phi$ w.r.t. each variable $x_i$ turns out to be positive. Therefore its maximum can occur at end points of the variables only. Being symmetric, we get $\phi(x_1,x_2,\cdots,x_r)\leq \max \{\phi(0,\cdots,0), \phi(a,0,\cdots,0),\cdots,\\\phi(a,a,\cdots,a)\}$, which will turn out to be at most zero for specific range of $a$. Hence we find that ${\rm max}~\phi(x_1,x_2,\cdots,x_r)\leq 0$. Similarly we get that ${\rm min} ~\phi(x_1,x_2,\cdots,x_r)\geq 2$, if second derivative of the symmetric function $\phi$ w.r.t. each variable $x_i$ turns out to be negative for $0\leq x_i\leq a$ and $\min \{\phi(0,\cdots,0), \phi(a,0,\cdots,0),\cdots,\phi(a,a,\cdots,a)\}\geq 2.$
\vspace{2mm}\\
 {$${\rm \bf Table~~~1}\vspace{4mm}$$
\begin{tabular}{cccccccccccc}
Case & A & B & C & D & E & F & G & H & I &Proposition& Inequalities \\ &&&&&&&&&&&\\

1&$>$&$>$&$>$&$>$&$>$&$>$&$\leq$&$\leq$&$>$ &2(i)&$(1,...,1,3,1)$ \\
2&$>$&$>$&$>$&$>$&$>$&$>$&$\leq$&$\leq$&$\leq$ &$24$ &$-$ \\

3&$>$&$>$&$>$&$>$&$>$&$\leq$&$>$&$\leq$&$>$ &1(i)&$(1,1,1,1,2,2,1)$ \\

4&$>$&$>$&$>$&$>$&$>$&$\leq$&$\leq$&$>$&$>$ &2(i)& $(1,1,1,1,3,1,1)$\\

5&$>$&$>$&$>$&$>$&$>$&$\leq$&$\leq$&$\leq$&$>$ &25&$ -$\\
6&$>$&$>$&$>$&$>$&$>$&$\leq$&$\leq$&$\leq$&$\leq$ &46& $-$\\

7&$>$&$>$&$>$&$>$&$\leq$&$>$&$>$&$\leq$&$>$ &1(i)&$(1,1,1,2,1,2,1)$ \\

8&$>$&$>$&$>$&$>$&$\leq$&$>$&$\leq$&$>$&$>$ &1(i)&$(1,1,1,2,2,1,1)$\\

9&$>$&$>$&$>$&$>$&$\leq$&$>$&$\leq$&$\leq$&$>$ &1(iv)&$(2,1,2,3,1)$ \\
10&$>$&$>$&$>$&$>$&$\leq$&$>$&$\leq$&$\leq$&$\leq$ &26& $-$\\
11&$>$&$>$&$>$&$>$&$\leq$&$\leq$&$>$&$>$&$>$ &2(i)&$(1,1,1,3,1,1,1)$ \\

12&$>$&$>$&$>$&$>$&$\leq$&$\leq$&$>$&$\leq$&$>$ &1(iv)&$(2,1,3,2,1)$\\

13&$>$&$>$&$>$&$>$&$\leq$&$\leq$&$\leq$&$>$&$>$ &27&$-$\\
 14&$>$&$>$&$>$&$>$&$\leq$&$\leq$&$\leq$&$\leq$&$>$ & 45&$-$\\
15&$>$&$>$&$>$&$>$&$\leq$&$\leq$&$\leq$&$\leq$&$\leq$&44&$-$\\

\end{tabular}}
\newpage

{\footnotesize \begin{tabular}{cccccccccccc}
 Case & A & B & C & D & E & F & G & H & I &Proposition& Inequalities \\ &&&&&&&&&&&\\

16&$>$&$>$&$>$&$\leq$&$>$&$>$&$>$&$\leq$&$>$&1(i)&$(1,1,2,1,1,2,1)$\\
17&$>$&$>$&$>$&$\leq$&$>$&$>$&$\leq$&$>$&$>$&1(i)&$(1,1,2,1,2,1,1)$\\
18&$>$&$>$&$>$&$\leq$&$>$&$>$&$\leq$&$\leq$&$>$&1(iv)&$(2,2,1,3,1)$\\
19&$>$&$>$&$>$&$\leq$&$>$&$>$&$\leq$&$\leq$&$\leq$&28&$-$\\
20&$>$&$>$&$>$&$\leq$&$>$&$\leq$&$>$&$>$&$>$&1(i)&$(1,1,2,2,1,1,1)$\\

21&$>$ &$>$&$>$&$\leq$&$>$&$\leq$&$>$&$\leq$&$>$&1(ii)&$(1,1,2,2,2,1)$\\
22&$>$&$>$&$>$&$\leq$&$>$&$\leq$&$\leq$&$>$&$>$&1(iv)&$(2,2,3,1,1)$\\
23&$>$&$>$&$>$&$\leq$&$>$&$\leq$&$\leq$&$\leq$&$>$&29& $-$\\
24&$>$&$>$&$>$&$\leq$&$>$&$\leq$&$\leq$&$\leq$&$\leq$&30&$-$\\
25&$>$&$>$&$>$&$\leq$&$\leq$&$>$&$>$&$>$&$>$&2(i)&$(1,1,3,1,1,1,1)$\\
26&$>$&$>$&$>$&$\leq$&$\leq$&$>$&$>$&$\leq$&$>$&1(iv)&$(2,3,1,2,1)$\\
27&$>$&$>$&$>$&$\leq$&$\leq$&$>$&$\leq$&$>$&$>$&1(iv)& $(2,3,2,1,1)$\\
28&$>$&$>$&$>$&$\leq$&$\leq$&$>$&$\leq$&$\leq$&$>$&1(vi)&$(2,3,3,1,1)$\\
29&$>$&$>$&$>$&$\leq$&$\leq$&$>$&$\leq$&$\leq$&$\leq$&31&$-$\\
30&$>$&$>$&$>$&$\leq$&$\leq$&$\leq$&$>$&$>$&$>$&32&$-$\\

31&$>$ &$>$&$>$&$\leq$&$\leq$&$\leq$&$>$&$\leq$&$>$&33&$-$\\
32&$>$&$>$&$>$&$\leq$&$\leq$&$\leq$&$\leq$&$>$&$>$&42&$-$\\
33&$>$&$>$&$>$&$\leq$&$\leq$&$\leq$&$\leq$&$\leq$&$>$ &46&$-$\\
34&$>$&$>$&$>$&$\leq$&$\leq$&$\leq$&$\leq$&$\leq$&$\leq$&34&$-$\\
35&$>$&$>$&$\leq$&$>$&$>$&$>$&$>$&$\leq$&$>$&1(i)&
$(1,2,1,1,1,2,1)$ \\

36&$>$&$>$&$\leq$&$>$&$>$&$>$&$\leq$&$>$&$>$&1(i)&
$(1,2,1,1,2,1,1)$ \\
37&$>$&$>$&$\leq$&$>$&$>$&$>$&$\leq$&$\leq$&$>$&1(iii)&
$(3,1,1,3,1)$ \\
38&$>$&$>$&$\leq$&$>$&$>$&$>$&$\leq$&$\leq$&$\leq$&19&
$-$ \\
39&$>$&$>$&$\leq$&$>$&$>$&$\leq$&$>$&$>$&$>$&1(i)&
$(1,2,1,2,1,1,1)$ \\

40&$>$&$>$&$\leq$&$>$&$>$&$\leq$&$>$&$\leq$&$>$&1(ii)&
$(1,2,1,2,2,1)$ \\
41&$>$&$>$&$\leq$&$>$&$>$&$\leq$&$\leq$&$>$&$>$&1(iii)&
$(3,1,3,1,1)$\\

42&$>$&$>$&$\leq$&$>$&$>$&$\leq$&$\leq$&$\leq$&$>$&15&
$-$ \\
43&$>$&$>$&$\leq$&$>$&$>$&$\leq$&$\leq$&$\leq$&$\leq$&14&
$-$ \\
44&$>$&$>$&$\leq$&$>$&$\leq$&$>$&$>$&$>$&$>$&1(i)&
$(1,2,2,1,1,1,1)$ \\

45&$>$&$>$&$\leq$&$>$&$\leq$&$>$&$>$&$\leq$&$>$&1(ii)&
$(1,2,2,1,2,1)$ \\

46&$>$&$>$&$\leq$&$>$&$\leq$&$>$&$\leq$&$>$&$>$&1(ii)&
$(1,2,2,2,1,1)$ \\
47&$>$&$>$&$\leq$&$>$&$\leq$&$>$&$\leq$&$\leq$&$>$&1(vi)&
$(3,2,3,1)$ \\
48&$>$&$>$&$\leq$&$>$&$\leq$&$>$&$\leq$&$\leq$&$\leq$&2(ii)&
$(3,1,...,1),(1,2,2,2,1,1)$ \\
49&$>$&$>$&$\leq$&$>$&$\leq$&$\leq$&$>$&$>$&$>$&1(iii)&
$(3,3,1,1,1)$ \\

50&$>$&$>$&$\leq$&$>$&$\leq$&$\leq$&$>$&$\leq$&$>$&1(vi)&

$(3,3,2,1)$\\

51&$>$&$>$&$\leq$&$>$&$\leq$&$\leq$&$\leq$&$>$&$>$&35&
$-$ \\

52&$>$&$>$&$\leq$&$>$&$\leq$&$\leq$&$\leq$&$\leq$&$>$&36&
$-$\\
53&$>$&$>$&$\leq$&$>$&$\leq$&$\leq$&$\leq$&$\leq$&$\leq$&2(ii)&
$(3,1,...,1),(1,2,2,1,..,1)$ \\
54&$>$&$>$&$\leq$&$\leq$&$>$&$>$&$>$&$>$&$>$&2(i)&
$(1,3,1,...,1)$ \\

55&$>$&$>$&$\leq$&$\leq$&$>$&$>$&$>$&$\leq$&$>$&13&
$-$\\
56&$>$&$>$&$\leq$&$\leq$&$>$&$>$&$\leq$&$>$&$>$&13&
$-$ \\

57&$>$&$>$&$\leq$&$\leq$&$>$&$>$&$\leq$&$\leq$&$>$&10&
$-$ \\
58&$>$&$>$&$\leq$&$\leq$&$>$&$>$&$\leq$&$\leq$&$\leq$&20&
$-$ \\
59&$>$&$>$&$\leq$&$\leq$&$>$&$\leq$&$>$&$>$&$>$&16&
$-$ \\

60&$>$&$>$&$\leq$&$\leq$&$>$&$\leq$&$>$&$\leq$&$>$&17&
$-$\\
61&$>$&$>$&$\leq$&$\leq$&$>$&$\leq$&$\leq$&$>$&$>$&18&
$-$ \\

62&$>$&$>$&$\leq$&$\leq$&$>$&$\leq$&$\leq$&$\leq$&$>$&21&
$-$ \\
63&$>$&$>$&$\leq$&$\leq$&$>$&$\leq$&$\leq$&$\leq$&$\leq$&2(ii)&
$(3,1,...,1),(1,2,1,2,1,1,1)$ \\
64&$>$&$>$&$\leq$&$\leq$&$\leq$&$>$&$>$&$>$&$>$&37&
$-$ \\65&$>$&$>$&$\leq$&$\leq$&$\leq$&$>$&$>$&$\leq$&$>$&22&
$-$ \\

66&$>$&$>$&$\leq$&$\leq$&$\leq$&$>$&$\leq$&$>$&$>$&22&
$-$ \\
\end{tabular}}
\newpage

{\footnotesize  \begin{tabular}{cccccccccccc}
 Case & A & B & C & D & E & F & G & H & I &Proposition& Inequalities \\ &&&&&&&&&&&\\

67&$>$&$>$&$\leq$&$\leq$&$\leq$&$>$&$\leq$&$\leq$&$>$&23&
$-$ \\
68&$>$&$>$&$\leq$&$\leq$&$\leq$&$>$&$\leq$&$\leq$&$\leq$&2(ii)&
$(3,1,...,1),(1,2,1,1,2,1,1)$ \\
69&$>$&$>$&$\leq$&$\leq$&$\leq$&$\leq$&$>$&$>$&$>$&38&
$-$ \\

70&$>$&$>$&$\leq$&$\leq$&$\leq$&$\leq$&$>$&$\leq$&$>$&39&
$-$ \\

71&$>$&$>$&$\leq$&$\leq$&$\leq$&$\leq$&$\leq$&$>$&$>$&41&
$-$ \\

72&$>$&$>$&$\leq$&$\leq$&$\leq$&$\leq$&$\leq$&$\leq$&$>$&40&
$-$ \\
73&$>$&$>$&$\leq$&$\leq$&$\leq$&$\leq$&$\leq$&$\leq$&$\leq$&2(ii)&
$(3,1,...,1),(1,2,1,...,1)$ \\

74&$>$&$\leq$&$>$&$>$&$>$&$>$&$>$&$\leq$&$>$ &1(i)&$(2,1,1,1,1,2,1)$\\

75&$>$&$\leq$&$>$&$>$&$>$&$>$&$\leq$&$>$&$>$&1(i)&$(2,1,1,1,2,1,1)$ \\

76&$>$&$\leq$&$>$&$>$&$>$&$>$&$\leq$&$\leq$&$>$ &1(v)&$(2,1,1,1,3,1)$ \\
77&$>$&$\leq$&$>$&$>$&$>$&$>$&$\leq$&$\leq$&$\leq$ &3 &(2,1,1,1,3,1),(2,1,1,1,2,2) \\
78&$>$&$\leq$&$>$&$>$&$>$&$\leq$&$>$&$>$&$>$ &1(i)&$(2,1,1,2,1,1,1)$ \\

79&$>$&$\leq$&$>$&$>$&$>$&$\leq$&$>$&$\leq$&$>$ &1(ii)&(2,1,1,2,2,1) \\

80&$>$&$\leq$&$>$&$>$&$>$&$\leq$&$\leq$&$>$&$>$ &1(v)& $(2,1,1,3,1,1)$\\
81&$>$&$\leq$&$>$&$>$&$>$&$\leq$&$\leq$&$\leq$&$>$ &3&(2,1,1,3,1,1),(2,1,1,2,2,1)\\
82&$>$&$\leq$&$>$&$>$&$>$&$\leq$&$\leq$&$\leq$&$\leq$ &11& $(2,1,1,2,2,1),(2,1,1,2,1,2)$,\\&&&&&&&&&&&$(2,1,1,3,1,1)$\\
83&$>$&$\leq$&$>$&$>$&$\leq$&$>$&$>$&$>$&$>$ &1(i)&$(2,1,2,1,1,1,1)$ \\

84&$>$&$\leq$&$>$&$>$&$\leq$&$>$&$>$&$\leq$&$>$ &1(ii)&(2,1,2,1,2,1)\\

85&$>$&$\leq$&$>$&$>$&$\leq$&$>$&$\leq$&$>$&$>$ &1(ii)&$(2,1,2,2,1,1)$\\

86&$>$&$\leq$&$>$&$>$&$\leq$&$>$&$\leq$&$\leq$&$>$ &1(iv)&$(2,1,2,3,1)$ \\
87&$>$&$\leq$&$>$&$>$&$\leq$&$>$&$\leq$&$\leq$&$\leq$ &2(ii)& $(3,1,...,1),(2,1,2,2,1,1)$\\
88&$>$&$\leq$&$>$&$>$&$\leq$&$\leq$&$>$&$>$&$>$ &1(v)&$(2,1,3,1,1,1)$ \\

89&$>$&$\leq$&$>$&$>$&$\leq$&$\leq$&$>$&$\leq$&$>$ &1(iv)&$(2,1,3,2,1)$\\
90&$>$&$\leq$&$>$&$>$&$\leq$&$\leq$&$\leq$&$>$&$>$ &3&$(2,1,3,1,1,1),(2,1,2,2,1,1)$\\

91&$>$&$\leq$&$>$&$>$&$\leq$&$\leq$&$\leq$&$\leq$&$>$ & 11&$(2,1,2,2,1,1),(2,1,2,1,2,1)$,\\&&&&&&&&&&&$(2,1,3,1,1,1)$\\

92 &$>$&$\leq$&$>$&$\leq$&$>$&$>$&$>$&$>$&$>$&1(i)&$(2,2,1,...,1)$\\

93&$>$&$\leq$&$>$&$\leq$&$>$&$>$&$>$&$\leq$&$>$&1(ii)&(2,2,1,1,2,1)\\
94&$>$
&$\leq$&$>$&$\leq$&$>$&$>$&$\leq$&$>$&$>$&1(ii)&$(2,2,1,2,1,1)$\\
95&$>$
&$\leq$&$>$&$\leq$&$>$&$>$&$\leq$&$\leq$&$>$&1(iv)&$(2,2,1,3,1)$\\
96&$>$
&$\leq$&$>$&$\leq$&$>$&$>$&$\leq$&$\leq$&$\leq$&2(ii)&$(3,1,...,1),(2,2,1,2,1,1)$\\
97&$>$
&$\leq$&$>$&$\leq$&$>$&$\leq$&$>$&$>$&$>$&1(ii)&$(2,2,2,1,1,1)$\\
98&$>$ &$\leq$&$>$&$\leq$&$>$&$\leq$&$>$&$\leq$&$>$&5&(2,2,2,2,1)\\
99&$>$
&$\leq$&$>$&$\leq$&$>$&$\leq$&$\leq$&$>$&$>$&1(iv)&$(2,2,3,1,1)$\\
100&$>$
&$\leq$&$>$&$\leq$&$>$&$\leq$&$\leq$&$\leq$&$>$&2(ii)& $(3,1,...,1),(2,2,2,1,1,1)$\\
101&$>$&$\leq$&$>$&$\leq$&$\leq$&$>$&$>$&$>$&$>$&1(v)&$(2,3,1,1,1,1)$\\
102&$>$
&$\leq$&$>$&$\leq$&$\leq$&$>$&$>$&$\leq$&$>$&1(iv)&$(2,3,1,2,1)$
\\
103&$>$&$\leq$&$>$&$\leq$&$\leq$&$>$&$\leq$&$>$&$>$&1(iv)& $(2,3,2,1,1)$\\
104&$>$
&$\leq$&$>$&$\leq$&$\leq$&$>$&$\leq$&$\leq$&$>$&1(vi)&
$(2,3,3,1)$\\
105&$>$
&$\leq$&$>$&$\leq$&$\leq$&$\leq$&$>$&$>$&$>$&3&$(2,3,1,1,1,1),(2,2,2,1,1,1)$
\\
106&$>$ &$\leq$&$>$&$\leq$&$\leq$&$\leq$&$>$&$\leq$&$>$&2(ii)&$(3,1,...,1),(2,2,1,1,2,1)$
\\
107&$>$
&$\leq$&$>$&$\leq$&$\leq$&$\leq$&$\leq$&$>$&$>$&11&$(2,2,2,1,1,1),(2,2,1,2,1,1)$,\\&&&&&&&&&&&$(2,3,1,1,1,1)$\\
108&$>$&$\leq$&$\leq$&$>$&$>$&$>$&$>$&$>$&$>$&2(i)&
$(3,1,...,1)$ \\

109&$>$&$\leq$&$\leq$&$>$&$>$&$>$&$>$&$\leq$&$>$&1(v)&
$(3,1,1,1,2,1)$ \\
110&$>$&$\leq$&$\leq$&$>$&$>$&$>$&$\leq$&$>$&$>$&1(v)&
$(3,1,1,2,1,1)$ \\
111&$>$&$\leq$&$\leq$&$>$&$>$&$>$&$\leq$&$\leq$&$>$&1(iii)&
$(3,1,1,3,1)$ \\
112&$>$&$\leq$&$\leq$&$>$&$>$&$>$&$\leq$&$\leq$&$\leq$&4&
$(3,1,1,3,1),(2,1,1,1,2,2)$ \\
113&$>$&$\leq$&$\leq$&$>$&$>$&$\leq$&$>$&$>$&$>$&1(v)&
$(3,1,2,1,1,1)$ \\
114&$>$&$\leq$&$\leq$&$>$&$>$&$\leq$&$>$&$\leq$&$>$&1(iv)&
$(3,1,2,2,1)$ \\
 \end{tabular}}
\newpage

{\footnotesize \begin{tabular}{cccccccccccc}
 Case & A & B & C & D & E & F & G & H & I &Proposition& Inequalities \\ &&&&&&&&&&&\\

115&$>$&$\leq$&$\leq$&$>$&$>$&$\leq$&$\leq$&$>$&$>$&1(iii)&
$(3,1,3,1,1)$\\

116&$>$&$\leq$&$\leq$&$>$&$>$&$\leq$&$\leq$&$\leq$&$>$&4&
$(3,1,3,1,1),(2,1,1,2,2,1)$ \\

117&$>$&$\leq$&$\leq$&$>$&$\leq$&$>$&$>$&$>$&$>$&1(v)&
$(3,2,1,1,1,1)$ \\

118&$>$&$\leq$&$\leq$&$>$&$\leq$&$>$&$>$&$\leq$&$>$&1(iv)&
$(3,2,1,2,1)$ \\

119&$>$&$\leq$&$\leq$&$>$&$\leq$&$>$&$\leq$&$>$&$>$&1(iv)&
$(3,2,2,1,1)$ \\

120&$>$&$\leq$&$\leq$&$>$&$\leq$&$>$&$\leq$&$\leq$&$>$&1(vi)&
(3,2,3,1) \\
121&$>$&$\leq$&$\leq$&$>$&$\leq$&$\leq$&$>$&$>$&$>$&1(iii)&
$(3,3,1,1,1)$ \\

122&$>$&$\leq$&$\leq$&$>$&$\leq$&$\leq$&$>$&$\leq$&$>$&1(vi)&
$(3,3,2,1)$\\

123&$>$&$\leq$&$\leq$&$>$&$\leq$&$\leq$&$\leq$&$>$&$>$&4&
$(3,3,1,1,1),(2,1,2,2,1,1)$ \\

124&$>$&$\leq$&$\leq$&$\leq$&$>$&$>$&$>$&$>$&$>$&8&
$(3,1,...,1),(4,1,...,1)$,\\&&&&&&&&&&&$(2,2,1,...,1)$ \\

125&$>$&$\leq$&$\leq$&$\leq$&$>$&$>$&$>$&$\leq$&$>$&7&
$(3,1,...,1),(4,1,1,2,1)$,\\&&&&&&&&&&&$(2,2,1,1,2,1)$\\

126&$>$&$\leq$&$\leq$&$\leq$&$>$&$>$&$\leq$&$>$&$>$&7&
$(3,1,...,1),(4,1,2,1,1)$,\\&&&&&&&&&&&$(2,2,1,2,1,1)$ \\

127&$>$&$\leq$&$\leq$&$\leq$&$>$&$>$&$\leq$&$\leq$&$>$&4&
$(3,1,1,3,1),(2,2,1,2,1,1)$ \\

128&$>$&$\leq$&$\leq$&$\leq$&$>$&$\leq$&$>$&$>$&$>$&7&
$(3,1,...,1),(4,2,1,1,1)$,\\&&&&&&&&&&&$(2,2,2,1,1,1)$ \\

129&$>$&$\leq$&$\leq$&$\leq$&$>$&$\leq$&$>$&$\leq$&$>$&2(ii)&
$(3,1,...,1),(2,1,1,2,2,1)$\\

130&$>$&$\leq$&$\leq$&$\leq$&$>$&$\leq$&$\leq$&$>$&$>$&4&
$(3,1,3,1,1),(2,2,2,1,1,1)$ \\

131&$>$&$\leq$&$\leq$&$\leq$&$\leq$&$>$&$>$&$>$&$>$&12&
$-$ \\

132&$>$&$\leq$&$\leq$&$\leq$&$\leq$&$>$&$>$&$\leq$&$>$&6&
$(2,1,...,1),(4,1,1,2,1)$,\\&&&&&&&&&&&$(2,1,2,1,2,1)$ \\

133&$>$&$\leq$&$\leq$&$\leq$&$\leq$&$>$&$\leq$&$>$&$>$&6&
$(2,1,...,1),(4,1,2,1,1)$,\\&&&&&&&&&&&$(2,1,2,2,1,1)$ \\

134&$>$&$\leq$&$\leq$&$\leq$&$\leq$&$\leq$&$>$&$>$&$>$&9&
$-$ \\

\end{tabular}}
\section{Easy Cases}
In this section we consider 111 easy cases. In many of these cases, the proof can be generalized for arbitrary $n$.\vspace{2mm}\\
{\noindent \bf Lemma 6.}  Let $X_{1},\cdots,X_{9}$ be positive real numbers, each $<2.1326324$ and satisfying
\begin{equation}X_{1}>1,~~~X_{1}X_{2}\cdots X_{9}=1~~~{\rm and}~~~x_i=|X_i-1|\end{equation}
Then the following hold :\vspace{2mm}\\
(i)~~~{If $X_{i}>1$ for $i=3,5,6,7,8,9$, then we have}\vspace{1mm}\\
\indent~~  $\mathfrak{S}_1=4X_1-\frac{2X_{1}^2}{X_2}+4X_3-\frac{2X_{3}^2}{X_4}+X_5+\cdots+X_9\leq9$.\vspace{2mm}\\
(ii)~~~{If $X_{i}>1$ for $i=3,5,7,8,9$, then we have}\vspace{1mm}\\
\indent~~  $\mathfrak{S}_2=4X_1-\frac{2X_{1}^2}{X_2}+4X_3-\frac{2X_{3}^2}{X_4}+4X_5-\frac{2X_{5}^2}{X_6}+X_7+X_8+X_9\leq9$.\vspace{2mm}\\
(iii)~~~{If $X_{i}>1$ for $i=4,7,8,9$ and $X_7\leq X_1$, $X_8\leq X_1$, $X_9\leq X_1$}\vspace{1mm}\\
\indent~~ then we have\vspace{1mm}\\
\indent~~  $\mathfrak{S}_3=4X_1-\frac{X_{1}^3}{X_{2}X_{3}}+4X_4-\frac{X_{4}^3}{X_{5}X_{6}}+X_7+X_8+X_9\leq9$.\vspace{2mm}\\
(iv)~~~{If $X_{i}>1$ for $i=3,5,8,9$ and $X_2\leq X_1$, $X_4\leq X_3$, $X_8\leq X_{1}X_{5}$,}\\
 \indent~~ $X_9\leq X_{1}X_{5}$, then we have\vspace{1mm}\\
\indent~~  $\mathfrak{S}_4=4X_1-\frac{2X_{1}^2}{X_2}+4X_3-\frac{2X_{3}^2}{X_4}+4X_5-\frac{X_{5}^3}{X_{6}X_{7}}+X_8+X_9\leq9$.\vspace{2mm}\\
(v)~~~{If $X_{i}>1$ for $i=4,6,7,8,9$ and $X_i\leq X_{1}X_{4}$ for $i=6,7,8,9$, $X_5\leq X_4$, }\vspace{1mm}\\
\indent~~  then we have\vspace{1mm}\\
\indent~~  $\mathfrak{S}_5=4X_1-\frac{X_{1}^3}{X_{2}X_{3}}+4X_4-\frac{2X_{4}^2}{X_{5}}+X_6+X_7+X_8+X_9\leq9$.\vspace{2mm}\\
(vi)~~~{If $X_{i}>1$ for $i=4,7,9$, $X_8\leq X_7$, $X_9\leq X_1X_{4}X_{7}$, then we have}\vspace{1mm}\\
\indent~~  $\mathfrak{S}_6=4X_1-\frac{X_{1}^3}{X_{2}X_{3}}+4X_4-\frac{X_{4}^3}{X_{5}X_{6}}+4X_7-\frac{2X_{7}^2}{X_8}+X_9\leq9$.\vspace{2mm}\\
\noindent  Proof is similar to that of Lemma 5 of [13] and is therefore omitted. \vspace{2mm}\\
{\noindent \bf Proposition 1.} The following cases do not arise:\vspace{2mm}\\
$\begin{array}{ll}
{\rm (i)}&(3),(7),(8),(16),(17),(20),(35),(36),(39),(44),(74),(75),(78),\\&(83),(92).\\
{\rm (ii)}&(21),(40),(45),(46),(79),(84),(85),(93),(94),(97).\\
{\rm (iii)}&(37),(41),(49),(111),(115),(121).\\
{\rm (iv)}&(9),(12),(18),(22),(26),(27),(86),(89),(95),(99),(102),(103),\\&(114),(118),(119).\\
{\rm (v)}&(76),(80),(88),(101),(109),(110),(113),(117).\\
{\rm (vi)}&(28),(47),(50),(104),(120),(122).\end{array}$
\vspace{2mm}\\
{\noindent \bf Proof.} Each part of Proposition 1 follows immediately from the corresponding part of Lemma 6, after selecting a suitable
inequality. The inequalities used are mentioned in Table 1. \vspace{2mm}\\
{\noindent \bf Lemma 7.} Let $X_{1},\cdots,X_{9}$ be positive real numbers, each $< 2.1326324$, satisfying $(4.1)$. Let
$$ \gamma={\displaystyle\sum_{4 \leq
i \leq 9\atop{X_{i} \leq 1}}}x_{i}~~~{\rm and} ~~~~\delta={\displaystyle\sum_{4 \leq i
 \leq 9\atop{X_{i}>1}}}x_{i}$$
Suppose that either\vspace{2mm}\\
$\begin{array}{ll}{\rm (i)}& X_{i}>1~ {\rm for~~ each} ~i~, ~4\leq i\leq 9~~ {\rm or}\\
{\rm (ii)}& \gamma\leq x_{1}\leq 0.5~~ {\rm or}\\
{\rm (iii)}& \gamma\leq \frac{2}{3}x_{1}~~ {\rm and}~~x_{1}\leq 1 ~~{\rm or}\\
{\rm (iv)}& \gamma\leq \delta/2~~{\rm and}~~\delta\leq 4x_{1}~~{\rm with}~~x_{1}\leq 0.226 ~~{\rm or}\\
{\rm (v)}& \delta\geq 2\gamma~~{\rm and}~~\gamma\leq 2x_{1}~~{\rm with}~~x_{1}\leq 0.226 ~~{\rm or}\\
{\rm (vi)}& \delta\geq \frac{4}{3}\gamma~~~{\rm and}~~\gamma\leq 2x_{1}~~{\rm with}~~x_{1}\leq 0.175, \end{array}$\\
then
$$\mathfrak{S}_7 = 4X_1-X_{1}^{4}X_{4}\ldots X_{9}+X_{4}+\ldots +X_{9}\leq 9.$$
\noindent  The simple proof  similar to that of Lemma 6 of [13]  is  omitted.\vspace{2mm} \\
{\noindent \bf Proposition 2.} The following cases do not arise:\\
$\begin{array}{ll}{\rm (i)}&(1),(4),(11),(25),(54),(108).\\
{\rm (ii)}&(48),(53),(63),(68),(73),(87),(96),(100),(106),(129).\end{array}$\\
{\noindent \bf Proof.}   We apply Lemma 7(i) for the cases in part (i)
and Lemma 7(ii) for the cases in part (ii). The inequalities used are mentioned in Table I.\vspace{2mm}\\
{\noindent \bf Lemma 8.} Let $X_{1},\cdots,X_{9}$ be positive real numbers satisfying $(4.1)$. Let\\
$1<X_{i}\leq X_{1}\leq \frac{4}{3}~~\mbox{for}~~ i=3,7,8,9,~~X_{i}\leq 1~~\mbox{for} ~~i=2,4,5,6 ~~\mbox{and}\\x_{6}<\frac{x_{7}+x_{8}+x_{9}}{2},$
then $\mathfrak{S}_{8}=4X_1-\frac{2X_{1}^2}{X_{2}}+4X_3-\frac{X_{3}^3}{X_{4}X_{5}}+X_6+X_7+X_8+X_9\leq9.$\vspace{2mm}\\
{\noindent \bf Proof.} Applying AM-GM inequality to $-\frac{X_{1}^{2}}{X_2}-\frac{X_{3}^3}{X_{4}X_{5}}$,  we get
$\mathfrak{S}_{8}\leq 4X_1-\frac{X_{1}^{2}}{X_2}+4X_{3}+X_{6}+X_{7}+X_{8}+X_{9}-2(X_{1}^{3}X_{3}^{4}X_{6}X_{7}X_{8}X_{9})^{1/2}.$
Further using $X_{2}\leq 1$ we have
$\mathfrak{S}_{8}\leq 4X_1-X_{1}^{2}+4X_{3}+X_{6}+X_{7}+X_{8}+X_{9}-2(X_{1}^{3}X_{3}^{4}X_{6}X_{7}X_{8}X_{9})^{1/2}.$
As right side is a decreasing function of $X_6$ for $X_6\leq 1$ and also we have $X_6=1-x_6>1-\frac{x_{7}+x_{8}+x_{9}}{2}$, we get\\
$\begin{array}{ll}\mathfrak{S}_{8}\leq &4X_1-X_{1}^{2}+4X_{3}+(1-\frac{x_{7}+x_{8}+x_{9}}{2})+X_{7}+X_{8}+X_{9}-2X_{1}^{\frac{3}{2}}X_{3}^{2}
\\&(1-\frac{x_{7}+x_{8}+x_{9}}{2})^{\frac{1}{2}} X_{7}^{\frac{1}{2}}X_{8}^{\frac{1}{2}}X_{9}^{\frac{1}{2}}. \end{array}$\\
Again right side is a decreasing function of $X_{3}$ for $x_7+x_8+x_9\leq 3x_1\leq 1$ and $X_3>1$. Replacing $X_{3}$ by 1 and simplifying we get
$\mathfrak{S}_{8}<11+2x_{1}-x_{1}^{2}+\frac{y}{2}-2(1+x_{1})^{\frac{3}{2}}(1-\frac{y}{2})^{\frac{1}{2}}(1+y)^{\frac{1}{2}}=\phi(x_1,y), {\rm say}$, where $y=x_{7}+x_{8}+x_{9}$.
One verifies that $\phi(x_1,y)$ is at most 9 for $0<y\leq 3x_1$ ~and~ $0<x_{1}\leq \frac{1}{3}$.
\vspace{2mm}\\
{\noindent \bf Proposition 3.} Cases (77), (81), (90), (105) do not arise.\vspace{2mm}\\
{\noindent \bf Proof.} It follows immediately from Lemma 8 after selecting suitable inequalities. The inequalities used are mentioned in Table I.\vspace{2mm}\\
{\noindent \bf Lemma 9.} Let $X_{1},\cdots,X_{9}$ be positive real numbers satisfying $(4.1)$. Let
$1<X_{i}\leq X_{1}\leq \frac{4}{3}~~ \mbox{for}~~ i=4,8,9,~~X_{i}\leq 1~~ \mbox{for} ~~i=2,3,5,6,7 ~~\mbox{and} ~~x_{7}<\frac{x_{8}+x_{9}}{2},$
then $\mathfrak{S}_{9}=4X_1-\frac{X_{1}^3}{X_{2}X_{3}}+4X_4-\frac{X_{4}^3}{X_{5}X_{6}}+X_7+X_8+X_9\leq9.$\vspace{2mm}\\
{\noindent \bf Proof.} Applying AM-GM inequality to $-\frac{X_{1}^{3}}{X_{2}X_{3}}-\frac{X_{4}^3}{X_{5}X_{6}}$ we get
$\mathfrak{S}_{9}\leq 4X_{1}+4X_{4}+X_{7}+X_{8}+X_{9}-2X_{1}^{2}X_{4}^{2}X_{7}^{\frac{1}{2}}X_{8}^{\frac{1}{2}}X_{9}^{\frac{1}{2}}$.
As right side is a decreasing function of $X_7$ for $X_7\leq 1$ and $X_7>1-\frac{x_{8}+x_{9}}{2}$. Replacing $X_7$ by $1-\frac{x_{8}+x_{9}}{2}$ we get
$\mathfrak{S}_{9}\leq 4X_{1}+4X_{4}+(1-\frac{x_{8}+x_{9}}{2})+X_{8}+X_{9}-2X_{1}^{2}X_{4}^{2}
(1-\frac{x_{8}+x_{9}}{2})^{\frac{1}{2}}X_{8}^{\frac{1}{2}}X_{9}^{\frac{1}{2}}.$
Again right side is a decreasing function of $X_{4}$ for $x_8+x_9\leq 2x_1\leq \frac{2}{3}$ and $X_4>1$. Replacing $X_{4}$ by 1 and simplifying we get
$\mathfrak{S}_{9}\leq11+4x_{1}+\frac{y}{2}-2(1+x_{1})^{2}(1-\frac{y}{2})^{\frac{1}{2}}(1+y)^{\frac{1}{2}}=\phi(x_1,y), ~{\rm say}$, where $y=x_{8}+x_{9}$.
One verifies that $\phi(x_1,y)$ is at most 9 for $0<y\leq 2x_1$ ~and~ $0<x_{1}\leq \frac{1}{3}$.
\vspace{2mm}\\
{\noindent \bf Proposition 4.} Cases (112), (116), (123), (127), (130) do not arise.\vspace{2mm}\\
{\noindent \bf Proof.} It follows immediately from Lemma 9 after selecting suitable inequalities. The inequalities used are mentioned in Table I.\vspace{2mm}\\
{\noindent \bf Lemma 10.} Let $X_{i}>1$ be real numbers for $1\leq i\leq 6$.\vspace{2mm}\\
(i) If $X_{1}^{5}\geq 2$, then\\
\indent$\mathfrak{S}_{10} = 4X_1-\frac{1}{2}X_{1}^{5}X_{2}X_{3}X_{4}X_{5}X_{6}+X_{2}+X_{3}+X_{4}+X_{5}+X_{6}<9$.\vspace{1mm}\\
(ii) If $X_{i}\leq X_{1}=A$ (say) for $2\leq i\leq 6$, ~$A< 2.1326324$, then\\
\indent $\mathfrak{S}_{11} = A+4X_2-\frac{1}{2}X_{2}^{5}X_{3}X_{4}X_{5}X_{6}A+X_{3}+X_{4}+X_{5}+X_{6}<9$.\vspace{2mm}\\
\noindent  The simple proof  similar to that of Lemma 7 of [13]  is  omitted. \vspace{2mm}\\
{\noindent \bf Proposition 5.} Case (98) i.e. $A>1, B\leq 1, C>1, D\leq 1, E>1, F\leq 1, G>1, H\leq1, I>1$ does not arise.\vspace{2mm}\\
{\noindent \bf Proof.} The inequality (2,2,2,2,1) holds here, i.e. $4A-\frac{2A^2}{B}+4C-\frac{2C^2}{D}+4E-\frac{2E^2}{F}+4G-\frac{2G^2}{H}+I>9$. Applying AM-GM,
we get $4A+4C+4E+4G+I-8(A^3C^3E^3G^3I)^{\frac{1}{4}}>9$. The left side is a decreasing function of $I$ for $1<I\leq A$, so we replace $I$ by 1 to get
$4A+4C+4E+4G-8(A^3C^3E^3G^3)^{\frac{1}{4}}>8$. The left side is a decreasing function of $G$ for $1<G\leq A$, so we can replace $G$ by 1 to get
$4A+4C+4E-8(A^3C^3E^3)^{\frac{1}{4}}>4$. Using similar argument with $E$ and $C$ we get $4A+4-8A^{\frac{3}{4}}>0$, which is not true for
$1<A<2.1326324.$\vspace{2mm}\\
{\noindent \bf Lemma 11.} All cases in which $B\leq1$ and any three out of $C,D,E,F,G,H,I$ are greater than 1 and $A<1.196$ do not arise.\vspace{2mm}\\
 The proof is similar to that of Proposition 3(ii) of [13].\vspace{2mm}\\
{\noindent \bf Proposition 6.} Cases (132) and (133) do not arise.\vspace{3mm}\\
{\noindent \bf Proof.} Firstly we consider Case (132) i.e. $A>1, B\leq 1, C\leq 1, D\leq 1, E\leq 1, F>1, G>1, H\leq 1, I>1$. Here $a\leq \frac{1}{3}$ by Lemma 3. If $A<1.196$ we get a contradiction by Lemma 11. So let $A\geq1.196$. Using the weak inequality (2,1,2,1,2,1), we have $-2b-c-2e+f-2h+i>0$. This gives $e+h<\frac{f+i}{2}=\frac{k}{2}$, say. Therefore
$EFHI\geq(1+f+i)(1-\frac{f+i}{2})>(1+k)(1-\frac{k}{2})\geq 1$, for $0<k\leq 2a\leq \frac{2}{3}$. Then $A^{4}EFGHI>A^{4}\geq2$, therefore
(4,1,1,2,1) holds, i.e. $4A-\frac{1}{2}\frac{A^4}{BCD}+E+F+4G-\frac{2G^2}{H}+I>9$. Using AM-GM inequality we get $4A+E+F+4G+I-2A^{\frac{5}{2}}G^{\frac{3}{2}}E^{\frac{1}{2}}F^{\frac{1}{2}}I^{\frac{1}{2}}>9$. Left side is a
decreasing function of $E$ and $E>1-\frac{k}{2}$, so we can replace E by $1-\frac{k}{2}$ and get
$4A+4G+\frac{k}{2}-2A^{\frac{5}{2}}G^{\frac{3}{2}}(1-\frac{k}{2})^{\frac{1}{2}}(1+k)^{\frac{1}{2}}>6$. Now left side is decreasing function of $G$ for $A\geq1.196$ and $G>1$, so replacing
$G$ by 1 and simplifying we get $\phi(a,k)=2+4a+\frac{k}{2}-2(a+1)^{\frac{5}{2}}(1-\frac{k}{2})^{\frac{1}{2}}(1+k)^{\frac{1}{2}}>0$. One verifies that $\phi(a,k)$
 is at most zero for $0<k\leq2a$ ~and~ $0<a\leq \frac{1}{3}$, giving thereby a contradiction.\vspace{1mm}\\
\indent Now consider Case (133). Using the weak inequality
(2,1,2,2,1,1) and the inequality (4,1,2,1,1) and proceeding as in
Case (132) we get a contradiction.\vspace{2mm}\\
{\noindent \bf Proposition 7.} Cases (125), (126) and (128) do not arise.\vspace{2mm}\\
{\noindent \bf Proof.} First consider Case (125), i.e. $A>1,~ B\leq1, ~C\leq1, ~D\leq1,~ E > 1, ~F > 1,~ G > 1, ~H \leq 1, ~I > 1.$
Here $a \leq \frac{1}{3}$. Using the weak inequality $(2,2,1,1,2,1)$ we get  $-2b-2d+e+f-2h+i >0 $, which gives $d+h<\frac{e+f+i}{2}=\frac{k}{2}$, say. Therefore
$EFHI\geq
(1-h)(1+k)>(1-\frac{k}{2})(1+k)\geq 1$ for $k\leq 3a\leq 1$. Suppose first that $A^4\geq 2$, then $A^{4}EFGHI> A^{4}\geq
2~$. Applying AM-GM inequality to $(4,1,1,2,1)$ we get $4A+4G+E+F+I-2(A^{5}G^{3}EFI)^{\frac{1}{2}}>9$. Left side is a decreasing function of $I$ for $1<I\leq A$, so we can replace $I$ by 1.
Similarly successively replacing $E$, $F$ and $G$ by 1, we get $4A-2A^{\frac{5}{2}}>2$, which is not true for $1<A\leq \frac{4}{3}$.\\
Hence we must have $A^4< 2$ i.e. $a<0.19$. Now we get a contradiction using Lemma 7(iv) with $\gamma = d+h, ~~\delta = e+f+g+i\leq4a$ and  $x_{1}=a<0.19$.
\vspace{1mm}\\
Proof of the Cases (126) and (128) is similar to the proof of Case (125) using the suitable inequalities. The inequalities used are mentioned in Table 1.
\vspace{2mm}\\
{\noindent \bf Proposition 8.} Case (124) i.e. $A>1, B\leq 1, C\leq 1, D\leq 1, E>1, F>1, G>1, H>1, I>1$ does not arise.\vspace{2mm}\\
{\noindent \bf Proof.} Suppose first that $A^{4}\geq2$. This gives $A^{4}EFGHI>A^{4}\geq 2$, therefore (4,1,1,1,1,1) holds. Using Lemma 10(i) with $X_{1}=A$,
$X_{2}=E$, $X_{3}=F$,
$X_{4}=G$, $X_{5}=H$, $X_{6}=I$  and noting that $X_{1}^5>X_{1}^4\geq 2$, we get a contradiction. So we have $A^{4}<2$, i.e. $a<0.19$. Now the inequality (3,1,1,1,1,1,1) gives
$1+4a-d+e+f+g+h+i-(1+a)^{4}(1-d)(1+e)(1+f)(1+g)(1+h)(1+i)>0$. Left side is increasing function of $d$. Also using the weak inequality
(2,2,1,1,1,1,1), we have $-2b-2d+e+f+g+h+i>0$. This gives $d<\frac{e+f+g+h+i}{2}$ . So replacing $d$ by $\frac{e+f+g+h+i}{2}$ we get
$\phi(e,f,g,h,i)=1+4a+\frac{e+f+g+h+i}{2}-(1+a)^{4}(1-\frac{e+f+g+h+i}{2})(1+e)(1+f)(1+g)(1+h)(1+i)>0$.
Following Remark 2, we find that ${\rm max}~\phi(e,f,g,h,i)\leq 0$ for $0<a<0.19$, giving thereby a contradiction.\vspace{2mm}\\
{\noindent \bf Proposition 9.} Case (134) i.e. $A>1, B\leq 1, C\leq 1, D\leq 1, E\leq1, F\leq1, G>1, H>1, I>1$ does not arise.\vspace{2mm}\\
{\noindent \bf Proof.} If $A<1.196$, we get a contradiction using Lemma 11. So let $A\geq1.196$. Also $A\leq \frac{4}{3}$. Using the weak inequalities (2,1,2,1,1,1,1)
and (2,2,2,1,1,1) we have $e<\frac{g+h+i}{2}-\frac{f}{2}$ and $f<\frac{g+h+i}{2}$ respectively. Now
$A^{4}EFGHI>(1+a)^{4}(1-\frac{g+h+i}{2}+\frac{f}{2})(1-f)(1+g)(1+h)(1+i) = \theta(f)>\theta(\frac{g+h+i}{2})$ = $(1+a)^{4}(1-\frac{g+h+i}{4})(1-\frac{g+h+i}{2})(1+g)(1+h)(1+i) = \phi(g,h,i)$, say. Following Remark 2, we find that ${\rm min}~\phi(g,h,i)\geq 2$ for $a\geq0.196$. Hence (4, 1, 1, 1, 1, 1) holds, i.e. $4A-\frac{1}{2}A^{5}EFGHI +E+F+G+H+I>9.$ As the
coefficient of $E$ on left side is negative, so replacing $E$ by its lower bound $ 1-\frac{g+h+i}{2}+\frac{f}{2}$
we get
\begin{equation}
\begin{array}{r}\varphi(f)=4a-\frac{f}{2}+\frac{g+h+i}{2}-{\frac{1}{2}}A^{5}(1-\frac{g+h+i}{2}+\frac{f}{2})(1-f)(1+g)(1+h)(1+i)>0.
\end{array}
\end{equation}
As $\varphi''(f)>0$ and $0\leq f<\frac{g+h+i}{2}$, we have $\varphi(f)\leq {\rm
max}\{\varphi(0),\varphi(\frac{g+h+i}{2})\}$. Let $\varphi(0)=\psi_{1}(g,h,i)~{\rm and}~ \varphi(\frac{g+h+i}{2})=\psi_{2}(g,h,i).$
Following Remark 2, we find that ${\rm max}~\psi_1(g,h,i)\leq 0$ and ${\rm max}~\psi_{2}(g,h,i)\leq0$ for $0< a \leq \frac{1}{3}$, which gives a contradiction to (4.2).
\vspace{2mm}\\
{\noindent \bf Proposition 10.} Case (57) i.e. $A>1, B>1, C\leq 1, D\leq 1, E>1, F>1, G\leq1, H\leq1, I>1$ does not arise.\vspace{2mm}\\
{\noindent \bf Proof.} Applying AM-GM inequality to $(1,3,1,3,1)$, we have $A+4B+E+4F-2B^{2}F^{2}(AEI)^{\frac{1}{2}}+I>9.$ Left side is a decreasing function of $E$ for $1<E\leq A$, so replacing $E$ by 1, we have
$A+4B+4F-2B^{2}F^{2}(AI)^{\frac{1}{2}}+I>8.$\vspace{1mm}\\
\textbf{Case(i)}~~~$A< B^{4}F^{4}I.$ Then left side is a decreasing
 function of $A$  and $A\geq I$. So we get $2I+4B+4F-2B^{2}F^{2}I>8.$ Again left side is a decreasing function of
 $B$  and $F$.  Successively replacing $B$ and $F$ by $1$ we get a contradiction.\\
\textbf{Case(ii)}~~~$A\geq B^{4}F^{4}I.$ This gives $a\geq4(b+f)+i> 2(b+f)+i$. From the weak inequality  $(2,2,2,1,1,1)$ we get $2b-2d+2f-g-h+i >0$
which gives $d+g+h < 2(b+f)+i <a$. Now using $(3,1,1,1,1,1,1)$ and applying Lemma 7(ii)  with $\gamma = d+g+h$ and $x_{1}=a$,
we get a contradiction.\vspace{2mm}\\
{\noindent \bf Proposition 11.} Cases (82), (91) and (107) do not arise.\vspace{2mm}\\
{\noindent \bf Proof.} First consider Case (82), i.e. $A>1,~ B\leq1, ~C>1, ~D>1,~ E>1, ~F\leq1,~ G\leq1, ~H \leq 1, ~I \leq 1.$ Here $a \leq \frac{1}{3}$.
Using the weak inequalities $(2,1,1,2,2,1)$ ~and~ $(2,1,1,2,1,2)$   we get
\begin{equation}-2b+c+d-2f-2h-i>0,\vspace{-2mm}\end{equation}
\begin{equation}-2b+c+d-2f-g-2i>0.\end{equation}
Therefore $h<\frac{c+d}{2}-\frac{i}{2}$~~and~~$i<\frac{c+d}{2}$.\\
The inequality (2,1,1,3,1,1) holds i.e. $4A-\frac{2A^{2}}{B}+C+D+4E-E^{4}ABCDHI+H+I>9$.
Applying  AM-GM inequality to $-\frac{A^{2}}{B}-E^{4}ABCDHI $ and using $B \leq 1$ we get
$4A-A^{2}+C+D+4E+H+I-2A^{\frac{3}{2}}E^{2}C^{\frac{1}{2}}D^{\frac{1}{2}}H^{\frac{1}{2}}I^{\frac{1}{2}}>9$. Left  side is a decreasing function of $E$ as $E>1$ and
$A^3CDHI>A^3CD(1-(h+i))>(1+a)^3(1+c)(1+d)(1-(c+d))\geq 1$, for $0<c\leq a$, $0<d\leq a$ and $a\leq\frac{1}{3}$. So replacing $E$ by $1$ we get $4A-A^{2}+C+D+H+I-2A^{\frac{3}{2}}C^{\frac{1}{2}}D^{\frac{1}{2}}H^{\frac{1}{2}}I^{\frac{1}{2}}>5$.
Again left side is a
decreasing function of $H$ as $CDI>(1+c+d)(1-\frac{c+d}{2})>1$ and  $1-\frac{c+d}{2}+\frac{i}{2}<H\leq 1$.  Replacing $H$ by $1-\frac{c+d}{2}+\frac{i}{2}$,
~we have
$4A-A^{2}+C+D+\left(1-\frac{c+d}{2}+\frac{i}{2}\right)+I-2A^{\frac{3}{2}}C^{\frac{1}{2}}D^{\frac{1}{2}}\left(1-\frac{c+d}{2}+\frac{i}{2}\right)^{\frac{1}{2}}I^{\frac{1}{2}}>5,$
i.e.\begin{equation}\begin{array}{ll}\phi(i) = & 3+4a-(1+a)^{2}+\frac{c+d}{2}-\frac{i}{2}-2(1+a)^{\frac{3}{2}}(1+c)^{\frac{1}{2}}(1+d)^{\frac{1}{2}}(1-i)^{\frac{1}{2}}\\&(1-\frac{c+d}{2}+
\frac{i}{2})^{\frac{1}{2}}>0.\end{array}\end{equation}
As $\phi''(i)>0$ and $0\leq i<\frac{c+d}{2}$,  we have $\phi(i)\leq {\rm
max}\{\phi(0),\phi(\frac{c+d}{2})\}$. Now let\\
$\phi(0)=\psi_{1}(c,d)$ and $\phi(\frac{c+d}{2})=\psi_{2}(c,d).$
Following Remark 2, we find that for $m=1,2$, ${\rm max}~\psi_m(c,d)\leq 0$, for $0<a\leq\frac{1}{3}$. This contradicts (4.5).
\vspace{1mm}\\
Proof of the Cases (91) and (107) is similar to the proof of Case (82) using the suitable inequalities. The inequalities used are mentioned in Table I.
\vspace{2mm}\\
{\noindent \bf Proposition 12.} Case (131) i.e. $A>1,~ B\leq1, ~C\leq1, ~D\leq1,~ E\leq1, ~F>1,~ G>1, ~H >1, ~I > 1$ does not arise.\vspace{2mm}\\
{\noindent \bf Proof.} Using the weak inequalities $(2,2,1,1,1,1,1)$  and $(2,1,2,1,1,1,1)$  we get\\
\begin{equation}-2b-2d-e+f+g+h+i>0,\vspace{-2mm}\end{equation}
\begin{equation}-2b-c-2e+f+g+h+i>0.\end{equation}
Therefore $d<\frac{f+g+h+i}{2}-\frac{e}{2}$~~and~~$e<\frac{f+g+h+i}{2}$.\\
\textbf{Case(i)}~~~$A\geq1.19$, then $A^{4}EFGHI> (1+a)^{4}(1-\frac{f+g+h+i}{2})(1+f)(1+g)(1+h)(1+i) = \phi(f,g,h,i)$. Following Remark 2, we find that ${\rm min}~\phi(f,g,h,i)\geq 2$ for $a\geq0.19$. Hence we have $A^{4}EFGHI>2$ and so (4, 1, 1, 1, 1, 1) holds, i.e.
$4A-\frac{1}{2}A^{5}EFGHI+E+F+G+H+I>9.$ As the coefficient of $E$ on left
side is negative, so we replace $E$ by $ 1-\frac{f+g+h+i}{2}$ and get that
$\varphi(f,g,h,i)=4a+\frac{f+g+h+i}{2}-\frac{1}{2}A^{5}(1-\frac{f+g+h+i}{2})(1+f)(1+g)(1+h)(1+i)>0.$\\
Following Remark 2, we find that ${\rm max}~\varphi(f,g,h,i)\leq0$ for $0<a\leq \frac{1}{3}$, a contradiction.\\
\textbf{Case(ii)~~~}$A<1.19$. The inequality (3,1,1,1,1,1,1) gives $4A-A^{4}DEFGHI+D+E+F+G+H+I>9.$ The coefficient of $D$ on the left side is
 negative, so replacing $D$ by $1-\frac{f+g+h+i}{2}+\frac{e}{2}$, we get\\
$
\begin{array}{ll}\phi(e)=&1+4a-\frac{e}{2}+\frac{f+g+h+i}{2}-(1+a)^{4}\left(1-\frac{f+g+h+i}{2}+\frac{e}{2}\right)(1-e)(1+f)\\&(1+g)(1+h)(1+i)>0.\end{array} $\\
Now $\phi^{''}(e)>0$~and~$0\leq e<\frac{f+g+h+i}{2}$, so we have $\phi(e)\leq {\rm
max}\{\phi(0),\phi(\frac{f+g+h+i}{2})\}$. Now let
$\phi(0)=\psi_{1}(f,g,h,i) ~{\rm and}~ \phi(\frac{f+g+h+i}{2})=\psi_{2}(f,g,h,i).$
Following Remark 2, we find that for $m=1,2$, ${\rm max}~\psi_{m}(f,g,h,i)\leq0$ for $0<a<0.19$, a contradiction.\vspace{2mm}\\
{\noindent \bf Proposition 13.} Cases (55) and (56) do not arise.\vspace{2mm}\\
{\noindent \bf Proof.}~~Consider Case (55) i.e. $A>1, B>1, C\leq 1, D\leq 1, E>1, F>1, G>1, H\leq 1, I>1.$\\
We can take $h\geq b$. For if $h<b$, then using $(1,3,1,1,1,1,1)$ and applying Lemma 7(ii) with $\gamma=h, ~x_{1}=b\leq \frac{1}{ 3}$, we get a contradiction. Now consider following two cases:
\vspace{2mm}\\
{\noindent \bf Case (i) $GH>1$}

As $B^{2}>CD$, we can use the weak inequality $(1,3,1,1,2,1)$ and get $A+4B-B^{4}EFGHIA+E+F+2H+I>9.$ Since $2H\leq 2$ and $-GH<-1$ we get
 $A+4B-B^{4}EFIA+E+F+I>7.$
As the coefficient of $I$ is negative we can replace $I$ by $1$. Similarly we can successively replace $E$ and $F$ by $1$ and get $A+4B-B^{4}A>4$
which is not possible for $A>1$ and $B>1$.
\vspace{2mm}\\
{\noindent \bf Case (ii) $GH \leq 1$}

Using the weak inequality $(2,2,1,1,2,1)$ we get  $2b-2d+e+f-2h+i >0 $, which gives $d<\frac{e+f+i}{2}$ as $h \geq b$.
Applying AM-GM inequality to (3,1,1,3,1) we get
 $4A+D+E+4F-2A^{2}F^{2}(DEI)^\frac{1}{2}+I>9.$ Left side is a decreasing function of $D$ and $D>1- \frac{e+f+i}{2}$, so
 replacing $D$ by $1- \frac{e+f+i}{2}$ we get
 $\phi(e,i)=2+4a+\frac{7f}{2}+\frac{e+i}{2}-2(1+a)^{2}(1+f)^{2}(1-\frac{e+f+i}{2})^\frac{1}{2}(1+e)^{\frac{1}{2}}(1+i)^{\frac{1}{2}}>0.$
Following Remark 2, we find that ${\rm max}~\phi(e,i)\leq0$, for $0<f\leq a\leq 0.5$, giving thereby a contradiction.\vspace{2mm}\\
Proof of the Case (56) is similar to the proof of Case (55) using the suitable inequalities.\vspace{2mm}\\
{\noindent \bf  Proposition 14.} Case (43) i.e. $A > 1,~ B > 1, ~C \leq 1, ~D > 1,~ E > 1, ~F \leq 1,~ G \leq 1, ~H \leq 1, ~I \leq 1$ does not arise.\vspace{2mm}\\
{\noindent \bf Proof.}  Using the weak inequalities (1,2,2,1,2,1) and $(1,2,2,2,2)$, we get
\begin{equation}a-2c+2e-f-2h-i>0,\vspace{-3mm}\end{equation}
\begin{equation}a-2c+2e-2g-2i>0.\end{equation}
This gives $h<e+\frac{a}{2}-\frac{i}{2}$~ and~ $i<e+\frac{a}{2}$.

{\noindent \bf Case (i) $E \geq 1.26$}

 Here $E^4ABCDI=\frac{E^{3}}{FGH}>2$. So (4*,4,1) holds, i.e. $4(ABCD)^{\frac{1}{4}}+4E-\frac{1}{2}E^{5}ABCDI+I > 9$. The left side is a decreasing function of $E$, so replacing $E$ by 1.26, we have
$\phi(x)=4x^{\frac{1}{4}}+4(1.26)-\frac{1}{2}(1.26)^{5}xI+I>9,$ where $x=ABCD$. Now
$\phi'(x)=0$ gives $x=(\frac{2}{1.26^{5}I})^{\frac{4}{3}}$ and $\phi''(x)<0$ at $x=(\frac{2}{1.26^{5}I})^{\frac{4}{3}}$.
So $\phi(x)\leq\phi((\frac{2}{1.26^{5}I})^{\frac{4}{3}})=4(1.26)+3\left(\frac{2}{1.26^5I}\right)^{\frac{1}{3}}+I$, which is at most 9
for $\frac{4}{9}\leq I\leq 1$, a contradiction.\vspace{2mm}\\
{\noindent \bf Case (ii) $E < 1.26$}

Applying AM-GM inequality to (3,1,3,1,1) we have $4A+D+4E+H+I-2A^{2}E^{2}D^{\frac{1}{2}}H^{\frac{1}{2}}I^{\frac{1}{2}}>9$. Left side is a decreasing function of $D$ as $A^3E^4HI>ABDEHI=\frac{1}{CFG}\geq 1$, so replacing $D$ by 1 we get $4A+4E+H+I-2A^{2}E^{2}H^{\frac{1}{2}}I^{\frac{1}{2}}>8$. Now left side is a decreasing function of $H$ as $A^4E^4I\geq ABDEI=\frac{1}{CFGH}\geq 1$ and $1-e-\frac{a}{2}+\frac{i}{2}<H\leq 1$.
So replacing $H$ by $1-e-\frac{a}{2}+\frac{i}{2}$ and on simplifying we get
$\theta(i)=2+\frac{7a}{2}+3e-\frac{i}{2}-2(1+a)^{2}(1+e)^{2}\sqrt{\left(1-e-\frac{a}{2}+\frac{i}{2}\right)(1-i)}>0.$
Now $\theta''(i)>0$ and $0\leq i<e+\frac{a}{2}$. Therefore $\theta(i)\leq{\rm max}\{\theta(0),\theta(e+\frac{a}{2})\}$,
which can be seen to be at most zero for $0<a\leq 0.5$ and $0<e\leq{\rm min}(a,0.26)$, giving thereby a contradiction.\vspace{2mm}\\
{\noindent \bf  Proposition 15.} Case (42) i.e. $A > 1,~ B > 1, ~C \leq 1, ~D > 1,~ E > 1, ~F \leq 1,~ G \leq 1, ~H \leq 1, ~I > 1$ does not arise.\vspace{2mm}\\
{\noindent \bf Proof.} Here $a\leq \frac{1}{2}$,~$b\leq \frac{1}{3}$. Using the weak inequality (1,2,2,1,2,1) we have $a-2c+2e-f-2h+i>0$, i.e.
$h<e+\frac{a+i}{2}<e+a$.\vspace{1mm}\\
{\noindent \bf Claim(i) $E^{4}ABCDI < 2$~ {\rm and} ~$FGH>\frac{1}{2}$}

Suppose $E^{4}ABCDI \geq 2$, then (4*,4,1) holds, i.e. $\phi(x)=4x^{\frac{1}{4}}+4E-\frac{1}{2}E^{5}xI+I> 9$, where $x = ABCD$. $\phi'(x) = x^{\frac{-3}{4}}-\frac{1}{2}E^{5}I = 0$ gives $x = (\frac{2}{E^{5}I})^{4/3}$.
Since $\phi''(x) < 0$ at $x = (\frac{2}{E^{5}I})^{4/3}$, so $\phi(x) \leq \phi((\frac{2}{E^{5}I})^{4/3}) = 4E+I+3(\frac{2}{E^{5}I})^{\frac{1}{3}}=\eta(E,I)$,
say. We find that $\eta(E,I)<9$ for $1<I\leq A\leq \frac{3}{2}$ and $1<E\leq \frac{4}{3}$. So $E^{4}ABCDI < 2$, i.e. $\frac{E^3}{FGH}<2$. It gives $FGH>\frac{1}{2}$.\vspace{3mm}

{\noindent \bf Claim(ii) $A^{4}EFGHI < 2$,~ $A < \sqrt{2}$,~ $E<1.189208$~{\rm and}~$I<1.31951$.}

Suppose $A^{4}EFGHI \geq 2$, then (4,1,1,1,1,1) gives  $\psi(E,I,y)=4A-\frac{1}{2}A^5EIy+E+I+2+y>9,$ where $y=FGH >1/2$ by Claim (i). The function $\psi(E,I,y)$ is symmetric in $E$ and $I$, also it is linear in $E$, $I$ and $y$. So for $1<E\leq A$, $1<I\leq A$ and $\frac{1}{2}<y\leq 1$, we have $\psi(E,I,y)\leq{\rm max}\{\psi(1,1,\frac{1}{2}),\psi(1,1,1),\psi(A,1,\frac{1}{2}),\psi(A,1,1),\psi(A,A,\frac{1}{2}),\psi(A,A,1)\}$, which can be easily seen to be less than $9$ for $1<A\leq 2$. So $A^4EFGHI< 2$.\vspace{1mm}\\Now
$A^{4}EFGHI < 2$~~and~~$E^{4}ABCDI < 2$ together gives $A^{4}E^{4}I<4$, which gives $A^4<4$, $E^8<A^4E^4<4$ and $I^5<A^4I<4$, i.e. $A< \sqrt{2}$,~ $E<1.189208$~and~$I<1.31951$.\vspace{3mm}

{\noindent \bf Final Contradiction:}

Applying AM-GM inequality to (3,1,3,1,1), we get  $4A+D+4E+H+I-2A^{2}E^{2}D^{\frac{1}{2}}H^{\frac{1}{2}}I^{\frac{1}{2}}>9$. Left side is a decreasing
function of $D$ as $A\geq D$ and $A^3E^4HI>A^3E^4(1-h)>(1+a)^3(1+e)^4(1-e-a)\geq1$, for $0<a<\sqrt{2}-1$ and $0<e<0.189208$. So we can replace $D$ by 1 and get that $4A+4E+H+I-2A^{2}E^{2}H^{\frac{1}{2}}I^{\frac{1}{2}}>8$.
Now left side is a decreasing function of $H$, replacing $H$ by $1-e-\frac{a+i}{2}$ and on simplifying we get\\
$\theta(i)=2+\frac{7a}{2}+3e+\frac{i}{2}-2(1+a)^{2}(1+e)^{2}\sqrt{(1-e-\frac{a+i}{2})(1+i)}>0.$\\
A simple calculation gives that $\theta''(i)>0$ and $0<i\leq{\rm min}(a,0.31951)$. Therefore $\theta(i)\leq{\rm max}\{\theta(0),\theta({\rm min}(a,0.31951))\}$, which
is at most zero for $0<a<\sqrt{2}-1$ and $0<e\leq{\rm min}(a,0.189208)$. This gives a contradiction\vspace{2mm}\\
{\noindent \bf  Proposition 16.} Case (59) i.e. $A > 1,~ B > 1, ~C \leq 1, ~D \leq 1,~ E > 1, ~F \leq 1,~ G > 1, ~H > 1, ~I > 1$ does not arise.\vspace{2mm}\\
{\noindent \bf Proof.} Here $a\leq \frac{1}{2}$,~$b\leq \frac{1}{3}$,~ $e\leq \frac{1}{3}$, ~$f\leq\frac{1}{4}$.
Using the weak inequalities (2,2,2,1,1,1) and (1,2,1,2,1,1,1) we have
\begin{equation}2b-2d-2f+g+h+i>0,\vspace{-3mm}\end{equation}
\begin{equation}a-2c-d-2f+g+h+i>0.\end{equation}
{\noindent \bf Claim (i) $f\geq b$ ~{\rm and} ~$d<\frac{g+h+i}{2}$}

Suppose $f<b$, then we get contradiction using Lemma 7(ii), with  $x_{1}=b\leq \frac{1}{3}$ and $\gamma=f<x_{1}$. Hence we must have $f\geq b$.
Now using (4.10) we have $d<\frac{g+h+i}{2}$.\vspace{3mm}

{\noindent \bf Claim (ii) $a\geq 0.359$}

Suppose $a<0.359$. The inequality (3,1,2,1,1,1) gives $4A-\frac{A^3}{BC}+D+4E-\frac{2E^2}{F}+G+H+I>9$. Now using that $4E-\frac{E^2}{F}<3E$
and also applying AM-GM inequality to
$\frac{-A^3}{BC}$ and $\frac{-E^2}{F}$, we get $4A+3E+D+G+H+I-2\sqrt{A^{4}E^{3}DGHI}>9$. Left side is a decreasing function of $E$ as $E>1$ and $A^4DGHI>(1+a)^4(1+g)(1+h)(1+i)(1-d)>(1+a)^4(1+g)(1+h)(1+i)(1-\frac{g+h+i}{2})\geq1$, for $0<g\leq a$, $0<h\leq a$, $0<i\leq a$ and $0<a\leq 0.5.$
So we replace $E$ by 1 to get $4A+3+D+G+H+I-2\sqrt{A^{4}DGHI}>9$. Again left side is a decreasing function of $D$,
so we replace $D$ by $1-\frac{g+h+i}{2}$ to get  $\phi(g,h,i)=2+4a+\frac{g+h+i}{2}-2\sqrt{(1+a)^{4}(1-\frac{g+h+i}{2})(1+g)(1+h)(1+i)}>0$.
 Following Remark 2, we find that ${\rm max}~\phi(g,h,i)\leq0$ for $0<a<0.359$, a contradiction. Hence we must have $a\geq0.359$.\vspace{3mm}

{\noindent \bf Final Contradiction:}

Let $g+h+i=k$, then from (4.11) we get $2f<a+k$. So $AFGHI>(1+a+k)(1-f)>(1+2f)(1-f)>1$, for $0\leq f<\frac{1}{4}$. It gives $A^{4}EFGHI>A^3>2$, so (4,2,1,1,1) holds. Applying AM-GM, we get $\psi(G,H,I)=4A+4E+G+H+I-2\sqrt{A^{5}E^{3}GHI}>9.$ Now $\psi(G,H,I)$ is a decreasing function of $G$, $H$ and $I$, so $\psi(G,H,I)<\psi(1,1,1)<9$ for $1.359\leq A\leq 1.5$ and $1<E\leq \frac{4}{3}$. Hence we get a contradiction.\vspace{2mm}\\
{\noindent \bf  Proposition 17.} Case (60) i.e. $A > 1,~ B > 1, ~C \leq 1, ~D \leq 1,~ E > 1, ~F \leq 1,~ G > 1, ~H \leq 1, ~I > 1$ does not arise.\vspace{2mm}\\
{\noindent \bf Proof.} Here $a\leq \frac{1}{2}$,~$b\leq \frac{1}{3}$,~ $e\leq \frac{1}{3}$. Using the weak inequalities (2,2,2,2,1) and (1,2,1,2,2,1) we have \vspace{-2mm}
\begin{equation}2b-2d-2f-2h+i>0,\vspace{-3mm}\end{equation}
\begin{equation}a-2c-d-2f-2h+i>0.\end{equation}
{\noindent \bf Claim (i) $A<1.226$}

Suppose $A\geq 1.226$. From (4.13) we get $f+h<\frac{a+i}{2}$, then $AFHI>(1+a)(1+i)(1-\frac{a+i}{2})>1$ for $0.226\leq a\leq 0.5$ and $0<i\leq a$.
Hence $A^4EFGHI>2$
and so (4,1,1,1,1,1) holds, i.e. $4A-\frac{1}{2}A^5EGI(1-(f+h))+E+G+I+2-(f+h)>9$. Left side is an increasing function
of $f+h$, so we replace $f+h$ by $\frac{a+i}{2}$ to get
$3.5a+e+g+\frac{i}{2}-\frac{1}{2}(1+a)^5(1+e)(1+g)(1+i)(1-\frac{a+i}{2})>0.$\\
Now left side is decreasing function of $e$ and $g$, so replacing $e$ by 0 and $g$ by 0, we get
$3.5a+\frac{i}{2}-\frac{1}{2}(1+a)^5(1+i)(1-\frac{a+i}{2})>0,$
which is not true for $0<a\leq 0.5$ and $0<i\leq a.$

{\noindent \bf Claim (ii) $f+h>2b$}

Suppose $f+h\leq2b$. Also $b\leq a<0.226$. Now we get a contradiction using Lemma 7(v) with $x_{1}=b$, $\delta=a+e+g+i$ and $\gamma=f+h$, as $b<0.226$, $f+h\leq2b$ and $f+h<\frac{a+i}{2}<\frac{a+e+g+i}{2}$, using (4.13). Hence we must have $f+h>2b$.\vspace{3mm}

{\noindent \bf Final Contradiction:}

From (4.12) we get $f+h<b+\frac{i}{2}$. Further using Claim (ii) we get $b<\frac{i}{2}$.  Also using (4.12) and Claim (ii), we have $d<\frac{i}{2}$. The inequality (1,2,1,1,1,1,1,1) gives
$A+4B+D+E+F+G+H+I-2B^{3}ADEFGHI>9$. The left side is a decreasing function of $D$, so we replace $D$ by $1-\frac{i}{2}$ and get that
$2+a+4b+e+g-(f+h)+\frac{i}{2}-2(1+b)^{3}(1+a)(1+e)(1+g)(1-(f+h))(1-\frac{i}{2})(1+i)>0$. Now the left side is an increasing function of $(f+h)$ as
$(1-\frac{i}{2})(1+i)>1$ for $0<i\leq a<0.226$, so replace $f+h$ by $b+\frac{i}{2}$ and get that $\phi(e)=2+a+3b+e+g-2(1+b)^{3}(1+a)(1+e)(1+g)(1-b-\frac{i}{2})(1-\frac{i}{2})(1+i)>0$.
One can check that $\phi'(e)<0$, so we replace $e$ by 0. Similarly $g$ can also be replaced by 0. So we have now
\begin{equation}\psi(i)=2+a+3b-2(1+b)^{3}(1+a)(1-b-\frac{i}{2})(1-\frac{i}{2})(1+i)>0\end{equation}
We find that $\psi''(i)>0$, therefore $\psi(i)\leq{\rm max}\{\psi(0),\psi(a)\}$, which can be seen to be at most zero for $0<b\leq a<0.226$, contradicting (4.14).\vspace{3mm}

{\noindent \bf  Proposition 18.} Case (61) i.e. $A > 1,~ B > 1, ~C \leq 1, ~D \leq 1,~ E > 1, ~F \leq 1,~ G \leq 1, ~H > 1, ~I > 1$ does not arise.\vspace{3mm}\\
{\noindent \bf Proof.} Here $a\leq \frac{1}{2}$,~$b\leq \frac{1}{3}$,~ $e\leq \frac{1}{3}$.
Using the weak inequalities (1,2,1,2,1,1,1),~\\ (1,2,1,1,2,1,1) and~ (2,2,2,1,1,1)  we have \vspace{-2mm}
\begin{equation}a-2c-d-2f-g+h+i>0,\vspace{-3mm}\end{equation}
\begin{equation}a-2c-d+e-2g+h+i>0,\vspace{-2mm}\end{equation}
\begin{equation}2b-2d-2f-g+h+i>0.\end{equation}
{\noindent \bf Claim (i) $B<1.05$}

Suppose $B\geq 1.05$. Applying AM-GM inequality to (1,3,3,1,1) we get $A+4B+4E+H+I-2\sqrt{B^4E^4AHI}>9.$ Left side is a decreasing function of $E$, $H$ and $I$, so replacing each of them by 1 we get $A+4B-2\sqrt{B^4A}>3$,
which is not true for $1<A\leq 1.5$, and $1.05\leq B\leq \frac{4}{3}$. Hence $B<1.05$.\vspace{3mm}

{\noindent \bf Claim (ii) $f+g\geq 2b$}

Suppose $f+g<2b$. From (4.15) and (4.16), we get $2f+g<a+h+i<a+e+h+i$ and $g<\frac{a+e+h+i}{2}$ respectively. Adding these two inequalities we get
$\gamma=f+g<\frac{3}{4}(a+e+h+i)=\frac{3}{4}\delta$. Now we get contradiction by Lemma 7(vi) with $x_1=b<0.05$, $\gamma=f+g<2x_1$ and $\delta=a+e+h+i$.\vspace{3mm}

{\noindent \bf Final Contradiction:}

Using $f+g\geq 2b$ in (4.17) we get $d<\frac{h+i}{2}$. Applying AM-GM inequality to (3,1,3,1,1) we get \\ $2+4a+4e-d+h+i-2\sqrt{(1+a)^{4}(1+e)^{4}(1-d)(1+h)(1+i)}>0.$\\  Left side is an increasing function of $d$, so we replace $d$ by $\frac{x}{2}$,
where $x=h+i$ and get that $2+4a+4e+\frac{x}{2}-2(1+a)^{2}(1+e)^{2}\sqrt{(1-\frac{x}{2})(1+x)}>0$. Now the left side is a decreasing function
of $e$ as $(1-\frac{x}{2})(1+x)\geq 1$ for $0<x=h+i\leq 2a\leq 1$, so we replace $e$ by 0 and get that $2+4a+\frac{x}{2}-2(1+a)^{2}\sqrt{(1-\frac{x}{2})(1+x)}>0$,
which is not true for $0<a\leq 0.5$ and $0<x\leq 2a$, giving a contradiction.\vspace{2mm}\\
{\noindent \bf Proposition 19.} Case (38) i.e. $A > 1,~ B > 1, ~C \leq 1, ~D > 1,~ E > 1, ~F > 1,~ G \leq 1, ~H \leq 1, ~I \leq 1$ does not arise.\vspace{3mm}\\
{\noindent \bf Proof.} Here $a\leq \frac{1}{2}$,~$b\leq \frac{1}{3}$. Using the weak inequality (1,2,1,2,1,2) we have \\
$a-2c+d+2f-g-2i>0$, i.e. $i<f+\frac{a+d}{2}$.\vspace{2mm}

{\noindent \bf Claim(i) $F^{4}ABCDE < 2$~ {\rm and}~$GHI>\frac{1}{2}$}

Suppose $F^{4}ABCDE \geq 2$, then (5*,4) holds, i.e. $5(ABCDE)^{\frac{1}{5}}+4F-\frac{1}{2}F^{5}ABCDE > 9$, i.e. $\phi(x)
=5x^{\frac{1}{5}}+4F-\frac{1}{2}F^{5}x>9$, where $x = ABCDE$.  Now $\phi'(x) = 0$ gives $x = (\frac{2}{F^{5}})^{5/4}$. Since $\phi''(x) < 0$ at $x = (\frac{2}{F^{5}})^{5/4}$, so
 $\phi(x) \leq \phi((\frac{2}{F^{5}})^{5/4}) < 9$, a contradiction. So $F^{4}ABCDE < 2$.
It gives $\frac{F^{3}}{GHI}<2$ and hence $GHI > \frac{1}{2}$.\vspace{3mm}

{\noindent \bf Claim(ii) $A^{4}EFGHI < 2$,~ $A < \sqrt{2}$~ {\rm and}~ $F<1.189208$}

Proof is similar to that of Claim (ii) of Case (42)(Proposition 15).\vspace{3mm}

{\noindent \bf Claim(iii) $A\leq1.26$}

Suppose $A>1.26$. Applying AM-GM inequality to (3,2,3,1) we get
$4A+4D+4F+I-3(2A^4D^3F^4I)^{\frac{1}{3}}>9$. Left side is a decreasing function of $I$, so replace $I$ by $1-f-\frac{a+d}{2}$. On simplifying we get\\
$\chi(d)=4+\frac{7a}{2}+\frac{7d}{2}+3f-3(2)^{\frac{1}{3}}(1+a)^{\frac{4}{3}}(1+f)^{\frac{4}{3}}(1+d)(1-\frac{a+d}{2}-f)^{\frac{1}{3}}>0.$\\
We find that $\chi''(d)>0$ and we have $0<d\leq a$. So $\chi(d)\leq {\rm max}\{\chi(0),\chi(a)\}$. Now $\chi(0)$ and $\chi(a)$ are functions in two variables $a$ and $f$ and can be checked to be less than zero for $0.26<a<\sqrt{2}-1$ and $0<f<0.189208$, giving thereby a contradiction. Hence we must have $a\leq0.26$.\vspace{3mm}

{\noindent \bf Final Contradiction:}\\
Applying AM-GM inequality to (3,1,1,3,1) we get
 $4A+4F+D+E+I-2A^{2}F^{2}D^{\frac{1}{2}}E^{\frac{1}{2}}I^{\frac{1}{2}}>9$. Left side is a decreasing function of $E$ as ~ $A^3F^4DI>(1+a)^3(1+f)^4(1-f-\frac{a+d}{2})\geq(1+a)^3(1+f)^4(1-f-a)\geq 1$, for $0<a\leq 0.26$ and $0<f<0.189208$. So replacing $E$ by 1 we get $4A+4F+D+I-2A^{2}F^{2}D^{\frac{1}{2}}I^{\frac{1}{2}}>8.$ Now the left side is a decreasing function of $I$
for $I\leq 1$. So we replace $I$ by $1-f-\frac{a+d}{2}$.
On simplifying we get $\psi(d)=2+\frac{7a}{2}+\frac{d}{2}+3f-2(1+a)^{2}(1+f)^{2}\sqrt{(1+d)(1-f-\frac{a+d}{2})}>0.$
We find that $\psi''(d)>0$ and we have $0<d\leq a$. So $\psi(d)\leq{\rm max}\{\psi(0),\psi(a)\}<0$ for $0<a\leq 0.26$ and $0<f<0.189208$, giving thereby a contradiction. \vspace{2mm}\\
{\noindent \bf  Proposition 20.} Case (58) i.e. $A > 1,~ B > 1, ~C \leq 1, ~D \leq 1,~ E > 1, ~F > 1,~ G \leq 1, ~H \leq 1, ~I \leq 1$ does not arise.\vspace{3mm}\\
{\noindent \bf Proof.} Here $a\leq \frac{1}{2}$,~$b\leq \frac{1}{3}$,~ $e\leq \frac{1}{2}$,~ $f\leq \frac{1}{3}$.
Using the weak inequalities (1,2,1,2,2,1) and (1,2,1,1,2,2) we have \vspace{-2mm}
\begin{equation}a-2c-d+2f-2h-i>0,\vspace{-3mm}\end{equation}
\begin{equation}a-2c-d+e-2g-2i>0.\end{equation}
{\noindent \bf Claim (i) $F^{4}ABCDE<2$~ {\rm and}~ $GHI>\frac{1}{2}$}

Proof is same as that of Claim (i) of Case (38)(Proposition 19). \vspace{2mm}

{\noindent \bf Claim (ii) $A^{4}EFGHI<2$~ {\rm and}~ $F<1.189208$}

Proof is similar to that of Claim (ii) of Case (42)(Proposition 15). \vspace{2mm}

{\noindent \bf Claim (iii) $c+d+i>f$~ {\rm and}~ $h<\frac{a+f}{2}$}

Suppose $c+d+i\leq f$. Now using Lemma 7(ii), with $x_1=f<0.189208$ and $\gamma=c+d+i\leq f=x_1$, we get a contradiction. Hence $c+d+i>f$.
Now using (4.18) we get $2h<(a+2f)-(c+d+i)$, which gives $h<\frac{a+f}{2}$.\vspace{2mm}

{\noindent \bf Claim (iv) $A<1.284$}

Suppose $A\geq 1.284$. From (4.19) we get $g+i<\frac{a+e}{2}$. Now $A^{4}EFGHI>(1+a)^{4}(1+e)(1+f)(1-\frac{a+e}{2})(1-\frac{a+f}{2})=\phi(e,f)$.
Following Remark 2 we find that ${\rm min}~\phi(e,f)\geq2$, for $0.284\leq a\leq 0.5$. Hence $A^{4}EFGHI>2$, contradicting Claim (ii).\vspace{2mm}

{\noindent \bf Claim (v) $c+h\geq \frac{a+f}{2}$}

Suppose $c+h<\frac{a+f}{2}$. Applying AM-GM inequality to (2,2,1,2,2) we get $4A+4C+E+4F+4H-8(A^{3}C^{3}F^{3}H^{3}E)^{\frac{1}{4}}>9$. Let $x=c+h$, then we get
$8+4a+4f-4x+e-8(1+a)^{\frac{3}{4}}(1+f)^{\frac{3}{4}}(1-x)^{\frac{3}{4}}(1+e)^{\frac{1}{4}}>0.$
Left side is an increasing function of $x$, so replacing $x$ by $\frac{a+f}{2}$,
we get $\psi(e)=8+2a+2f+e-8(1+a)^{\frac{3}{4}}(1+f)^{\frac{3}{4}}(1-\frac{a+f}{2})^{\frac{3}{4}}(1+e)^{\frac{1}{4}}>0$.
Now $\psi'(e)<0$, so replace $e$ by 0, i.e.
$\psi(e)\leq \psi(0)=8+2a+2f-8(1+a)^{\frac{3}{4}}(1+f)^{\frac{3}{4}}(1-\frac{a+f}{2})^{\frac{3}{4}}$, which is at most zero for $0<a<0.284$~ and~ $0<f<0.189208$; a contradiction.\vspace{2mm}

{\noindent \bf Final Contradiction:}

Using (4.18) and Claim (v), we have $d+i<f$. Applying AM-GM inequality  to (3,1,1,3,1) we get
$2+4a+4f+e-(d+i)-2(1+a)^{2}(1+f)^{2}\sqrt{(1+e)(1-(d+i))}>0.$
Left side is an increasing function of $(d+i)$, so replacing $(d+i)$ by $f$ we get $2+4a+3f+e-2(1+a)^{2}(1+f)^{2}\sqrt{(1+e)(1-f)}>0$. Now left side is
decreasing function of $e$, so we replace $e$ by 0 and get that
$2+4a+3f-2(1+a)^{2}(1+f)^{2}\sqrt{(1-f)}>0$, which is not true for $0<f<0.189208$ and $f\leq a<0.284$.\vspace{3mm}\\
{\noindent \bf  Proposition 21.} Case (62) i.e. $A > 1,~ B > 1, ~C \leq 1, ~D \leq 1,~ E > 1, ~F \leq 1,~ G \leq 1, ~H \leq 1, ~I > 1$ does not arise.\vspace{3mm}\\
{\noindent \bf Proof.} Here $a\leq \frac{1}{2}$,~$b\leq \frac{1}{3}$,~ $e\leq \frac{1}{3}$.
Using the weak inequalities (1,2,1,2,1,1,1),~ \\(1,2,1,1,2,1,1) and~ (1,2,1,2,2,1)  we have \vspace{-2mm}
\begin{equation}a-2c-d-2f-g-h+i>0,\vspace{-3mm}\end{equation}
\begin{equation}a-2c-d+e-2g-h+i>0,\vspace{-2mm}\end{equation}
\begin{equation}a-2c-d-2f-2h+i>0.\end{equation}
{\noindent \bf Claim (i) $E^{4}ABCDI<2$~{\rm and}~ $FGH>\frac{1}{2}$}

Proof is same as that of Claim (i) of Case (42)(Proposition 15).\vspace{2mm}

{\noindent \bf Claim (ii) $A^{4}EFGHI<2$~ {\rm and}~ $E<1.189208$}

Proof is same as that of Claim (ii) of Case (42)(Proposition 15). \vspace{2mm}

{\noindent \bf Claim (iii) $A<1.301$}

Suppose $A\geq1.301$. Using (4.20), ~(4.21) and (4.22) respectively we get $f<\frac{a+i}{2}-\frac{g+h}{2}$, ~$g<\frac{a+e+i}{2}-\frac{h}{2}$ ~and ~$h<\frac{a+i}{2}\leq a$.
Therefore  $A^{4}EFGHI=(1+a)^{4}(1+e)(1+i)(1-f)(1-g)(1-h)>(1+a)^{4}(1+e)(1+i)(1-\frac{a+i}{2}+\frac{g+h}{2})(1-g)(1-h)=\phi(g)$, say. As $\phi'(g)<0$, $\phi(g)>\phi(\frac{a+e+i}{2}-\frac{h}{2})$.
Now $\phi(\frac{a+e+i}{2}-\frac{h}{2})=(1+a)^{4}(1+e)(1+i)(1-h)(1-\frac{a+i}{4}+\frac{h}{4}+\frac{e}{4})(1-\frac{a+e+i}{2}+\frac{h}{2})=\eta(h).$
We find that $\eta''(h)<0$. So $\eta(h)\geq{\rm min}\{\eta(0),\eta(a)\}$. Let $\eta(0)=\theta_{1}(i)$ and  $\eta(a)=\theta_{2}(i).$
It is easy to check that for $m=1,2$, $\theta_{m}''(i)<0$, therefore $\theta_{m}(i)\geq {\rm min}\{\theta_{m}(0),\theta_{m}(a)\}$. But for $m=1,2$, $\theta_{m}(0)$ and $\theta_{m}(a)$ are greater than 2 for $0.301\leq a\leq 0.5$ and $0<e<0.189208$ i.e.
$A^4EFGHI>2$,
contradicting Claim (ii). Hence we must have $a<0.301$.\vspace{2mm}

{\noindent \bf Claim (iv) $c+g>\frac{a+e}{2}$}

Suppose $c+g\leq\frac{a+e}{2}$. Applying AM-GM inequality to (2,2,2,2,1) we get $4A+4C+4E+4G+I-8(A^{3}C^{3}E^{3}G^{3}I)^{\frac{1}{4}}>9$, i.e.
$8+4a+4e+i-4(c+g)-8(1+a)^{\frac{3}{4}}(1+e)^{\frac{3}{4}}(1+i)^{\frac{1}{4}}(1-(c+g))^{\frac{3}{4}}>0$. Here the left side is an increasing function of $c+g$,
so we replace $c+g$ by $\frac{a+e}{2}$ to get $8+2a+2e+i-8(1+a)^{\frac{3}{4}}(1+e)^{\frac{3}{4}}(1+i)^{\frac{1}{4}}(1-\frac{a+e}{2})^{\frac{3}{4}}>0.$
Now the left side is a decreasing function of $i$ as $2(1+e)^{\frac{3}{4}}(1-\frac{a+e}{2})^{\frac{3}{4}}>1$, so on replacing $i$ by 0, we get
$8+2a+2e-8(1+a)^{\frac{3}{4}}(1+e)^{\frac{3}{4}}(1-\frac{a+e}{2})^{\frac{3}{4}}>0.$ But this is not true for $0<a<0.301$ and $0<e<0.189208$. Hence we must have $c+g>\frac{a+e}{2}$.\vspace{3mm}

{\noindent \bf Final Contradiction:}

We have $c+g>\frac{a+e}{2}$, using it with (4.21), we get $d+h<i$. Applying AM-GM  inequality to (3,1,3,1,1) we have  $2+4a+4e-d-h+i-2\sqrt{(1+a)^{4}(1+e)^{4}(1-(d+h))(1+i)}>0.$  Left side is an increasing function of $(d+h)$, so we replace $(d+h)$ by $i$ and get that
$2+4a+4e-2(1+a)^{2}(1+e)^{2}\sqrt{(1-i)(1+i)}>0.$ Now the left side is an increasing function of $i$, so on replacing $i$ by $a$, we get
$2+4a+4e-2(1+a)^{2}(1+e)^{2}\sqrt{(1-a)(1+a)}>0.$ But this is not true for $0<a<0.301$ and $0<e<0.189208$.\vspace{3mm}

{\noindent \bf  Proposition 22.} Case (65) and Case (66) do not exist.\vspace{3mm}\\
{\noindent \bf Proof.}
Consider first Case (65) i.e.   $A > 1,~ B > 1, ~C \leq 1, ~D \leq 1,~ E \leq 1, ~F > 1,~ G > 1, ~H \leq 1, ~I > 1$\vspace{3mm}\\
 Here $0<a\leq \frac{1}{2}$ and $0<b\leq \frac{1}{3}$. Using the weak inequalities (1,2,2,1,2,1), (2,2,1,1,2,1) and (2,1,2,1,2,1), we have \vspace{-2mm}
\begin{equation}a-2c-2e+f-2h+i>0,\vspace{-3mm}\end{equation}
\begin{equation}2b-2d-e+f-2h+i>0,\vspace{-2mm}\end{equation}
\begin{equation}2b-c-2e+f-2h+i>0.\end{equation}
{\noindent \bf  Claim (i) $A<1.226$}

Suppose $A\geq 1.226$. Using (4.23), we have $A^4EFGHI\geq (1+a)^{4}(1+f)(1+i)(1-(e+h))>(1+a)^{4}(1+f)(1+i)(1-\frac{a+f+i}{2})=\phi(f,i)$. Following Remark 2 we find that ${\rm min}~\phi(f,i)>2$, for $0.226\leq a\leq 0.5$. Hence $A^{4}EFGHI>2$ and so (4,1,1,1,1,1) holds,
i.e. $4a-\frac{1}{2}(1+a)^{5}(1-(e+h))(1+f)(1+g)(1+i)-(e+h)+f+g+i>0$. Left side is an increasing function of
$e+h$ for $A\geq 1.226$, so we replace $e+h$ by $\frac{a+f+i}{2}$ and get that $\theta(g)=\frac{7a}{2}+\frac{f+i}{2}+g-\frac{1}{2}(1+a)^{5}(1+f)(1+i)(1+g)(1-\frac{a+f+i}{2})>0$.
Now $\theta'(g)<0$ for $0<f\leq a$ and $0<i\leq a$. So we can replace $g$ by 0 and get that
 $\psi(f,i)=\frac{7a}{2}+\frac{f+i}{2}-\frac{1}{2}(1+a)^{5}(1+f)(1+i)(1-\frac{a+f+i}{2})>0.$ Following Remark 2 we find that ${\rm max}~\psi(f,i)\leq 0$ for $0<a\leq 0.5$, a contradiction.\vspace{2mm}

{\noindent \bf  Claim (ii) $e+h>2b$ {\rm and} $d<\frac{f+i}{2}$}

Suppose $e+h\leq 2b$. We get a contradiction by Lemma 7(v) with $x_{1}=b\leq a<0.226$,
$\gamma=e+h\leq 2x_{1}$ and $\gamma<\frac{a+f+g+i}{2}=\frac{\delta}{2}$,
using (4.23). Hence we must have $e+h>2b$. Together with (4.24), we get $d<\frac{f+i}{2}$.\vspace{2mm}

{\noindent \bf  Claim (iii) $A\geq 1.21$ {\rm and} $f+i>a$}

Suppose $A<1.21$. The inequality (1,2,1,1,1,1,1,1) gives $A+4B-2B^{3}ADEFGHI+D+E+F+G+H+I>9.$ Coefficients of $G$ and $D$ are negative on the left side for $1<G\leq \frac{4}{3}$ and $D\leq 1$ as $2B^3AEFHI=\frac{2B^2}{CDG}\geq\frac{2}{4/3}>1$, so we replace $G$ by 1 and $D$ by $1-\frac{f+i}{2}$ and get that $2+a+4b-(e+h)+\frac{f+i}{2}-2(1+b)^{3}(1+a)(1+f)(1+i)(1-(e+h))(1-\frac{f+i}{2})>0.$ Now the left side is an increasing
function of $(e+h)$ as $(1+f)(1+i)(1-\frac{f+i}{2})\geq 1$, so we replace $(e+h)$ by $b+\frac{f+i}{2}$, using (4.25) and get that
\begin{equation}\eta(f,i)=2+a+3b-2(1+b)^{3}(1+a)(1-b-\frac{f+i}{2})(1-\frac{f+i}{2})(1+f)(1+i)>0.\end{equation}
Following Remark 2 we find that ${\rm max}~\eta(f,i)<0$ for $0<a<0.21$ and $0<b\leq a$. It contradicts (4.26). Hence we must have $a\geq0.21$.\\
 Further if $f+i\leq a$, then we have $\eta(f,i)<{\rm max}\{\eta(0,0),\eta(a,0)\}\leq0$, in full range of $a$ i.e.
in $0<a<0.226$ and $0<b\leq a$. Hence $f+i>a$.
\vspace{2mm}

{\noindent \bf  Final Contradiction:}

We have $1.21\leq A<1.226$, $e+h<\frac{a+f+i}{2}$ and $f+i>a$.\\
$A^{4}EFGHI>(1+a)^{4}(1+f+i)(1-(e+h))>(1+a)^{4}(1+x)(1-\frac{a+x}{2})=\varphi(x)$, where $x=f+i$. Now $\varphi''(x)<0$, so for $a<x\leq 2a$, we have $\varphi(x)\geq{\rm min}\{\varphi(a),\varphi(2a)\}>2$, for $0.21\leq a<0.226$. Hence $A^4EFGHI>2$. Now we use (4,1,1,1,1,1)  and get contradiction proceeding as
in Claim (i).\vspace{2mm}

Case (66) can be dealt with working exactly same as in Case (65) using the suitable inequalities.
\vspace{2mm}\\
{\noindent \bf  Proposition 23.} Case (67) i.e. $A > 1,~ B > 1, ~C \leq 1, ~D \leq 1,~ E \leq 1, ~F > 1,~ G \leq 1, ~H \leq 1, ~I > 1$ does not arise.\vspace{2mm}\\
{\noindent \bf Proof.} Here $a\leq \frac{1}{2}$,~$b\leq \frac{1}{3}$,~$f\leq \frac{1}{3}$.
Using the weak inequalities (1,2,2,1,2,1), (1,2,1,1,2,1,1)~and (1,2,2,2,1,1)  we have
\begin{equation}a-2c-2e+f-2h+i>0,\vspace{-4mm}\end{equation}
\begin{equation}a-2c-d-e-2g-h+i>0,\vspace{-2mm}\end{equation}
\begin{equation}a-2c-2e-2g-h+i>0.\end{equation}
{\noindent \bf Claim (i)~~$A<1.35$}

Suppose $A\geq 1.35$.  $A^4EFGHI>\frac{3}{4}(1+a)^{4}(1+f)(1+i)(1-(e+h))>\frac{3}{4}(1+a)^{4}(1+f)(1+i)(1-\frac{a+f+i}{2})=\eta(f,i)$, using (4.26) and $G>\frac{3}{4}$.
Following Remark 2 we find that ${\rm min}~\eta(f,i)>2$, for $0.35\leq a\leq 0.5$. Hence $A^4EFGHI>2$. Applying AM-GM inequality to (4,1,2,1,1) we get $2+4a+4f-(e+h)+i-2\sqrt{(1+a)^{5}(1+f)^{3}(1-(e+h))(1+i)}>0$. Now the left side is an increasing
function of $(e+h)$, so we replace $(e+h)$ by $\frac{a+f+i}{2}$ and get
$\eta(i)=2+\frac{7a}{2}+\frac{7f}{2}+\frac{i}{2}-2\sqrt{(1+a)^{5}(1+f)^{3}(1-\frac{a+f+i}{2})(1+i)}>0.$
Now $\eta''(i)>0$, therefore $\eta(i)\leq{\rm max}\{\eta(0),\eta(a)\},$ which can be seen to be at most zero for $0.35\leq a\leq 0.5$ and $0<f\leq \frac{1}{3}$, giving thereby a contradiction.\vspace{2mm}

{\noindent \bf Claim (ii)~~$d+e>\frac{a+i}{2}$}

Suppose $d+e\leq\frac{a+i}{2}$. Applying AM-GM inequality to (3,1,1,3,1) we get $2+4a+4f+i-(d+e)-2\sqrt{(1+a)^4(1+f)^4(1+i)(1-(d+e))}>0.$
Left side is an increasing function of $(d+e)$, so we replace $(d+e)$ by $\frac{a+i}{2}$ and get that $2+\frac{7a}{2}+4f+\frac{i}{2}-2\sqrt{(1+a)^4(1+f)^4(1+i)(1-\frac{a+i}{2})}>0$.
Now the left side is a decreasing function of $f$ as $(1+a)(1+i)(1-\frac{a+i}{2})\geq1$, so we replace $f$ by 0 and get that
$2+\frac{7a}{2}+\frac{i}{2}-2\sqrt{(1+a)^4(1+i)(1-\frac{a+i}{2})}>0$, which is not true for $0<a<0.35$ and $0<i\leq a$. Hence $d+e>\frac{a+i}{2}$.\vspace{2mm}

{\noindent \bf  Final Contradiction:}

Using $d+e>\frac{a+i}{2}$ together with (4.28), we get $c+g<\frac{a+i}{4}$. Applying AM-GM inequality to (2,2,1,1,2,1) we get $4A+4C+4G+E+F+I-6(A^3C^3G^3EFI)^{\frac{1}{3}}>9$. Here the left side is a decreasing function of $E$ as
$ACG>(1+a)(1-\frac{a+i}{4})>2(1+a)(1-\frac{a}{2})>1$.
From (4.29) we have $E=1-e>1-\frac{a+i}{2}+(c+g)$.
 So we replace  $E$ by $1-\frac{a+i}{2}+(c+g)$ and get that
$\phi(x)=6+\frac{7a}{2}-3x+f+\frac{i}{2}-6(1+a)(1-x)(1-\frac{a+i}{2}+x)^{\frac{1}{3}}(1+f)^{\frac{1}{3}}(1+i)^{\frac{1}{3}}>0,$\\
where $x=c+g$.
It is easy to check that $\phi''(x)>0$, so $\phi(x)<{\rm max}\{\phi(0),\phi(\frac{a+i}{4})\}$.\\
Let $\phi(0)=\eta_{1}(f)$ and $\phi(\frac{a+i}{4})=\eta_{2}(f).$
We find that $\eta_{1}(f)$ and $\eta_{2}(f)$ are decreasing functions of $f$. Therefore
$\eta_1(f)\leq\eta_{1}(0)$ and $\eta_2(f)\leq\eta_{2}(0)$. But $\eta_{1}(0)$ and $\eta_{2}(0)$ are at most zero for $0<a<0.35$ and $0<i\leq a$. Hence we have a contradiction.
\section{Difficult Cases}
{\noindent \bf  Proposition 24.} Case (2) i.e. $ A>1,~ B>1, ~C>1, ~D>1,~ E>1, ~ F>1, ~G\leq 1,~ H\leq 1,~ I\leq 1$ does not
arise.\vspace{2mm}\\ {\noindent \bf Proof.}\vspace{2mm}

\noindent{\bf Claim(i)} $F^4ABCDE< 2, ~F<1.1487$ ~and~ $ GHI > 1/2$

Suppose $F^4ABCDE\geq2$. We use the inequality $(1,1,1,1,1,4)$ and get contradiction applying Lemma 10(ii) with $X_{2}=F,~X_{3}=B,~X_{4}=C,~X_{5}=D,~X_{6}=E$. Therefore $F^4ABCDE< 2$. This implies $F^5\leq 2$ which gives $F<1.1487$. Also $ GHI=\frac{1}{ABCDEF}> \frac{F^{3}}{2} >\frac{1}{2}.$\vspace{2mm}

\noindent{\bf Claim(ii) $A^4EFGHI< 2 {\rm ~and} ~A<\sqrt{2}$}

Proof is similar to that of Claim(ii) of Case (42)(Proposition 15).\vspace{2mm}

\noindent{\bf Claim(iii)} $B< 1.31951$

Suppose  $B\geq1.31951$. Then using $GHI >1/2$, we get $B^4AFGHI>B^5\times \frac{1}{2}>2$. So the inequality (1,4,1,1,1,1) holds. i.e. $\chi(y,F)=A+4B-\frac{1}{2}B^5yAF+2+y+F>9,$ where $y=GHI>\frac{1}{2}$ by Claim(i). $\chi(y,F)$ is a decreasing function of $y$ and a decreasing function of $F$ as well for $A\geq B\geq1.31951$. So we have $\chi(y,F)\leq\chi(\frac{1}{2},1)$, which is at most 9 for $1<B\leq A\leq \sqrt{2}$. This gives a contradiction. \vspace{2mm}

\noindent{\bf Claim(iv)} $C<1.2174$

Suppose $C\geq 1.2174$. Using $(2,2,1,4^{*})$ and  the AM-GM inequality we get
$\varphi(A,C,E,x)=4A+4C+E-4A^{\frac{3}{2}}C^{\frac{3}{2}}E^{\frac{1}{2}}x^{\frac{1}{2}}+4x^{\frac{1}{4}}>9,$ where $x=FGHI$. Since $x > GHI >
\frac{1}{2}$ and $\varphi$ is a decreasing function of $x$, we get
$\varphi(A,C,E,x)\leq \varphi(A,C,E,\frac{1}{2}).$  Further $\varphi(A,C,E,\frac{1}{2})$ is a decreasing function of $E$ and  $E>1$, so $\varphi(A,C,E,\frac{1}{2})<\varphi(A,C,1,\frac{1}{2})$  which
can be easily  verified to be less than $9$ for $A\geq C\geq 1.2174.$ This gives a contradiction.\vspace{1mm}\\
Also we see that $\varphi(A,C,1,\frac{1}{2})<9$ for $A>1.27$ and $C>1.171$. Hence we can take
\begin{equation} C\leq1.171 ~{\rm if} ~A>1.27.\end{equation}
\noindent{\bf Claim(v)} $D<1.2174$

Suppose $D\geq 1.2174$. Using $(2,1,2,4^{*})$ and  AM-GM inequality we get
$\varphi(A,D,C,x)=4A+C+4D-4A^{\frac{3}{2}}D^{\frac{3}{2}}C^{\frac{1}{2}}x^{\frac{1}{2}}+4x^{\frac{1}{4}}>9,$ where $x=FGHI$. Working as in Claim(iv), we get $\varphi(A,D,C,x)\leq\varphi(A,D,1,\frac{1}{2})<9$ for $A\geq D\geq 1.2174.$ This gives a contradiction.\vspace{1mm}\\
Also we see that $\varphi(A,D,1,\frac{1}{2})<9$ for $A>1.27$ and $D>1.171$. Hence we can take
\begin{equation}D\leq1.171~ {\rm  if}~ A>1.27.\end{equation}
Further if $C>1.182$, then ~$\varphi(A,D,C,\frac{1}{2})<\varphi(A,D,1.182,\frac{1}{2})\leq 9$ for $A\geq C> 1.182$ and $D>1.151.$ This gives a contradiction. So we can take
\begin{equation}C\leq1.182~ {\rm  if}~ D>1.151.\end{equation}
\noindent{\bf Claim(vi)} $B\leq1.182$ if $D>1.151$

Suppose $D>1.151$. Using $(1,2,2,4^{*})$ and  AM-GM inequality we get
$\varphi(B,D,A,x)=A+4B+4D-4B^{\frac{3}{2}}D^{\frac{3}{2}}A^{\frac{1}{2}}x^{\frac{1}{2}}+4x^{\frac{1}{4}}>9,$ where $x=FGHI$. Working as in Claim(iv), we get  $\varphi(B,D,A,x)\leq\varphi(B,1.151,A,\frac{1}{2})\leq9$ for $A\geq B> 1.182,$ a contradiction. So we can take
$B\leq1.182$ if $D>1.151$.\vspace{1mm}\\
Also if we have $D>1.136$, then $\varphi(B,D,A,x)\leq\varphi(B,1.136,A,\frac{1}{2})\leq9$ for $B>1.16$ and $A>1.27$, a contradiction. So we can take
\begin{equation}B\leq1.16~ {\rm  if}~A>1.27~ {\rm  and}~D>1.136.\end{equation}
Further if we have $D>1.168$, then $\varphi(B,D,A,x)\leq\varphi(B,1.168,A,\frac{1}{2})\leq9$ for $B> 1.168$ and $A>1.173$, a contradiction. So we can take
\begin{equation}D\leq1.168~ {\rm  if}~B>1.168~ {\rm  and}~A>1.173.\end{equation}
\noindent{\bf Final Contradiction}

\indent Using $(3,1,1,3,1)$ and  AM-GM inequality we get $4A+4F+D+E+I-2A^2F^2\sqrt{D}\sqrt{E}\sqrt{I}$ $> 9$.
Left side of this inequality is  a quadratic in $\sqrt{I}.$ Since $A^4F^4DE-4A-4F-D-E+9>0$, we have
\begin{equation}\sqrt{I}<A^2F^2\sqrt{D}\sqrt{E}-(A^4F^4DE-4A-4F-D-E+9)^{\frac{1}{2}}=\alpha {\rm ~~ (say).}\end{equation} Using AM-GM inequality in $(1,2,2,2,2)$, we
get $A+2C+4D+4F+4H-6DFHA^{\frac{1}{3}}B^{\frac{1}{3}}C^{\frac{1}{3}}$ $>9$, which gives
$H<(A+2C+4D+4F-9)(6DFA^{\frac{1}{3}}B^{\frac{1}{3}}C^{\frac{1}{3}}-4)^{-1}$. Substituting this upper bound of $H$ in the weak inequality
$(2,2,2,2,1)$, we get
\begin{equation}I>9-2B-2D-2F-2\left\{\frac{A+2C+4D+4F-9}{6DFA^{\frac{1}{3}}B^{\frac{1}{3}}C^{\frac{1}{3}}-4}\right\}=\beta {\rm ~~
(say).}\end{equation} From (5.6) and (5.7) we have $\beta<\alpha^2$. On simplifying we
get\begin{equation}\begin{array}{l}\eta(B)=\left\{A^4F^4DE-2A-F+B+\frac{D}{2}-\frac{E}{2}+\frac{A+2C+4D+4F-9}{6DFA^{\frac{1}{3}}B^{\frac{1}{3}}C^{\frac{1}{3}}-4}
\right\}^2~~\\~~-A^4F^4\{A^4F^4D^2E^2-4ADE-4FDE-D^2E-E^2D+9DE\}>0.\end{array}\end{equation}
Now
$$\begin{array}{l}\eta^{'}(B)=2\left\{A^4F^4DE-2A-F+B+\frac{D}{2}-\frac{E}{2}+\frac{A+2C+4D+4F-9}{6DFA^{\frac{1}{3}}B^{\frac{1}{3}}C^{\frac{1}{3}}-4}
\right\}~~~~~~~~~~~~~~~\\~~~~~~~~~~~~~~ \times
\left\{1-\frac{2DFA^{\frac{1}{3}}C^{\frac{1}{3}}B^{\frac{-2}{3}}(A+2C+4D+4F-9)}{(6DFA^{\frac{1}{3}}B^{\frac{1}{3}}C^{\frac{1}{3}}-4)^2}\right
\}.\end{array}$$One finds that $\eta^{'}(B)>0$ if $\psi(A, B, C, D, F)=(6DFA^{\frac{1}{3}}B^{\frac{1}{3}}C^{\frac{1}{3}}-4)^2-2DFA^{\frac{1}{3}}C^{\frac{1}{3}}B^{\frac{-2}{3}}(A+2C+4D+4F-9)>0$. As $\psi(A, B, C, D, F)$ is an increasing function of $A,~B, ~C, ~D$ and $F$, we get $\psi(A, B, C, D, F)>\psi(1, 1, 1, 1, 1)>0$. Hence $\eta(B)$ is an
increasing function of $B$, for all $B>1$. Take

 \begin{equation}\lambda=\left \{ \begin{array}{l}A~~~~~~{\rm if}~~1<A\leq 1.173\\A~~~~~~{\rm if}~~1.173<A\leq 1.27~{\rm and}~1<D\leq 1.151\\1.168~~{\rm if}~~1.173<A\leq 1.182, ~D>1.151~{\rm and}~B\leq1.168\\A~~~~~~{\rm if}~~ 1.173<A\leq 1.182, ~D>1.151~{\rm and}~B>1.168\\1.168~~{\rm if}~~1.182<A\leq 1.27, ~D>1.151~{\rm and}~B\leq1.168\\1.182~~{\rm if}~~ 1.182<A\leq 1.27, ~D>1.151~{\rm and}~B>1.168\\{\rm min}(A,1.31951)~~{\rm if} ~~1.27<A<\sqrt{2}~{\rm and}~D\leq1.136\\1.16~~~~{\rm if}~~1.27<A<\sqrt{2}~ {\rm and}~ D>1.136. \end{array}\right.{}\end{equation}
From Claim(iii), Claim (vi) and (5.4) we have $B\leq \lambda$.
So $\eta(B)\leq\eta(\lambda)$, which gives:
\begin{equation}\begin{array}{l}\zeta(C)=\left\{A^4F^4DE-2A-F+\lambda+\frac{D}{2}-\frac{E}{2}+\frac{A+2C+4D+4F-9}{6DFA^{\frac{1}{3}}\lambda^{\frac{1}{3}}C^{\frac{1}{3}}-4}
\right\}^2~~\\~~-A^4F^4\{A^4F^4D^2E^2-4ADE-4FDE-D^2E-E^2D+9DE\}>0.\end{array}\end{equation}
Let
$$f(C)=A^4F^4DE-2A-F+\lambda+\frac{D}{2}-\frac{E}{2}+\frac{A+2C+4D+4F-9}{6DFA^{\frac{1}{3}}\lambda^{\frac{1}{3}}C^{\frac{1}{3}}-4}. $$
So we have $\zeta^{''}(C)=2(f^{'}(C))^2+2f(C)f^{''}(C)$, as the second summand in $\zeta(C)$ is independent of $C$. One can easily show that $\zeta^{''}(C)>0$ by proving that $f^{''}(C)>0$ and $f(C)>0$, for all $C$, $1<C\leq {\rm min}(A,1.2174)$. \\Take  $\mu_1=1$ and
\begin{equation}\mu_2=\left \{ \begin{array}{ll}A&{\rm if}~~1<A\leq 1.173\\1.2174&{\rm if}~~1.173<A\leq 1.27~{\rm and}~D\leq1.151\\A&{\rm if}~~ 1.173<A\leq 1.182~{\rm and}~D>1.151\\&({\rm whether}~ B>1.168 ~{\rm or}~ B\leq1.168)\\1.182&{\rm if}~~1.182<A\leq 1.27~{\rm and}~D>1.151\\&({\rm whether}~ B>1.168 ~{\rm or}~ B\leq1.168)\\1.171&{\rm if} ~~1.27<A<\sqrt{2}.\end{array}\right.{} \end{equation}
 From Claim(iv), (5.1) and (5.3), we have $\mu_1 \leq C \leq \mu_2$. This implies $\zeta(C)\leq{\rm max}\{\zeta(\mu_1),\zeta(\mu_2)\}$. Write $\zeta(\mu_i)=\psi_{i}(D)$, where $i=1, 2$. So\\
$\psi_{i}(D)=\left\{A^4F^4DE-2A-F+\lambda+\frac{D}{2}-\frac{E}{2}+\frac{A+2\mu_i+4D+4F-9}{6DFA^{\frac{1}{3}}\lambda^{\frac{1}{3}}\mu_i^{\frac{1}{3}}-4}
\right\}^2~~\\~~~~~~~~-A^4F^4\{A^4F^4D^2E^2-4ADE-4FDE-D^2E-E^2D+9DE\}.$\\
$~~~~~~~~~=(\chi_1(D))^2-\chi_2(D)$, say.\vspace{1mm}\\ Therefore
$\psi_{i}^{''}(D)=2(\chi_{1}'(D))^2+2\chi_{1}(D)\chi_{1}''(D)-\chi_{2}''(D),$ where
$$\chi_1(D)=A^4F^4DE-2A-F+\lambda+\frac{D}{2}-\frac{E}{2}+\frac{A+2\mu_i+4D+4F-9}{6DFA^{\frac{1}{3}}\lambda^{\frac{1}{3}}\mu_i^{\frac{1}{3}}-4} $$
$$\chi_2(D)=A^4F^4\{A^4F^4D^2E^2-4ADE-4FDE-D^2E-E^2D+9DE\} $$
One can show
that $\psi_{i}^{''}(D)>0$ by proving  that
$~2(\chi_{1}'(D))^2-\chi_{2}''(D)>0, ~\chi_{1}(D)>0~{\rm and~}\chi_{1}''(D)>0$, for all $D$, $1<D\leq {\rm min}(A,1.2174)$.\\ Take
\begin{equation}(\nu_1,\nu_2)=\left \{ \begin{array}{l}(1,A)~~~~~~~~{\rm if}~~1<A\leq 1.173\\(1,1.151)~~~{\rm if}~~1.173<A\leq 1.27~{\rm and}~1<D\leq1.151\\(1.151,A)~~{\rm if}~~1.173<A\leq 1.182,~D>1.151~{\rm and}~B\leq1.168\\(1.151,1.168)~~{\rm if}~~1.173<A\leq 1.182,~D>1.151~{\rm and}~B>1.168\\(1.151,{\rm min}(A,1.2174))~~{\rm if}~~1.182<A\leq 1.27, D>1.151\\~~~~~~~~~~~~~~~~~~~~~~~~~~~~~~~~~~~~~{\rm and}~~B\leq1.168\\(1.151,1.168)~{\rm if}~~1.182<A\leq 1.27, D>1.151~{\rm and}~B>1.168\\(1,1.136)~~~~~~{\rm if} ~~1.27<A<\sqrt{2}~~{\rm and}~1<D\leq1.136\\(1.136,1.171)~{\rm if} ~~1.27<A<\sqrt{2}~~{\rm and}~D>1.136 .\end{array}\right.{} \end{equation}
 From Claim(v), (5.2) and (5.5), we have $\nu_1 \leq D \leq \nu_2$. This implies $\psi_{i}(D)\leq{\rm max}\{\psi_{i}(\nu_1),\psi_{i}(\nu_2)\}$. Let  $\psi_{i}(\nu_j)=\psi_{ij}(E)$, for $1\leq i \leq 2$, ~$1\leq j \leq 2$, where\\
 $\psi_{ij}(E)=\left\{A^4F^4\nu_jE-2A-F+\lambda+\frac{\nu_j}{2}-\frac{E}{2}+\frac{A+2\mu_i+4\nu_j+4F-9}{6\nu_jFA^{\frac{1}{3}}\lambda^{\frac{1}{3}}\mu_i^{\frac{1}{3}}-4}
\right\}^2~~\\~~-A^4F^4\{A^4F^4\nu_j^2E^2-4A\nu_jE-4F\nu_jE-\nu_j^2E-E^2\nu_j+9\nu_jE\}.$\\
It is easy to check that $\psi_{ij}^{''}(E)>0$. So $\psi_{ij}^{'}(E)$ is an increasing function of $E$. So for $E\leq A$, $\psi_{ij}^{'}(E)\leq \psi_{ij}^{'}(A)$, which can be seen to be negative for all $A$ and $F$, $1<A<\sqrt{2}$, $1<F<1.1487$. Hence $\psi_{ij}(E)$ is a decreasing function of $E$. So for $E>1$, $\psi_{ij}(E)<\psi_{ij}(1)=\theta_{ij}(A,F)$, say, where
\begin{equation}\begin{array}{l}\theta_{ij}(A,F)=\left\{A^4F^4\nu_j-2A-F+\lambda+\frac{\nu_j}{2}-\frac{1}{2}+\frac{A+2\mu_i+4\nu_j+4F-9}{6\nu_jFA^{\frac{1}{3}}\lambda^{\frac{1}{3}}\mu_i^{\frac{1}{3}}-4}
\right\}^2~~~~\\~~~~-A^4F^4\{A^4F^4\nu_{j}^2-4A\nu_j-4F\nu_j-\nu_{j}^2-\nu_j+9\nu_j\}.\end{array}\end{equation}
Also using $F<1.1487$ and $F^4\leq \frac{2}{A}$ (from Claim (i)), we can take
\begin{equation}F<\left \{ \begin{array}{lll}1.1487&{\rm if}&1<A\leq 1.27\\1.121&{\rm if}&1.27<A<\sqrt{2}. \end{array}\right.{} \end{equation}
Now for $1\leq i \leq 2$, $1\leq j \leq 2$, $\theta_{ij}(A,F)$ are functions in two variables $A$ and $F$ and can be seen to be negative for $A$ and $F$ lying in the corresponding ranges as defined in (5.14) and for~ $\lambda$, $\mu_i$, $\nu_j$ as defined in (5.9), (5.11) and (5.12) for each $i$ and $j$. This contradicts (5.10) and hence (5.8).\vspace{3mm}\\
{\noindent \bf  Proposition 25.} Case (5) i.e. $ A>1,~ B>1, ~C>1, ~D>1,~ E>1, ~ F\leq1, ~G\leq 1,~ H\leq 1,~ I> 1$ does not
arise.\vspace{2mm}\\ {\noindent \bf Proof.}\vspace{2mm}

\noindent{\bf Claim(i)} $E^4ABCDI<2, ~E<1.1487$ ~{\rm and}~ $ FGH > 1/2$

Suppose $E^4ABCDI\geq2$. We use the inequality $(1,1,1,1,4,1)$ and get contradiction applying Lemma 10(ii) with $X_{2}=E,~X_{3}=B,~X_{4}=C,~X_{5}=D,~X_{6}=I$. Therefore $E^4ABCDI< 2$. This implies $E^5\leq 2$ which gives $E<1.1487$. Also $ FGH=\frac{1}{ABCDEI}> \frac{E^{3}}{2} >\frac{1}{2}.$\vspace{2mm}

\noindent{\bf Claim(ii) $A^4EFGHI< 2, ~A<\sqrt{2}  {\rm ~and} ~I<1.31951$}

Proof is similar to Claim(ii) of Case (42)(Proposition 15).\vspace{2mm}

\noindent{\bf Claim(iii)} $C<1.2174$

Suppose $C\geq 1.2174$. Using $(2,2,4^{*},1)$ and  AM-GM inequality we get
$\varphi(A,C,I,x)=4A+4C+I-4A^{\frac{3}{2}}C^{\frac{3}{2}}I^{\frac{1}{2}}x^{\frac{1}{2}}+4x^{\frac{1}{4}}>9,$ where $x=EFGH$. Now working similarly as in Claim(iv) of Case (2)(Proposition 24), we get $\varphi(A,C,I,x)\leq \varphi(A,C,1,\frac{1}{2})\leq9$ for $A\geq C\geq 1.2174.$ This gives a contradiction.\vspace{1mm}\\
Further if we have $I>1.086$, then $\varphi(A,C,I,x)\leq \varphi(A,C,1.086,\frac{1}{2})\leq 9$ for $A\geq C>1.191$. So we must have
\begin{equation}C\leq1.191 ~{\rm if}~ I>1.086,\end{equation}
Similarly we get \vspace{-2mm}
\begin{equation}I\leq1.165 ~{\rm if}~ A>1.182 ~{\rm and}~ C>1.158;\vspace{-2mm}\end{equation}
\begin{equation}C\leq 1.062 ~{\rm if}~ A>1.27 ~{\rm and}~ I>1.24;\vspace{-2mm}\end{equation}
\begin{equation}C\leq1.171 ~{\rm if}~ A>1.27;\vspace{-2mm}\end{equation}
\begin{equation}C\leq1.135 ~{\rm if}~ A>1.31951;\vspace{-2mm}\end{equation}
\begin{equation}C\leq1.102 ~{\rm if}~ A>1.373.\end{equation}
\noindent{\bf Final Contradiction}

Using $(3,1,3,1,1)$ and  AM-GM inequality we get $4A+4E+D+H+I-2A^2E^2\sqrt{D}\sqrt{H}\sqrt{I}$ $> 9$.
Left side of this inequality is  a quadratic in $\sqrt{H}.$ Since $A^4E^4DI-4A-4E-D-I+9>0$, we have
\begin{equation}\sqrt{H}<A^2E^2\sqrt{D}\sqrt{I}-(A^4E^4DI-4A-4E-D-I+9)^{\frac{1}{2}}=\alpha {\rm ~~ (say).}\end{equation}  Using AM-GM inequality in $(2,2,2,2,1)$, we
get $G<(4A+2D+4E+I-9)(6AEC^{\frac{1}{3}}D^{\frac{1}{3}}I^{\frac{1}{3}}-4)^{-1}$. Substituting this upper bound of $G$ in the inequality
$(1,2,2,2,1,1)$, we get
\begin{equation}H>9-A-2C-2E-I-2\left\{\frac{4A+2D+4E+I-9}{6AEC^{\frac{1}{3}}D^{\frac{1}{3}}I^{\frac{1}{3}}-4}\right\}=\beta {\rm ~~
(say).}\end{equation} From (5.21) and (5.22) we have $\beta<\alpha^2$. On simplifying we
get\begin{equation}\begin{array}{l}\eta(C)=A^4E^4DI-\frac{3A}{2}-E-\frac{D}{2}+C+\frac{4A+2D+4E+I-9}{6AEC^{\frac{1}{3}}D^{\frac{1}{3}}I^{\frac{1}{3}}-4}
~~\\~~-A^2E^2\sqrt{A^4E^4D^2I^2-4ADI-4EDI-D^2I-I^2D+9DI}>0.\end{array}\end{equation}
Working as in Case(2)(Proposition 24) we find that $\eta(C)$ is an
increasing function of $C$, for all $C>1$.\\ Take
 \begin{equation}\lambda=\left \{ \begin{array}{l}A~~~~~~~{\rm if}~~1<A\leq 1.182\\{\rm min}(A,1.2174)~~{\rm if}~~1.182<A\leq 1.27 ~{\rm and}~I\leq1.086\\1.158~~{\rm if}~~ 1.182<A\leq 1.27, I>1.086~{\rm and}~C\leq 1.158\\{\rm min}(A,1.191)~~{\rm if} ~~1.182<A\leq1.27, I>1.086~{\rm and}~C>1.158\\1.171~~~~{\rm if}~~1.27<A\leq1.31951~{\rm and}~I\leq 1.24\\1.062~~~~{\rm if}~~1.27<A\leq1.31951~{\rm and}~I>1.24\\1.135~~~~{\rm if}~~1.31951<A\leq1.373\\1.102~~~~{\rm if}~~1.373<A<\sqrt{2}. \end{array}\right.{} \end{equation}
Claim(iii) and (5.15), (5.17), (5.18), (5.19), (5.20) imply  $C\leq \lambda$. So $\eta(C)\leq
\eta(\lambda)$, which gives
\begin{equation}\begin{array}{l}\zeta(I)=A^4E^4DI-\frac{3A}{2}-E-\frac{D}{2}+\lambda+\frac{4A+2D+4E+I-9}{6AE\lambda^{\frac{1}{3}}D^{\frac{1}{3}}I^{\frac{1}{3}}-4}
~~\\~~-A^2E^2\sqrt{A^4E^4D^2I^2-4ADI-4EDI-D^2I-I^2D+9DI}>0.\end{array}\end{equation}
Write $\zeta(I)=\phi_1(I)-\sqrt{\phi_2(I)}$, so that
$\zeta^{''}(I)=\phi_1^{''}(I)+\frac{(\phi_2^{'}(I))^2-2\phi_2(I)\phi_2^{''}(I)}{4(\phi_2(I))^{3/2}}. $
One can easily show that $\zeta^{''}(I)>0$, by proving that $\phi_1^{''}(I)>0$,  $(\phi_2^{'}(I))^2-2\phi_2(I)\phi_2^{''}(I)>0$  and  $\phi_2(I)>0$.\\
Take
\begin{equation}(\mu_1,\mu_2)=\left \{ \begin{array}{ll}(1, A)&{\rm if}~~1<A\leq 1.182\\(1,1.086)&{\rm if}~~1.182<A\leq 1.27 ~{\rm and}~I\leq1.086\\(1.086, A)&{\rm if}~~1.182<A\leq 1.27, I>1.086~{\rm and}~C\leq 1.158\\(1.086, 1.165)&{\rm if} ~~1.182<A\leq1.27, I>1.086~{\rm and}~C>1.158\\(1,1.24)&{\rm if}~~1.27<A\leq1.31951~{\rm and}~I\leq 1.24\\(1.24, A)&{\rm if}~~1.27<A\leq1.31951~{\rm and}~I>1.24\\(1,1.31951)&{\rm if}~~1.31951<A<\sqrt{2}. \end{array}\right.{} \end{equation}

Claim(ii) and (5.16) imply $\zeta(I)\leq{\rm max}\{\zeta(\mu_1),\zeta(\mu_2)\}$. For $j= 1, 2$, let
\begin{equation*}\begin{array}{l}\zeta(\mu_j)=\eta_{j}(D)=A^4E^4D\mu_j
-\frac{3A}{2}-E-\frac{D}{2}+\lambda+\frac{4A+2D+4E+\mu_j-9}{6AED^{\frac{1}{3}}\lambda^{\frac{1}{3}}\mu_j^{\frac{1}{3}}-4}
~~~~~~~~~~~~~~~\\~~~~~~~~~~~~~~~-A^2E^2\sqrt{A^4E^4D^2\mu_j^2-4AD\mu_j-4ED\mu_j-D^2\mu_j-\mu_j^2D+9D\mu_j}.\end{array}\end{equation*}
$$\begin{array}{l}\eta_{j}'(D)=A^4E^4\mu_j-\frac{1}{2}+\frac{8EA\lambda^{\frac{1}{3}}D^{\frac{1}{3}}\mu_{j}^{\frac{1}{3}}-8-
2(4A+4E+\mu_j-9)EA\lambda^{\frac{1}{3}}D^{\frac{-2}{3}}\mu_{j}^{\frac{1}{3}}}{(6EA\lambda^{\frac{1}{3}}D^{\frac{1}{3}}\mu_{j}^{\frac{1}{3}}-4)^2}
~~~~~\\~~~~~~~~~
-\frac{A^2E^2\{2A^4E^4D\mu_j^2-4A\mu_j-4E\mu_j-2D\mu_j-\mu_j^2+9\mu_j\}}{2\sqrt{A^4E^4D^2\mu_j^2-4AD\mu_j-4ED\mu_j-\mu_j^2D-D^2\mu_j+9D\mu_j}}.
\end{array}$$
Now we prove that $\eta_{j}^{'}(D)<0$ for $1<D\leq A$, $1<A<\sqrt{2}$ and $1<E<1.1487$. Let \begin{equation}(\alpha_1, \alpha_2)= \left \{ \begin{array}{lll}(0.438,0.062) & {\rm if} & A\leq 1.31951\\(0.45,0.05) & {\rm if} &  1.31951<A<\sqrt{2}.
\end{array}\right.{} \end{equation}

Let $\eta_{j}'(D)=P+Q $, where   $$\begin{array}{l}P= A^4E^4\mu_j-\alpha_1
-\frac{A^2E^2\{2A^4E^4D\mu_j^2-4A\mu_j-4E\mu_j-2D\mu_j-\mu_j^2+9\mu_j\}}{2\sqrt{A^4E^4D^2\mu_j^2-4AD\mu_j-4ED\mu_j-\mu_j^2D-D^2\mu_j+9D\mu_j}}.
\end{array}$$ and
$$\begin{array}{l}
 Q =\frac{8EA\lambda^{\frac{1}{3}}D^{\frac{1}{3}}\mu_{j}^{\frac{1}{3}}-8-
2(4A+4E+\mu_j-9)EA\lambda^{\frac{1}{3}}D^{\frac{-2}{3}}\mu_{j}^{\frac{1}{3}}}{(6EA\lambda^{\frac{1}{3}}D^{\frac{1}{3}}\mu_{j}^{\frac{1}{3}}-4)^2} -\alpha_2.\end{array}$$ To prove $\eta_{j}'(D)<0$, we  show that $P<0$ and $Q<0$
for $1<D\leq A,~ 1<A<\sqrt{2}$ and $1<E<1.1487$. \vspace{2mm}\\
Now $P<0$ ~if~  $\theta_{1}(D)=4(A^4E^4\mu_j-\alpha_1)^2(A^4E^4D^2\mu_j^2-4AD\mu_j-4ED\mu_j-\mu_j^2D-D^2\mu_j+9D\mu_j)- A^4E^4\{2A^4E^4D\mu_j^2-4A\mu_j-4E\mu_j-2D\mu_j-\mu_j^2+9\mu_j\}^2<0.$ As
$\theta_{1}(D)$ is an increasing function of $D$, so $\theta_{1}(D)\leq \theta_{1}(A)$ which is a function in two
variables $A$ and $E$ and can be shown to be less than zero for $1<E<1.1487$, $1<A<\sqrt{2}$ and for $\alpha_1$, $(\mu_1,\mu_2)$ as given in (5.27) and (5.26) respectively.\\
Further $Q<0$ if $\theta_{2}(D,E)=8EA\lambda^{\frac{1}{3}}D^{\frac{1}{3}}\mu_{j}^{\frac{1}{3}}-8-
2(4A+4E+\mu_j-9)EA\lambda^{\frac{1}{3}}D^{\frac{-2}{3}}\mu_{j}^{\frac{1}{3}}-\alpha_2(6EA\lambda^{\frac{1}{3}}D^{\frac{1}{3}}\mu_{j}^{\frac{1}{3}}-4)^2<0.$ One finds that
$\theta_{2}(D,E)$ is a decreasing function of $E$ and an increasing function of $D$, therefore $\theta_{2}(D,E)\leq
\theta_{2}(A,1)$, which is less than zero for $1<A<\sqrt{2}$ and for $\alpha_1$, $(\mu_1,\mu_2)$ as given in (5.27) and (5.26) respectively. \\Thus $\eta_{j}(D)$ is a decreasing function of $D$, therefore
$\eta_{j}(D)\leq \eta_{j}(1)$, where
$$\begin{array}{l}\eta_{1}(1)=A^4E^4\mu_1-\frac{3A}{2}-E-\frac{1}{2}+\lambda+\frac{4A+2+4E+\mu_1-9}{6AE\lambda^{\frac{1}{3}}\mu_1^{\frac{1}{3}}-4}
~~~~~~~~~~~~~~~\\~~~~~~~~~~~~~~~-A^2E^2\sqrt{A^4E^4\mu_1^2-4A\mu_1-4E\mu_1-\mu_1-\mu_1^2+9\mu_1}=\chi_1(A,E), {\rm say}.\end{array}$$

$$\begin{array}{l}\eta_{2}(1)=A^4E^4\mu_2-\frac{3A}{2}-E-\frac{1}{2}+\lambda+\frac{4A+2+4E+\mu_2-9}{6AE\lambda^{\frac{1}{3}}\mu_2^{\frac{1}{3}}-4}
~~~~~~~~~~~~~~~\\~~~~~~~~~~~~~~~-A^2E^2\sqrt{A^4E^4\mu_2^2-4A\mu_2-4E\mu_2-\mu_2-\mu_2^2+9\mu_2}=\chi_2(A,E), {\rm say}.\end{array}$$
Also using $E<1.1487$ and $E<(\frac{2}{AI})^{1/4}$, from Claim (i), we can take
\begin{equation}E<\left \{ \begin{array}{lll}1.1487&{\rm if}&1<A\leq 1.182\\1.1173&{\rm if}&A>1.182~{\rm and}~I>1.086\\1.1203&{\rm if}&A>1.27\\1.1096&{\rm if}&A>1.31951.\end{array}\right.{} \end{equation}
 Now $\chi_{1}(A,E)$ and $\chi_{2}(A,E)$ are functions in two variables $A$ and $E$ and can be seen to be negative for $\lambda$, $\mu_1$, $\mu_2$ as defined in (5.24) and (5.26) and for $A$, $E$ lying in the corresponding ranges as defined in (5.28). This contradicts (5.25) and hence (5.23).\vspace{2mm}\\
{\noindent \bf  Proposition 26.} Case (10) i.e. $A > 1,~ B > 1, ~C > 1, ~D > 1,~ E \leq 1, ~F > 1,~ G \leq 1, ~H \leq 1, ~I \leq 1$ does not arise.\vspace{3mm}\\
{\noindent \bf Proof.} Here $A<2.1326324$, ~$B\leq 2$, ~$C\leq \frac{3}{2}$ and ~$D\leq \frac{4}{3}$.
Using the weak inequalities (1,2,2,1,2,1) ~and~  (1,2,2,2,2) we have
\begin{equation}a+2c-2e+f-2h-i>0,\end{equation}
\begin{equation}a+2c-2e-2g-2i>0.\end{equation}

\noindent{\bf Claim(i) $F^{4}ABCDE<2$, ~$GHI>\frac{1}{2}$ ~{\rm and} ~$F<1.22$}

Proof is same as that of Claim (i) of Case (38)(Proposition 19).\\
Further $2>F^4ABCDE>F^5E\geq F^{5}(\frac{3}{4})$ gives $F<1.22$.\vspace{2mm}

{\noindent \bf Claim(ii) $A < 1.52$}

Suppose $A \geq 1.52$, then using $GHI>\frac{1}{2}$ and $E\geq\frac{3}{4}$ we have $A^4EFGHI>(1.52)^{4}\times \frac{3}{4}\times \frac{1}{2}>2$. Therefore (4,1,1,1,1,1) holds, i.e. $\phi(y,E,F)=4A-\frac{1}{2}A^5EFy+E+F+2+y>9,$ where $y=GHI >1/2$. One easily checks that $\phi(y,E,F)$ is a decreasing function of each of the variables $y$, $E$ and $F$, so $\phi(y,E,F)\leq\phi(\frac{1}{2},\frac{3}{4},1)\leq9$ for $1.52\leq A<2.1326324$, a contradiction.\vspace{2mm}

{\noindent \bf Claim(iii) $C \geq 1.19$}

Suppose $C<1.19$. Using (5.29) and (5.30) we have $h<\frac{a+f}{2}+c-\frac{e+i}{2}$ and $e+i<\frac{a}{2}+c$, respectively. Applying AM-GM inequality to (2,2,1,2,1,1) we get $4A+4C+E+4F+H+I-6ACF(EHI)^{\frac{1}{3}}>9$. Left side is a decreasing function of $H$, as $A^3C^3F^3EI\geq ABCDEFI=\frac{1}{GH}\geq 1$, so we replace $H$ by $1-(\frac{a}{2}+c+\frac{f}{2}-\frac{e+i}{2})$ and get that $\phi(x)=6+\frac{7a}{2}+3c+\frac{7f}{2}-\frac{x}{2}-6(1+a)(1+c)(1+f)(1-x)^{\frac{1}{3}}(1-\frac{a}{2}-\frac{f}{2}-c+\frac{x}{2})^{\frac{1}{3}}>0,$ where $x=e+i.$ We find that $\phi''(x)>0$, therefore $\phi(x)\leq{\rm max}\{\phi(0),\phi(\frac{a}{2}+c)\}$. Let $\phi(0)=\varphi_{1}(f)$ and $\phi(\frac{a}{2}+c)=\varphi_{2}(f).$
Now for $m=1,2$, ~$\varphi''_{m}(f)>0$, so $\varphi_{m}(f)<{\rm max}\{\varphi_{m}(0),\varphi_{m}(0.22)\}$, which can be seen to be less than $0$ for $0<c\leq{\rm min}(a,0.19)$ and $0<a<0.52$. Hence we must have $C\geq1.19$.
\vspace{3mm}

{\noindent \bf Claim(iv) $C^4ABGHI<2$ ~{\rm and} ~$C < 1.32$}

Suppose $C^4ABGHI \geq 2$, then (1,1,4,1,1,1) holds, i.e. $\phi(y,B)=A+B+4C-\frac{1}{2}C^5ABy+2+y>9,$ where $y=GHI >1/2$ by Claim(i). $\phi(y,B)$ is a linear function of $y$ and $B$, so $\phi(y,B)\leq{\rm max}\{\phi(\frac{1}{2},1),\phi(\frac{1}{2},A),\phi(1,1),\phi(1,A)\}$, which can be seen to be less than  9 for $1.19\leq C\leq \frac{3}{2}$ and $1< A<1.52$. Hence $C^4ABGHI<2$.\vspace{1mm}\\
Further $2>C^4ABGHI>C^5GHI>C^{5}\times \frac{1}{2}$ gives $C< 1.32$.\vspace{2mm}

{\noindent \bf Claim(v) $C > D$}

Suppose $C\leq D$. Using (5.30) we have $i<\frac{a}{2}+c$. The inequality (3,2,3,1) holds, i.e. $4A-\frac{A^3}{BC}+4D-\frac{2D^2}{E}+4F-\frac{F^3}{GH}+I>9$. Using $E\leq 1$ and applying AM-GM inequality to $-\frac{A^3}{BC},-\frac{D^2}{E},-\frac{F^3}{GH}$, we get $4A+4D-D^2+4F+I-3(A^4D^3F^4I)^{\frac{1}{3}}>9.$ Note that left side is a decreasing function of $F$ and $D$ as $A^4I>(1+a)^4(1-\frac{a}{2}-c)>(1+a)^4(1-\frac{3a}{2})\geq1$, for $0<a<0.52$. So we replace $F$ by 1 and $D$ by $C$ to get that $\eta(I)=4A+4C-C^2+I-3(A^4C^3I)^{\frac{1}{3}}>5$. Now $\eta'(I)<0$, therefore $\eta(I)<\eta(1-\frac{a}{2}-c)$, which can be verified to be at most five for $0.19\leq c<0.32$ and $0<a<0.52$.\vspace{2mm}

{\noindent \bf Claim(vi) $i > c$ ~{\rm and} ~$e+g<\frac{a}{2}$}

Suppose $i\leq c$. Applying AM-GM inequality to (1,1,3,3,1) we have $A+B+4C+4F+I-2C^2F^2(ABI)^{\frac{1}{2}}>9$. Left side is a decreasing function of $F$ and $B$ as $C^4I>(1+c)^4(1-c)>1$, for
$0.19\leq c<0.32$, so we replace $F$ and $B$  by 1 to get that $A+4C+I-2C^2(AI)^{\frac{1}{2}}>4.$ Now the left side is a decreasing function of $I$, so we replace $I$ by $1-c$ and get that $2+a+3c-2(1+c)^{2}(1+a)^{\frac{1}{2}}(1-c)^{\frac{1}{2}}>0$, which is not true for $0.19\leq c\leq{\rm min}(a,0.32)$ and $0.19\leq a<0.52$. Hence we have $i\leq c$.\vspace{1mm}\\
Now using (5.30), we get $e+g<\frac{a}{2}$.\vspace{2mm}

{\noindent \bf Claim(vii) $A>1.38$ ~{\rm and} ~$B>1.035$}

Suppose $A\leq 1.38$. Applying AM-GM to the inequality (2,2,2,2,1) we get  $4A+4C+4E+4G+I-8A^{\frac{3}{4}}C^{\frac{3}{4}}E^{\frac{3}{4}}G^{\frac{3}{4}}I^{\frac{1}{4}}>9$. Here left side is a decreasing function of $I$ for $e+g<\frac{a}{2}$, so we replace $I$ by $1-\frac{a}{2}-c+e+g$, using (5.30) and get that $\phi(x)=8+\frac{7a}{2}+3c-3x-8(1+a)^{\frac{3}{4}}(1+c)^{\frac{3}{4}}(1-x)^{\frac{3}{4}}(1-\frac{a}{2}-c+x)^{\frac{1}{4}}>0$, where $e+g=x$. Now $\phi''(x)>0$ and $0\leq x<\frac{a}{2}$. Therefore $\phi(x)\leq{\rm max}\{\phi(0),\phi(\frac{a}{2})\}$.
Let $\phi(0)=\psi_{1}(a,c)$ and $\phi(\frac{a}{2})=\psi_{2}(a,c).$ We find that $\psi_{1}(a,c)<0$  for $0.19\leq c<0.32$ and $c\leq a\leq0.52$ and $\psi_{2}(a,c)<0$  for $0.19\leq c<0.32$ and $0.19\leq a\leq0.38$. Hence we must have $A>1.38$.\vspace{1mm}\\
As $B\geq\frac{3}{4}A$, we get $B>1.035$.\vspace{2mm}

{\noindent \bf Claim(viii) $C < 1.294$}

If $C\geq 1.294$, then using $GHI>\frac{1}{2}$, we have $C^4ABGHI>(1.294)^{4}\times1.38\times 1.035\times \frac{1}{2}>2$, a contradiction to Claim(iv). This proves $C< 1.294$.\vspace{2mm}

{\noindent \bf Claim(ix) $i > 1.157c$ ~{\rm and} ~$e+g<\frac{a}{2}-0.157c$}

Suppose $i\leq 1.157c$. We proceed as in Claim(vi). Here we still have $C^4I>1$, so we replace $F$ by 1, $B$ by 1.035 and $I$ by $1-1.157c$ to get that $2.035+a+2.843c-2\sqrt{1.035}(1+c)^{2}(1+a)^{\frac{1}{2}}(1-1.157c)^{\frac{1}{2}}>0$, which is not true for $0.19\leq c<0.294$ and $0.38<a<0.52$. Hence we have $i > 1.157c$.\vspace{1mm}\\
Now using (5.30), we get $e+g<\frac{a}{2}-0.157c$.\vspace{2mm}

{\noindent \bf Claim(x) $A>1.483$ ~{\rm and} ~$B>1.112$}

Suppose $A\leq 1.483$. Now we proceed as in Claim(vii) and use $x\leq \frac{a}{2}-0.157c$ instead of $x\leq \frac{a}{2}$. In place of $\psi_{2}(a,c)$ we have $\psi_{3}(a,c)=8+2a+3.471c-8(1+a)^{\frac{3}{4}}(1+c)^{\frac{3}{4}}(1-\frac{a}{2}+0.157c)^{\frac{3}{4}}(1-1.157c)^{\frac{1}{4}}$ and we find that $\psi_{3}(a,c)<0$  for $0.19\leq c<0.294$ and $0.38<a\leq0.483$. Hence we must have $A>1.483$ and $B\geq\frac{3}{4}A>1.112$.\vspace{3mm}

{\noindent \bf Claim(xi) $C < 1.248$}

If $C\geq 1.248$, then $C^4ABGHI>(1.248)^{4}\times1.483\times 1.112\times \frac{1}{2}>2$, a contradiction to Claim(iv).\vspace{3mm}

{\noindent \bf Claim(xii) $i > 1.35c$ ~{\rm and} ~$e+g<\frac{a}{2}-0.35c$}

Suppose $i\leq 1.35c$. Working as  in Claim (vi) and using that $B>1.112$ we get a contradiction for $0.19\leq c<0.248$ and $0.483<a<0.52$.\vspace{1mm}\\
Now using (5.30), we get $e+g<\frac{a}{2}-0.35c$.\vspace{2mm}

{\noindent \bf Final Contradiction}

We proceed as in Claim (vii) and use $x\leq \frac{a}{2}-0.35c$ instead of $x\leq \frac{a}{2}$. In place of $\psi_{2}(a,c)$ we have $\psi_{4}(a,c)=8+2a+4.05c-8(1+a)^{\frac{3}{4}}(1+c)^{\frac{3}{4}}(1-\frac{a}{2}+0.35c)^{\frac{3}{4}}(1-1.35c)^{\frac{1}{4}}$ and we find that $\psi_{4}(a,c)<0$  for $0.19\leq c<0.248$ and $0.483<a<0.52$. Hence we get a contradiction.\vspace{3mm}\\
{\noindent \bf  Proposition 27.} Case (13) i.e. $ A>1,~ B>1, ~C>1, ~D>1,~ E\leq1, ~ F\leq1, ~G\leq 1,~ H> 1,~ I> 1$ does not
arise.\vspace{2mm}\\ {\noindent \bf Proof.}~
The proof is similar to that of Case (2)(Proposition 24). We give here an outline of the proof omitting the details.

\noindent{\bf Claim(i)} $D^4ABCHI< 2, ~D<1.1487$ ~and~ $ EFG > 1/2$\\$~~~~~~~~~~~~~~~$(inequality used is (1,1,1,4,1,1).)\vspace{1mm}\\
\noindent{\bf Claim(ii) $A^4EFGHI<2, ~A<\sqrt{2}  {\rm ~and} ~I< 1.31951$}\\$~~~~~~~~~~~~~~~$(inequality used is (4,1,1,1,1,1).)\vspace{1mm}\\
\noindent{\bf Claim(iii)} ~$H<1.2174$\\$~~~~~~~~~~~~~~~$(inequality used is $(2,1,4^{*},2)$.)\vspace{1mm}\\
\noindent{\bf Claim(iv)} $B^4FGHIA<2$, ~$B^5HI\leq 4$  and $~B< 1.31951$\\$~~~~~~~~~~~~~~~$(inequality used is (1,4,1,1,1,1); other bounds used are $\frac{1}{2}<FG\leq 1$, $1<H<1.2174$, $1<I<1.31951$ and $1<B\leq A<\sqrt{2}$. $B^4FGHIA<2$ further gives ~$B^5HI\leq 4$  and $~B< 1.31951$ )\vspace{2mm}

\noindent{\bf Final Contradiction}

Applying AM-GM to the inequality $(3,3,1,1,1)$ we get a quadratic inequality in $\sqrt{G},$ which gives an upper bound $\alpha$ on  $\sqrt{G}$. Using AM-GM inequality in $(1,2,2,2,1,1)$, we
find an upper bound on $F$ and substituting this upper bound of $F$ in the weak inequality
$(2,2,2,1,1,1)$, we get a lower bound $\beta$ on $G$. Comparing these lower and upper bounds on $G$, we get that
\begin{equation}
\begin{array}{l}\eta(B)=\left\{A^4D^4HI-2A-D+B+\frac{A+4B+4D+H+I-9}{6BDA^{\frac{1}{3}}H^{\frac{1}{3}}I^{\frac{1}{3}}-4}
\right\}^2~~\\~~-A^4D^4\{A^4D^4H^2I^2-4AHI-4DHI-H^2I-I^2H+9HI\}>0.\end{array}\end{equation}
 We find that $\eta(B)$ is an increasing function of $B$ and $B\leq \lambda$, where

 \begin{equation}\lambda=\left \{ \begin{array}{l}A~~{\rm if}~~1<A\leq 1.213\\A~~{\rm if}~~1.213<A\leq 1.254, 1<H<1.2174~{\rm and}~1<I\leq1.181\\A~~{\rm if}~~ 1.213<A\leq 1.254, 1<H\leq1.14~{\rm and}~1.181<I\leq A\\{\rm min}(A,1.2247)~{\rm if}~1.213<A\leq 1.254, H>1.14~{\rm and}~1.181<I\leq1.199\\1.207~~{\rm if}~1.213<A\leq 1.254, 1.14<H<1.2174~{\rm and}~1.199<I\leq A\\ (\frac{4}{HI})^{\frac{1}{5}}~{\rm if}~1.254<A<\sqrt{2}. \end{array}\right.{} \end{equation}
So $\eta(B)\leq\eta(\lambda)$. Now consider $\eta(\lambda)$ as a function of $H$, say $\zeta(H)$. One finds that $\zeta^{''}(H)>0$ and $\mu_1 \leq H \leq \mu_2$, where
\begin{equation}(\mu_1,\mu_2)=\left \{ \begin{array}{ll}(1,A)&{\rm if}~~1<A\leq 1.213\\(1,1.2174)&{\rm if}~~1.213<A\leq 1.254~{\rm and}~1<H< 1.2174\\(1,1.14)&{\rm if}~~ 1.213<A\leq 1.254~{\rm and}~1<H\leq1.14\\(1.14,1.2174)&{\rm if}~~ 1.213<A\leq 1.254~{\rm and}~1.14<H<1.2174\\(1,1.2174)&{\rm if}~~ 1.254<A<\sqrt{2}.\end{array}\right.{} \end{equation}
 Therefore $\zeta(H)\leq{\rm max}\{\zeta(\mu_1),\zeta(\mu_2)\}$. Let $\zeta(\mu_i)=\psi_{i}(I)$, for $i=1,2$. We find that $\psi_{i}^{''}(I)>0$ for all $A$ and $D$, $1<A<\sqrt{2}$ and $1<D<1.1487$. Therefore $\psi_{i}(I)\leq{\rm max}\{\psi_{i}(\nu_1),\psi_{i}(\nu_2)\}$, where
\begin{equation}(\nu_1,\nu_2)=\left \{ \begin{array}{ll}(1,A)&{\rm if}~~1<A\leq 1.213\\(1,1.181)&{\rm if}~~1.213<A\leq 1.254~{\rm and}~1<I\leq1.181\\(1.181,A) &{\rm if}~~ 1.213<A\leq 1.254~{\rm and}~1.181<I\leq A\\(1.181,1.199)&{\rm if}~~ 1.213<A\leq 1.254~{\rm and}~1.181<I\leq1.199\\(1.199,A)&{\rm if}~~ 1.213<A\leq 1.254~{\rm and}~1.199<I\leq A\\(1,1.31951)&{\rm if}~~ 1.254<A< \sqrt{2}. \end{array}\right.{} \end{equation}
Let $\psi_{i}(\nu_j)=\theta_{ij}(A,D)$. So
\begin{equation}
\begin{array}{l}\theta_{ij}(A,D)=\left\{A^4D^4\mu_i\nu_j-2A-D+\lambda+\frac{A+4\lambda+4D+\mu_i+\nu_j-9}{6\lambda DA^{\frac{1}{3}}\mu_i^{\frac{1}{3}}\nu_j^{\frac{1}{3}}-4}
\right\}^2~~\\~~-A^4D^4\{A^4D^4\mu_i^2\nu_j^2-4A\mu_i\nu_j-4D\mu_i\nu_j-\mu_i^2\nu_j-\nu_j^2\mu_i+9\mu_i\nu_j\}.\end{array}
\end{equation}
Also using $D<1.1487$ and $D<(\frac{2}{AI})^{1/4}$, from Claim (i), we can take
\begin{equation}D<\left \{ \begin{array}{lll}1.1487&{\rm if}&1<A\leq 1.213\\1.1332&{\rm if}&A>1.213\\1.088&{\rm if}&A>1.213~{\rm and}~I>1.181\\1.1238&{\rm if}&A>1.254.\end{array}\right.{} \end{equation}
For~ $\lambda$, $\mu_i$, $\nu_j$ defined in (5.32), (5.33) and (5.34) the functions $\theta_{ij}(A,D)$ in two variables can be seen to be negative for $A$ and $D$ lying in the corresponding ranges given in (5.36). This contradicts (5.31).\vspace{3mm}\\
{\noindent \bf Proposition 28.} Case (19) i.e. $A > 1,~ B > 1, ~C > 1, ~D \leq 1,~ E > 1, ~F > 1,~ G \leq 1, ~H \leq 1, ~I \leq 1$ does not arise.\vspace{3mm}\\
{\noindent \bf Proof.} Here $a\leq 1$,~$b\leq \frac{1}{2}$,~$c\leq\frac{1}{3}$,~$e\leq \frac{1}{2}$,~$f\leq \frac{1}{3}$. Using the weak inequalities
(1,1,2,1,2,2),~(1,1,2,2,2,1),~(1,2,1,2,1,2),~(2,2,2,1,2)~and~(2,2,1,2,2), we have
\begin{equation}a+b-2d+e-2g-2i>0,\vspace{-3mm}\end{equation}
\begin{equation}a+b-2d+2f-2h-i>0,\vspace{-2mm}\end{equation}
\begin{equation}a+2c-d+2f-g-2i>0,\vspace{-2mm}\end{equation}
\begin{equation}2b-2d+2f-g-2i>0,\vspace{-2mm}\end{equation}
\begin{equation}2b-2d+e-2g-2i>0.\end{equation}
{\noindent \bf Claim(i) $F^{4}ABCDE < 2$~ {\rm and}~$GHI>\frac{1}{2}$}

Proof is same as that Claim (i) of Case (38)(Proposition 19).\vspace{2mm}

{\noindent \bf Claim(ii) $A^{4}EFGHI < 2$,~ $A < \sqrt{2}$~ {\rm and}~ $F<1.189208$}

Proof is similar to that of Claim (ii) of Case (42)(Proposition 15).\vspace{2mm}\\
{\noindent \bf Claim(iii) $d+i>2f$~ {\rm and}~ $h<\frac{a+b}{2}$}

If $d+i\leq2f$, we get a contradiction using (5.37) and Lemma 7(v) with $\gamma=d+i$, $\delta=a+b+c+e$ and $x_{1}=f<0.189208$.
\vspace{2mm}

{\noindent \bf Claim(iv) $B>C$}

Suppose $B\leq C$. The inequality (2,2,1,3,1) holds, i.e. $4A-\frac{2A^{2}}{B}+4C-\frac{2C^{2}}{D}+E+4F-\frac{F^{3}}{GH}+I>9$.
Applying AM-GM inequality to $\frac{-A^2}{B},\frac{-C^2}{D},\frac{-F^3}{GH}$ and using $B\leq A$,~ $D\leq 1$ we get \begin{equation}3A+4C-C^2+E+4F+I-3(A^{3}C^{3}F^{4}EI)^{\frac{1}{3}}>9.\end{equation}
Using (5.40) we have $i<b+f\leq a+f$. So left side of (5.42) is a decreasing function of $C$ as $A^3F^4I>(1+a)^3(1+f)^4(1-a-f)\geq 1$. Therefore we replace $C$ by $B$ and get that $3A+4B-B^2+E+4F+I-3(A^{3}B^{3}F^{4}EI)^{\frac{1}{3}}>9.$ Now
left side is a decreasing function of $E$ as $A\geq E$ and $B^3F^4I>B^3F^4(1-b-f)\geq1$, so replace $E$ by 1. Therefore we have
$3A+4B-B^2+4F+I-3(A^{3}B^{3}F^{4}I)^{\frac{1}{3}}>8.$
Now left side is a decreasing function of $I$ for $I\leq 1$. So we can replace $I$ by $1-b-f$ and get
$\phi(f)=4+3a+3b+3f-(1+b)^{2}-3(1+a)(1+b)(1+f)^{\frac{4}{3}}(1-b-f)^{\frac{1}{3}}>0.$
As $\phi''(f)>0$ and $0<f\leq{\rm{min}}(a,0.189208)$, we have $\phi(f)\leq{\rm{max}}\{\phi(0),\phi(\rm{min}(a,0.189208))\}$, which is at most zero for $0<a<\sqrt{2}-1$
and $0<b\leq a$, a contradiction. Hence we must have $B>C$.\vspace{2mm}

{\noindent \bf Claim(v) $B< 1.31951$}

Proof is same as that of Claim (iii) of Case (2)(Proposition 24).
\vspace{3mm}

{\noindent \bf Claim(vi) $B< 1.068$}

Suppose $B\geq1.068$. As $B>C$ and $D\leq 1$, we have $B^{2}>CD$, so (1,3,1,3,1) holds. Applying AM-GM inequality to (1,3,1,3,1) we have $A+4B+E+4F+I-2B^{2}F^{2}A^{\frac{1}{2}}E^{\frac{1}{2}}I^{\frac{1}{2}}>9$.
Left side is a decreasing function of $I$, so replacing $I$ by $1-b-f$, we get $\phi(e)=2+a+3b+3f+e-2(1+b)^{2}(1+f)^{2}(1+a)^{\frac{1}{2}}(1+e)^{\frac{1}{2}}(1-b-f)^{\frac{1}{2}}>0$.
We find that $\phi'(e)<0$. Therefore $\phi(e)<\phi(0)=2+a+3b+3f-2(1+b)^{2}(1+f)^{2}(1+a)^{\frac{1}{2}}(1-b-f)^{\frac{1}{2}}=\psi(f)$.
As $\psi''(f)>0$ and $0<f<0.189208$, $\psi(f)<\rm{max}\{\psi(0),\psi(0.189208)\}$, which is at most zero for
$0<a<\sqrt{2}-1$ and $0.068 \leq b \leq \rm{min}(a,0.31951)$, a contradiction. Hence we must have $b<0.068$.\vspace{2mm}

 \noindent{\bf Final Contradiction}:

The inequality (1,2,1,1,1,1,1,1) gives $A+4B+D+E+F+G+H+I-2B^3ADEFGHI>9$.
Coefficients of $F$ and $H$ on the left side are negative for $B>C$, ~$1<F<1.189208$ and $H\leq 1$, so replacing $F$ by 1 and $H$ by $1-\frac{a+b}{2}$ (using Claim (iii)), we get,\\
$2+\frac{a}{2}+\frac{7b}{2}-d-g-i+e-2(1+b)^{3}(1+a)(1+e)(1-(d+g+i))(1-\frac{a+b}{2})>0. $\\
Now the left side is an increasing function of $(d+g+i)$, as $a+b<1$. Also $d+g+i<b+\frac{e}{2}$ using (5.41). So replacing $(d+g+i)$ by $b+\frac{e}{2}$, we have,
 $\phi(e)=2+\frac{a}{2}+\frac{5b}{2}+\frac{e}{2}-2(1+b)^{3}(1+a)(1+e)(1-b-\frac{e}{2})(1-\frac{a+b}{2})>0$. As $\phi''(e)>0$ and $0<e\leq a$, we have $\phi(e)\leq \rm{max}\{\phi(0),\phi(a)\}$, which is at most zero for $0<a<\sqrt{2}-1$~ and~ $0<b<0.068$. Hence we get a contradiction.\vspace{2mm}\\
{\noindent \bf Remark 3.} Suppose the function $\phi(x_1,x_2,\cdots,x_r)$ of Remark 2 is not symmetric in the variables $x_1,x_2,\cdots,x_r$ and the variables $x_i$ satisfy $0\leq x_i\leq\eta_i(x_{i+1},\cdots,x_r)$ for $1\leq i\leq r-1$, where $\eta_i$ are linear functions in variables $x_{i+1},\cdots,x_r$. As in Remark 2 we need to maximize or minimize the function $\phi$. Here the order of the variables to be considered in succession is very vital, so we call it $\phi_{ord}(x_1,x_2,\cdots,x_r)$. We find that second derivative of $\phi_{ord}$ w.r.t the variable $x_1$ is positive in the given ranges of $x_1,x_2,\cdots,x_r$. So $\phi_{ord}(x_1,x_2,\cdots,x_r)\leq {\rm max}\{\phi_{ord}(0,x_2,\cdots,x_r),\phi_{ord}(\eta_1(x_{2},\cdots,x_r),x_2,\cdots,x_r)\}$. Further we find second derivatives of the functions $\phi_{ord}(0,x_2,\cdots,x_r)$ and $\phi_{ord}(\eta_1(x_{2},\cdots,x_r),x_2,\cdots,x_r)$ w.r.t. $x_2$ is positive within the given ranges of $x_2,\cdots,x_r$. Therefore $\phi_{ord}(x_1,x_2,\cdots,x_r)\leq {\rm max}\{\phi_{ord}(0,0,x_3,\cdots,x_r),\\ \phi_{ord}(0,\eta_2(x_3,\cdots,x_r),x_3,\cdots,x_r),\phi_{ord}(\eta_1(0,x_3,\cdots,x_r),0,x_3,\cdots,x_r),\\
\phi_{ord}(\eta_1(\eta_2(x_3,\cdots,x_r),x_3,\cdots,x_r),\eta_2(x_3,\cdots,x_r),x_3,\cdots,x_r)\}$.
Continuing like this we get that $\phi_{ord}(x_1,x_2,\cdots,x_r)$ is bounded by a number of functions in two variables $x_{r-1}$, $x_r$ (sometimes in just one variable), all of which can be verified to be at most zero within specific ranges of $x_{r-1}$ and $x_r$ using the software Mathematica. Hence we find that ${\rm max}~\phi_{ord}(x_1,x_2,\cdots,x_r)\leq 0$. Similarly we get that ${\rm min}~\phi_{ord}(x_1,x_2,\cdots,x_r)\geq 2$, if the second derivative of the relative functions w.r.t corresponding variable turn out to be negative.\vspace{2mm}\\
{\noindent \bf  Proposition 29.} Case (23) i.e. $A > 1,~ B > 1, ~C > 1, ~D \leq 1,~ E > 1, ~F \leq 1,~ G \leq 1, ~H \leq 1, ~I > 1$ does not arise.\vspace{3mm}\\
{\noindent \bf Proof.} Here $a\leq 1$, ~$b\leq \frac{1}{2}$, $c\leq \frac{1}{3}$~ and~  $e\leq \frac{1}{3}$.
Using the weak inequalities (2,2,2,2,1)~ and (2,2,1,2,1,1) we have
\begin{equation}2b-2d-2f-2h+i>0,\vspace{-3mm}\end{equation}
\begin{equation}2b-2d+e-2g-h+i>0.\end{equation}
{\noindent \bf Claim(i)~~$e<0.183$ ~{\rm or}~ $i<0.7$}

Suppose $e\geq 0.183$ and $i \geq 0.7$, then $E^{4}ABCDI > E^{4}I^{2}D \geq (1.183)^{4}(1.7)^{2}\times(\frac{3}{4})>2$. So (1,1,1,1,4,1) holds, i.e.
$\eta(B,C,D,I)=A+B+C+D+4E-\frac{1}{2}E^{5}ABCDI+I>9$. Now $\eta(B,C,D,I)$ is a decreasing function in each of the variables $B$, $C$, $D$ and $I$ for $B>1$, $C>1$, $D>\frac{3}{4}$, $E\geq1.183$ and $A\geq I\geq 1.7$. So $\eta(B,C,D,I)<\eta(1,1,\frac{3}{4},1.7)<9$ for $1.183\leq E\leq\frac{4}{3}$ and  $1<A\leq 2$. Hence we must have $e<0.183$ ~or~ $i<0.7$.\vspace{2mm}

{\noindent \bf Claim(ii)~~$E^4ABCDI<2$ ~{\rm and}~ $FGH>\frac{E^3}{2}$}

Suppose $E^4ABCDI\geq2$, then (4*,4,1) holds, i.e. $\phi(x)=4x^{\frac{1}{4}}+4E-\frac{1}{2}E^{5}xI+I-9>0$, where $x=ABCD$. Working as in Claim(i) of
Case (42)(Proposition 15) we get a contradiction when either $1<E<1.183$, $1<I\leq 2$ or  $1<E\leq \frac{4}{3}$, $1<I<1.7$, which gives a
contradiction. Hence $E^4ABCDI<2$, i.e. $\frac{E^{3}}{FGH}<2$, which gives $FGH>\frac{E^3}{2}$.\vspace{2mm}

{\noindent \bf Claim(iii) $A^{4}EFGHI < 2$,~ $A < \sqrt{2}$~ {\rm and}~ $E<1.189208$}

Proof is same as that of Claim (ii) of Case (42)(Proposition 15).\vspace{2mm}

{\noindent \bf Claim(iv) $B<1.31951$}

Proof is similar to Claim (iii) of Case (2)(Proposition 24).


 {\noindent \bf Claim(v) $d+f>0.38b+0.2i$}

 Suppose $d+f\leq 0.38b+0.2i$. Applying AM-GM inequality to (1,2,2,2,1,1) we get $\phi_{ord}(h,d+f,i,b,a)=6+a+4b-4(d+f)-h+i-6(1+b)(1-(d+f))(1+a)^{\frac{1}{3}}(1-h)^{\frac{1}{3}}(1+i)^{\frac{1}{3}}>0,$ where $0\leq h<b+\frac{i}{2}-(d+f)$(using (5.43)), $0\leq d+f\leq 0.38b+0.2i$, $0<i\leq a$ and $0<b\leq{\rm min}(a,0.31951)$.
 Following Remark 3 we find that ${\rm max}~\phi_{ord}(h,d+f,i,b,a)\leq 0$ for $0<a<\sqrt{2}-1$. Hence we must have $d+f>0.38b+0.2i$.\vspace{2mm}


{\noindent \bf  Final Contradiction}

Using (5.43) and Claim (v) we have $h<b+\frac{i}{2}-(d+f)<0.62b+0.3i$. Also using (5.44), we get $d+g<b+\frac{e+i}{2}-\frac{h}{2}.$  Applying AM-GM  inequality to (1,2,1,2,1,1,1) we get $4+a+4b+4e-(d+g)-h+i-4(1+b)^{\frac{3}{2}}(1+e)^{\frac{3}{2}}(1+a)^{\frac{1}{2}}(1-(d+g))^{\frac{1}{2}}(1-h)^{\frac{1}{2}}(1+i)^{\frac{1}{2}}>0.$
Left side is an increasing function of $d+g$, so replacing $d+g$ by
$b+\frac{e+i}{2}-\frac{h}{2}$ we get that $\phi_{ord}(h,e,i,b,a)=4+a+3b+\frac{7e}{2}+\frac{i}{2}-\frac{h}{2}-4(1+b)^{\frac{3}{2}}(1+e)^{\frac{3}{2}}(1+a)^{\frac{1}{2}}(1-b-\frac{e+i}{2}+\frac{h}{2})^{\frac{1}{2}}
(1-h)^{\frac{1}{2}}(1+i)^{\frac{1}{2}}>0,$ where $0\leq h<0.62b+0.3i$, $0<e<0.189208$, $0<i\leq a$ and $0<b\leq{\rm min}(a,0.31951)$.  Following Remark 3 we find that ${\rm max}~\phi_{ord}(h,e,i,b,a)\leq 0$ for $0<a<\sqrt{2}-1$, a contradiction.\vspace{3mm}\\
{\noindent \bf Proposition 30.} Case (24) i.e. $A > 1,~ B > 1, ~C > 1, ~D \leq 1,~ E > 1, ~F \leq 1,~ G \leq 1, ~H \leq 1, ~I \leq 1$ does not arise.\vspace{3mm}\\
{\noindent \bf Proof.} Here $a\leq 1$,~$b\leq \frac{1}{2}$,~$c\leq \frac{1}{3}$,~$e\leq \frac{1}{3}$.
 Using the weak inequalities (2,2,2,1,1,1),~(2,2,1,2,1,1),~(2,2,2,2,1) ~and~ (2,2,2,1,2), we have
\begin{equation}2b-2d-2f-g-h-i>0,\vspace{-3mm}\end{equation}
\begin{equation}2b-2d+e-2g-h-i>0,\vspace{-2mm}\end{equation}
\begin{equation}2b-2d-2f-2h-i>0,\vspace{-2mm}\end{equation}
\begin{equation}2b-2d-2f-g-2i>0.\end{equation}

{\noindent \bf Claim(i) $E<1.26$}

Proof is same as that of Case(i) of Case (43)(Proposition 14).\vspace{2mm}

{\noindent \bf Claim(ii) $B\geq 1.357$}

Suppose $B<1.357$.  The inequality (1,2,1,1,1,1,1,1) gives $A+4B+DF+E+G+H+I-2B^3ADEFGHI>8$.
Coefficient of $DF$ on left side is negative as $2B^3AEGHI>\frac{2}{C}>1$ and from (5.45) we have
$DF>1-(b-\frac{g+h+i}{2})$, so replacing $DF$ by $1-(b-\frac{g+h+i}{2})$
 and simplifying we get\\
$\begin{array}{ll}\phi_{ord}(g,e,h,i,b,a)=&2+a+3b+e-\frac{g+h+i}{2}-2(1+b)^3(1+a)(1+e)(1-g)\\&(1-h)(1-i)\left(1-b+\frac{g+h+i}{2}\right)>0.\end{array}$\\
Using (5.46), (5.47) and (5.48) respectively, we have $0\leq g<b+\frac{e}{2}-\frac{h+i}{2}$, $0\leq h<b-\frac{i}{2}$ and $0\leq i<b$. Also $0<e<0.26$. Therefore following Remark 3 we find that ${\rm max}~\phi_{ord}(g,e,h,i,b,a)\leq 0$ for $0<b\leq {\rm{min}}(a,0.357)$~and~$0<a\leq1$. Hence we get a contradiction.\vspace{2mm}

{\noindent \bf Claim(iii) $B\leq 1.39$}

Suppose $B>1.39$. We first prove that $B^4AFGHI>2$.\\
$B^4AFGHI>(1+b)^5(1-f)(1-g)(1-h)(1-i)>(1+b)^5(1-b+\frac{g+h+i}{2})(1-g)(1-h)(1-i)=\phi(g)$, using (5.45). One finds that $\phi'(g)<0$, and
$g<b+\frac{e}{2}-\frac{h+i}{2}<b+\frac{0.26}{2}-\frac{h+i}{2}$, using (5.46) and Claim (i). So $\phi(g)>\phi(b+\frac{0.26}{2}-\frac{h+i}{2})=
(1+b)^5(1-\frac{b}{2}+\frac{h+i}{4}+\frac{0.13}{2})(1-b+\frac{h+i}{2}-0.13)(1-h)(1-i)=\psi_{ord}(h,i,b),$ where $0\leq h<b-\frac{i}{2}$ and $0\leq i<b.$
 Following Remark 3 we find that ${\rm min}~\psi_{ord}(g,h,i,b,a)\geq 2$ for $0.39<b\leq 0.5$. Hence $B^4AFGHI>2$ and so the inequality (1,4,1,1,1,1) holds, i.e.
$\chi(F)=A+4B-\frac{1}{2}B^5AFGHI+F+G+H+I>9.$ Now the left side is a decreasing function of $F$ for $B>1.39$, $0\leq g<b+\frac{0.26}{2}-\frac{h+i}{2}$, $0\leq h<b-\frac{i}{2}$ and $0\leq i<b$.
Therefore we can replace $F$ by $1-b+\frac{g+h+i}{2}$ to get\\
$\varphi_{ord}(g,h,i,b,a)=a+3b-\frac{g+h+i}{2}-\frac{1}{2}(1+b)^5(1+a)(1-g)(1-h)(1-i)(1-b+\frac{g+h+i}{2})>0,$ where $0\leq g<b+\frac{0.26}{2}-\frac{h+i}{2}$, $0\leq h<b-\frac{i}{2}$ and $0\leq i<b$.\\
Following Remark 3 we find that ${\rm max}~\varphi_{ord}(g,h,i,b,a)\leq 0$ for $0.39<b\leq\frac{1}{2}$ and $0.39<a\leq 1$. Hence we must have $B\leq 1.39$.\vspace{2mm}

{\noindent \bf Claim(iv) $E\leq 1.142$}

Suppose $E>1.142$. From (5.48) we have $d+i<b$. Therefore $E^4ABCDI>(1+e)^4(1+b)^2(1-d)(1-i)>(1+e)^4(1+b)^2(1-(d+i))>(1+e)^4(1+b)^2(1-b)>2$, for $e>0.142$ and $0.357<b<0.39$. Hence $E^4ABCDI>2$. So (4*,4,1) holds. Now working as in Case (43)(Proposition 14) and using lower bound of $E$ as 1.142, we get a contradiction.\vspace{2mm}

{\noindent \bf Final Contradiction:}

We have $0.357<b<0.39$ and $e<0.142$. Here we have $g<b+\frac{0.142}{2}-\frac{h+i}{2}$ instead of $g<b+\frac{0.26}{2}-\frac{h+i}{2}$. Proceeding as in Claim (iii),
 we get $B^4AFGHI>2$, for $0.357<b<0.39$. So the inequality (1,4,1,1,1,1) holds
 and we again get a contradiction as in Claim (iii).\vspace{3mm}\\
 {\noindent \bf  Proposition 31.} Case (29) i.e. $A > 1,~ B > 1, ~C > 1, ~D \leq 1,~ E \leq 1, ~F > 1,~ G \leq 1, ~H \leq 1, ~I \leq 1$ does not arise.\vspace{3mm}\\
{\noindent \bf Proof.} Here $a\leq 1$, ~$b\leq \frac{1}{2}$, $c\leq \frac{1}{3}$~ and~  $f\leq \frac{1}{3}$.
Using the weak inequalities (2,2,1,2,1,1),~ (2,2,1,1,2,1), ~(2,2,1,2,2) and (2,1,2,2,2) we have
\begin{equation}2b-2d-e-2g-h-i>0,\vspace{-3mm}\end{equation}
\begin{equation}2b-2d-e+f-2h-i>0,\vspace{-2mm}\end{equation}
\begin{equation}2b-2d-e-2g-2i>0,\vspace{-2mm}\end{equation}
\begin{equation}2b+c-2e-2g-2i>0.\end{equation}
{\noindent \bf Claim(i) $F^4ABCDE<2$ ~{\rm and}~ $GHI>\frac{1}{2}$}

Proof is same as that of Claim (i) of Case (38)(Proposition 19),

{\noindent \bf Claim(ii) $B\leq 1.29$}

Suppose $B>1.29$, we prove that $B^4AFGHI \geq 2$. Using (5.49), (5.50) and (5.51) respectively we have $g<b-\frac{h+i}{2}$, ~$h<b+\frac{f}{2}-\frac{i}{2}$
and $i<b$. So $B^4AFGHI > (1+b)^5(1+f)(1-b+\frac{h+i}{2})(1-h)(1-i)=\eta(h)$. Now $\eta'(h)<0$, so $\eta(h)>\eta(b+\frac{f}{2}-\frac{i}{2})$.
Therefore $B^4AFGHI > (1+b)^5(1+f)(1-\frac{b}{2}+\frac{f}{4}+\frac{i}{4})(1-b-\frac{f}{2}+\frac{i}{2})(1-i)=\phi(i)$. We find that $\phi''(i)<0$, so
$\phi(i)\geq{\rm min}\{\phi(0),\phi(b)\}$, which can be seen to be greater than 2 for $0.29<b\leq 0.5$ and $0<f\leq \frac{1}{3}$. So we have $B^4AFGHI \geq 2$ and
therefore (1,4,1,1,1,1) holds. Now we proceed as in Claim (iii) of Case (2)(Proposition 24) and get a contradiction.
\vspace{2mm}

{\noindent \bf Claim(iii) $B>C$}

Suppose $B\leq C$. The inequality (2,2,1,3,1) holds. Now working as in Claim(iv) of Case (19)(Proposition 28), we have (5.42), i.e $4+3a+4c-(1+c)^2+4f-(e+i)-3(1+a)(1+c)(1+f)^{\frac{4}{3}}(1-(e+i))^{\frac{1}{3}}>0.$ A simple calculation shows that left side is a decreasing function of $f$ first, so we replace $f$ by 0; an increasing function of $e+i$ so we replace $e+i$ by $b+\frac{c}{2}$ (using (5.52)); then a decreasing of $c$, so we replace $c$ by $b$ to get $4+3a+\frac{5b}{2}-(1+b)^2-3(1+a)(1+b)(1-\frac{3b}{2})^{\frac{1}{3}}>0,$ but this is not true for ~$0<a\leq 1$ and $0<b\leq 0.29$. Hence $B> C$.\vspace{2mm}

{\noindent \bf Claim(iv) $a\leq3b$}

If $b<\frac{a}{3}$, then we get a contradiction using (5.49) and Lemma 7(iii) with $x_1=a\leq 1$ and $\gamma=d+e+g+h+i$.\vspace{3mm}

{\noindent \bf Claim(v) $e\geq \frac{b}{2}$}

Suppose $e<\frac{b}{2}$. Applying AM-GM to the inequality (1,3,1,3,1) we get
 $A+4B+E+4F+I-2\sqrt{B^4F^4AEI}>9$. Left side is a decreasing function of $F$ as $F>1$ and $B^4AEI>B^5(1-e)(1-i)>B^5(1-\frac{b}{2})(1-b)\geq1$ for $0<b\leq 0.29$, so replacing $F$ by 1 we get $A+4B+E+I-2\sqrt{B^4AEI}>5$. Now left side is a decreasing function of $I$ as $B^4AE>1$, so replacing $I$ by $1-b+\frac{e}{2}$, we get $\theta(e)=2+a+3b-\frac{e}{2}-2(1+b)^{2}\sqrt{(1+a)(1-e)(1-b+\frac{e}{2})}>0$. Now $\theta''(e)>0$ and
we have $0\leq e<\frac{b}{2}$, so $\theta(e)\leq{\rm max}\{\theta(0),\theta(\frac{b}{2})\}$.
We find that $\theta(0)$ and $\theta(\frac{b}{2})$ are at most zero for $b<a\leq 3b$ and $0<b\leq 0.29$. Hence we must have $e\geq \frac{b}{2}$.\vspace{2mm}

\noindent{\bf Final Contradiction:}

Using (5.51) and Claim(v) we have $d+g<b-\frac{e}{2}-i<\frac{3b}{4}-i$. Applying AM-GM to the inequality (1,2,2,1,2,1), we have $A+4B+4D+F+4G+I-6BDG(AFI)^{\frac{1}{3}}>9$. Left side is a
decreasing function of $F$ as $\frac{2BDG(AI)^{\frac{1}{3}}}{F^{\frac{2}{3}}}>2\times(\frac{3}{4})^{\frac{1}{3}}B(1-(d+g))(1-i)^{\frac{1}{3}}>
2\times(\frac{3}{4})^{\frac{1}{3}}B(1-\frac{3b}{4})(1-b)^{\frac{1}{3}}>1$, so we replace $F$ by 1 and get that
$6+a+4b-4(d+g)-i-6(1+b)(1-(d+g))(1+a)^{\frac{1}{3}}(1-i)^{\frac{1}{3}}>0$. Now the left side is an increasing function of $d+g$ for $i<b$,
so we replace $d+g$ by $\frac{3b}{4}-i$ to get that $\eta(i)=6+a+b+3i-6(1+b)(1-\frac{3b}{4}+i)(1+a)^{\frac{1}{3}}(1-i)^{\frac{1}{3}}>0$. Now $\eta''(i)>0$,
so $\eta(i)\leq{\rm max}\{\eta(0),\eta(b)\}.$
We find that $\eta(0)$ and $\eta(b)$ are at most zero for $0<b\leq0.29$ and $b\leq a\leq 3b$. Hence we have a contradiction.\vspace{3mm}

{\noindent \bf  Proposition 32.} Case (30) i.e. $ A>1,~ B>1, ~C>1, ~D\leq 1,~ E\leq 1, ~F\leq 1,~ G>1,~ H> 1,~ I> 1$ does not
arise.\vspace{2mm}\\ {\noindent \bf Proof.} Here $a\leq 1, ~b\leq \frac{1}{ 2},~ c\leq \frac{1}{ 3}$. Using  the weak inequalities $(1,1,2,2,1,1,1), \\ (2,2,2,1,1,1),(1,2,1,2,1,1,1) $ and
$(1,2,2,1,1,1,1)$ we get
\begin{equation}a+b-2d-2f+g+h+i >0,\vspace{-3mm}\end{equation}
\begin{equation}2b-2d-2f+g+h+i>0,\vspace{-2mm}\end{equation}
\begin{equation}a+2c-d-2f+g+h+i>0,\vspace{-2mm}\end{equation}
\begin{equation}a+2c-2e-f+g+h+i >0.\end{equation}
\noindent{\bf Claim(i)} $a<0.588$

Suppose $a\geq 0.588.$ Then $A^6GHI>A^6>16.$
 Therefore the inequality $(6,1,1,1)$ holds. That is $4A-\frac{1}{16}A^7GHI+G+H+I>9$. As the left hand side is a decreasing
 function of $G$, $H$ and of $I$, we can replace each of $G$, $H$ and $I$ by $1$ to get $4A-\frac{1}{16}A^7>6$, which is not true for $a\geq 0.588.$ \vspace{2mm}

\noindent{\bf Claim(ii)} $C^4GHIAB<2$ and  $c<0.149$

Suppose  $C^4GHIAB\geq 2$. Therefore $(1,1,4,1,1,1)$ holds i.e.  $A+B+4C-\frac{1}{2}C^5GHIAB+G+H+I>9,$ which is not true, by
Lemma 10(ii) with $X_{2}=C,~X_{3}=B,~X_{4}=G,~X_{5}=H,~X_{6}=I$. \\Now $C^4GHIAB<2$ implies $C^5<2$ and so $c<0.149$.\vspace{2mm}

\noindent{\bf Claim(iii) $b<0.202$}

Suppose $b\geq 0.202$. We first show that $B^4FGHIA>2$. Let $g+h+i=k.$ We consider the following cases:\\
{\bf Case(i)}
$k < a.$ Here using $f<b+\frac{k}{2}$ from (5.54) we have $B^4FGHIA>(1+b)^4(1+a)(1-b-\frac{k}{2})(1+k)=\phi(k)\geq{\rm
min}\{\phi(0),\phi(a)\}>2$ for $ b\leq a<0.588$ and $0.202\leq b\leq0.5$.\\ {\bf Case(ii)} $ a\leq k<2a, ~ a\geq 0.4.$
Here using that $F\geq0.46873B$, we have $B^4FGHIA\geq(0.46873)B^5A(1+a)\geq0.46873(1.202)^5(1.4)^2>2.$\\{\bf Case(iii)} $ a\leq k<2a, ~ a<
0.4$. $B^4FGHIA>(1+b)^4(1+a)(1-b-\frac{g+h+i}{2})(1+g)(1+h)(1+i)=\phi(g,h,i)\geq{\rm min}\{\phi(a,0,0),\phi(a,a,0)\}>2$ for $ 0.202\leq b\leq a<0.4$ .\\ {\bf Case(iv)} $ 2a\leq k<3a, ~ a\geq 0.3.$
Here using that $F\geq0.46873B$, we have $B^4FGHIA>(0.46873)B^5A(1+2a)>0.46873(1.202)^5(1.3)(1+2\times 0.3)>2$.\\{\bf Case(v)} $ 2a\leq k<3a, ~ a<
0.3$. $B^4FGHIA>(1+b)^4(1+a)(1-b-\frac{g+h+i}{2})(1+g)(1+h)(1+i)=\phi(g,h,i)\geq{\rm min}\{\phi(a,a,0),\phi(a,a,a)\}>2$ for $0.202\leq b\leq a<0.3$.\\
Therefore the inequality $(1,4,1,1,1,1)$ holds i.e.  $A+4B-\frac{1}{2}B^5FGHIA+F+G+H+I>9$. The
left side of this inequality is a decreasing function of $F$, for $B>1.13$. We have $2F\geq C$ and $F\geq \frac{2}{3}D$, i.e.
$1-c\geq 2f$ and $\frac{d}{2}\geq \frac{3f-1}{4}$. Also from (5.55), $c+\frac{a+g+h+i}{2}>f+\frac{d}{2}$. Adding all
these we get $f<\frac{1}{3}+\frac{2(a+g+h+i)}{15}$. So we replace $F$ by $1-\frac{1}{3}-\frac{2(a+g+h+i)}{15}$ to get\\
 $\begin{array}{ll}
\chi(g,h,i)=&4b-\frac{1}{3} +\frac{13(a+g+h+i)}{15}-\frac{1}{2}(1+b)^5(1+a)(1+g)(1+h)(1+i)\\&\left(\frac{2}{3} -\frac{2(a+g+h+i)}{15}\right)
>0.\end{array}$\\ Following Remark 2 we find that ${\rm max}~\chi(g, h, i)\leq0$ for $0.202\leq b\leq\frac{1}{2}$ and $0< a < 0.588$. This gives a contradiction.\vspace{3mm}

\noindent{\bf Claim(iv) $a>0.276$}

Suppose  $a\leq 0.276$.\vspace{2mm}

\noindent{\bf Case(i)} $f<0.45b+0.395(g+h+i)$, i.e. $F>1-0.45b-0.395(g+h+i).$\\ Using the inequality $(2,3,1,1,1,1)$ we have $2B+4C-C^4FGHIAB+F+G+H+I>9$. As $ABFGHI>AB(1-0.45b-0.395(g+h+i))GHI>1$, the left side is a decreasing function of $C$ and $F$, we can replace $C$ by 1 and
$F$ by $1-0.45b-0.395(g+h+i)$ to get $$\begin{array}{l}\phi(g,h,i)=1+1.55b+
0.605(g+h+i)-(1+a)(1+b)\times\\~~~~~~~~~~\{1-0.45b-0.395(g+h+i)\}(1+g)(1+h)(1+i)>0.\end{array}$$ Following Remark 2 we find that ${\rm max}~\phi(g, h, i)\leq0$ for $0<b\leq{\rm min}(a,0.202)$ and $0<a\leq 0.276$. This gives a contradiction.\vspace{2mm}

\noindent{\bf Case(ii)} $f\geq 0.45b+0.395(g+h+i)$.\\ From (5.54) we get $d<b+\frac{g+h+i}{2}-f<b+\frac{g+h+i}{2}-(0.45b+0.395(g+h+i))$, i.e. $d<0.55b+0.105(g+h+i)$. Applying AM-GM to the inequality
$(1,2,2,1,1,1,1)$ we have $A+4B+4D+F+G+H+I-4\sqrt{B^3D^3AFGHI}>9.$ The left hand
side is a decreasing function of $F$ so replacing $F$ by $1-b-\frac{g+h+i}{2}+d$ and simplifying we get\\
$
\begin{array}{ll}
\eta_{ord}(d,g,h,i,b,a)&=4+a+3b+\frac{g+h+i}{2}-3d-4(1+b)^{\frac{3}{2}}(1-d)^{\frac{3}{2}}\\ &(1-b- \frac{g+h+i}{2}+d)^{\frac{1}{2}}(1+a)^{\frac{1}{2}}(1+g)^{\frac{1}{2}}
(1+h)^{\frac{1}{2}}(1+i)^{\frac{1}{2}}>0, \end{array}
$\\
where $0\leq d<0.55b+0.105(g+h+i)$, $0<g\leq a$, $0<h\leq a$ and $0<i\leq a$.
Following Remark 3 we find that ${\rm max}~\eta_{ord}(d,g,h,i,b,a)\leq0$ for $0<b\leq {\rm min}(a,0.202)$ and $0<a\leq 0.276$, which gives a contradiction. Hence $A>1.276$.\vspace{3mm}

\noindent{\bf Claim(v)} $b<0.174$

Suppose $b\geq 1.174$. We first show that $B^4FGHIA>2$. We consider the following cases:\\ {\bf Case(i)}
$k < a.$ Here using $f<b+\frac{k}{2}$ we have $B^4FGHIA>(1+b)^4(1+a)(1-b-\frac{k}{2})(1+k)=\phi(k)\geq{\rm
min}\{\phi(0),\phi(a)\}>2$ for $ 0.276< a<0.588$ and $0.174\leq b<0.202$.\\ {\bf Case(ii)} $ a\leq k<2a, ~ a\geq 0.53.$
Here using that $F\geq0.46873B$, we have $B^4FGHIA\geq0.46873(1.174)^5(1.53)^2>2.$\\{\bf Case(iii)} $ a\leq k<2a, ~ a<
0.53$. $B^4FGHIA>(1+b)^4(1+a)(1-b-\frac{g+h+i}{2})(1+g)(1+h)(1+i)=\phi(g,h,i)\geq{\rm min}\{\phi(a,0,0),\phi(a,a,0)\}>2$ for $ 0.276\leq a<0.53$ and
$0.174\leq b<0.202$.\\{\bf Case(iv)} $ 2a\leq k<3a$.~ Here using that $F\geq0.46873B$, we have $B^4FGHIA>(0.46873)B^5A(1+2a)>(0.46873)(1.174)^5(1.276)(1+2\times0.276)>2$.\\   Now working as in Claim(iii) and using lower bound of $b$ as 0.174, the inequality $(1,4,1,1,1,1)$ gives a contradiction.\vspace{3mm}

\noindent{\bf Claim(vi) $a>0.3$}

Suppose  $a\leq 0.3$.\vspace{2mm}

\noindent{\bf Case(i)} $f<0.5b+0.425(g+h+i).$\\We use the inequality $(2,3,1,1,1,1)$ and work as in Case(i) of Claim(iv). We still have $ABFGHI>1$. So we replace $C$ by 1 and $F$ by $1-0.5b-0.425(g+h+i)$ to get a contradiction for $0<b<0.174$ and $0.276<a<0.3$.\vspace{2mm}

\noindent{\bf Case(ii)} $f\geq 0.5b+0.425(g+h+i)$.\\ Using (5.54) we get $d<0.5b+0.075(g+h+i)$. We use the inequality $(1,2,2,1,1,1,1)$ and proceed as in Case(ii) of Claim(iv) to get $\eta(d)\leq{\rm max}\{\eta(0),\eta(0.5b+0.075(g+h+i))\}$, which is again at most zero for $0.276<a\leq 0.3$ and $0<b<0.174$, giving thereby a contradiction. Hence $A>1.3$.\vspace{3mm}

\noindent{\bf Claim(vii)} $f>2c$, ~$c<0.076$ ~and~ $g+h+i<2a$

Suppose $f\leq 2c$. From (5.53) we have $f<\frac{1}{2}(a+b+g+h+i)$. We use $(1,1,3,1,1,1,1)$ and apply Lemma 7(v) with $X_1=C$, $\gamma=f$ and $\delta=a+b+g+h+i$. So we have $\gamma<\frac{\delta}{2}$, ~$x_1<0.149$ ~and~ $\gamma\leq 2x_1$. This gives a contradiction. Hence we must have $f>2c.$\\
 $~~~~$
 If  $c\geq 0.076$ then $C^4GHIAB>C^4A(1+b+g+h+i)>C^4(1.3)(1+2c)>2$, as $b+g+h+i\geq b+\frac{g}{2}+\frac{h}{2}+\frac{i}{2}>f>2c$
from inequality (5.54). This contradicts Claim (ii).\\
$~~~~$
Also if $g+h+i \geq 2a$, then $C^4ABGHI>(1+c)^4(1+a)(1+g+h+i)\geq (1+a)(1+2a)>2$, for $a>0.3$. Again a contradiction to Claim (ii).\vspace{3mm}

\noindent{\bf Claim(viii) $A^4EFGHI<2$ ~{\rm and} ~$A<1.362$}

Suppose $A^4EFGHI\geq 2$. Applying AM-GM to the inequality $(4,2,1,1,1)$ we get $4A+4E-2A^{\frac{5}{2}}E^{\frac{3}{2}}G^{\frac{1}{2}}H^{\frac{1}{2}}I^{\frac{1}{2}}+2+GHI>9.$ The left
 side is a decreasing function of $E$ for $E> \frac{2}{3}$ and $A>1.3$.  So replacing $E$  by
$\frac{2}{3}$  we get $\phi(x)=4A+\frac{8}{3}-2A^{\frac{5}{2}}(\frac{2}{3})^{\frac{3}{2}}x^{\frac{1}{2}}+2+x>9,$ where $1<x=GHI\leq A^3$. As $\phi''(x)>0$, we have $\phi(x)\leq {\rm max}\{\phi(1),\phi(A^3)\}$. One can easily check that $\phi(1)$ and $\phi(A^3)$ are less than $9$  for $1.3<A<1.588$; a contradiction. So $A^4EFGHI<2$.\vspace{1mm}\\
Now if $A\geq 1.362$, then $A^4EFGHI>\frac{A^3}{BC}>\frac{(1.362)^3}{(1.174)(1.076)}>2$, a contradiction.\vspace{3mm}

\noindent{\bf Claim(ix)} $f \geq 1.6c+0.57(g+h+i)$

Suppose  $f< 1.6c+0.57(g+h+i)$\\ Using the inequality $(2,3,1,1,1,1)$ we have $2B+4C-C^4FGHIAB+F+G+H+I>9$. As
the coefficient of $B$ namely  $2-C^4FGHIA> 0$ by Claim (ii), we can replace $B$ by 1.174 and
then $F$ by $1-1.6c-0.57(g+h+i)$ to get $$\begin{array}{l}\phi_{ord}(g,h,i,c,a)=1+2(0.174)+2.4c+
0.43(g+h+i)-(1+c)^4(1+a)\times\\~~~~~~~~~~\{1-1.6c-0.57(g+h+i)\}(1+g)(1+h)(1+i)(1+0.174)>0.\end{array}$$
Following Remark 3 one finds that ${\rm max}~\phi_{ord}(g,h,i,c,a)\leq0$ for $0<g<2a-h-i$, $0<h\leq a$, $0<i\leq a$, $0<c<0.076$ and $0.3<a<0.362$. This gives a contradiction. So we must have $f\geq 1.6c+0.57(g+h+i)$.\vspace{3mm}

\noindent{\bf Claim(x) $g+h+i>a$}

Suppose if possible $k=g+h+i\leq a$. From (5.56) and Claim (ix) we get $e<\frac{a+2c+k-f}{2}<0.5a+0.2c+0.215k$. Also from (5.54) we have
$d+f<b+\frac{k}{2}<0.174+\frac{k}{2},$ i.e. $DF>1-(0.174+\frac{k}{2})$. Using  inequality $(2,1,1,1,1,1,1,1)$ we have $4A-2A^3CDEFGHI+C+DF+E+G+H+I>8$. Left side is a decreasing function of $DF$ and $E$, replacing $DF$ by  $1-(0.174+\frac{k}{2})$ and $E$ by $1-(0.5a+0.2c+0.215k)$ we get\\
$\begin{array}{l}\psi(k)=1.826+3.5a+0.8c+0.285k-2(1+a)^3(1+c)\times \\~~~~~~~~~~~~~~~~~~~~~
(0.826-\frac{k}{2})(1-0.5a-0.2c-0.215k)(1+k)>0.\end{array}$\\ Again one finds that $\psi''(k)>0$, therefore $\psi(k)\leq
{\rm max}\{\psi(0),~\psi(a)\},$ which is at most $0$ for  $0.3< a<0.362$ and $0<c<0.076$. This gives a
contradiction.\vspace{3mm}

\noindent{\bf Claim(xi) $b<0.134$}

Suppose $b\geq 0.134$. Here using $f<b+\frac{g+h+i}{2}$, we have
$B^4FGHIA>(1+b)^4(1+a)(1-b-\frac{g+h+i}{2})(1+g)(1+h)(1+i)=\phi(g,h,i)>{\rm min}\{\phi(a,0,0),~\phi(a,a,0)\}>2 $, for $0.134\leq b<0.174$ and $0.3<a<0.362$. Therefore $B^4FGHIA>2$. Now working as in Claim (iii) and using lower bound of $b$ as 0.134, the inequality $(1,4,1,1,1,1)$ gives a contradiction.\vspace{3mm}

\noindent{\bf Claim(xii)} $c<0.044$ and $a<0.3329$

If $c\geq 0.044$ then $C^4GHIAB>C^4(1+a)^2>2$ for $a>0.3,$ which contradicts Claim (ii). So we have $c<0.044$.\\ $~~~~$Further if $a\geq 0.3329$,
then $A^4EFGHI>\frac{A^3}{BC}>\frac{A^3}{(1.134)(1.044)}>2,$ a contradiction to Claim (viii).\vspace{3mm}

\noindent{\bf Final contradiction}

If $f<0.8b+0.5(g+h+i)$, we use the inequality $(2,3,1,1,1,1)$ and work as in Case (i) of Claim (iv) to get a contradiction
for  $b<0.134$ and $0.3<a<0.3329$. If  $f\geq 0.8b+0.5(g+h+i)$, we find $d<0.2b$. Again using
the inequality $(1,2,2,1,1,1,1)$ and working as in Case (ii) of Claim (iv), we get a contradiction for $0.3< a<
0.3329$ and $0<b<0.134$.\vspace{3mm}\\
{\noindent \bf  Proposition 33.} Case (31) i.e. $A > 1,~ B > 1, ~C > 1, ~D \leq 1,~ E \leq 1, ~F \leq 1,~ G > 1, ~H \leq 1, ~I > 1$ does not arise.\vspace{3mm}\\
{\noindent \bf Proof.} Here $a\leq 1$, ~$b\leq \frac{1}{2}$, $c\leq \frac{1}{3}$~ and~  $g\leq \frac{1}{3}$.
Using the weak inequalities (2,2,2,2,1),~ (1,2,1,2,2,1) and (1,2,2,1,2,1) we have
\begin{equation}2b-2d-2f-2h+i>0,\vspace{-3mm}\end{equation}
\begin{equation}a+2c-d-2f-2h+i>0,\vspace{-2mm}\end{equation}
\begin{equation}a+2c-2e-f-2h+i>0.\end{equation}

{\noindent \bf Claim(i) $C<1.217$}

Suppose $C\geq 1.217$, then $C^4ABGHI> C^{5}H\geq (1.217)^{5}(\frac{3}{4})>2$. Applying AM-GM to the inequality $(1,1,4,2,1)$ we get $A+B+4C+4G+I-2\sqrt{C^5G^3ABI}>9$.
Left side is a decreasing function of $I$ for $1<I\leq A$, we replace $I$ by 1 to get that $A+B+4C+4G-2\sqrt{C^5G^3AB}>8$. We can successively replace $B$ and $G$ by 1 with the similar argument to get that $A+4C-2\sqrt{C^5A}>3$. But this is not true for $1.217\leq C\leq \frac{4}{3}$ and $1<A\leq 2$.\vspace{3mm}

{\noindent \bf Claim(ii) $A<1.54$}

Suppose $A\geq 1.54$, then $A^4EFGHI=\frac{A^3}{BCD}>\frac{(1.54)^3}{1.5\times1.217}>2$. Also $2E>1\geq F$, so (4,2,2,1) holds, i.e.
$4A-\frac{1}{2}A^{5}EFGHI+4E-\frac{2E^2}{F}+4G-\frac{2G^2}{H}+I>9$. Using $G>H$ and applying AM-GM to $-\frac{1}{2}A^{5}EFGHI$, $-\frac{2E^2}{F}$
and $-\frac{G^2}{H}$ we get $4A+4E+3G+I-3(A^{5}E^{3}G^{3}I)^{\frac{1}{3}}>9$. Left side is a decreasing function of $E$ for $A\geq 1.2$,
so we can replace $E$ by $0.46873A$ to get $4A+4(0.46873)A+3G+I-3(0.46873)(A^{8}G^{3}I)^{\frac{1}{3}}>9$. Again left side is a decreasing function of $G$ for $A>1.33$; so we can replace $G$ by 1 to get $4A+4(0.46873)A+I-3(0.46873)(A^{8}I)^{\frac{1}{3}}>6$, which is not true for $1<I\leq A\leq 2$. So $A<1.54$.\vspace{3mm}

{\noindent \bf Claim(iii) $B<1.21$}

Suppose $B\geq 1.21$. Using (5.57) we have $f+h<b+\frac{i}{2}$. Then $B^4AFGHI>B^4AI(1-(f+h))>(1+b)^{4}(1+a)(1+i)(1-b-\frac{i}{2})=\eta(i)$, say.
We find that $\eta''(i)<0$, so for $0<i\leq a$, $\eta(i)\geq{\rm min}\{\eta(0),\eta(a)\}>$2 for $0.21\leq b\leq a<0.54$. Applying AM-GM to (1,4,1,2,1) we get $A+4B+F+4G+I-2\sqrt{B^5G^3AFI}>9$. The left side is a
decreasing function of $F$. Replacing $F$ by $0.46873B$ we get $A+4B+0.46873B+4G+I-2\sqrt{(0.46873)B^6G^3AI}>9$. Now for $B\geq 1.21$, the left side
is decreasing function of $G$ and $I$, so
replacing $G$ and $I$ by 1, we get $A+4B+0.46873B-2\sqrt{(0.46873)B^6A}>4$, which is not true for $1.21\leq B\leq A\leq 1.54$. Hence we must have
$B<1.21$.\vspace{3mm}

{\noindent \bf Claim(iv) $C<1.122$}

Suppose $C\geq 1.122$. As (2,3,1,1,1,1) holds i.e. $4A-\frac{2A^2}{B}+4C-\frac{C^3}{DE}+F+G+H+I>9$. Now using $A\geq B$ and applying AM-GM to
$-\frac{A^2}{B}-\frac{C^3}{DE}$, we get $3A+4C+F+G+H+I-2\sqrt{C^4A^3FGHI}>9$. Now left side is a decreasing function of $C$ as $A^3FGHI>ABCFGHI=\frac{1}{DE}\geq 1$, so we replace $C$ by 1.122 and get that $2+3a+4(0.122)+g+i-(f+h)-2(1.122)^{2}\sqrt{(1+a)^{3}(1+g)(1+i)(1-(f+h))}>0$.
As the left side is an increasing function of $f+h$ and $f+h<b+\frac{i}{2}<0.21+\frac{i}{2}$, we have\\
$\eta(i)=2.278+3a+g+\frac{i}{2}-2(1.122)^{2}\sqrt{(1+a)^{3}(1+g)(1+i)(0.79-\frac{i}{2})}>0$. Now $\eta''(i)>0$, so
$\eta(i)\leq{\rm max}\{\eta(0),\eta(a)\}$, which can be verified to be at most zero for $0.122\leq a\leq 0.54$ and $0<g\leq \frac{1}{3}$. Hence we must have $C<1.122$.\vspace{3mm}

{\noindent \bf Claim(v) $A>1.336$}

Suppose $A\leq 1.336$. Applying AM-GM to the inequality (2,2,1,1,1,1,1) we have
$4A+4C+E+F+G+H+I-4\sqrt{A^3C^3EFGHI}>9$. Left side is a decreasing function of $G$ as $2\sqrt{A^2C^3EFHI}>2\sqrt{\frac{1}{DG}}\geq2\sqrt{\frac{3}{4}}>1$, so we replace $G$ by 1 and use $E+H<1+EH$ to get $4A+4C+EH+F+I-4\sqrt{A^3C^3EHFI}>7$. Now the left side is a decreasing function of $EH$, so we replace $EH$
by $1-(\frac{a+i}{2}+c-\frac{f}{2})$, using (5.59) and get
$\phi_{ord}(f,i,c,a)=4+\frac{7a}{2}+\frac{i}{2}+3c-\frac{f}{2}-4\sqrt{(1+a)^3(1+c)^3(1+i)(1-f)(1-\frac{a+i}{2}-c+\frac{f}{2})}>0,$ where $0\leq f<\frac{a+i}{2}+c$ (using (5.58)) and $0< i\leq a$. Following Remark 3 we find that ${\rm max}~\phi_{ord}(f,i,c,a)\leq 0$ for $0<a\leq 0.336$ and $0<c<0.122$. This gives a contradiction.\vspace{3mm}

{\noindent \bf Claim(vi) $C \leq 1.085$}

Suppose $C>1.085$. Now proceeding as in Claim (iv) and replacing $C$ by 1.085 instead of 1.122 we get a contradiction for $0.336<a<0.54$ and $0<g\leq \frac{1}{3}$.\vspace{3mm}

{\noindent \bf Final Contradiction:}

We have $1.336<A<1.54$, $C\leq1.085$ and $B<1.21$. We prove that $A^4EFGHI>2$. Note that $e+h<\frac{a+i}{2}+c-\frac{f}{2}<\frac{a+i}{2}+0.085-\frac{f}{2}$. So
$A^4EFGHI>(1+a)^4(1+i)(1-(e+h))(1-f)>(1+a)^4(1+i)(1-\frac{a+i}{2}-0.085+\frac{f}{2})(1-f)=\eta(f)$. Now $\eta''(f)<0$ and
$0\leq f<b+\frac{i}{2}<0.21+\frac{i}{2}$. Therefore $\eta(f)\geq{\rm min}\{\eta(0),\eta(0.21+\frac{i}{2})\}$, which can be verified to be greater than 2 for $0.336<a<0.54$ and $0<i\leq a$. Hence $A^4EFGHI>2$. So (4,2,2,1) holds. Working as in Claim (ii) we get a contradiction.\vspace{2mm}\\
{\noindent \bf  Proposition 34.} Case (34) i.e. $A > 1,~ B > 1, ~C > 1, ~D \leq 1,~ E \leq 1, ~F \leq 1,~ G \leq 1, ~H \leq 1, ~I \leq 1$ does not arise.\vspace{3mm}\\
{\noindent \bf Proof.} Here $a\leq 1$,~$b\leq \frac{1}{2}$ ~and~ $c\leq \frac{1}{3}$. Using the weak inequalities
(2,2,2,1,1,1),~ (2,2,2,1,2),~ (2,2,1,1,1,1,1),~ (2,2,1,2,1,1), ~(2,2,1,2,2),~ (2,2,2,2,1)~\\ and~  (1,2,2,1,1,1,1) we have \vspace{-2mm}
\begin{equation}2b-2d-2f-g-h-i>0,\vspace{-3mm}\end{equation}
\begin{equation}2b-2d-2f-g-2i>0,\vspace{-2mm}\end{equation}
\begin{equation}2b-2d-e-f-g-h-i>0,\vspace{-2mm}\end{equation}
\begin{equation}2b-2d-e-2g-h-i>0,\vspace{-2mm}\end{equation}
\begin{equation}2b-2d-e-2g-2i>0,\vspace{-2mm}\end{equation}
\begin{equation}2b-2d-2f-2h-i>0,\vspace{-2mm}\end{equation}
\begin{equation}a+2c-2e-f-g-h-i>0.\end{equation}

\noindent{\bf Claim(i) $b<0.313$}

 Suppose $b \geq 0.313.$ From (5.60), ~(5.63), ~(5.65) and ~(5.61) we have  $f<b-\frac{g+h+i}{2}$, ~$g<b-\frac{h+i}{2}$, ~$h<b-\frac{i}{2}$~ and ~$i<b$, respectively.
 Therefore $B^{4}FGHIA \geq (1+b)^{5}(1-b+\frac{g+h+i}{2})(1-g)(1-h)(1-i)=\phi(g)$. As  $\phi'(g)<0$, we have $\phi(g)> \phi(b-\frac{h+i}{2})=
(1+b)^{5}(1-\frac{b}{2}+\frac{h+i}{4})(1-b+\frac{h+i}{2})(1-h)(1-i)=\psi(h).$ Now $\psi''(h)<0$, so  $\psi(h)\geq  {\rm min}\{\psi(0),\psi(b-\frac{i}{2})\}$, which is greater than 2 for $0.313\leq b\leq 0.5$ and $0\leq i<b$.
  Hence we have $B^4AFGHI>2$, so $(1,4,1,1,1,1)$ holds  i.e.
 $A+4B-\frac{1}{2}B^{5}FGHIA+F+G+H+I>9$. Since $B^5FGHI>2$, the coefficient of $A$ is negative, we can replace $A$ by $B$ to get
 $5B-\frac{1}{2}B^{6}FGHI+F+G+H+I>9$. Now the coefficient of $F$  is negative  and $F> 1-b+\frac{g+h+i}{2}$, therefore we
 get\\ $\phi_{ord}(g,h,i,b)=4b-\frac{1}{2}(1+b)^{6}\left(1-b+\frac{g+h+i}{2}\right)(1-g)(1-h)(1-i)-\frac{g+h+i}{2}>0,$\\where $0\leq g<b-\frac{h+i}{2}$, $0\leq h<b-\frac{i}{2}$ and $0\leq i<b$.
 Following Remark 3 we find that ${\rm max}~\phi_{ord}(g,h,i,b)\leq 0$ for $0<b\leq 0.5$. Hence we get a contradiction.\vspace{3mm}

\noindent{\bf Claim(ii) $B>C$}

Suppose $B\leq C$. Using (5.64) and (5.65) we have $g+i<b$ and $f+h<b$, respectively, i.e. $GI>1-b$ and  $FH>1-b$ respectively.
The inequality (1,1,3,1,1,1,1) holds, i.e. $A+B+4C-C^{4}ABFGHI+F+G+H+I>9$. Left side is a decreasing function of $C$ as $C^3ABFGHI>B^5(1-b)^2\geq1$, for $0<b<0.313$,
so we replace $C$ by $B$. Also using $F+H<1+FH$ and $G+I<1+GI$, we get $A+5B-B^{5}A(FH)(GI)+FH+GI>7.$ Now the left side is a decreasing function of $FH$ as $B^5AGI>B^6(1-b)>1$ for $0<b<0.313$, so we replace $FH$ by $1-b$. Similarly replacing $GI$ by $1-b$, we get $1+a+3b-(1+b)^{5}(1+a)(1-b)^{2}>0$. Now the left side is a decreasing function of $a$, so replacing $a$ by $b$, we get $1+4b-(1+b)^{6}(1-b)^{2}>0$, which is not true for
$0<b<0.313$. Hence we must have $B>C$.

\noindent{\bf Claim(iii) $C\leq 1.215$}

Suppose $C>1.215$, then $b>0.215$ using Claim (ii). $C^{4}ABGHI>C^{4}B^{2}(1-b+\frac{h+i}{2})(1-h)(1-i)=\phi(h)>\phi(b-\frac{i}{2})=(1+0.215)^{4}(1+b)^{2}(1-\frac{b}{2}+\frac{i}{4})(1-b+\frac{i}{2})(1-i)=\eta(b,i)$.
Now $\eta(b,i)>2$ for $0.215<b<0.313$ and $0\leq i<b$. Hence (1,1,4,1,1,1) holds, i.e. $A+B+4C-\frac{1}{2}C^{5}ABGHI+G+H+I>9$.
The left side is a decreasing function of $C$ and $G$, so we replace $C$ by 1.215 and $G$ by $1-b+\frac{h+i}{2}$ to get $\phi_{ord}(h,i,b,a)=a+4\times 0.215-\frac{h+i}{2}-\frac{1}{2}(1.215)^{5}(1+a)(1+b)(1-b+\frac{h+i}{2})(1-h)(1-i)>0$, where $0\leq h<b-\frac{i}{2}$ and $0\leq i<b$. Following Remark 3 we find that ${\rm max}~\phi_{ord}(h,i,b,a)\leq 0$ for $0<b<0.313$ and $0<a\leq 1$. Hence $C\leq 1.215$.\vspace{3mm}

\noindent{\bf Claim(iv) $d+f > 0.472b$ ~{\rm and} ~$g+h+i<1.056b$}

Suppose $d+f \leq 0.472b$. Applying AM-GM inequality to $(1,2,2,2,1,1)$ we have
$6+a+4b-4(d+f)-h-i-6(1+b)(1-d)(1-f)(1+a)^{\frac{1}{3}}(1-h)^{\frac{1}{3}}(1-i)^{\frac{1}{3}}>0.$ Left side is an increasing function of $h$ as $2(1+b)(1-d-f)(1-i)^{\frac{1}{3}}>
2(1+b)(1-0.472b)(1-b)^{\frac{1}{3}}>1$
and we have $h<b-d-f-\frac{i}{2}$, using (5.65). Therefore replacing  $h$ by $b-d-f-\frac{i}{2}$ we get\\
$\phi_{ord}(i,d+f,b,a)=6+a+3b-3(d+f)-\frac{i}{2}-6(1+b)(1-d-f)(1+a)^{\frac{1}{3}}(1-b+d+f+\frac{i}{2})^{\frac{1}{3}}(1-i)^{\frac{1}{3}}>0,$ where
$0\leq i\leq b-d-f$ (using (5.61)) and $0\leq d+f\leq 0.472b$. Following Remark 3 we find that ${\rm max}~\phi_{ord}(i,d+f,b,a)\leq 0$ for $0<b<0.313$ and $b\leq a\leq 1$.
Hence we have $d+f > 0.472b$. Now using (5.60), we get $g+h+i<2b-2(d+f)<2b-2(0.472b)<1.056b$.\vspace{3mm}

\noindent{\bf Claim(v) $B\leq 1.255$}

Suppose $B>1.255$, then $B^{4}AFGHI>(1+b)^{5}(1-b+\frac{g+h+i}{2})(1-(g+h+i))>(1+b)^{5}(1-b+\frac{1.056b}{2})(1-1.056b)>2$, for $0.255<b<0.313$.
Now we get a contradiction using (1,4,1,1,1,1) and proceeding as in Claim (i).\vspace{3mm}

\noindent{\bf Claim(vi) $C\leq 1.157$}

Suppose $C>1.157$, then $C^{4}ABGHI>(1.157)^{4}(1+b)^{2}(1-(g+h+i))>(1.157)^{4}(1+b)^{2}(1-1.056b)>2$, for $0.157<c<b\leq0.255$.
Now working as in Claim (iii) and using lower bound of $c$ as 0.157, the inequality (1,1,4,1,1,1) gives a contradiction. Hence $C\leq 1.157$.\vspace{3mm}

\noindent{\bf Claim(vii) $A\leq 1.43$}

Suppose $A>1.43$, then $A^4EFGHI=\frac{A^3}{BCD}>\frac{(1.43)^3}{1.255\times1.157}>2$. So (4,1,1,1,1,1) holds. Using $E+F+G+H+I<4+EFGHI$, we have $4A-\frac{1}{2}A^5EFGHI+EFGHI>5$. Coefficient of $EFGHI$ is negative and $EFGHI=\frac{1}{ABCD}>\frac{1}{(1.157)(1.255)A}$.
Therefore
$4A-\frac{A^{4}}{2\times1.157\times1.255}+\frac{1}{(1.157)(1.255)A}>5$, which is not true for $1.43<A\leq 2$. Hence we have $A\leq 1.43$.\vspace{3mm}

\noindent{\bf Claim(viii) $f+g+h+i>2c$ {\rm and} $e<\frac{a}{2}$}

Suppose $f+g+h+i \leq 2c$.  Using  inequality  $(1,1,3,1,1,1,1)$, we have $ A+B+4C-C^4FGHIAB+F+G+H+I>9$. Now the coefficients of $A$ and $B$ are negative on the left side, as $C^4FGHI>(1+c)^4(1-2c)>1$. So we can replace both $A$ and $B$ by
$C$ to get $6C-C^6FGHI+F+G+H+I>9.$ This implies
$1+6c-(f+g+h+i)-(1+c)^6(1-f-g-h-i)>0$, which is not true for $f+g+h+i\leq 2c$ and $c<0.157$. Now using (5.66) we have $e<\frac{a}{2}+c-(f+g+h+i)<\frac{a}{2}$.\vspace{3mm}

\noindent{\bf Final contradiction}

The inequality $(1,2,1,1,1,1,1,1)$ gives $A+4B-2B^{3}DEFGHIA+D+E+F+G+H+I >9$.  As the coefficient
of $D$  is negative, we can replace $D$ by $1-b+\frac{e+f+g+h+i}{2}$ to get $\phi_{ord}(e,f,x,b,a)=2+a+3b-\frac{e+f+x}{2}-2(1+b)^{3}(1+a)(1-b+\frac{e+f+x}{2})(1-e)(1-f)(1-x)>0$, where $g+h+i=x$, $0\leq e\leq\frac{a}{2}$, $0\leq f\leq b-\frac{x}{2}$ and $0\leq x\leq1.056b$. Following Remark 3 we find that ${\rm max}~\phi_{ord}(e,f,x,b,a)\leq 0$ for $0<b\leq \min(a,0.255)$ and $0<a<0.43$.\vspace{2mm}\\
{\noindent \bf  Proposition 35.} Case (51) i.e. $A > 1,~ B > 1, ~C \leq 1, ~D > 1,~ E \leq 1, ~F \leq 1,~ G \leq 1, ~H > 1, ~I > 1$ does not arise.\vspace{3mm}\\
{\noindent \bf Proof.} Here $a\leq \frac{1}{2}$,~$b\leq \frac{1}{3}$,~$d\leq \frac{1}{3}$. Using the weak inequality (1,2,2,2,1,1), we get
\begin{equation}a-2c-2e-2g+h+i>0\end{equation}
{\noindent \bf Claim(i) $D^{4}ABCHI < 2$~ {\rm and}~ $EFG>\frac{1}{2}$}

Suppose $D^4ABCHI\geq 2$, then (3,4,1,1) holds. Applying AM-GM inequality we get
$4A+4D+H+I-\sqrt{2D^5A^4HI}-9>0,$ i.e. $\psi(h,i)=1+4a+4d+h+i-\sqrt{2(1+d)^5(1+a)^4(1+h)(1+i)}>0.$ Following Remark 2 we find that ${\rm max}~\psi(h,i)\leq 0$ for $0<a\leq \frac{1}{2}$ and $0<d\leq \frac{1}{3}$, giving thereby a contradiction.\\
Hence we have $D^{4}ABCHI < 2$, which implies $\frac{D^{3}}{EFG}<2$, i.e. $EFG>\frac{1}{2}$.\vspace{3mm}

{\noindent \bf Claim(ii) $A^{4}EFGHI < 2$,~ {\rm and}~  $A < \sqrt{2}$}

Proof is similar to that of Claim(ii) of Case (38)(Proposition 19).\vspace{3mm}

{\noindent \bf Claim(iii) $g\geq 0.406(a+h+i)$ ~{\rm and}~ $c+e<0.094(a+h+i)$}

Suppose $g<0.406(a+h+i)$. Applying AM-GM inequality to (3,3,1,1,1) we have $4A+4D+G+H+I-2A^{2}D^{2}\sqrt{GHI}>9$. Left side is a decreasing function of $D$ as $D>1$ and $A^4GHI>(1+a)^4(1-0.406(a+h+i))(1+h)(1+i)\geq 1$, for $0<h\leq a$, $0<i\leq a$ and $0<a<\sqrt{2}-1$. So we replace $D$ by 1 and get that $4A+G+H+I-2A^{2}\sqrt{GHI}>5$. Now left side is a decreasing function of $G$ as
$A^4HI>1$ and $G\leq 1$. So replacing $G$ by $1-0.406(a+h+i)$, we get
$\phi(h,i)=2+3.594a+0.594(h+i)-2(1+a)^{2}\sqrt{(1-0.406(a+h+i))(1+h)(1+i)}>0.$
Following Remark 2 we find that ${\rm max}~\phi(h,i)\leq0$ for $0<a<\sqrt{2}-1$. Hence we must have $g\geq 0.406(a+h+i)$. Now using (5.67) we get $c+e<\frac{a+h+i}{2}-g<0.094(a+h+i)$.\vspace{3mm}

{\noindent \bf Claim(iv) $A\geq 1.388$}

Suppose $A<1.388$. Applying AM-GM inequality to (2,2,2,1,1,1) we get
$4A+4C+4E+G+H+I-6(A^{3}C^{3}E^{3}GHI)^{\frac{1}{3}}>9$. Left side is a decreasing function of $G$ as $ACE(HI)^{\frac{1}{3}}>A(1-0.094(a+h+i))(HI)^{\frac{1}{3}}>\frac{1}{2}$,
for $0<a\leq 0.5$, $0<h\leq a$ and $0<i\leq a$. So replacing $G$ by $1-\frac{a+h+i}{2}+(c+e)$, we get\\
$\eta(x)=6+\frac{7a}{2}+\frac{h+i}{2}-3x-6(1+a)(1-x)\left(1-\frac{a+h+i}{2}+x\right)^{\frac{1}{3}}(1+h)^{\frac{1}{3}}(1+i)^{\frac{1}{3}}>0,$\\
where $x=c+e$. Now $\eta''(x)>0$ and we have $0\leq x<0.094(a+h+i)$. Therefore $\eta(x)\leq{\rm max}\{\eta(0),\eta(0.094(a+h+i))\}$. Let
$\eta(0)=\eta_{1}(h,i)$ and $\eta(0.094(a+h+i))=\eta_{2}(h,i)$. Following Remark 2 we find that for $m=1,2$, ${\rm max}~\eta_m(h,i)\leq0$ for $0<a<0.388$. Hence we must have $A\geq 1.388$.\vspace{3mm}

{\noindent \bf Claim(v) $G\leq 0.73$}

Suppose $G>0.73$. Applying AM-GM inequality to (3,3,1,1,1) we get $4A+4D+G+H+I-2A^{2}D^{2}\sqrt{GHI}>9$. Left side is a decreasing function of $D$ and $G$. Replacing $D$ by 1 and $G$ by 0.73, we get
$$\phi(h,i)=1.73+4a+h+i-2\sqrt{0.73}(1+a)^{2}\sqrt{(1+h)(1+i)}>0.$$
Following Remark 2 we find that ${\rm max}~\phi(h,i)\leq0$ for $0.388\leq a<0.415$. Hence we must have $G\leq 0.73$.\vspace{3mm}

{\noindent \bf Final Contradiction:}

Using $G\leq 0.73$ and $EFG>\frac{1}{2}$, we get $F>\frac{1}{2EG}\geq \frac{1}{2G}>0.684$.\\
Now $A^{4}EFGHI>0.684(1+a)^{4}(1-e)(1-\frac{a+h+i}{2}+e)(1+h)(1+i)=\chi(e)$.\\
Using that $e<0.094(a+h+i)$, we find that $\chi'(e)>0$. So we have $\chi(e)>\chi(0),$ for $e>0$, i.e. $A^{4}EFGHI>0.684(1+a)^{4}(1-\frac{a+h+i}{2})(1+h)(1+i)=\psi(h,i)$.
Following Remark 2 we find that ${\rm min}~\psi(h,i)\geq2$ for  $0.388\leq a<\sqrt{2}-1$.
Hence $A^{4}EFGHI>2$, contradicting Claim (ii).\vspace{2mm}

{\noindent \bf  Proposition 36.} Case (52) i.e. $A > 1,~ B > 1, ~C \leq 1, ~D > 1,~ E \leq 1, ~F \leq 1,~ G \leq 1, ~H \leq 1, ~I > 1$ does not arise.\vspace{3mm}\\
{\noindent \bf Proof.} Here $a\leq \frac{1}{2}$,~$b\leq \frac{1}{3}$. Using the weak inequalities (1,2,2,2,1,1),(1,2,1,2,1,1,1),
\\(1,2,2,1,1,1,1) and (1,2,2,1,2,1) we get \vspace{-2mm}
\begin{equation}a-2c-2e-2g-h+i>0,\vspace{-3mm}\end{equation}
\begin{equation}a-2c+d-2f-g-h+i>0,\vspace{-2mm}\end{equation}
\begin{equation}a-2c-2e-f-g-h+i>0,\vspace{-2mm}\end{equation}
\begin{equation}a-2c-2e-f-2h+i>0.\end{equation}
{\noindent \bf Claim(i) $c+e+g\geq 0.279(a+i)$~ {\rm and}~ $h<0.442(a+i)$}

Suppose $c+e+g<0.279(a+i)$. Applying AM-GM inequality to (2,2,2,2,1) we get $4A+4C+4E+4G+I-8(A^{3}C^{3}E^{3}G^{3}I)^{\frac{1}{4}}>9.$
Let $c+e+g=x$, then we have $\phi(x)=8+4a-4x+i-8((1+a)^{3}(1-x)^{3}(1+i))^{\frac{1}{4}}>0.$ Now $\phi(x)$ is an increasing function of $x$,
 so replacing $x$ by $0.279(a+i)$, we get\\
$8+4a-4(0.279)(a+i)+i-8(1+a)^{\frac{3}{4}}(1-0.279(a+i))^{\frac{3}{4}}(1+i)^{\frac{1}{4}}>0. $\\
But this is not true for $0<a\leq 0.5$ and $0<i\leq a$. Hence we must have $c+e+g\geq 0.279(a+i)$. Now using (5.68), we get $h<(a+i)-2(c+e+g)=0.442(a+i)$.\vspace{3mm}

{\noindent \bf Claim(ii) $A\leq 1.35$.}

Suppose $A>1.35$. We prove that $A^{4}EFGHI>2$. Now $F\geq\frac{2D}{3}$ gives $\frac{d}{2}-\frac{1}{4}+\frac{3f}{4}<0$ and from (5.69) we have
$f-\frac{a+d+i}{2}+\frac{g+h}{2}<0$.
Adding these two we get $f\leq\frac{1}{7}+\frac{2}{7}(a+i)-\frac{2}{7}(g+h).$
As $e<\frac{a+i}{2}-\frac{f+g+h}{2}$ from (5.70), we have $A^{4}EFGHI>(1+a)^{4}(1-\frac{a+i}{2}+\frac{f+g+h}{2})(1-f)(1-g)(1-h)(1+i)$.\\
Now right side is a decreasing function of $f$ for $a\leq \frac{1}{2}$, so replacing $f$ by $\frac{1}{7}+\frac{2}{7}(a+i)-\frac{2}{7}(g+h)$, we get
$A^{4}EFGHI>(1+a)^{4}\left(\frac{15}{14}-\frac{5}{14}(a+i)+\frac{5}{14}(g+h)\right)\\ \left(\frac{6}{7}-\frac{2}{7}(a+i)+\frac{2}{7}(g+h)\right)(1-g)(1-h)(1+i)
=\phi_{ord}(g,h,i,a)$, say, where $0\leq g<\frac{a+i}{2}-\frac{h}{2}$(using (5.68)) and $0\leq h<0.442(a+i)$ and $0<i\leq a$. Following Remark 3 we find that ${\rm min}~\phi_{ord}(g,h,i,a)\geq2$ for $0.35<a\leq 0.5$. Hence $A^{4}EFGHI>2$.
\\So (4,2,1,1,1) holds. Applying AM-GM inequality we get $\psi(E)=4A+4E+G+H+I-2\sqrt{A^{5}E^{3}GHI}>9.$
$\psi'(E)=4-3\sqrt{A^{5}EGHI}<0$, for $A>1.35$, $E>\frac{3}{4}$, $0\leq g<\frac{a+i}{2}-\frac{h}{2}$ and $0\leq h<0.442(a+i)$.
So replacing $E$ by $\frac{3}{4}$ we get $4A+3+G+H+I-2\sqrt{\frac{27}{64}A^{5}GHI}>9$. Now left side is a decreasing function of $G$ for  $A>1.35$, $G\leq 1$ and
$0\leq h<0.442(a+i)$, so replacing $G$ by $1-\frac{a+i}{2}+\frac{h}{2}$, we get\\
$\theta(h)=1+\frac{7a}{2}+\frac{i}{2}-\frac{h}{2}-2\sqrt{\frac{27}{64}(1+a)^{5}(1-\frac{a+i}{2}+\frac{h}{2})(1-h)(1+i)}>0.$
Now $\theta''(h)>0$, therefore $\theta(h)\leq{\rm max}\{\theta(0),\theta(0.442(a+i))\}$, which turns out to be less than zero for $0<a\leq 0.5$ and $0<i\leq a$.
Hence $A\leq 1.35$.\vspace{3mm}

{\noindent \bf Claim(iii) $c+e\geq 0.18(a+i)$~ {\rm and}~ $f+g+h<0.64(a+i)$}

Suppose $c+e<0.18(a+i)$. Applying AM-GM inequality to (2,2,2,1,1,1) we get
 $\phi_{ord}(g,h,c+e,i,a)=6+4a-4(c+e)-g-h+i-6(1+a)(1-(c+e))(1-g)^{\frac{1}{3}}(1-h)^{\frac{1}{3}}(1+i)^{\frac{1}{3}}>0,$ where $0\leq g <\frac{a+i}{2}-\frac{h}{2}-(c+e)$ (using (5.68)), $0\leq h<\frac{a+i}{2}-(c+e)$ (using (5.71)), $0\leq c+e \leq 0.18(a+i)$ and $0<i\leq a$. Following Remark 3 we find that ${\rm max}~\phi_{ord}(g,h,c+e,i,a)\leq 0$ for $0<a\leq 0.35$. Hence we must have $c+e\geq 0.18(a+i)$.\\
Now using (5.70), we get $f+g+h<0.64(a+i)$.
\vspace{3mm}

{\noindent \bf Claim(iv) $A\leq 1.29$.}

Suppose $A>1.29$. We prove that $A^{4}EFGHI>2$.\\
Using (5.70) we have $A^{4}EFGHI>(1+a)^{4}(1-\frac{a+i}{2}+\frac{f+g+h}{2})(1-f)(1-g)(1-h)(1+i)$. Let $f+g+h=y$, then $A^{4}EFGHI>(1+a)^{4}(1-\frac{a+i}{2}+\frac{y}{2})(1-y)(1+i)=\phi(y).$
It is easy to check that $\phi'(y)<0$. So using Claim(iii) we have $\phi(y)>\phi(0.64(a+i))>2$ for $0.29<a\leq 0.35$ and $0<i\leq a$. Hence $A^{4}EFGHI>2$ and therefore (4,2,1,1,1) holds. Now proceeding as in Claim (ii) we get a contradiction .
\vspace{3mm}

{\noindent \bf Claim(v) $c+e+g\geq 0.377(a+i)$~ {\rm and}~ $h<0.246(a+i)$}

Suppose $c+e+g<0.377(a+i)$. Now proceeding as in Claim (i) and replacing $c+e+g$ by $0.377(a+i)$ in place of $0.279(a+i)$, we get
a contradiction for $0<a\leq 0.29$ and $0<i\leq a$. Hence we must have $c+e+g\geq 0.377(a+i)$. Now using (5.68), we get $h<0.246(a+i)$.\vspace{3mm}

{\noindent \bf Final Contradiction:}

Applying AM-GM inequality to (3,3,1,1,1) we get $4A+4D+G+H+I-2\sqrt{A^{4}D^{4}GHI}>9$. Left side is a decreasing function of $D$ as $A^4GHI>1$, for $g<\frac{a+i}{2}-\frac{h}{2}$ and $0\leq h<0.246(a+i)$. So replacing $D$ by 1 we get $4A+4+G+H+I-2\sqrt{A^{4}GHI}>9$. Now left side is a decreasing function of $G$.
So we replace $G$ by $1-\frac{a+i}{2}+\frac{h}{2}$ and get that
\begin{equation}\phi(h)=2+\frac{7a}{2}+\frac{i}{2}-\frac{h}{2}-2(1+a)^{2}\sqrt{(1-\frac{a+i}{2}+\frac{h}{2})(1-h)(1+i)}>0.\end{equation}
But $\phi''(h)>0$, so $\phi(h)\leq{\rm max}\{\phi(0),\phi(0.246(a+i))\}$, which can be easily checked to be at most zero for $0<a\leq 0.29$ and $0<i\leq a$, which contradicts (5.72).\vspace{2mm}\\
{\noindent \bf Proposition 37.}  Case (64) i.e. $ A>1,~ B>1, ~C\leq1, ~D\leq1,~ E\leq 1, ~F> 1,~ G>1, ~H>1, ~I>1$ does not
arise.\vspace{2mm}\\
 {\noindent \bf Proof.} Here $ a\leq \frac{1}{2}, ~b\leq \frac{1}{3}$. Using  the weak inequalities
$(1,2,2,1,1,1,1)$, $(2,1,2,1,1,1,1)$ and  $(2,2,1,1,1,1,1)$ we get \vspace{-2mm}
\begin{equation}a-2c-2e+f+g+h+i>0,\vspace{-3mm}\end{equation}
\begin{equation}2b-c-2e+f+g+h+i>0,\vspace{-2mm}\end{equation}
\begin{equation}2b-2d-e+f+g+h+i>0.\end{equation}

\noindent{\bf Claim(i)} $B^{4}FGHIA \leq 2$ and  $B<1.149$

Suppose $B^{4}FGHIA >2$. Therefore $(1,4,1,1,1,1)$ holds i.e.  $A+4B-\frac{1}{2}B^{5}FGHIA $ $+F+G+H+I>9.$ This is not true, by Lemma 10(ii) with $X_{2}=B,~X_{3}=F,~X_{4}=G,~X_{5}=H,~X_{6}=I$. Further $B^{4}FGHIA \leq 2$ implies $B^5\leq 2$ i.e. $B<1.149$.\vspace{3mm}

\noindent{\bf Claim(ii)} $e>2b$ ~and ~$f+g+h+i>2b$

Suppose $e\leq 2b$. Applying  Lemma 7(v) with
$x_{1}=b, ~\gamma=e $ and $~\delta=a+f+g+h+i$, we get
a contradiction as $\gamma=e<\frac{\delta}{2}$ (using (5.73)), ~$\gamma=e\leq 2b=2x_{1}$ and $x_1<0.149$. So we must have $e>2b.$ Now (5.74) gives that $f+g+h+i>2b$.\vspace{3mm}

\noindent{\bf Claim(iii)} $b<0.106$

Suppose $b \geq 0.106$, then $B^{4}FGHIA >B^{5}(1+f+g+h+i)>B^{5}(1+2b)>2$. This contradicts Claim(i).\vspace{3mm}

\noindent{\bf Claim(iv)} $a<0.223$

Suppose $a\geq 0.223$. We prove that $A^4EFGHI>2$. Consider following cases:\\
\textbf{Case(i)} $0.223\leq a<0.36$\\
From the inequality (5.74) we have $e<b+\frac{f+g+h+i}{2}<0.106+\frac{f+g+h+i}{2}$. Therefore $A^4EFGHI>(1+a)^4(1-0.106-\frac{f+g+h+i}{2})(1+f)(1+g)(1+h)(1+i)=\phi(f,g,h,i).$ Following Remark 2 we find that ${\rm min}~\phi(f,g,h,i)>2$, for $0.223\leq a<0.36$.\\
\textbf{Case(ii)} $0.36\leq a<0.5$\\
Here $A^4EFGHI=\frac{A^3}{BCD}>\frac{1.36^3}{1.106}>2.$ \vspace{2mm}\\
Hence (4,1,1,1,1,1) holds, i.e. $4A-\frac{1}{2}A^{5}EFGHI+E+F+G+H+I>9$. As for $A\geq1.223$, the coefficient of $E$ namely
$1-\frac{1}{2}A^{5}FGHI$ is negative and $E\geq 0.46873A$, so replacing $E$ by $0.46873A$ we get $\eta(f,g,h,i)=4.46873a+f+g+h+i-\frac{1}{2}(0.46873)(1+a)^{6}(1+f)(1+g)(1+h)(1+i)-0.53127>0$. Following Remark 2 we find that ${\rm max}~\eta(f,g,h,i)\leq 0$, for $0\leq a<0.5$. Hence we have $A<1.223$.  \vspace{3mm}

\noindent{\bf Claim(v)} $e>2b+0.165k$, ~$d<0.4175k$ ~and ~$k>\frac{2b}{0.67}$, where $k=f+g+h+i$

Suppose $e \leq 2b+0.165k$. Using $(1,3,1,1,1,1,1)$ we have $A+4B-B^{4}EFGHIA +E+F+G+H+I>9$. As the coefficient of $E$ is
negative and $E>1-2b-0.165k$, we have $A+4B-B^{4}(1-2b-0.165k)FGHIA +1-2b-0.165k+F+G+H+I>9$. This gives
$\theta(k)=1+a+2b+0.835k-(1+b)^{4}(1-2b-0.165k)(1+a)(1+k)>0$. As $\theta''(k)>0$ and $2b < k \leq 4a$, we get
$\theta(k) \leq {\rm max}\{\theta(2b),\theta(4a)\}<0$ for $0<a<0.223$ and $0<b\leq{\rm min}(a,0.106)$. This gives a
contradiction.\vspace{1mm}\\
Now using (5.75), we have $2d<2b-e+k$, i.e. $d<0.4175k$. Also (5.74) gives $k>2e-2b>2(2b+0.165k)-2b$. It gives $k>\frac{2b}{0.67}$.\vspace{3mm}

\noindent{\bf Claim(vi)} $A>1.16$

Suppose $A\leq 1.16$. Consider following two cases:\\
\noindent{\bf Case I.} $1<A\leq 1.118$\\
Applying AM-GM inequality to (1,2,2,1,1,1,1), we get $A+4B+4D+F+G+H+I-4\sqrt{B^3D^3AFGHI}>9.$ Left side is a decreasing function of $D$, as $BDAFGHI=\frac{1}{CE}\geq1$ and we have $D>1-0.4175(f+g+h+i)$. So replacing $D$ by $1-0.4175(f+g+h+i)$, we get\\
$\begin{array}{ll}\eta(f,g,h,i)=&4+a+4b-0.67(f+g+h+i)-4(1+b)^{\frac{3}{2}}(1-0.4175(f+g+\\&h+i))^{\frac{3}{2}}(1+a)^{\frac{1}{2}}(1+f)^{\frac{1}{2}}(1+g)^{\frac{1}{2}}
(1+h)^{\frac{1}{2}}(1+i)^{\frac{1}{2}}>0.\end{array}$ Following Remark 2 we find that ${\rm max}~\eta(f,g,h,i)\leq0$ for $0<a\leq0.118$ and $0<b<0.106$.\\

\noindent{\bf Case II.} $1.118<A\leq 1.16$\\
Applying AM-GM inequality to (3,2,1,1,1,1) we get $4A+4D+F+G+H+I-2\sqrt{2}(A^4D^3FGHI)^{\frac{1}{2}}>9$. Left side is a decreasing function of $D$, so we replace $D$ by $1-0.4175(f+g+h+i)$, i.e.
$\phi(f,g,h,i)=3+4a-0.67(f+g+h+i)-2\sqrt{2}(1+a)^{2}(1-0.4175(f+g+h+i))^{\frac{3}{2}}(1+f)^{\frac{1}{2}}(1+g)^{\frac{1}{2}}(1+h)^{\frac{1}{2}}(1+i)^{\frac{1}{2}}>0$.
Following Remark 2 we find that ${\rm max}~\phi(f,g,h,i)\leq0$ for $0.118<a\leq 0.16$, giving thereby a contradiction.\vspace{3mm}

\noindent{\bf Claim (vii)} $b<0.084$

Suppose $b \geq 0.084$. Using  $k=f+g+h+i>\frac{2b}{0.67}$ we have $B^{4}FGHIA
>1.16B^{4}(1+k)>1.16B^{4}(1+\frac{2b}{0.67})>2$ for $b \geq 0.084$. This contradicts Claim(i).\vspace{3mm}

\noindent{\bf Claim (viii)} $e>2b+0.33k$ ~and~ $d<0.335k$

Suppose $e \leq 2b+0.33k$. Working as in Claim(v) and replacing $E$ by $1-2b-0.33k$, we get
$\theta(f,g,h,i)=1+a+2b+0.67(f+g+h+i)-(1+b)^{4}(1-2b-0.33(f+g+h+i))(1+a)(1+f)(1+g)(1+h)(1+i)>0$. Following Remark 2 we find that ${\rm max}~ \theta(f,g,h,i)\leq0$ for $0<b<0.084$ and $0.16<a<0.223$, giving thereby a contradiction.\\
Now using (5.75), we have $2d<2b-e+k$, i.e. $d<0.335k$.\vspace{3mm}

\noindent{\bf Final contradiction}

We proceed as in Case II of Claim (vi) and replace $d$ by $0.335(f+g+h+i)$ in place of $0.4175(f+g+h+i)$. Then we get $\phi(f,g,h,i)=3+4a-0.34(f+g+h+i)-2\sqrt{2}(1+a)^{2}(1-0.335(f+g+h+i))^{\frac{3}{2}}(1+f)^{\frac{1}{2}}(1+g)^{\frac{1}{2}}(1+h)^{\frac{1}{2}}(1+i)^{\frac{1}{2}}>0$.
Following Remark 2 we find that ${\rm max}~\phi(f,g,h,i)\leq0$ for $0.16<a<0.223$, giving thereby a contradiction.\vspace{3mm}\\
{\noindent \bf Proposition 38.} Case (69) i.e. $ A>1,~ B>1, ~C\leq 1, ~D\leq1,~ E\leq 1, ~F\leq 1,~ G > 1, ~H > 1, ~I > 1$
 does not arise.\vspace{2mm}

{\noindent \bf Proof.} Here  $ a\leq \frac{1}{2},~ b\leq \frac{1}{ 3}$. Using the weak inequalities
$(1,2,2,1,1,1,1),(2,2,2,\\1,1,1)$, ~$(2,2,1,1,1,1,1)$ ~and ~(1,2,1,2,1,1,1)  we get \vspace{-2mm}
\begin{equation}a-2c-2e-f+g+h+i>0,\vspace{-3mm}\end{equation}
\begin{equation}2b-2d-2f+g+h+i>0,\vspace{-2mm}\end{equation}
\begin{equation}2b-2d-e-f+g+h+i>0,\vspace{-2mm}\end{equation}
\begin{equation}a-2c-d-2f+g+h+i>0.\end{equation}
Using (5.76), we have \begin{equation}e<\frac{a+g+h+i}{2}-\frac{f}{2}\end{equation}
Using $F\geq\frac{2}{3}D$ and (5.77) respectively we have $1-d\leq\frac{3}{2}-\frac{3f}{2}$ and $d+f<\frac{2b+g+h+i}{2}$. Adding these two we get
\begin{equation}f<\frac{1+2b+g+h+i}{5}\end{equation}
 Now $E\geq\frac{2}{3}C$ and (5.76) respectively implies  $1-c<\frac{3}{2}-\frac{3e}{2}$ and $c+e<\frac{a+g+h+i}{2}-\frac{f}{2}$. Adding these two we get
\begin{equation}e<\frac{1+a+g+h+i-f}{5}\end{equation}

 \noindent{\bf Claim (i)} $A < 1.313$

Suppose $A\geq 1.313$. First we prove that $A^4EFGHI\geq 2$. We have $A^4EFGHI=\frac{A^3}{BCD}\geq \frac{3}{4}A^3>2$, for $A\geq 1.39$. So now consider $1.313\leq A<1.39$. Using (5.80) we have $A^4EFGHI>A^4(1+g)(1+h)(1+i)(1-\frac{a+g+h+i}{2}+\frac{f}{2})(1-f)=\eta(f)$. Now $\eta''(f)<0$ and $0\leq f<\frac{1}{3}+\frac{g+h+i}{5}$, using (5.81) and that $b\leq\frac{1}{3}.$  Therefore $\eta(f)\geq{\rm min}\{\eta(0),\eta(\frac{1}{3}+\frac{g+h+i}{5})\}$. \\ Let
$\eta(0)=\phi_1(g,h,i)$ and $\eta(\frac{1}{3}+\frac{g+h+i}{5})=\phi_2(g,h,i).$\\
Following Remark 2 we find that for $m=1,2$, ${\rm min}~\phi_m(g,h,i)\geq2$ for $0.313\leq a<0.39$. Hence we have $A^4EFGHI>2$ for $A\geq1.313$. So the inequality (4,1,1,1,1,1) gives $4a-e-f+g+h+i-\frac{1}{2}(1+a)^5(1-e)(1-f)(1+g)(1+h)(1+i)>0.$ Coefficient of $e$ on left side is positive for $f<\frac{1}{3}+\frac{g+h+i}{5}$, $0<g\leq a$, $0<h\leq a$, $0<i\leq a$ and $0.246\leq a\leq \frac{1}{2}$. So using (5.82) we replace $e$ by $\frac{1+a+g+h+i-f}{5}$  and get that   $\psi_{ord}(f,g,h,i,a)=\frac{-1}{5}+\frac{19a}{5}+\frac{4}{5}(g+h+i)-\frac{4}{5}f-\frac{1}{2}(1+a)^5(1+g)(1+h)(1+i)(1-f)
(\frac{4}{5}-\frac{a+g+h+i}{5}+\frac{f}{5})>0,$ where $0\leq f<\frac{1}{3}+\frac{g+h+i}{5}$, $0<g\leq a$, $0<h\leq a$ and $0<i\leq a$.
Following Remark 3 we find that ${\rm max}~\psi_{ord}(f,g,h,i,a)\leq0$ for $0<a\leq \frac{1}{2}$. Hence we have a contradiction. So we must have $a<0.313$.\vspace{2mm}\\
\noindent{\bf Claim (ii)} $B < 1.17$

 Suppose $B\geq 1.17$. Using (5.79) we have $B^4AFGHI>(1+b)^4(1+a)(1+g)(1+h)(1+i)(1-\frac{a+g+h+i}{2})=\phi(g,h,i)$. But $\phi(g,h,i)\geq{\rm min}\{\phi(0,0,0), \phi(a,0,0),\\ \phi(a,a,0), \phi(a,a,a)\}>2$ for $0.17\leq b\leq a<0.313.$ Hence $B^4AFGHI>2$. So (1,4,1,1,1,1) holds. That is $A+4B-\frac{1}{2}B^5AFGHI+F+G+H+I>9$. Coefficient of $F$, namely $1-\frac{1}{2}B^5AGHI$ is negative for $B>1.124$, so replacing $F$ by  $1-\frac{1+2b+g+h+i}{5}$ (using (5.81)) we get\\ $\eta(g,h,i)=\frac{-1}{5}+a+\frac{18}{5}b+\frac{4}{5}(g+h+i)-\frac{1}{2}(1+b)^5(1+a)(1+g)(1+h)(1+i)(\frac{4}{5}-\frac{2b+g+h+i}{5})>0$. But following Remark 2 we find that ${\rm max}~\eta(g,h,i)\leq0$ for $0<b\leq a<0.313$. Hence we must have $B<1.17$.\vspace{3mm}

 \noindent{\bf Claim (iii)} $A<1.29$

 Suppose $A\geq 1.29$. We work as in Claim(i) and use $f<\frac{1+2b+g+h+i}{5}<\frac{1.34+g+h+i}{5}$ instead of $f<\frac{1}{3}+\frac{g+h+i}{5}$. We find that $A^4EFGHI>2$ for $A\geq1.29$ and then (4,1,1,1,1,1) gives a contradiction.\vspace{3mm}

 \noindent{\bf Claim (iv)} $e+f>1.5b+0.351(g+h+i)$ and $d<\frac{b}{4}+0.3245(g+h+i)$

 Suppose $e+f\leq 1.5b+0.351(g+h+i).$ As (1,3,1,1,1,1,1) holds, we have $1+a+4b-(e+f)+g+h+i-(1+b)^4(1+a)(1+g)(1+h)(1+i)(1-(e+f))>0$. The coefficient of $e+f$ is clearly positive, so we replace $e+f$ by $1.5b+0.351(g+h+i)$ to get $\phi(g,h,i)=1+a+4b+g+h+i-(1.5b+0.351(g+h+i))-(1+b)^4(1+a)(1+g)(1+h)(1+i)(1-1.5b-0.351(g+h+i))>0$. But following Remark 2 we find that  ${\rm max}~\phi(g,h,i)\leq0$ for $0<b\leq{\rm min}(a,0.17)$ and $0<a<0.29$. This gives a contradiction. Hence we must have
 $e+f>1.5b+0.351(g+h+i)$. Now using (5.78) we have $2d<2b+g+h+i-(e+f)$. It gives $d<\frac{b}{4}+0.3245(g+h+i)$.\vspace{3mm}

 \noindent{\bf Claim (v)}  $A> 1.176$ and $g+h+i>1.97a$

 Suppose $A\leq 1.176$. Applying AM-GM to (1,2,2,1,1,1,1) we get $A+4B+4D+F+G+H+I-4(B^3D^3AFGHI)^{\frac{1}{2}}>9.$ Left side is a decreasing function of $F$ as $B^3D^3AGHI>\frac{1}{4}$, for $d<\frac{b}{4}+0.3245(g+h+i)$. Also $F=1-f>1-(b+\frac{g+h+i}{2}-d)$, using (5.77). So replacing $F$ by $1-b-\frac{g+h+i}{2}+d$, we get
 $\phi_{ord}(d,g,h,i,b,a)=4+a+3b-3d+\frac{g+h+i}{2}-4(1+b)^{\frac{3}{2}}(1-d)^{\frac{3}{2}}(1+a)^{\frac{1}{2}}(1+g)^{\frac{1}{2}}(1+h)^{\frac{1}{2}}
 (1+i)^{\frac{1}{2}}(1-b-\frac{g+h+i}{2}+d)^{\frac{1}{2}}>0,$ where $0\leq d<\frac{b}{4}+0.3245(g+h+i)$, $0<g\leq a$, $0<h\leq a$ and $0<i\leq a$. Following Remark 3 we find that  ${\rm max}~\phi_{ord}(d,g,h,i,b,a)\leq0$ for $0<b<0.17$ and  $0<a\leq0.176$. So $a>0.176$.\vspace{1mm}\\
Further if $g+h+i\leq1.97a$, then we have $0\leq g \leq {\rm min}(1.97a-h-i, a)$, ~$0\leq h \leq {\rm min}(1.97a-i, a)$, ~$0\leq i\leq a$. Now we find that ${\rm max}~\phi_{ord}(d,g,h,i,b,a)\leq0$ for $0<b<0.17$ and  $0<a\leq0.29$. Hence we must have $g+h+i>1.97a$.\vspace{3mm}

 \noindent{\bf Claim (vi)} $c+e>0.281(a+g+h+i)$ and $f<0.438(a+g+h+i)$

 Suppose $c+e\leq 0.281(a+g+h+i)$. Applying AM-GM inequality to (2,2,2,1,1,1) we get
 $6+4a-4(c+e)+g+h+i-6(1+a)(1-(c+e))(1+g)^{\frac{1}{3}}(1+h)^{\frac{1}{3}}(1+i)^{\frac{1}{3}}>0$. The coefficient of $c+e$ is positive, so we replace $c+e$ by
 $0.281(a+g+h+i)$ to get that $\phi(g,h,i)=6+4a+g+h+i-4\times 0.281(a+g+h+i)-6(1+a)(1-0.281(a+g+h+i))(1+g)^{\frac{1}{3}}(1+h)^{\frac{1}{3}}(1+i)^{\frac{1}{3}}>0$. Following Remark 2 we find that ${\rm max}~\phi(g,h,i)\leq0$ for $0<a<0.29$. This gives a contradiction. Hence we must have
 $c+e>0.281(a+g+h+i)$. Now using (5.76) we have $f<a+g+h+i-2(c+e)$. It gives $f<0.438(a+g+h+i)$.\vspace{3mm}

\noindent{\bf Claim (vii)} $b < 0.126$

Suppose $b\geq 0.126$, then $B^4AFGHI>(1+0.126)^4(1+a)(1+g)(1+h)(1+i)(1-0.438(a+g+h+i))=\phi(g,h,i)$. Now $1.97a<g+h+i\leq 3a$ implies ${\rm max}(1.97a-h-i,0)\leq g \leq a$, ~${\rm max}(0.97a-i,0)\leq h \leq a$ ~and~ $0\leq i \leq a$, so $\phi(g,h,i)$ being symmetric function in $g$, ~$h$ ~and~ $i$, we get that $\phi(g,h,i)>{\rm min}\{\phi(a,0.97a,0), \phi(a,a,0), \phi(a,a,0.97a), \phi(a,a,a)\}$, which can be verified to be greater than 2 for $0.176<a<0.29$. Hence $B^4AFGHI>2$ for $b\geq 0.126$ and so (1,4,1,1,1,1) holds. Now we work as in Claim (ii) and get a contradiction.\vspace{3mm}

\noindent{\bf Claim (viii)} $d > 0.171$

Suppose $d\leq 0.171$. We use (1,2,2,1,1,1,1) and proceed as in Claim(v). We find that ${\rm max}~\phi_{ord}(d,g,h,i,b,a)\leq0$ for $0<d\leq0.171$, ${\rm max}(1.97a-h-i,0)\leq g \leq a$, ~${\rm max}(0.97a-i,0)\leq h \leq a$, $0\leq i \leq a$, $0<b<0.126$ and  $0.176<a\leq0.29$. This gives a contradiction. Hence we must have $d > 0.171$.\vspace{3mm}

\noindent{\bf Claim (ix)} $A<1.232$

Suppose $A\geq 1.232$, then $A^4EFGHI=\frac{A^3}{BCD}>\frac{1.232^3}{1.126\times 0.829}>2$. Now working as in Claim(i) we get a contradiction with refined bounds $1.97a<g+h+i\leq3a$ and $a\geq0.232$.\vspace{3mm}

 \noindent{\bf Claim (x)} $c+e>0.312(a+g+h+i)$ and $f<0.376(a+g+h+i)$

 We use (2,2,2,1,1,1) and work as in Claim(vi). Here we replace $c+e$ by $0.312(a+g+h+i)$ and get a contradiction for $0.176<a<0.232$. Hence $c+e>0.312(a+g+h+i)$. Now (5.76) gives  $f<0.376(a+g+h+i)$.\vspace{3mm}

\noindent{\bf Claim (xi)} $b < 0.114$

Suppose $b\geq 0.114$. Now working as in Claim(vii) and using refined bounds on $f$, namely $f<0.376(a+g+h+i)$ we get contradiction. \vspace{3mm}

\noindent{\bf Claim (xii)} $e+f>1.5b+0.41(g+h+i)$ and $d<\frac{b}{4}+0.295(g+h+i)$

 Suppose $e+f\leq 1.5b+0.41(g+h+i).$ We proceed as in Claim (iv) and use refined bound on $e+f$, namely $1.5b+0.41(g+h+i)$ to get a contradiction for $0<b<0.114$ and $0.176<a<0.232$. Hence we must have
 $e+f>1.5b+0.41(g+h+i)$. Now using (5.78) we have $d<\frac{b}{4}+0.295(g+h+i)$.\vspace{3mm}

\noindent{\bf Claim (xiii)}  $A\geq 1.224$

 Suppose $A<1.224$. We use (1,2,2,1,1,1,1) and proceed as in Claim (v).
 We find that ${\rm max}~\phi_{ord}(d,g,h,i,b,a)\leq0$ for $0\leq d<\frac{b}{4}+0.295(g+h+i)$, $0\leq g \leq a$, ~$0\leq h \leq a$, $0\leq i \leq a$, $0<b<0.114$ and  $0<a\leq0.224$. This gives a contradiction. \vspace{3mm}

\noindent{\bf Claim (xiv)} $d > 0.19$

Suppose $d\leq 0.19$. We use (1,2,2,1,1,1,1) and proceed as in Claim(viii). Here we get a contradiction for $0\leq d\leq0.19$, $0.224<a<0.232$ and $0<b<0.114$.\vspace{3mm}

\noindent{\bf Final Contradiction}

We have $A^4EFGHI=\frac{A^3}{BCD}>\frac{1.224^3}{1.114\times 0.81}>2$. Now we get a contradiction proceeding as in Claim (i).\vspace{3mm}


{\noindent \bf  Proposition 39.} Case (70) i.e. $A > 1,~ B > 1, ~C \leq 1, ~D \leq 1,~ E \leq 1, ~F \leq 1,~ G > 1, ~H \leq 1, ~I > 1$ does not arise.\vspace{3mm}\\
{\noindent \bf Proof.} Here $a\leq \frac{1}{2}$, ~$b\leq \frac{1}{3}$ ~ and~  $g\leq \frac{1}{3}$.
Using the weak inequalities (1,2,2,1,2,1), (1,2,1,2,2,1), (2,2,2,2,1) and (2,2,1,1,2,1) we have \vspace{-2mm}
\begin{equation}a-2c-2e-f-2h+i>0,\vspace{-3mm}\end{equation}
\begin{equation}a-2c-d-2f-2h+i>0,\vspace{-2mm}\end{equation}
\begin{equation}2b-2d-2f-2h+i>0,\vspace{-2mm}\end{equation}
\begin{equation}2b-2d-e-f-2h+i>0.\end{equation}
\noindent{\bf Claim (i)} $B<1.17$

Suppose $B\geq 1.17$. Then using (5.84) we have $B^4AFGHI>(1+b)^4(1+a)(1+i)(1-(f+h))>1.17^4(1+a)(1+i)(1-\frac{a+i}{2})>2$ for $0.17\leq b\leq a\leq 0.5$ and $0<i\leq a$. Therefore (1,4,1,2,1) holds. Applying AM-GM inequality to (1,4,1,2,1) we have $A+4B+F+4G+I-2(B^5G^3AFI)^{\frac{1}{2}}>9$. Left side is a decreasing function of $F$ and $F>1-(b+\frac{i}{2})$, using (5.85). So we replace $F$ by $1-(b+\frac{i}{2})$ and  get that $2+a+3b+4g+\frac{i}{2}-2(1+b)^{\frac{5}{2}}(1+g)^{\frac{3}{2}}(1+a)^{\frac{1}{2}}(1+i)^{\frac{1}{2}}\left(1-b-\frac{i}{2}\right)^{\frac{1}{2}}>0.$ Now left side can be verified to be decreasing function of $g$, therefore replacing $g$ by $0$ we get
$\eta(i)=2+a+3b+\frac{i}{2}-2(1+b)^{\frac{5}{2}}(1+a)^{\frac{1}{2}}(1+i)^{\frac{1}{2}}\left(1-b-\frac{i}{2}\right)^{\frac{1}{2}}>0.$
But $\eta''(i)>0$ and $0<i\leq a$. So $\eta(i)\leq{\rm max}\{\eta(0),\eta(a)\},$ which can be seen to be at most zero for $0.17\leq b\leq a\leq 0.5$ Hence we get  a contradiction.\vspace{3mm}

\noindent{\bf Claim (ii)} $e+f+h>2b$ ~and ~$d<\frac{i}{2}$

Suppose $e+f+h\leq 2b$. We use the inequality  $(1,3,1,1,1,1,1)$ and apply  Lemma 7(vi) with
$X_{1}=B, ~\gamma=e+f+h $ and $~\delta=a+g+i$. Using (5.83) and (5.84) we get $2e+f+2h<a+i$ and $f<\frac{a+i}{2}$ respectively. Adding these two inequalities we get $e+f+h<\frac{3}{4}(a+i)<\frac{3}{4}(a+g+i)$, i.e $\gamma<\frac{3}{4}\delta$. Also $x_1=b<0.175$ and $\gamma=e+f+h\leq 2x_1$. Hence we have a contradiction.\vspace{1mm}\\
Now using (5.86), we have $d<b+\frac{i}{2}-\frac{e+f+h}{2}<\frac{i}{2}.$\vspace{3mm}

\noindent{\bf Claim (iii)} $A>1.3$

Suppose $A\leq 1.3$. We have $CE>1-(c+e)>1-\frac{a+i}{2}+\frac{f}{2}$ and $f<\frac{a+i}{2}-\frac{d}{2}$, using (5.83) and (5.84) respectively. Also $d<\frac{i}{2}$ from Claim (ii). Applying AM-GM inequality to (2,1,1,1,1,2,1) and using $C+E<1+CE$, we get $4A+CE+D+F+4G+I-4(A^3G^3CDEFI)^{\frac{1}{2}}>8.$ Left side is a decreasing function of $G$ as $G>1$ and~ $9A^3CDEFI>9A^3I(1-c-e)(1-d)(1-f)>9A^3I(1-\frac{a+i}{2})^2(1-\frac{i}{2})>4$. So we can replace $G$ by 1 to get $4A+CE+D+F+I-4(A^3CDEFI)^{\frac{1}{2}}>4.$ The left side is a decreasing function of $CE$ as $A^3DFI\geq ABGDFI=\frac{1}{CEH}\geq1$, so we replace $CE$ by $1-\frac{a+i}{2}+\frac{f}{2}$ and get that $\phi_{ord}(f,d,i,a)=4+\frac{7a}{2}+\frac{i}{2}-d-\frac{f}{2}-4(1+a)^{\frac{3}{2}}(1-\frac{a+i}{2}+\frac{f}{2})^{\frac{1}{2}}(1-d)^{\frac{1}{2}}(1-f)^{\frac{1}{2}}
(1+i)^{\frac{1}{2}}>0$, where $0\leq f<\frac{a+i}{2}-\frac{d}{2}$ (using 5.84), $0\leq d<\frac{i}{2}$ and $0<i\leq a$. Following Remark 3 we find that ${\rm max}~\phi_{ord}(f,d,i,a)\leq 0$ for $0<a\leq 0.3$, giving thereby a contradiction.\vspace{3mm}

\noindent{\bf Final Contradiction}

We have $A^4EFGHI>(1+a)^4(1-(e+h))(1-f)(1+i)>(1+a)^4(1-\frac{a+i}{2}+\frac{f}{2})(1-f)(1+i)=\eta(f)\geq{\rm min}\{\eta(0),\eta(\frac{a+i}{2})\}>2$, for $0.3<a\leq \frac{1}{2}$ and $0<i\leq a$. Applying AM-GM inequality to (4,1,1,2,1) we have $4A+4G+E+F+I-2(A^5G^3EFI)^{\frac{1}{2}}>9$. Now left side is a decreasing function of $G$ as $G>1$  and $9A^5EFI>9(1+a)^5(1+i)(1-(e+f))>9(1+a)^5(1+i)(1-\frac{3}{4}(a+i))>16$, for $0.3<a\leq 0.5$ and $0<i\leq a$. So we replace $G$ by 1 and get that $4A+4+E+F+I-2(A^5EFI)^{\frac{1}{2}}>9$. Now the left side is a decreasing function of $E$. We replace $E$ by $1-\frac{a+i}{2}+\frac{f}{2}$. So we have
$\eta(f)=2+\frac{7a}{2}+\frac{i}{2}-\frac{f}{2}-2(1+a)^{\frac{5}{2}}(1+i)^{\frac{1}{2}}(1-f)^{\frac{1}{2}}\left(1-\frac{a+i}{2}+\frac{f}{2}\right)^{\frac{1}{2}}>0.$\\
Now $\eta''(f)>0$ and $0\leq f<\frac{a+i}{2}$, so $\eta(f)\leq{\rm max}\{\eta(0),\eta(\frac{a+i}{2})\}$.
 But $\eta(0)$ and $\eta(\frac{a+i}{2})$ are functions in two variables $a$ and $i$ and can be seen to be at most zero for $0<a\leq 0.5$ and $0<i\leq a$, a contradiction.\vspace{3mm}\\
%
{\noindent \bf Proposition 40.} Case (72) i.e. $ A>1,~ B>1, ~C\leq 1, ~D\leq1,~ E\leq 1, ~F\leq 1,~ G \leq 1, ~H \leq 1, ~I > 1$
 does not arise.\vspace{2mm}

{\noindent \bf Proof.} Here  $ a\leq \frac{1}{2},~ b\leq \frac{1}{ 3}$.  Also
$~2F\geq1\geq G$,  for if $~F  < \frac{1}{2}$ then $~E \leq \frac{4F}{3} < \frac{2}{3}$ so  that $~EF < \frac{1}{3}$ which implies
$ ABCDGHI > 3$. But $ ~ABCDGHI\leq ABI\leq A^2B \leq \frac{9}{4}\times \frac{4}{3}~= 3$, a contradiction.  Similarly we have $2G\geq H$.
Using the weak inequalities
$(1,2,2,1,1,1,1),(1,2,1,2,1,1,1)$, $(1,2,2,2,1,1)$, ~$(1,2,2,1,2,1)$ and ~(2,2,1,1,1,1,1)  we get \vspace{-2mm}
\begin{equation}a-2c-2e-f-g-h+i>0,\vspace{-3mm}\end{equation}
\begin{equation}a-2c-d-2f-g-h+i>0,\vspace{-2mm}\end{equation}
\begin{equation}a-2c-2e-2g-h+i>0,\vspace{-2mm}\end{equation}
\begin{equation}a-2c-2e-f-2h+i>0,\vspace{-2mm}\end{equation}
\begin{equation}2b-2d-e-f-g-h+i>0.\end{equation}

 \noindent{\bf Claim (i)} $A < 1.313$

Suppose $A\geq 1.313$. We prove that $A^4EFGHI>2$.\\
From (5.87), (5.88), (5.89) and (5.90) respectively we get  $0\leq e< \frac{a+i}{2}-\frac{f+g+h}{2}, ~0\leq f <
\frac{a+i}{2}-\frac{g+h}{2}$, $0\leq g< \frac{a+i}{2}-\frac{h}{2}$ and $0\leq h<\frac{a+i}{2}$. Now $A^4EFGHI>A^4(1-\frac{a+i}{2}+\frac{f+g+h}{2})(1-f)(1-g)(1-h)(1+i)=\phi_{ord}(f,g,h,i,a),$ say. Following Remark 3 we find that ${\rm min}~\phi_{ord}(f,g,h,i,a)\geq 2$ for $0.313\leq a\leq \frac{1}{2}$ and $0<i\leq a$.
Hence $ A^4EFGHI>2$ for $A\geq 1.313$. So $(4,1,1,1,1,1)$ holds i.e. $4A-\frac{1}{2}A^{5}EFGHI +E+F+G+H+I>9.$ It gives $4a-e-f-g-h+i-\frac{1}{2}(1+a)^5(1-e)(1-f)(1-g)(1-h)(1+i)>0$. Coefficient of $e$ on left hand side is positive for $0\leq f<\frac{a+i}{2}-\frac{g+h}{2}$, $0\leq g<\frac{a+i}{2}-\frac{h}{2}$, $0\leq h<\frac{a+i}{2}$ and $a\geq0.21$. So replacing $e$ by $ \frac{a+i}{2}-\frac{f+g+h}{2}$ and simplifying we get\\
$\begin{array}{ll}
 \zeta_{ord}(f,g,h,i,a)=&\frac{7a}{2}-\frac{f+g+h}{2}+\frac{i}{2}-\frac{1}{2}\{A^{5}(1-\frac{a+i}{2}+\frac{f+g+h}{2})(1-f)\\&(1-g)(1-h)(1+i)\} > 0.
 \end{array}$\\ Following Remark 3 we find that ${\rm max}~\zeta_{ord}(f,g,h,i,a)\leq 0$ for $0<a\leq 0.5$ and $0<i\leq a$. Hence $A<1.313$.\vspace{3mm}

 \noindent{\bf Claim (ii)} $c+e+g >0.36(a+i)$ and $h<0.28(a+i)$

  Suppose that  $x=c+e+g \leq 0.36(a+i)$. Applying  AM-GM inequality to (2,2,2,2,1) we get
  $4A+4C+4E+4G+I-8A^{\frac{3}{4}}C^{\frac{3}{4}}E^{\frac{3}{4}}G^{\frac{3}{4}}I^{\frac{1}{4}}>9 $ which implies
  $8+4a-4x+i-8(1+a)^{\frac{3}{4}}(1-x)^{\frac{3}{4}}(1+i)^{\frac{1}{4}}>0$. Since the left side is an increasing  function of $x$, therefore  replacing $x$ by  $ 0.36(a+i)$ we
  get
$
8+2.56a-0.44i-8(1+a)^{\frac{3}{4}}(1-0.36(a+i))^{\frac{3}{4}}(1+i)^{\frac{1}{4}}>0.
$ But this is not true for $0<a<0.313$ and $0<i\leq a$. \vspace{1mm}\\
Now using (5.89) we get $h<a+i-2\times 0.36(a+i)<0.28(a+i)$.\vspace{3mm}

 \noindent{\bf Claim (iii)} $a < 0.29$

 Suppose that $a \geq 0.29$. We prove that $A^4EFGHI>2$. Working as in Claim (i) and using $h<0.28(a+i)$ in place of $h<\frac{a+i}{2}$ we find that ${\rm min}~\phi_{ord}(f,g,h,i,a)\geq 2$ for $0.29\leq a<0.313$ and $0<i\leq a$. Hence $A^4EFGHI>2$ for $A>1.29$. So (4,1,1,1,1,1) holds. Again working as in Claim (i) we get a contradiction. Hence $a<0.29$.\vspace{3mm}

\noindent{\bf Claim (iv)} $c+e+g >0.377(a+i)$ and $h<0.246(a+i)$

  Assume that  $x=c+e+g \leq 0.377(a+i)$. We proceed as in Claim (ii) and replace $x$ by  $ 0.377(a+i)$ to
  get a contradiction for $0<a<0.29$ and $0<i\leq a$. \vspace{1mm}\\
Now using (5.89) we get $h<0.246(a+i)$.\vspace{3mm}

\noindent{\bf Claim (v)} $c+e >0.26(a+i)$ and $f+g+h<0.48(a+i)$

 Suppose that  $y=c+e \leq 0.26(a+i)$. Also using (5.89) and (5.90) respectively we have $g<\frac{a+i}{2}-y-\frac{h}{2}$ and $h<\frac{a+i}{2}-y$. Applying  AM-GM inequality  to $(2,2,2,1,1,1)$ we have
  $4A+4C+4E-6ACEG^{\frac{1}{3}}H^{\frac{1}{3}}I^{\frac{1}{3}}+G+H+I>9 $ which implies
  $6+4a-4y-g-h+i-6(1+a)(1-y)(1-g)^{\frac{1}{3}}(1-h)^{\frac{1}{3}}(1+i)^{\frac{1}{3}}>0$, where $y=c+e$. The left side is an increasing  function of $g$ as $\frac{2(1+a)(1-y)(1-h)^{\frac{1}{3}}(1+i)^{\frac{1}{3}}}{(1-g)^{\frac{2}{3}}}>2(1+a)(1-0.26(a+i))(1-\frac{a+i}{2})^{\frac{1}{3}}(1+i)^{\frac{1}{3}}>1$. Therefore  replacing $g$ by  $ \frac{a+i}{2}-y-\frac{h}{2}$ we
  get\\
  $
\begin{array}{ll}
\varphi_{ord}(h,y,i,a)=&6+\frac{7a}{2}+\frac{i}{2}-3y-\frac{h}{2}-6\{(1+a)(1-y)(1-\frac{a+i}{2}+y+\frac{h}{2})^{\frac{1}{3}}\\&(1-h)^{\frac{1}{3}}(1+i)^{\frac{1}{3}}\}>0.
\end{array}
$\\ Following Remark 3 we find that ${\rm max}~\varphi_{ord}(h,y,i,a)\leq 0$ for $0\leq h<\frac{a+i}{2}-y$,  $0\leq y\leq 0.26(a+i)$, $0<i\leq a$ and $0<a<0.29$. Hence we have a contradiction. Therefore we must have $c+e >0.26(a+i)$. \vspace{1mm}\\
Now using (5.87) we get $f+g+h<(a+i)-2(c+e)<0.48(a+i).$\vspace{3mm}

\noindent{\bf Claim (vi)} $a < 0.2482$

Suppose $a\geq 0.2482$. We have $A^4EFGHI>(1+a)^4(1-\frac{a+i}{2}+\frac{f+g+h}{2})(1-(f+g+h))(1+i)=\phi(f+g+h)$. As $\phi''(f+g+h)<0$ and $0\leq f+g+h<0.48(a+i)$, we have $\phi(f+g+h)\geq{\rm min}\{\phi(0),\phi(0.48(a+i))\}$, which can be verified to be greater than 2 for $0<i\leq a$ and $0.2482\leq a< 0.29$. Hence we have $A^4EFGHI>2$. Now we get contradiction by working as in Claim (i).\vspace{3mm}

\noindent{\bf Claim (vii)} $c+e >0.324(a+i)$ and $f+g+h<0.352(a+i)$

 Using the refined bound on $a$, namely $0<a<0.2482$, and working as in Claim(v), we get the desired result.\vspace{3mm}

\noindent{\bf Claim(viii)} $b < 0.17$

Suppose that $b\geq 0.17$. We have $B^4AFGHI>(1.17)^4(1+a)(1+i)(1-0.352(a+i))>2$, for $0<i\leq a$ and $0.17\leq b\leq a<0.2482$. Therefore $(1,4,1,1,1,1)$ holds, which gives $\phi(z)=a+4b-z+i-\frac{1}{2}(1+b)^{5}(1+a)(1-z)(1+i) >
0,$ where $z=f+g+h$. The coefficient of $z$ on left side is positive and $z < 0.352(a+i)$. Therefore
$\phi(z)<\phi(0.352(a+i))=0.648a+4b+0.648i-\frac{1}{2}(1+b)^{5}(1+a)(1-0.352(a+i))(1+i)=\psi(i)$, say.
  As $\psi''(i)>0$, we have  $\psi(i)\leq {\rm max}\{\psi(0),\psi(a)\}$ which can be easily seen to be less than zero
  for $0.17 \leq b \leq a < 0.2482$. This gives a contradiction .\vspace{3mm}

\noindent{\bf Claim(ix)} $e+f+g+h>2b$ and $d<\frac{i}{2}$

 Suppose $e+f+g+h \leq 2b $. From (5.87) we have $2e+f+g+h<a+i$. Also  $ f+g+h < 0.352(a+i)$. Adding these two we get $ e+f+g+h <
 0.676(a+i)$. Now we get a
 contradiction from Lemma 7(vi) with $\gamma =e+f+g+h<\frac{3}{4}\delta$,  where  $\delta = a+i$ and $\gamma <2b$,
 $b=x_1<0.17$. \vspace{1mm}\\
 Now using (5.91) we have $2d<2b+i-(e+f+g+h)<i$, i.e. $d<\frac{i}{2}$.

\noindent{\bf Final contradiction}

Using $(3,1,1,1,1,1,1)$ we have $4A -A^{4}DEFGHI+D
 +E+F+G+H+I > 9$. Coefficient of $D$ is negative on the left side, as $A^4EFGHI>A^4(1-(e+f+g+h))(1+i)>(1+a)^4(1-0.676(a+i))(1+i)>1$, so we replace $D$ by $1-\frac{i}{2}$. This gives $1+4a+\frac{i}{2}-e-z -(1+a)^{4}(1-\frac{i}{2})(1+i)(1-e)(1-z) > 0,$ where $z=f+g+h$. As the coefficient of $e$ is
 positive for $z<0.352(a+i)$ and  we have $e < \frac{a+i}{2}-\frac{z}{2}$, so we can replace $e$ by $\frac{a+i}{2}-\frac{z}{2}$ to get $\eta(z)=1+\frac{7a}{2}-\frac{z}{2} -(1+a)^{4}(1-\frac{i}{2})(1+i)(1-\frac{a+i}{2}+\frac{z}{2})(1-z) > 0$. Now $\eta''(z)>0$ and $0\leq z<0.352(a+i)$. Therefore $\eta(z)\leq{\rm max}\{\eta(0),\eta(0.352(a+i))\}<0$, for $0<i\leq a$ and $0<a<0.2482$. This gives a contradiction.
 \vspace{3mm}

{\noindent \bf Proposition 41.} Case (71) i.e. $ A>1,~ B>1, ~C\leq 1, ~D\leq1,~ E\leq 1, ~F\leq 1,~ G \leq 1, ~H > 1, ~I > 1$
 does not arise.\vspace{2mm}

{\noindent \bf Proof.}Here  $ a\leq \frac{1}{2},~ b\leq \frac{1}{ 3}$. Using the weak inequalities
(1,2,2,1,1,1,1), ~$(1,2,2,2,1,1)$, ~$(2,2,1,1,1,1,1)$, ~(2,1,2,1,1,1,1), ~(2,1,2,2,1,1), ~$(2,2,2,1,1,1)$ ~and ~(1,2,1,2,1,1,1)  we get \vspace{-2mm}
\begin{equation}a-2c-2e-f-g+h+i>0,\vspace{-3mm}\end{equation}
\begin{equation}a-2c-2e-2g+h+i>0\vspace{-2mm}\end{equation}
\begin{equation}2b-2d-e-f-g+h+i>0,\vspace{-2mm}\end{equation}
\begin{equation}2b-c-2e-f-g+h+i>0\vspace{-2mm}\end{equation}
\begin{equation}2b-2d-2f-g+h+i>0\vspace{-2mm}\end{equation}
\begin{equation}a-2c-d-2f-g+h+i>0\end{equation}
Using $G>\frac{2}{3}E$ and (5.93) we have $1-e<\frac{3}{2}-\frac{3}{2}g$ and $e+g<\frac{a+h+i}{2}$, respectively. Adding these two we get  $g<\frac{1+a+h+i}{5}$. Similarly using $E>\frac{2C}{3}$ and (5.93) we get $e<\frac{1+a+h+i-2g}{5}$. \\

 \noindent{\bf Claim (i)} $A < 1.331$

Suppose $A\geq 1.331$. We prove first that $A^4EFGHI>2$.\\
If $A\geq 1.387$ then $A^4EFGHI=\frac{A^3}{BCD}\geq(1.387)^3\times \frac{3}{4}>2$.\\ Let now $0.331\leq a\leq 0.387$.
From (5.92)  we have $e<\frac{a+h+i}{2}-\frac{f+g}{2}$.
 Therefore $A^4EFGHI>A^4(1-\frac{a+h+i}{2}+\frac{f+g}{2})(1-f)(1-g)(1+h)(1+i)=\phi_{ord}(f,g,h,i,a),$ say. Following Remark 3 we find that ${\rm min}~\phi_{ord}(f,g,h,i,a)\geq 2$ for $f<\frac{a+h+i}{2}-\frac{g}{2}$ (from (5.97)), $g<\frac{1+a+h+i}{5}$,  $0<h\leq a$ and $0<i\leq a$ and $0.331\leq a\leq 0.387$. As $2E>1>F$, the inequalities (4,2,1,1,1) and (4,1,1,1,1,1) hold. We consider following cases:

 \noindent{\bf Case I}~~ $0.331\leq a<0.45$

 The inequality (4,2,1,1,1), using AM-GM  gives
  \begin{equation}2+4a-4e-g+h+i-2\sqrt{(1+a)^5(1-e)^3(1-g)(1+h)(1+i)}>0\end{equation}

  \par Now left side of (5.98) is an increasing function of $e$.  So we replace $e$ by $\frac{1+a+h+i-2g}{5}$  to get that $\psi_{ord}(g,h,i,a)=\frac{6}{5}+\frac{16a}{5}+\frac{h+i}{5}+\frac{3g}{5}-2(1+a)^{\frac{5}{2}}(\frac{4}{5}-\frac{a+h+i}{5}+\frac{2g}{5})^{\frac{3}{2}}(1-g)^{\frac{1}{2}}(1+h)^{\frac{1}{2}}
  (1+i)^{\frac{1}{2}}>0$. Following Remark 3 we find that ${\rm max}~\psi_{ord}(g,h,i,a)< 0$ for $0<g<\frac{1+a+h+i}{5}$, $0<h\leq a$, $0< i \leq a$ and  $0.331\leq a\leq 0.45$. This contradicts (5.98). We notice that $\psi_{ord}(g,h,i,a)<0$  for $0.058<g<\frac{1+a+h+i}{5}$  also if $a>0.45$.

\noindent{\bf Case II}~~ $0.45\leq a \leq 0.5$ and $g\leq 0.058$  \\
We have $G=1-g\geq 0.942$, ~$E\geq0.46873A>0.67$,  ~$F\geq\frac{3}{4}E>0.5$.\\
We use now (4,1,1,1,1,1). Following Remark 3 we find that $\chi_{ord}(G,F,E,H,I,A)$ $=4A-\frac{1}{2}A^5EFGHI+E+F+G+H+I<9$ for $0.942 \leq G\leq 1$,  $0.67 \leq E\leq 1$, $0.5 \leq F\leq 1$, $1<H\leq A$, $1<I\leq A$ and $1<A \leq 1.5$.

\noindent{\bf Claim (ii)} $c+e>0.195(a+h+i)$ and $g<0.305(a+h+i)$.

Suppose $c+e\leq 0.195(a+h+i)$. The inequality (2,2,2,1,1,1) with AM-GM gives
 $6+4a-4x-g+h+i-6(1+a)(1-x)(1-g)^{\frac{1}{3}}(1+h)^{\frac{1}{3}}(1+i)^{\frac{1}{3}}>0$, where $x=c+e$. Left side is an increasing function of $g$, So we replace $g$ by $\frac{a+h+i}{2}-(c+e)$, using (5.93). But then we find that $\phi_{ord}(x, h,i,a)=6+\frac{7a}{2}+\frac{h+i}{2}-3x-6(1+a)(1-x)(1-\frac{a+h+i}{2}+x)^{\frac{1}{3}}(1+h)^{\frac{1}{3}}(1+i)^{\frac{1}{3}}<0$ for $0\leq x\leq 0.195(a+h+i)$, $0<h\leq a$, $0< i \leq a$ and $0<a<0.331$.
  Hence we have $c+e>0.195(a+h+i)$.\vspace{1mm}\\
 Now using (41.2) get $g<0.305(a+h+i)$.\vspace{3mm}

\noindent{\bf Claim (iii)} $A<1.289$

 We proceed as in Claim (i) and use $0\leq g<0.305(a+h+i)$ in place of $0<g<\frac{1+a+h+i}{5}$ and find that $A^4EFGHI>2$, for $A\geq 1.289$. Now (4,2,1,1,1) gives a contradiction.\vspace{3mm}

\noindent{\bf Claim (iv)} $c+e>0.233(a+h+i)$ and $g<0.267(a+h+i)$.

If $c+e\leq 0.233(a+h+i)$, we proceed as in Claim (ii) to get a contradiction for $0<a<0.289$.\vspace{1mm}\\
 Now using (5.93) we have $g<0.267(a+h+i)$.\vspace{3mm}

 \noindent{\bf Claim (v)} $A<1.281$

   We proceed as in Claims (i)and (iii) and use $g<0.267(a+h+i)$ to prove $A^4EFGHI>2$ for $a\geq 0.281$. Then (4,2,1,1,1) gives a contradiction.\vspace{3mm}

 \noindent{\bf Claim (vi)} $c+e>0.24(a+h+i)$, ~$f+g<0.52(a+h+i)$  ~and  ~$g<0.26(a+h+i)$.

 If $c+e\leq 0.24(a+h+i)$, we proceed as in Claim (ii) to get a contradiction.\vspace{1mm}\\
 Now using (5.92) and (5.93) respectively we get ~$f+g<0.52(a+h+i)$  ~and  ~$g<0.26(a+h+i)$.\vspace{3mm}

\noindent{\bf Claim (vii)} $A>1.2$

Suppose $A\leq1.2$. Using (5.94)  we have $d<b+\frac{h+i}{2}-\frac{e+f+g}{2}$. As (1,2,1,1,1,1,1,1) holds, we have $A+4B+D+E+F+G+H+I-2B^3ADEFGHI>9$. Left side is a decreasing function of $D$, as $BAEFGHI=\frac{1}{CD}\geq 1$. So we replace $D$ by $1-(b+\frac{h+i}{2}-\frac{e+f+g}{2})$ and get $\eta_{ord}(e,y,h,i, a,b)=2+a+3b+\frac{h+i}{2}-\frac{e+y}{2}-2(1+b)^3(1+a)(1+h)(1+i)(1-b-\frac{h+i}{2}+\frac{e+y}{2})(1-e)(1-y)>0$ where $y=f+g$. Following Remark 3, we find that $\eta_{ord}(e,y,h,i, a,b)<0$ for $0<e<b+\frac{h+i}{2}-\frac{y}{2}$ (using (5.95)), $0<y<0.52(a+h+i)$, $0<h\leq a$, $0< i \leq a$ and $0<b\leq a\leq0.2$.
This gives a contradiction. Hence we have $a>0.2$.\vspace{3mm}

\noindent{\bf Claim (viii)} $B<1.168$

Suppose $B\geq 1.168$, then $B^4AFGHI>(1+0.168)^4(1+a)(1+h)(1+i)(1-(f+g))>2 $ for $f+g<0.52(a+h+i)$, $0<h\leq a$, $0< i \leq a$ and $0.2<a<0.281$. So (1,4,1,1,1,1) holds. It gives $\eta_{ord}(f+g,h,i,a,b)=a+b-(f+g)+h+i-\frac{1}{2}(1+b)^5(1+a)(1+h)(1+i)(1-(f+g))>0$. Left side is an increasing function of $f+g$ for $b>0.11$, so we replace $f+g$ by $0.52(a+h+i)$ and find that $\eta_{ord}(f+g,h,i,a,b)<0$ for $0<b\leq a\leq 0.281$. Hence we have $B<1.168$.\vspace{3mm}

\noindent{\bf Claim (ix)} $h+i>1.415a$

Suppose $h+i\leq 1.415a$. We proceed as in Claim (vii). Here we have
 $0<h\leq 1.1415a-i$. Following Remark 3, we find that $\eta_{ord}(e,y,h,i, a,b)<0$ for $0<e<b+\frac{h+i}{2}-\frac{y}{2}$, $0<y<0.52(a+h+i)$, $0<h\leq 1.1415a-i$, $0< i \leq a$ and for $0<b<0.168$ and $0.2<a<0.281$. This gives a contradiction. Hence we have $h+i>1.415a$.\vspace{2mm}

\noindent{\bf Claim (x)} $B<1.145$

Suppose $B\geq 1.145$, then $B^4AFGHI>2$ for for $f+g<0.52(a+h+i)$, $1.1415a-i<h\leq a$, $0< i \leq a$ and $0.2<a<0.281$. So (1,4,1,1,1,1) holds and we work as in Claim (viii) to get a contradiction.\vspace{3mm}

\noindent{\bf Claim (xi)} $e+f+g>1.59b+0.4(h+i)$ and $d<\frac{0.41b}{2}+0.3(h+i)$

Suppose $e+f+g\leq 1.59b+0.4(h+i)$. As (1,3,1,1,1,1,1) holds, we have  $a+4b-(e+f+g)+h+i-(1+b)^4(1+a)(1+h)(1+i)(1-(e+f+g))>0$.  We can replace $e+f+g$ by $1.59b+0.4(h+i)$ to get
$\phi_{ord}(h,i,a,b)=1+a+2.41b+0.6(h+i)-(1+b)^4(1+a)(1+h)(1+i)(1-1.59b-0.4(h+i))>0$.
Following Remark 3, we see that $\phi_{ord}(h,i,a,b)<0$ for $1.415a-i<h<a$, $0<i \leq a$ and for $0.2<a<0.281$, ~$0<b<0.145$. Hence we have $e+f+g>1.59b+0.4(h+i)$.\vspace{1mm}\\
Now  using (5.94) we get  $d<\frac{0.41b}{2}+0.3(h+i)$.\vspace{3mm}

\noindent{\bf Claim (xii)} $A\geq 1.25$ ~and~ $g>0.16(a+h+i)$

Suppose $A<1.25$. We work as in claim (vii), first  replace $E$ by $1-(b+\frac{h+i}{2}-\frac{f+g}{2})$  to get $\eta_{ord}(f,d,g,h,i,a,b)=2+a+3b+\frac{h+i}{2}-d-\frac{f+g}{2}-2(1+b)^3(1+a)(1+h)(1+i)(1-d)(1-f)(1-g)(1-b-\frac{h+i}{2}+\frac{f+g}{2})>0$. Following Remark 3, we find that $\eta_{ord}(f,d,g,h,i,a,b)<0$ for $0<f< b+\frac{h+i}{2}-\frac{g}{2}-d$ (using (5.96)), $0\leq d<\frac{0.41b}{2}+0.3(h+i)$, $0 \leq g<0.26(a+h+i)$, $1.415a-i<h<a$, $0<i \leq a$  and for $0.2<a<0.25$, $0<b<0.145$. Hence we must have $a\geq 0.25$.\vspace{2mm}

Now if we use $g\leq 0.16(a+h+i)$ instead of $g<0.26(a+h+i)$,, then we get $\eta_{ord}(f,d,g,h,i,a,b)<0$ in full range of $a$ and $b$, i.e. in $0.2<a<0.281$ and $0<b<0.145$. Hence we have $g>0.16(a+h+i)$.\vspace{3mm}

\noindent{\bf Claim (xiii)} $c>0.04(a+h+i)$

Suppose $c\leq 0.04(a+h+i).$ As (2,2,1,1,1,1,1) holds, applying AM-GM we get $4A+4C+E+F+G+H+I-4\sqrt{A^3C^3EFGHI}>9.$ Now left side is a decreasing function of $E$, so we replace $E$ by $1-\frac{a+h+i}{2}+\frac{f+g}{2}+c$ (using (5.92)) to get $\phi_{ord}(f,g,c,h,i,a)=4+\frac{7a}{2}-3c+\frac{h+i}{2}-\frac{f+g}{2}-4(1+a)^{\frac{3}{2}}(1-c)^{\frac{3}{2}}(1-f)^{\frac{1}{2}}(1-g)^{\frac{1}{2}}(1+h)^{\frac{1}{2}}
(1+i)^{\frac{1}{2}}(1-\frac{a+h+i}{2}+\frac{f+g}{2}+c)^{\frac{1}{2}}>0$. Following Remark 3, we find that $\phi_{ord}(f,g,c,h,i,a)<0$ for $0\leq f<\frac{a+h+i}{2}-\frac{g}{2}-c$ (using (5.97)), $0.16(a+h+i)<g<0.26(a+h+i)$, $0<c< 0.04(a+h+i)$,  $0<h\leq a$, ~$0<i\leq a$ and  ~$0.25<a<0.281$.
 Hence we must have $c>0.04(a+h+i)$.\vspace{3mm}

\noindent{\bf Final Contradiction}

Using $c>0.04(a+h+i)$, (5.92) and (5.97) we have  $e<0.46(a+h+i)-\frac{f}{2}-\frac{g}{2}$ and $f<0.46(a+h+i)-\frac{g}{2}$.  We find that
$A^4EFGHI>(1+a)^4(1+h)(1+i)(1-0.46(a+h+i)+\frac{f}{2}+\frac{g}{2})(1-f)(1-g)=\eta_{ord}(f,g,h,i,a)>2$ for  $0\leq f<0.46(a+h+i)-\frac{g}{2}$, $0.16(a+h+i)<g<0.26(a+h+i)$, $1.415a-i<h<a$, $0<i \leq a$  and for $0.25<a<0.281$.  So (4,2,1,1,1) holds. Now we work as in Claim (i) to get a contradiction.\vspace{3mm}
\section{Most Difficult Cases.}
In this section we will discuss the cases 6, 14, 15, 32 and 33.
We shall use  the following inequalities repeatedly, whenever applicable. \vspace{2mm}

\noindent 1. (3,4,2) with AM-GM, $\phi_1=4A+4D+4H-3H(D^5A^4)^{\frac{1}{3}}-9>0.$

\noindent 2. (1,4,4) with AM-GM, $\phi_2=A+4B+4F-(B^5F^5A)^{\frac{1}{2}}-9>0.$

\noindent 3. (1,4,1,1,1,1), $\phi_3=A+4B-\frac{1}{2}B^5AFGHI+F+G+H+I-9>0.$

\noindent 4. (4,2,2,1) with AM-GM, $\phi_4=4A+4E+4G+I-3EG(2A^5I)^{\frac{1}{3}}-9>0$.

\noindent 5. (4,1,2,2) with AM-GM, $\phi_5=4A+E+4F+4H-3FH(2A^5E)^{\frac{1}{3}}-9>0$.

\noindent 6. (1,4,2,2) with AM-GM, $\phi_6=A+4B+4F+4H-3FH(2B^5A)^{\frac{1}{3}}-9>0$.

\noindent 7. (1,4,1,2,1) with AM-GM, $\phi_7=A+4B+F+4G+I-2(B^5G^3AFI)^{\frac{1}{2}}-9>0$.

\noindent 8. (1,4,3,1) with AM-GM, $\phi_8=A+4B+4F+I-2(0.5B^5F^4AI)^{\frac{1}{2}}-9>0$.

\noindent 9. (2,1,4,2) with AM-GM, $\phi_9=4A-\frac{2A^2}{B}+C+4D+4H-2(D^5H^3ABC)^{\frac{1}{2}}-9>0$.

\noindent 10. (2,2,4,1), $\phi_{10}=4A-\frac{2A^2}{B}+4C-\frac{2C^2}{D}+4E-\frac{1}{2}E^5ABCDI+I-9>0$.\\
We note that $\phi_{10}$ is a decreasing function of $A$, a decreasing function of $C$ if $C>D$, a decreasing function of $E$,
 whenever $E^4ABCDI>2$ and a decreasing function of $I$ whenever $E^5ABCD>2$.\vspace{2mm}

\noindent 11. (2,4,2,1) with AM-GM, $\phi_{11}=4A-\frac{2A^2}{B}+4C+4G-2(C^5G^3ABI)^{\frac{1}{2}}+I-9>0$. We note that
 $\phi_{11}$ is a decreasing function of $A$. \\

\noindent 12. (1,1,4,2,1) with AM-GM, $\phi_{12}=A+B+4C+4G-2(C^5G^3ABI)^{\frac{1}{2}}+I-9>0$.\vspace{2mm}

\noindent 13. (1,2,4,2) with AM-GM,\\ $\phi_{13}=A+4B-\frac{2B^2}{C}+4D+4H-2(D^5H^3ABC)^{\frac{1}{2}}-9>0$.\\
Note that $\phi_{13}$  is a decreasing function of $B$ whenever $B>C$, a decreasing function of $H$ whenever $D^5ABCH>1$.\vspace{2mm}

\noindent 14. (1,2,1,1,4)  : $\phi_{14}=A+4B-\frac{2B^2}{C}+D+E+4F-\frac{1}{2}F^5ABCDE-9>0$.\\
Note that $\phi_{14}$ is linear in $A,~D,~E$ and is a decreasing function of $B$ whenever $B>C$.\vspace{2mm}

\noindent 15. (2,2,1,4), $\phi_{15}=4A-\frac{2A^2}{B}+4C-\frac{2C^2}{D}+E+4F-\frac{1}{2}F^5ABCDE-9>0$.\\
Note that $\phi_{15}$ is a decreasing function of $A$ always, a decreasing function of $C$ whenever $C>D$ and is  a decreasing
function of $F$ and $E$ whenever $F^4ABCD>2.$\vspace{2mm}

\noindent 16. (1,2,4,1,1), $\phi_{16}=A+4B-\frac{2B^2}{C}+4D-\frac{1}{2}D^5ABCHI +H+I-9>0$.\\
Note that $\phi_{16}$ is  a decreasing function of $B$ whenever $B>C$ and is  a decreasing function of $H$
 and $I$ whenever $D^5ABCI>2$ or $D^5ABCH>2$.\vspace{2mm}

\noindent 17. (2,1,4,1,1),  $\phi_{17}=4A-\frac{2A^2}{B}+C+4D-\frac{1}{2}D^5ABCHI+H+I-9>0$.\\
Note that $\phi_{17}$ is linear function of $C, H$ and $I$ and a decreasing function of $A$.\vspace{2mm}

\noindent 18. (1,2,1,4,1), $\phi_{18}=A+4B-\frac{2B^2}{C}+D+4E-\frac{1}{2}E^5ABCDI+I-9>0$.\\
We note that $\phi_{18}$ is a linear function of $A$, $D$ and $I$ and a decreasing function of $B$ if $B>C$.\vspace{2mm}

\noindent 19. (2,2,2,2,1) with AM-GM, \\
$\phi_{19}=4A-\frac{2A^2}{B}+4C+4E+4G+I-6CEG(ABI)^{\frac{1}{3}}-9>0$.\\
We note that $\phi_{19}$ is a decreasing  function of $A$. i.e.\\
$\phi_{19}^{(1)}=8+4a-\frac{2(1+a)^2}{1+b}+4c-4(e+g)+i-6(1+c)(1-e-g)(ABI)^{\frac{1}{3}}>0$.\\
$\phi_{19}^{(2)}=4A-\frac{2A^2}{B}+4C-\frac{2C^2}{D}+4E+4G+I-4(E^3G^3ABCDI)^{\frac{1}{2}}-9>0$.\\\vspace{2mm}
$\phi_{19}^{(3)}=4A+4C-\frac{2C^2}{D}+4E+4G+I-6(A^3E^3G^3CDI)^{\frac{1}{2}}-9>0$.\\\vspace{2mm}

\noindent 20. (2,1,1,1,4), $\phi_{20}=4A-\frac{2A^2}{B}+C+D+E+4F-\frac{1}{2}F^5ABCDE-9>0$.
Note that $\phi_{20}$ is a decreasing function of $A$, and is linear in $C,D,E$.\vspace{2mm}

\noindent 21. (1,1,4,1,1,1), $\phi_{21}=A+B+4C-\frac{1}{2}C^5ABGHI+G+H+I-9>0.$\vspace{2mm}

\noindent 22. (1,2,2,2,1,1) with AM-GM, \\$\phi_{22}=A+4B-\frac{2B^2}{C}+4D+4F+H+I-4(D^3F^3ABCHI)^{\frac{1}{2}}-9>0$\\
i.e. $\phi_{22}^{(1)}= 6+a+4b-\frac{2(1+b)^2}{1+c}+4d-4f-h+i-4\{(1+d)^3(1-f)^3(1+a)$ $~~~~~~~~~~~~~~~~~~~~~~~~~~~~~~~~~~~~~~~
~~~~~~~~~~~~~~~~~~~~~(1+b)(1+c)(1+i)(1-h)\}^{\frac{1}{2}}>0$.\\
$\phi_{22}^{(2)}=A+4B+4D+4F+H+I-6BDF(AHI)^{\frac{1}{3}}-9>0$.\vspace{2mm}\\
$\phi_{22}^{(3)}=A+4B-\frac{2B^2}{C}+4D-\frac{2D^2}{E}+4F+H+I-2(F^3ABCDHI)-9>0$\vspace{2mm}

\noindent 23. (2,2,2,1,1,1) with AM-GM,\\ $\phi_{23}=4A-\frac{2A^2}{B}+4C+4E+G+H+I-4(C^3E^3ABGHI)^{\frac{1}{2}}-9>0$.\\
 $\phi_{23}^{(1)}=4A+4C+4E+G+H+I-6ACE(GHI)^{\frac{1}{3}}-9>0$.\vspace{2mm}

\noindent 24. (2,3,1,1,1,1),  $\phi_{24}=4A-\frac{2A^2}{B}+4C-C^4ABFGHI+F+G+H+I-9>0$.\\
$\phi_{24}^{(1)}=2B+4C-C^4ABFGHI+F+G+H+I-9>0$.
\vspace{2mm}

\noindent 25. (3,3,1,1,1) with AM-GM, $\phi_{25}=4A+4D+G+H+I-2(A^4D^4GHI)^{\frac{1}{2}}-9>0$. \vspace{2mm}

\noindent 26. (4,4,1) with AM-GM, $\phi_{26}=4A+4E+I-(A^5E^5I)^{\frac{1}{2}}-9>0.$\vspace{2mm}

\noindent 27. (4,1,2,1,1) with AM-GM, $\phi_{27}=4A+E+4F+H+I-2(A^5F^3EHI)^{\frac{1}{2}}-9>0$.\vspace{2mm}

\noindent 28. (4,3,1,1) with AM-GM, $\phi_{28}=4A+4E+H+I-2(0.5A^5E^4HI)^{\frac{1}{2}}-9>0$.\vspace{2mm}

\noindent 29. (1,4,2,1,1) with AM-GM : $\phi_{29}=A+4B+4F+H+I-2(B^5F^3AHI)^{\frac{1}{2}}-9>0$.\vspace{2mm}

\noindent 30. (1,3,4,1) with AM-GM, $\phi_{30}=A+4B+4E+I-2(0.5E^5B^4AI)^{\frac{1}{2}}-9>0$.\vspace{2mm}

\noindent 31. (1,2,2,1,1,1,1) with AM-GM, \\$\phi_{31}=A+4B+4D+F+G+H+I-4(B^3D^3AFGHI)^{\frac{1}{2}}-9>0$.\vspace{2mm}

\noindent 32. (2,4,1,1,1), $\phi_{32}=2B+4C-\frac{1}{2}C^5ABGHI+G+H+I-9>0$.\\
Note that $\phi_{32}$  is linear in $G,H,I$ and also in $B$.\\
$\phi_{32}^{(1)}=4A-\frac{2A^2}{B}+4C-\frac{1}{2}C^5ABGHI+G+H+I-9>0$.
\vspace{2mm}

\noindent 33. (4,1,1,1,1,1), $\phi_{33}=4A-\frac{1}{2}A^5EFGHI+E+F+G+H+I-9>0$.\vspace{2mm}

\noindent 34. (1,2,1,2,1,1,1) with AM-GM, \\$\phi_{34}=A+4B+D+4E+G+H+I-4(B^3E^3ADGHI)^{\frac{1}{2}}-9>0$.\vspace{2mm}

\noindent 35. (1,3,2,1,1,1) with AM-GM, $\phi_{35}=A+4B+4E+G+H+I-2(2B^4E^3AGHI)^{\frac{1}{2}}-9>0$.\vspace{2mm}

\noindent 36. (1,1,1,4,1,1)~$\phi_{36}=A+B+C+4D-\frac{1}{2}D^5ABCHI+H+I-9>0$.\vspace{2mm}

\noindent 37. (4,2,1,1,1) with AM-GM : $\phi_{37}=4A+4E+G+H+I-2(A^5E^3GHI)^{\frac{1}{2}}-9>0$.\vspace{2mm}

\noindent 38. (3,4,1,1) with AM-GM, $\phi_{38}=4A+4D+H+I-2(0.5A^4D^5HI)^{\frac{1}{2}}-9>0$.\vspace{2mm}

\noindent 39. (3,1,3,1,1) with AM-GM, $\phi_{39}=4A+D+4E+H+I-2(A^4E^4DHI)^{\frac{1}{2}}-9>0$.\vspace{2mm}

\noindent 40. (1,2,2,2,2) with AM-GM, \\
$\phi_{40}^{(1)}=A+4B+4D+4F+4H-8(B^3D^3F^3H^3A)^{\frac{1}{4}}-9>0$.\vspace{2mm}\\
$\phi_{40}^{(2)}=A+4B+4D-\frac{2D^2}{E}+4F+4H-6(B^3F^3H^3ADE)^{\frac{1}{3}}-9>0$.\vspace{2mm}\\
$\phi_{40}^{(3)}=A+4B-\frac{2B^2}{C}+4D-\frac{2D^2}{E}+4F+4H-4(F^3H^3ABCDE)^{\frac{1}{2}}-9>0$.\vspace{2mm}

\noindent 41. (2,1,1,4,1) with AM-GM, $\phi_{41}=4A-\frac{2A^2}{B}+C+D+4E-\frac{1}{2}E^5ABCDI+I-9>0$.
\vspace{2mm}

\noindent 42. (2,4,1,2) with AM-GM, $\phi_{42}=4A-\frac{2A^2}{B}+4C-\frac{1}{2}C^5ABGHI+G+4H-\frac{2H^2}{I}-9>0$.
\vspace{2mm}

\noindent 43. (2,4,3) with AM-GM, $\phi_{43}=4A-\frac{2A^2}{B}+4C+4G-\sqrt{2C^5G^4AB}-9>0$.
\vspace{2mm}

\noindent 44. (1,2,1,1,2,2) with AM-GM, $\phi_{44}=A+4B-\frac{2B^2}{C}+D+E+4F-\frac{2F^2}{G}+4H-\frac{2H^2}{I}>0$.
\vspace{2mm}

\noindent 45. (1,2,1,2,2,1) with AM-GM, $\phi_{45}=A+4B-\frac{2B^2}{C}+D+4E-\frac{2E^2}{F}+4G-\frac{2G^2}{H}+I>0$.
\vspace{2mm}

\noindent ($1^*$) (2,1,$6^*$), $\phi_1^*=4A-\frac{2A^2}{B}+C+6(\frac{1}{ABC})^{1/6}-9>0$.\vspace{2mm}

\noindent ($2^*$) (1,2,$6^*$), $\phi_2^*=A+4B-\frac{2B^2}{C}+6(\frac{1}{ABC})^{1/6}-9>0$.\vspace{2mm}

\noindent ($3^*$) (2,$7^*$), $\phi_3^*=4A-\frac{2A^2}{B}+7(\frac{1}{AB})^{1/7}-9>0$.\vspace{2mm}

\noindent ($4^*$) (4,$5^*$), $\phi_4^*=4A-\frac{1}{2}A^5x+5(x)^{1/5}-9>0$, where $x=EFGHI=\frac{1}{ABCD}$.\vspace{2mm}

\noindent ($5^*$) (2,$6^*$,1),  $\phi_5^*=4A-\frac{2A^2}{B}+6(\frac{1}{ABI})^{1/6}+I-9>0$.\vspace{2mm}

\noindent ($6^*$)~~(2,$5^*$,1,1),  $\phi_6^*=4A-\frac{2A^2}{B}+5(\frac{1}{ABHI})^{\frac{1}{5}}+H+I-9>0$.\vspace{2mm}

\noindent ($7^*$)~~(2,1,$5^*$,1),  $\phi_7^*=4A-\frac{2A^2}{B}+C+5(\frac{1}{ABCI})^{\frac{1}{5}}+I-9>0$.\vspace{2mm}

 \noindent ($8^*$) ($4^*$,4,1), $\phi_8^*=4(x)^{1/4}+4E-\frac{1}{2}E^5xI+I-9>0$, where $x=ABCD$.\vspace{2mm}

 \noindent ($9^*$)~~(1,2,$5^*$,1),  $\phi_9^*=A+4B-\frac{2B^2}{C}+5(\frac{1}{ABCI})^{\frac{1}{5}}+I-9>0$.\vspace{2mm}

 \noindent ($10^*$)~~(2,2,$5^*$),  $\phi_{10}^*=4A-\frac{2A^2}{B}+4C-\frac{2C^2}{D}+5(\frac{1}{ABCD})^{\frac{1}{5}}-9>0$.\vspace{2mm}

  \noindent ($11^*$)~~(1,2,1,$5^*$),  $\phi_{11}^*=A+4B-\frac{2B^2}{C}+D+5(\frac{1}{ABCD})^{\frac{1}{5}}-9>0$.\vspace{2mm}

   \noindent ($12^*$)~~($5^*,4$),  $\phi_{12}^*=5(x)^{1/5}+4F-\frac{1}{2}F^5x-9>0$, where $x=ABCDE$.\vspace{2mm}
\\
{\noindent \bf Proposition 42.} Case (33) i.e. $ A>1,~ B>1, ~C> 1, ~D\leq 1,~ E\leq 1, ~F\leq 1,~ G \leq 1, ~H \leq 1, ~I >1$
 does not arise.\vspace{2mm}

{\noindent \bf Proof.} Here  $a\leq 1, ~ b\leq \frac{1}{2},~ c\leq \frac{1}{ 3}$.\vspace{2mm}
Using the weak inequalities $(2,2,2,2,1)$, $(2,2,1,2,1,1),(2,2,2,1,1,1)$, $(1,2,2,1,1,1,1)$ we have
\begin{equation}2b-2d-2f-2h+i>0\vspace{-3mm}\end{equation}
\begin{equation}2b-2d-e-2g-h+i>0\vspace{-2mm}\end{equation}
\begin{equation}2b-2d-2f-g-h+i>0\vspace{-2mm}\end{equation}
\begin{equation}a+2c-2e-f-g-h+i>0\end{equation}

\noindent In some of cases or claims that follow, we will need to prove $\phi_{29}<0$ for different bounds on $F,~H$. We note that
$\phi_{29}$ is a decreasing function  of $F$. If we have  $F \geq \mu_2$, where $\mu_2$ a real number we use $H>\frac{2}{3}F$
(provided it is a decreasing function of $H$). We can also use $H\geq \frac{2}{3}F$ and $F>0.46873B$ to get \vspace{-2mm}
$$\phi_{29} <\phi_{29}^{(1)}=A+4B+\frac{14}{3}\times 0.46873B+I-2(0.46873)^2(\frac{2}{3}B^9AI)^{1/2}-9. \vspace{-2mm}$$

 When   $F<\mu_2$. The weak inequality (2,2,2,2,1) gives $2B+2D+2F+2H+I>9$ which further gives $F> \frac{7-I}{2}-B-H$.
 We also have $H>\frac{7-I}{2}-B-F> \frac{7-I}{2}-B-\mu_2$. \\ Sometimes we use the maximum of the two lower bounds on $F$
 namely $$F>\left\{\begin{array}{lll}0.46873B & {\rm if }& H\geq \frac{7-I}{2}-1.46873B\\\frac{7-I}{2}-B-H
& {\rm if }& \frac{7-I}{2}-B-\mu_2<H<\frac{7-I}{2}-1.46873B\end{array}\right. $$
In the first case we replace $F$ by $0.46873B$ and $H$ by $\frac{7-I}{2}-1.46873B$ (whenever $\phi_{29}$ is decreasing function
of $H$) to get $\phi_{29}\leq$\vspace{2mm} \\$ \phi_{29}^{(2)}=A+3(1.46873)B-2(0.46873)^{\frac{3}{2}}B^4 (AI)^{\frac{3}{2}}
(\frac{7-I}{2}-1.46873B)^{\frac{1}{2}}-\frac{7-I}{2}$\vspace{2mm} \\
In the second case we replace $F$ by $\frac{7-I}{2}-B-H$ to get $\phi_{29}\leq
\phi_{29}^{(3)}$ and then consider it as a function of $H$, where $\frac{7-I}{2}-B-\mu_2<H<\frac{7-I}{2}-1.46873B$. \vspace{2mm}

\noindent{\bf Claim (i)} $A< 1.685$\vspace{2mm}

Suppose $A\geq 1.685$, then $\frac{A^3}{BCD}>2$. When $\frac{E^3}{FGH}>2$ we find $\phi_{26}<0$ for $1<I\leq A$, $E>0.46873A$
and $1.685\leq  A<2$. Let now $\frac{E^3}{FGH}\leq 2$. If $G\geq 0.6082$, we find $\phi_4<0$ for $I>1$, $E>0.46873A$ and $1.685\leq A<2$.
 So we can take $G<0.6082$ which gives $H\geq \frac{E^3}{2\times 0.6082F}$. If $1>F\geq 0.691$, we find  $\phi_{27}<0$ for $1<I\leq A$,
 $E>0.46873A$ and $1.685\leq  A<2$. So we can take $F<0.691$ which implies $E^2>FG$. But then $\phi_{28}<0$ for
 $H\geq \frac{E^3}{2\times 0.6082\times 0.691}$, $1<I\leq A$, $E>0.46873A$ and $1.685\leq  A<2$.\vspace{3mm}

\noindent{\bf Claim (ii)} $C< 1.26$\vspace{2mm}

Suppose $C\geq 1.26$, then $\frac{C^3}{DEF}>2$ and we find that  $\phi_{11}<0$ for $G>0.46873C$, $1<I\leq A$, $A\geq B$,
$1<B\leq \frac{3}{2}$, $1.26<C\leq \frac{4}{3}$.\vspace{3mm}

\noindent{\bf Claim (iii)} $B< 1.375$\vspace{2mm}

Suppose $B\geq 1.375$, then $\frac{B^3}{CDE}>2$. We note that $\phi_{29}$ is a decreasing function  of $F$. If $F \geq 0.77$,
 we find that  $\phi_{29}<0$ for $H>\frac{2}{3}F$, $F \geq 0.77$, $1<I\leq A$, $B\leq A<1.685$, $1.375<B\leq \frac{3}{2}$.\vspace{2mm}

When
 $F<0.77=\mu_2$. We find that $\phi_{29}^{(2)}<0, ~\phi_{29}^{(3)}<0$  for $1<I\leq A$, $B\leq A<1.685$, $1.375<B\leq \frac{3}{2}$.\vspace{2mm}

\textbf{We divide the range of $A$ into several
 subintervals and  arrive at a contradiction in each.}\vspace{2mm}

\noindent \textbf{Case I: 1.64 $\leq A < $ 1.685}\vspace{2mm}

\noindent If $B\leq 1.327$, we find that $\phi_5^*<0$, so can take $B>1.327$. Here $\frac{A^3}{BCD}>2$.\vspace{2mm}

If $G\geq 0.625$, we find $\phi_4<0$ for $1<I\leq A$, $E>0.46873A$ and $1.64\leq A<1.685$. So we can take $G<0.625$.
If $F\geq 0.77$, we find  $\phi_{27}<0$ for $H\geq \frac{2}{3}F$, $1<I\leq A$, $E>0.46873A$ and $1.64\leq  A<1.685$.
So we can take $F<0.77=\mu_2$ which implies $E^2>(0.46873A)^2>FG$.\vspace{2mm}

Next we find that $\phi_{28}<0$ for $1<I\leq A$, $H\geq 0.615$, $E>0.46873A$ and $1.64\leq  A<1.685$. Also $\phi_{28}<0$ for
$1<I\leq A$, $E\geq 0.888$, $H>\frac{1}{2}E$ and $1.64\leq  A<1.685$. So we can assume that $E<0.888$ and $H<0.615$. Then
$\frac{B^3}{CDE}>2$ and $\phi_{29}<0$ for $F\geq \frac{7-I}{2}-B-H$, $\frac{7-I}{2}-B-\mu_2<H<0.615$, $1<I\leq A$, $1.64\leq A<1.685$,
$1.327<B\leq 1.375$. This gives a contradiction.\vspace{2mm}

\noindent \textbf{Case II: 1.58 $\leq A < $ 1.64}\vspace{2mm}

 If $B\leq 1.304$, we find that $\phi_5^*<0$, so can take $B>1.304$. Here $\frac{A^3}{BCD}>2$. \\
Working as in Case I
 we can take $G<0.642$, $F<0.788$ which implies $E^2>(0.46873A)^2>FG$. Further considering $\phi_{28}$, can take
 $E<0.888$ and $H<0.704$. \vspace{2mm}

Suppose $B\geq 1.345$, then $\frac{B^3}{CDE}>2$. We work as in Claim (iii) and find that  $\phi_{29}<0$ for $H>\frac{2}{3}F$,
$F \geq 0.785$, $1<I\leq A$, $1.58\leq A<1.64$, $1.345<B\leq 1.375$. When $F<0.785=\mu_2$, we find that $\phi_{29}^{(2)}<0, ~\phi_{29}^{(3)}<0$
 for $1<I\leq A$, $1.58\leq A<1.64$, $1.345<B\leq 1.375$. Therefore we can take $B<1.345$.\vspace{2mm}

Suppose $C\geq 1.1186$, then $\frac{C^3}{DEF}>2$ and we find that  $\phi_{11}<0$ for $G>0.46873C$, $1<I\leq A$, $A\geq 1.58$,
$1.304<B\leq 1.345$ and $1.1186<C\leq 1.26$. Therefore we can take $C<1.1186$.\vspace{2mm}

 Finally for $B>1.304$ we find $\frac{B^3}{CDE}>2$. Again working as in Claim (iii) and find that  $\phi_{29}<0$ for
 $F>\frac{7-I}{2}-B-H$,  $\frac{7-I}{2}-B-0.788<H<0.704$, $1<I\leq A$, $1.58\leq A<1.64$, $1.304<B\leq 1.345$. This gives a
 contradiction.\vspace{2mm}

\noindent \textbf{Case III: 1.51 $\leq A < $ 1.58}\vspace{2mm}

 If $B\leq 1.276$, we find that $\phi_5^*<0$, so can take $B>1.276$.\vspace{2mm}

Suppose $B\geq 1.361$, then $\frac{B^3}{CDE}>2$ for $C<1.26$. We work as in Claim (iii) and find that   $\phi_{29}^{(2)}<0, ~\phi_{29}^{(3)}<0$
with $\mu_2=1$ for $1<I\leq A$, $1.51\leq A<1.58$, $1.361<B\leq 1.375$. Therefore we can take $B<1.361$.\\
Now $\frac{A^3}{BCD}>2$. If $E\geq 0.893$, we find $\phi_4<0$ for  $G\geq \frac{2}{3}E$,$1<I\leq A$ and $1.51\leq A<1.58$. So we
can take $E<0.893$.
If $F\geq 0.812$, we find  $\phi_{27}<0$ for $H\geq \frac{2}{3}F$, $1<I\leq A$, $E>0.46873A$ and $1.51\leq  A<1.58$. So we can
take $F<0.812$\vspace{2mm}

Suppose $C\geq 1.152$, then $\frac{C^3}{DEF}>2$ and we find that  $\phi_{11}<0$ for $G>0.46873C$, $1<I\leq A$, $A\geq 1.51$,
$1.276<B\leq 1.361$ and $1.152<C\leq 1.26$. Therefore we can take $C<1.152$. \vspace{2mm}

Suppose $B\geq 1.312$, then $\frac{B^3}{CDE}>2$ for $C<1.152$. Again working as in Claim (iii) and find that
$\phi_{29}^{(2)}<0, ~\phi_{29}^{(3)}<0$ with $\mu_2=1$ for $1<I\leq A$, $1.51\leq A<1.58$, $1.312<B\leq 1.361$.
Therefore we can take $B<1.312$.\vspace{2mm}

Suppose $d+f<0.35b+0.14i$. From (6.1) we have $h<b+\frac{i}{2}-d-f$. Then we find that $\phi_{22}^{(1)}<0$
for $1<I\leq A$, $1.276<B \leq 1.312$ and $1.51\leq A<1.58$. Therefore
we can take $d+f\geq 0.35b+0.14i$ which together with (6.1) gives $h<0.65b+0.36i$. \vspace{2mm}

Finally for $B>1.276$, we find $\frac{B^3}{CDE}>2$ for $C<1.152$ and $E<0.893$. If $H>0.915$, we find $\phi_{29}<0$. If $H<0.915$
i.e. $h>0.085$
we find $\phi_{3}<0$ for $0<f<b+\frac{i}{2}-\frac{g+h}{2}$, $0<g<b+\frac{i}{2}-\frac{e+h}{2}<b+\frac{i}{2}-\frac{0.107+h}{2}$ as
$E<0.893$, $0.085<h<0.65b+0.36i$, $0<i\leq a$, $0.51\leq a<0.58$, $0.276<b\leq 0.312$. This gives a contradiction.\vspace{2mm}

\noindent \textbf{Case IV: 1.41 $\leq A < $ 1.51}\vspace{2mm}

\noindent{\bf Subcase (i)} $B\geq 1.35$\vspace{2mm}

Suppose $\frac{B^3}{CDE}>2$. Then we find that $\phi_{29}<0$ for .......So we can take $\frac{B^3}{CDE}\leq 2$ which gives
$C\geq \frac{B^3}{2}>1.23,~ D\geq \frac{B^3}{2C}>0.976, ~E\geq \frac{B^3}{2C}>0.976.$ Then $E^4ABCDI>2$ and $B^2>CD$. But
$\phi_{30}<0$ for $E>0.976$, $1<I\leq A$, $1.35\leq B \leq 1.375$ and $1.41\leq A<1.51$. This gives a contradiction.\vspace{2mm}

\noindent{\bf Subcase (ii)} $1.267\leq B< 1.35$\vspace{2mm}

Suppose $d+f<0.329b+0.12i$. From (6.1) we have $h<b+\frac{i}{2}-d-f$. Then we find that $\phi_{22}^{(1)}<0$
for $1<I\leq A$, $1.267<B \leq 1.35$ and $1.41\leq A<1.51$. Therefore
we can take $d+f\geq 0.329b+0.12i$ which together with weak inequalities (2,2,2,2,1) gives $h<0.671b+0.38i$. Now $B^4AFGHI> 2$ for
$f<b-\frac{g+h}{2}+\frac{i}{2}, ~g<b-\frac{h}{2}+\frac{i}{2},~h<0.671b+0.38i$, $A>1.41$, $1<I\leq A<1.51$ and  $1.267\leq B \leq 1.35$
but $\phi_{29}<0$. This gives a contradiction.\vspace{2mm}

\noindent{\bf Subcase (iii)} $1 \leq B< 1.267$\\
If $B\leq 1.237$, we find that $\phi_5^*<0$, so can take $B>1.237$. \vspace{2mm}

Suppose $d+f<0.45b+0.208i$. As before we find that $\phi_{22}^{(1)}<0$ for
 $1<I\leq A$, $1.237<B \leq 1.267$ and $1.41\leq A<1.51$. Therefore
we can take $d+f\geq 0.45b+0.208i$ which gives $h<0.55b+0.292i$. Now $B^4AFGHI> 2$ but $\phi_{3}<0$ for $f<b-\frac{g+h}{2}+\frac{i}{2},
~0<g<b-\frac{h}{2}+\frac{i}{2},~0<h<0.55b+0.292i$, $1<I\leq A$, $1.41 \leq A <  1.51$ and  $1.237\leq B \leq 1.267$ . This gives a contradiction.\vspace{2mm}

\noindent \textbf{Case V: 1.32 $\leq A < $ 1.41}\vspace{2mm}

\noindent{\bf Subcase (i)} $B\geq 1.303$\vspace{2mm}

Working as in  Subcase (iii) of Case IV, we can take $d+f\geq 0.3b+0.103i$ which gives $h<0.7b+0.397i$. Now $B^4AFGHI> 2$  for
$f<b-\frac{g+h}{2}+\frac{i}{2}, ~0<g<b-\frac{h}{2}+\frac{i}{2},~0<h<0.7b+0.397i$, $1<I\leq A$, $1.32 \leq A <  1.41$ and
 $1.303\leq B \leq 1.375$. And we find that $\phi_{3}<0$ for $f<b-\frac{g+h}{2}+\frac{i}{2}, ~0<g<b-\frac{h}{2}+\frac{i}{2},~0<h<0.7b+0.397i$,
  $1<I\leq A$, $1.32 \leq A <  1.41$ and  $1 \leq B \leq 1.375$. This gives a contradiction.\vspace{2mm}

\noindent{\bf Subcase (ii)} $1.269\leq B< 1.303$\vspace{2mm}

Working as in  Subcase (iii) of Case IV, we can take $d+f\geq 0.45b+0.187i$ which gives $h<0.55b+0.313i$. Now $B^4AFGHI> 2$ but
$\phi_{3}<0$  for $f<b-\frac{g+h}{2}+\frac{i}{2}, ~0<g<b-\frac{h}{2}+\frac{i}{2},~0<h<0.55b+0.313i$, $1<I\leq A$, $1.32 \leq A <  1.41$ and
 $1.269\leq B \leq 1.303$.  This gives a contradiction.\vspace{2mm}

\noindent{\bf Subcase (iii)} $1 \leq B< 1.269$\vspace{2mm}

If $B\leq 1.198$, we find that $\phi_5^*<0$, so can take $B>1.198$. \\
Working as in  subcase(i) we can take $d+f\geq 0.5b+0.218i$ which gives $h<0.5b+0.282i$ and $g+h<b+0.564i$.\\
Suppose $C\geq 1.133$. Then $C^4ABGHI>2$ for $g+h<b+0.564i$, $1<I\leq A$, $1.32 \leq A <  1.41$ and  $1.198\leq B \leq 1.269$.
But then $\phi_{21}<0$  for $g+h<b+0.564i$, $B<1.269$, $1.32 \leq A <  1.41$ and  $1\leq C \leq 1.26$. So we can take $C<1.133$.\vspace{2mm}

Suppose $e<0.34b+0.196i$. Then $\phi_{23}<0$ for $g<\frac{a+i}{2}+c-e-\frac{h}{2}$, $0<h<0.5b+0.282i$, $0<c<0.133$, $0<e<0.34b+0.196i$,
$0<i\leq a$, $0.32\leq a<0.41$ and $0<b<0.269$. Therefore we can assume $e\geq 0.34b+0.196i$ which together with (6.2)
gives $g<0.83b+0.402i-\frac{h}{2}$.\vspace{2mm}

Suppose $B^4AFGHI> 2$. We find that $\phi_{3}<0$ for $f<b-\frac{g+h}{2}+\frac{i}{2}, ~0<g<b-\frac{h}{2}+\frac{i}{2},~0<h<0.5b+0.282i$,
$0<i\leq a$, $0.32 \leq a <  0.41$ and  $0 \leq b \leq 0.375$. Therefore we can assume that $B^4AFGHI< 2$.\vspace{2mm}

If $B \geq 1.2161$, we find that $B^4AFGHI> 2$ for $f<b-\frac{g+h}{2}+\frac{i}{2}, ~0<g<0.83b+0.402i-\frac{h}{2},~0<h<0.5b+0.282i$,
$0<i\leq a$, $0.32 \leq a <  0.41$. Therefore we can assume $B<1.2161$.\vspace{2mm}

Now when $B$ is further reduced, working as above we find that $d+f\geq 0.6b+0.252i$, $h<0.4b+0.248i$, $e\geq 0.48b+0.3i$, and
$g<0.76b+0.35i-\frac{h}{2}$. But then $B^4AFGHI> 2$ for $B>1.198$ which gives a contradiction.\vspace{2mm}

\noindent \textbf{Case VI: 1 $\leq A <$1.17.}\vspace{2mm}

 If $f+g+h \leq 1.6c+0.39(a+i)$, we find that $\phi_{24}^{(1)}<0$ for  $1<B\leq A$, $1<I\leq A$, $1< A\leq 1.17$. When $f+g+h > 1.6c+0.39(a+i)$,
 we get $e<0.2c+0.305(a+i)$ and then find that $\phi_{23}^{(1)}<0$ for $g<\frac{a+i}{2}+c-e-\frac{h}{2},~0<h<\frac{a+i}{2}+c-e,$
 $e<0.2c+0.305(a+i),~0<i\leq a $, $0<c \leq a\leq 0.17$. \vspace{2mm}

\noindent \textbf{Case VII: 1.17 $\leq A <$1.32.}\vspace{2mm}

\noindent{\bf Claim (i)} $B>1.125$\vspace{2mm}

If $B\leq 1.125$, we find that $\phi_5^*<0$, so can take $B>1.125$. \\

\noindent{\bf Claim (ii)} $h<0.584b+0.3i$ and $g+h<1.168b+0.6i$ \vspace{2mm}

Suppose $d+f<0.416b+0.2i$. From (6.1) we have $h<b+\frac{i}{2}-d-f$. we find that $\phi_{22}^{(1)}<0$
 for $1<I\leq A$, $1<B \leq A\leq 1.32$ . Therefore
we can take $d+f\geq 0.416b+0.2i$ which together with (6.1) and (6.3) gives $h<0.584b+0.3i$
and $g+h<1.168b+0.6i$ respectively. \vspace{2mm}

\noindent{\bf Claim (iii)} $B^4AFGHI< 2$ and $B<1.286$ \vspace{2mm}

 Suppose $B^4AFGHI>2$, then we find that $\phi_{3}<0$ for $f<b-\frac{g+h}{2}+\frac{i}{2}$,   $0<g<b-\frac{h}{2}+\frac{i}{2},
 ~0<h<0.584b+0.3i$, $0<i\leq a$, and $0 \leq b\leq a <  0.32$. Therefore we can assume that $B^4AFGHI< 2$.\vspace{2mm}

If $B \geq 1.286$, we find that $B^4AFGHI> 2$ for $f<b-\frac{g+h}{2}+\frac{i}{2}$, $0<g<b-\frac{h}{2}+\frac{i}{2},~0<h<0.584b+0.3i$,
 $0<i\leq a$, and $0.286 \leq b\leq a <  0.32$. Therefore we can assume $B<1.286$.\vspace{2mm}

\noindent{\bf Claim (iv)} $h<0.5b+0.294i$ \vspace{2mm}

With $B$ reduced, working as above we can further take $d+f\geq 0.5b+0.206i$ which  gives $h<0.5b+0.294i$. \vspace{2mm}

\noindent{\bf Claim (v)} $g<0.875b+0.395i-d-\frac{h}{2}$ \vspace{2mm}

Suppose $e<0.25b+0.21i$. Then $\phi_{23}<0$ for $g<\frac{a+i}{2}+c-e-\frac{h}{2}$, $0<h<0.5b+0.294i$, $0<c<0.26$, $0<e<0.25b+0.21i$,
$0<i\leq a$, $0.17\leq a<0.32$ and $0<b< \min \{a,0.286\}$. Therefore we can assume $e\geq 0.25b+0.21i$ which together with (6.2) gives $g<0.875b+0.395i-d-\frac{h}{2}$.\vspace{2mm}

\noindent{\bf Claim (vi)}  $B<1.253$ and $h<0.45b+0.26i$ \vspace{2mm}

If $B \geq 1.253$, we find that $B^4AFGHI> 2$ for $f<b-\frac{g+h}{2}+\frac{i}{2}, ~0<g<0.875b+0.395i-\frac{h}{2},~0<h<0.5b+0.294i$,
$0<i\leq a$, $0.253 \leq b \leq a <  0.32$. Therefore we can assume $B<1.253$.
With $B$ further reduced, working as above we can  take $d+f\geq 0.55b+0.24i$ which  gives $h<0.45b+0.26i$. \vspace{2mm}

\noindent{\bf Claim (vii)} $d\geq 0.36b+0.188i$ \vspace{2mm}

Suppose $d<0.36b+0.188i$. We find that $\phi_{31}<0$ for $f<b-d-\frac{g+h}{2}+\frac{i}{2}$, $0<g<b-d-\frac{e+h}{2}+\frac{i}{2}<
b-d-\frac{h}{2}-\frac{0.25b+0.21i}{2}+\frac{i}{2}$, $0<h<0.45b+0.26i$, $0<d<0.36b+0.188i$, $0<i\leq a$, $0.17\leq a<0.32$ and
$0<b< \min \{a,0.253\}$. Therefore we can assume that $d\geq 0.36b+0.188i$ which gives $f<0.64b+0.312i-\frac{g+h}{2}$ and
$0<g<0.515b+0.207i-\frac{h}{2}$. \vspace{2mm}

\noindent{\bf Claim (viii)} $C<1.137$ \vspace{2mm}

Suppose $C\geq 1.137$ Then $C^4ABGHI>2$ for $g<0.515b+0.207i-\frac{h}{2}$, $0<h<0.45b+0.26i$, $0<i\leq a$, $0.17\leq a<0.32$ and $0.125<b< 0.253$.
But then we find that  $\phi_{21}<0$ for $g+h<1.168b+0.6i$, $1<I\leq A$, $1<C<1.26$, $1<B \leq A\leq 1.32$. Hence we can take $C<1.137$.\vspace{2mm}

\noindent{\bf Claim (ix)} $i\geq 0.6a$ \vspace{2mm}

Suppose $i<0.6a$. If $f+g+h<1.64c+0.442(a+i),$ we find that $\phi_{24}<0$ for $b< \min\{a, 0.253\}$, $0<i<0.6a$, $0<c<0.137$ and $0.17<a<0.32$.
If $f+g+h\geq 1.64c+0.442(a+i),$ we get (using weak inequality (1,2,2,1,1,1,1)) that $e<0.18c+0.279(a+i)$. But then $\phi_{23}^{(1)}<0$ for
 $g<\frac{a+i}{2}+c-e-\frac{h}{2},~0<h<\frac{a+i}{2}+c-e,$ $e<0.18c+0.279(a+i),~0<i\leq 0.6a $, $0<c \leq 0.137$, $0.17<a \leq 0.32$.\vspace{2mm}

\noindent{\bf Claim (x)} $A\geq 1.223$ and $B>1.153$ \vspace{2mm}

We work as in Claim (vi). If $f+g+h<1.58c+0.43(a+i),$ we find that $\phi_{24}<0$  If $f+g+h\geq 1.58c+0.43(a+i),$ so that  $e<0.21c+0.285(a+i)$.
Then $\phi_{23}^{(1)}<0$ $0.6a<i\leq a $, $0<c \leq 0.137$, $0.17<a \leq 0.223$. Therefore we can take $A \geq 1.223$.
 But then $\phi_5^*<0$ for $B\leq 1.153$.\vspace{2mm}

\noindent{\bf Claim (xi)}  $B<1.164$ \vspace{2mm}

Suppose $B\geq 1.164$. Using $f<0.64b+0.312i-\frac{g+h}{2}$ and $0<g<0.515b+0.207i-\frac{h}{2}$, $0<h<0.45b+0.26i$, $0.6a\leq i\leq a$,
$0.223\leq a<0.32$ and $0.164\leq b< 0.253$ we find that $B^4AFGHI> 2$. Therefore we must have $B<1.164$.\vspace{2mm}

\noindent{\bf Claim (xii)} $C<1.103$ \vspace{2mm}

With improved bounds on $i,~a$ and $b$ we work as in Claim (viii) and get the result.\vspace{2mm}

\noindent{\bf Final contradiction} \vspace{2mm}

Suppose $A\leq 1.27$. If $f+g+h<1.6c+0.51(a+i),$ we find that $\phi_{24}<0$  If $f+g+h\geq 1.6c+0.51(a+i),$ so that  $e<0.2c+0.245(a+i)$.
 Then $\phi_{23}^{(1)}<0$. Therefore we can take $A \geq 1.27$. But then $\phi_5^*<0$ for $B\leq 1.164$. This gives a contradiction.\vspace{2mm}

\noindent \textbf{Proposition 43.} Case (32) i.e. $A>1$, $B>1$, $C>1$, $D\leq 1$, $E\leq 1$, $F\leq1$, $G\leq1$, $H>1$, $I>1$ does not arise.\vspace{2mm}

{\noindent \bf Proof.} Here $A\leq2$, $B\leq1.5$, $C\leq\frac{4}{3}$.
 Using the weak inequalities $(2,2,1,2,1,1),(1,2,2,1,1,1,1),(2,2,2,1,1,1)$ we have
\begin{equation}2b-2d-e-2g+h+i>0\vspace{-3mm}\end{equation}
\begin{equation}a+2c-2e-f-g+h+i>0\vspace{-2mm}\end{equation}
\begin{equation}2b-2d-2f-g+h+i>0\end{equation}

\noindent \textbf{ Note I :} We find that $\phi_{32}<0$ using  $G\geq 0.46873C$, $1<H\leq A$, $1<I\leq A$ when \\
(i) $1<B\leq A\leq 2,~ 1.155<C\leq \frac{4}{3}$ or\\
(ii) $1<B\leq A\leq 1.465,~ 1.11<C\leq \frac{4}{3}$ or\\
(iii) $1<B\leq A\leq 1.3,~ 1.05<C\leq \frac{4}{3}$ or\\
(iv) $1<B\leq 1.27,~1.02<A\leq 1.465,~ 1.01<C\leq \frac{4}{3}$.\vspace{2mm}

\noindent \textbf{ Note II :} If $B^4AFGHI>2$, then we find that $\phi_{29}<0$ for $F\geq 0.46873B$, $1<H\leq A$ and $1<I\leq A$ in the
 following cases \\
(i) $1.245<B\leq A\leq 1.7325$ or\\
(ii) $1<B\leq A\leq 1.3$ or\\
(iii) $1.2<B\leq 1.25,~1.29\leq A<1.32$ or \\
(iv) $1<B\leq A,~1.32\leq A<1.35$, $h+i<1.48a$ or \\
(v) $1.2<B\leq A,~1.35\leq A<1.38$, $h+i<1.48a$ or\\
(vi) $1.213<B\leq A,~1.38\leq A<1.41$, $h+i<1.48a$. \vspace{2mm}

\noindent\textbf{Claim(i) $C\leq1.155$ and $A\leq1.7325$}\vspace{2mm}

Suppose $C>1.155$. Using \eqref{1} we have $g<b+\frac{h+i}{2}$. Also $G>\frac{4}{9}C$. So we have
$$g<\left\{\begin{array}{lll}b+\frac{h+i}{2} & {\rm if }& 0<h+i<\frac{10}{9}-\frac{8}{9}c-2b\\\frac{5}{9}-\frac{4}{9}c
& {\rm if }& h+i>\frac{10}{9}-\frac{8}{9}c-2b. \end{array}\right. $$

We find that $C^4ABGHI>(1+c)^4(1+a)(1+b)(1+h+i)(1-g)>2$ for $g,~h+i$ as above and $0<b< \min\{a,0.5\}$, $0.155<c\leq a\leq 1$.
 Also  $\phi_{32}<0$ by Note I.  This gives a contradiction. Now $A\leq \frac{3}{2}C<1.7325$. \vspace{2mm}

\noindent\textbf{Claim(ii) $B\leq1.322$ }\vspace{2mm}

Suppose $B>1.322$, then $\frac{B^3}{CDE}>2$. But $\phi_{29}<0$ by Note II, which gives a contradiction.\vspace{2mm}

\noindent\textbf{Claim(iii) $A\leq1.465$}\vspace{2mm}

Suppose $A>1.465$, then $\frac{A^3}{BCD}>2$. Now if $E^2>FG$, then $\phi_{28}<0$ for $E\geq 0.46873A$, $1<H\leq A$ and $1<I\leq A$.
If $E^2\leq FG$, then $\phi_{33}<0$ using $FG>E^2$,  $E\geq0.46873A$, $1<H\leq A$ and $1<I\leq A$. Hence $A\leq1.465$.\vspace{2mm}

\noindent\textbf{Claim(iv) $A\geq 1.198$}\vspace{2mm}

Suppose $A<1.198$. If $f+g<1.6c+0.373(a+h+i)$, we find that
 $\phi_{24}^{(1)}<0$  for $0<b\leq a$, $0<h\leq a$, $0<i\leq a$ and $0<c\leq \min\{c,0.155\}$.
If $f+g\geq1.6c+0.373(a+h+i)$,
using (6.6) we have $e<0.2c+0.3135(a+h+i)$. But then
$\phi_{23}^{(1)}<0$, for $0<g<\frac{a+h+i}{2}+c-e$, $0<e<0.2c+0.3135(a+h+i)$, $0<h<a$, $0<i<a$ and $0<c<0.155$ and $0<a<0.198$,
a contradiction.\vspace{2mm}

\textbf{We divide the range of $A$ into several
 subintervals and  arrive at a contradiction in each.}\vspace{2mm}

 \noindent  \textbf{Case I: 1.198 $\leq A < $ 1.23}\vspace{2mm}


Firstly if $B\leq 1.1128$, we find that $\phi_6^*<0$, so can take $B>1.1128$. \vspace{2mm}


Next we find that for $C>1.1482$, $C^4ABGHI>2$ for $g<b+\frac{h+i}{2}$, $0<h\leq a,~0<i\leq a$ but $\phi_{32}<0$ by Note I, giving
thereby a  contradiction. So we can take $C\leq 1.1482$. \vspace{2mm}


Suppose now $h+i\leq1.73a$. If $f+g<1.55c+0.3958(a+h+i)$, we find that
 $\phi_{24}^{(1)}<0$  for $0.1128<b\leq a$, $0<h\leq 1.73a-i$, $0<i\leq a$ and $0<c\leq 0.1482$.
If $f+g\geq1.55c+0.3958(a+h+i)$,
 we have $e<0.225c+0.3021(a+h+i)$. But then
$\phi_{23}^{(1)}<0$, for $0<g<\frac{a+h+i}{2}+c-e$, $0<e<0.225c+0.3021(a+h+i)$, $0<h<1.73a-i$, $0<i<a$ and $0<c<0.1482$ and
$0.198<a<0.23$, a contradiction. Therefore we can take $h+i>1.73a$. \vspace{2mm}


 We find that  $\phi_{22}^{(1)}<0$ for $d+f\leq0.735b+0.26(h+i)$, $0<h\leq a$, $0<i\leq a$ and $0<b\leq a<0.23$. Hence
 $d+f>0.735b+0.26(h+i)$. Then using (6.7) we get $g<0.53b+0.48(h+i)$. Now suppose $C>1.0879$, then $C^4ABGHI >2$ for
 $1.73a-i<h\leq a$, $0<i\leq a$, $0.1128<b\leq a$ and $0.198<a<0.23$. But  $\phi_{32}<0$ by Note I; giving thereby a  contradiction.
 Hence we must have $C\leq 1.0879$. \vspace{2mm}


Suppose  $e\leq0.1884$, then working as in Claim (iii), we find that $\phi_{23}^{(1)}<0$, for $0<g<\frac{a+h+i}{2}+c-e$, $0<e<0.1884$,
 $1.73a-i<h\leq a$, $0<i\leq a$, $0<c\leq0.0879$ and $0.198<a<0.23$. Therefore we can take $e>0.1884$, i.e. $E<0.8116$. Now
 $\frac{B^3}{CDE}>2$  for $B\geq1.2088$ and $\phi_{29}<0$ by Note II. Hence we can take $B<1.2088$. \vspace{2mm}


Finally if $f+g<1.5c+0.4299(a+h+i)$, we find that
 $\phi_{24}^{(1)}<0$  for $0.1128<b\leq a$, $0<h\leq a$, $0<i\leq a$ and $0<c\leq 0.0879$.
If $f+g\geq1.5c+0.4299(a+h+i)$,
 we have $e<0.225c+0.28505(a+h+i)$. But then
$\phi_{23}^{(1)}<0$, for $0<g<\frac{a+h+i}{2}+c-e$, $0<e<0.225c+0.28505(a+h+i)$, $0<h\leq a$, $0<i\leq a$ and $0<c<0.1482$ and
 $0.198<a<0.23$, a contradiction.\vspace{2mm}

 \noindent  \textbf{Case II: 1.23 $\leq A < $ 1.26}\vspace{2mm}

Firstly if $B\leq 1.1128$, we find that $\phi_6^*<0$, so can take $B>1.1128$. \vspace{2mm}

As in Case I, for $C>1.14468$ we find that $C^4ABGHI>2$  but $\phi_{32}<0$ by Note I, giving thereby a  contradiction. So we
can take $C\leq 1.14468$. \vspace{2mm}


Suppose $h+i\leq1.49a$. We consider the cases $f+g<1.55c+0.3989(a+h+i)$ and $f+g\geq1.55c+0.3989(a+h+i)$,
 and work as in  Case I to get a contradiction. Therefore we can take $h+i>1.49a$. \vspace{2mm}

 We can further take $E<0.8061$ by showing $\phi_{23}^{(1)}<0$, then take  $B<1.2266$  by proving $\frac{B^3}{CDE}>2$
  for $B\geq1.2266$ and $\phi_{29}<0$ by Note II. \vspace{2mm}



We work as in  Case I to get $d+f>0.754b+0.25(h+i)$ which gives $g<0.492b+0.5(h+i)$.  Now $C^4ABGHI >2$ for $C>1.0834$,
$1.49a-i<h\leq a$, $0<i\leq a$, $0.1228<b\leq \min \{a,0.2266\}$ and $0.23<a<0.26$. But  $\phi_{32}<0$ by Note I. So we must
have $C\leq 1.0834$. Now $\frac{B^3}{CDE}>2$  for $B\geq1.20431$ and $\phi_{29}<0$ by Note II. So we can take $B<1.20431$. \vspace{2mm}


With $B$ and $C$ reduced, we can take  $h+i>1.78a$ by considering the cases $f+g<1.5c+0.4416(a+h+i)$ and $f+g\geq1.5c+0.4416(a+h+i)$.\vspace{2mm}


Suppose next $d+f\leq 0.783b+0.26(h+i)$. We find that  $\phi_{22}^{(1)}<0$ for $d+f\leq0.783b+0.26(h+i)$, $1.78a-i<h\leq a$,
$0<i\leq a$ and $0.1228<b\leq 0.20431,~0.23\leq a<0.26$. Hence $d+f>0.783b+0.26(h+i)$. Then using (6.7) we get
$g<0.434b+0.48(h+i)$. \vspace{2mm}

Suppose $d\leq 0.097$. We find that $\phi_{31}<0$ for $0<f<b+\frac{h+i}{2}-d-\frac{g}{2}$, $0<g<0.434b+0.48(h+i)$,
$0<d<0.097$, $0<h\leq a$, $0<i\leq a$, $0.1228<b<0.20431$ and $0.23<a<0.26$. Hence $d>0.097$.\vspace{2mm}


Suppose $B \geq 1.15$, then we find that $B^4AFGHI>2$ for $0<f<b+\frac{h+i}{2}-d-\frac{g}{2}<b+\frac{h+i}{2}-0.097-\frac{g}{2}$,
$0<g<b+\frac{h+i}{2}-d-\frac{e}{2}<b+\frac{h+i}{2}-0.097-\frac{0.1939
}{2}$, $1.78a-i<h\leq a$, $0<i\leq a$ and $0.23<a<0.26$. But then we get a contradiction using  Note II. Hence we can take
$B\leq 1.15$. \vspace{2mm}


 We find that  $\phi_{22}^{(1)}<0$ for $d+f\leq0.833b+0.297(h+i)$, $1.78a-i<h\leq a$, $0<i\leq a$ and $0.1228<b\leq 0.20431,~0.23\leq a<0.26$.
  Hence we must have $d+f>0.783b+0.26(h+i)$. Then using (6.7) we get $g<0.334b+0.406(h+i)$. Suppose $d\leq 0.147$. We find that $\phi_{31}<0$
  for $0<f<b+\frac{h+i}{2}-d-\frac{g}{2}$, $0<g<0.334b+0.406(h+i)$, $0<d<0.147$, $0<h\leq a$, $0<i\leq a$, $0.1228<b<0.15$ and $0.23<a<0.26$.
  Hence we must have $d>0.147$.\vspace{2mm}


Finally we find that $B^4AFGHI>2$ for $0<f<b+\frac{h+i}{2}-d-\frac{g}{2}<b+\frac{h+i}{2}-0.147-\frac{g}{2}$,
$0<g<b+\frac{h+i}{2}-d-\frac{e}{2}<b+\frac{h+i}{2}-0.147-\frac{0.1939
}{2}$, $1.78a-i<h\leq a$, $0<i\leq a$, $0.1228<b<0.15$ and $0.23<a<0.26$. Then we get contradiction as $\phi_{29}<0$ by  Note II.\vspace{1mm}

 \noindent  \textbf{Case III: 1.26 $\leq A < $ 1.29}\vspace{2mm}

 As in case I we can first take $B>1.1313$ by proving $\phi_6^*<0$, then can take $C\leq 1.14227$.
 We can further take $h+i>1.25a$ by considering the  cases $f+g<1.55c+0.4011(a+h+i)$ and $f+g\geq1.55c+0.4011(a+h+i)$.

%

 We can further take $E<0.8053$ by showing $\phi_{23}^{(1)}<0$, then take  $B<1.2254$
 by proving $\frac{B^3}{CDE}>2$  for $B>1.2254,~C<1.14227,~D<1,~E<0.8053$ and $\phi_{29}<0$ by Note II. \vspace{2mm}

%
We work as in  Case I to get $d+f>0.742b+0.25(h+i)$ which gives $g<0.516b+0.5(h+i)$.
 Now $C^4ABGHI >2$ for $C>1.0796$, $1.25a-i<h\leq a$, $0<i\leq a$, $0.1313<b\leq 0.2254$ and $0.26<a<0.29$.
  But  $\phi_{32}<0$ by Note I. So we must have $C\leq 1.0796$. Now $\frac{B^3}{CDE}>2$  for $B\geq1.2025$
  and $\phi_{29}<0$ by Note II. So we can take $B<1.2025$. \vspace{2mm}

With $B$ and $C$ reduced, we can take  $h+i>1.58a$ by considering the cases $f+g<1.5c+0.451(a+h+i)$ and $f+g\geq1.5c+0.451(a+h+i)$.\vspace{2mm}

We can further take $E<0.798$ by showing $\phi_{23}^{(1)}<0$, then take  $B<1.1989$  by proving $\frac{B^3}{CDE}>2$
 for $B>1.1989,~C<1.0796,~D<1,~E<0.798$ and $\phi_{29}<0$ by Note II. \vspace{2mm}

With $B$ and $C$ further reduced, we can take $g<0.42b+0.488(h+i)$ and $d>0.08$.\vspace{2mm}
%

%


Now suppose $B \geq 1.156$, then we find that $B^4AFGHI>2$ for $0<f<b+\frac{h+i}{2}-d-\frac{g}{2}<b+\frac{h+i}{2}-0.08-\frac{g}{2}$,
$0<g<b+\frac{h+i}{2}-d-\frac{e}{2}<b+\frac{h+i}{2}-0.08-\frac{0.202
}{2}$, $1.58a-i<h\leq a$, $0<i\leq a$ and $0.26<a<0.29$. So we get a
 contradiction using  Note II. Hence we can take $B\leq1.156$. \vspace{2mm}

With $B$ and $C$ further reduced, we can take $g<0.28b+0.45(h+i)$ and $d>0.13$ \vspace{2mm}
%

Finally We find that $B^4AFGHI>2$ for $0<f<b+\frac{h+i}{2}-0.13-\frac{g}{2}$, $0<g<b+\frac{h+i}{2}-0.13-\frac{0.202
}{2}$, $1.58a-i<h\leq a$, $0<i\leq a$, $0.1313<b<0.156$ and $0.26<a<0.29$. So we get a contradiction using  Note II.\vspace{1mm}

 \noindent \textbf{Case IV: 1.29 $\leq A < $ 1.32}\vspace{2mm}

\noindent{\bf Claim (i)} $B>1.1393$\\
\noindent \textbf{Claim(ii)}  $C\leq 1.14065$\\
\noindent{\bf Claim (iii)} $h+i>0.99a$ \vspace{2mm}

We consider the cases $f+g<1.55c+0.4006(a+h+i)$ and $f+g\geq1.55c+0.4006(a+h+i)$,
 and work as in Claim (iii) of Case I to get a contradiction.\vspace{2mm}

\noindent\textbf{Claim(iv)} $E<0.8103$ and $B<1.2273$ \vspace{2mm}

\noindent For the second result we use  $\frac{B^3}{CDE}>2$ for $B>1.2273,~C<1.14065,~D<1,~E<0.8103$ and $\phi_{29}<0$ by Note II.\vspace{2mm}

\noindent \textbf{Claim(v)} $d+f>0.7197b+0.25(h+i)$ and $C<1.0796$\\
\noindent \textbf{Claim(vi)}  $B\leq1.20498$\vspace{2mm}

We use here $\frac{B^3}{CDE}>2$ for $B>1.20498,~C<1.0796,~D<1,~E<0.8103$ and see that $\phi_{29}<0$  by Note II.\vspace{2mm}

\noindent \textbf{ Claim(vii)}  $f+g>1.5c+0.4586(a+h+i)$ and $h+i>1.37a$\\
\noindent\textbf{Claim(viii)} $E<0.803$ and $B<1.2014$ \vspace{2mm}

\noindent For the second result we use  $\frac{B^3}{CDE}>2$ for $B>1.2014,~C<1.0796,~D<1,~E<0.803$ and see that  $\phi_{29}<0$  by Note II.\vspace{2mm}

\noindent \textbf{Claim(ix)} $g<0.504b+0.48(h+i)$, $d>0.0342$ and $B\leq1.1875$ \vspace{2mm}

For the last result we use $\frac{B^3}{CDE}>2$ for $B>1.1875,~C<1.0796,~D<0.9658,~E<0.803$ and $\phi_{3}<0$ for $f<b+\frac{h+i}{2}-\frac{g}{2}$,
$0<g<b+\frac{h+i}{2}$, $0<h\leq a$ and $0<i\leq a$; a contradiction.\vspace{2mm}

\noindent \textbf{Claim(x)}  $d+f>0.785b+0.26(h+i)$ and $C<1.0619$ and $B\leq1.18097$\vspace{2mm}

For the last result we use  $\frac{B^3}{CDE}>2$ for $B>1.18097,~C<1.0619,~D<0.9658,~E<0.803$ and $\phi_{3}<0$ for $f<b+\frac{h+i}{2}-\frac{g}{2}$,
$0<g<b+\frac{h+i}{2}$, $0<h\leq a$ and $0<i\leq a$; a contradiction.\vspace{2mm}

\noindent \textbf{Claim(xi)} $g<0.408b+0.476(h+i)$ and $d>0.076$ \\
\noindent \textbf{Claim(xii)}  $B\leq1.1495$ \vspace{2mm}

Suppose $B \geq 1.1495$, then we find that $B^4AFGHI>2$ but $\phi_{3}<0$ for $0<f<b+\frac{h+i}{2}-d-\frac{g}{2}<b+\frac{h+i}{2}-0.076-\frac{g}{2}$,
 $0<g<b+\frac{h+i}{2}-d-\frac{e}{2}<b+\frac{h+i}{2}-0.076-\frac{0.197
}{2}$, $\max\{0.99a-i,0\}<h\leq a$, $0<i\leq a$ and $0.29<a<0.32$. So we get a contradiction.\vspace{2mm}

\noindent \textbf{Claim(xiii)} $g<0.304b+0.45(h+i)$ \\
\noindent \textbf{Final contradiction}\vspace{2mm}

We find that $B^4AFGHI>2$ but $\phi_{3}<0$ for $0<f<b+\frac{h+i}{2}-d-\frac{g}{2}<b+\frac{h+i}{2}-0.076-\frac{g}{2}$,
$0<g<b+\frac{h+i}{2}-d-\frac{e}{2}<b+\frac{h+i}{2}-0.076-\frac{0.197
}{2}$, $\max\{0.99a-i,0\}<h\leq a$, $0<i\leq a$, $0.1393<b<0.1495$ and $0.29<a<0.32$. This gives a contradiction.\vspace{1mm}

\noindent  \textbf{Case V: 1.32 $\leq A < $ 1.35}\vspace{2mm}

As in case I we can first take $B>1.14696$ by proving $\phi_6^*<0$, then can take $C\leq 1.14$\vspace{2mm}


Now we  distinguish the cases  $h+i<1.48a$ or~~$h+i\geq1.48a$.

\noindent  \textbf{Subcase I:~~~$ h+i<1.48a$}\vspace{2mm}

We can  take $E<0.836$ by showing $\phi_{23}^{(1)}<0$, then take  $B<1.2399$  by proving $\frac{B^3}{CDE}>2$
 for $B>1.2399,~C<1.14,~D<1,~E<0.836$ and see that  $\phi_{29}<0$  by Note II. \vspace{2mm}

We work as in  Case I to get $d+f>0.742b+0.25(h+i)$ which gives $g<0.538b+0.47(h+i)$ and then $d>0.096$.\vspace{2mm}

%


Finally we find that $B^4AFGHI>2$ for $0<f<b+\frac{h+i}{2}-d-\frac{g}{2}<b+\frac{h+i}{2}-0.096-\frac{g}{2}$,
$0<g<b+\frac{h+i}{2}-d-\frac{e}{2}<b+\frac{h+i}{2}-0.096-\frac{0.164
}{2}$, $0<h\leq \min\{a, 1.48a-i\}$, $0<i\leq a$, $0.14696<b<0.2399$ and $0.32<a<0.35$. But then $\phi_{29}<0$ by Note II, a contradiction.

\noindent {\center \textbf{Subcase II:~~~$ h+i\geq 1.48a$}}\vspace{2mm}

\noindent\textbf{Claim(i)} $E<0.929$ \vspace{2mm}

Suppose $E\geq0.929$. Here $B^2>CD$. We find that $\phi_{35}<0$ for $G\geq\frac{2}{3}E$, $E\geq0.929$, $1.48a<h+i<2a$,
 $0.14696<b<0.322$ and $0.32<a<0.35$. \vspace{2mm}

\noindent\textbf{Claim(ii)} $D<0.9651$  and $B<1.2692$\vspace{2mm}

If $D\geq 0.9651$, then $D^4ABCHI\geq D^4ABC(1+1.48a-i)(1+i)>2$ for $0.48a\leq i\leq a$, $D>0.9651$, $A>1.32$ and $B>1.14696$.
But then $\phi_{36}<0$
for $1<B<1.322$, $1<C\leq A$, $1<H\leq A$, $1<I\leq A$, $0.9651<D\leq 1$ and $1.32<A<1.35$. \\
For the second result we use  $\frac{B^3}{CDE}>2$ for $B>1.2692,~C<1.14,~D<1,~E<0.929$ and see that  $\phi_{29}<0$  by Note II. \vspace{2mm}

\noindent\textbf{Claim(iii)} $g<0.754b+0.5(h+i)$ and $C\leq1.0948$\vspace{2mm}

We work as in  Case I to get $d+f>0.623b+0.25(h+i)$ which gives $g<0.754b+0.5(h+i)$.  Now $C^4ABGHI >2$ for $C>1.0948$,
$1.48a-i\leq h\leq a$, $0.48a<i\leq a$, $0.14696\leq b\leq a$ and $0.32<a\leq 0.35$. But  $\phi_{32}<0$ by Note I.\vspace{2mm}

\noindent\textbf{Claim(iv)} $E<0.8183$ and  $B\leq1.2003$\vspace{2mm}

Suppose first $e\leq0.1817$, then  we find that $\phi_{23}^{(1)}<0$, for $0<g<\frac{a+h+i}{2}+c-e$, $0<e<0.1817$,
$1.48a-i<h\leq a$, $0<i\leq a$, $0<c\leq0.0948$ and $0.32<a<0.35$. Therefore we can take $e>0.1817$, i.e. $E<0.8183$.
Now $\frac{B^3}{CDE}>2$  for $B\geq1.2003$ and $\phi_{3}<0$ for $f<b+\frac{h+i}{2}-d-\frac{g}{2}<b+\frac{h+i}{2}-0.0349-\frac{g}{2}$,
$g<0.754b+\frac{h+i}{2}$, $0<h\leq a$, $0<i \leq a$ when $0<b\leq 0.25$. When $b>0.25$, $\phi_{29}<0$ by Note II.\vspace{2mm}

We repeat the cycle and get that\vspace{2mm}

\noindent\textbf{Claim(v)} $g<0.484b+0.5(h+i)$, $C<1.0582$, $B\leq 1.1868$ (proving here $\phi_{3}<0$) \\
\noindent\textbf{Claim(vi)} $A\leq 1.3434$ \vspace{2mm}

Suppose $A>1.3434$, we get $\frac{A^3}{BCD}>2$.  If $G>\frac{6.46}{9}$, we find that  $\phi_{37}<0$ for  $G>\frac{6.46}{9}$,
$E\geq0.46873A$,  $1<H,I\leq A$ and $1<A\leq 2$. Let now $G\leq\frac{6.46}{9}$, i.e. $\frac{2.54}{9}<g<0.53127$.
Also $E\geq\frac{2}{3}C$, i.e. $\frac{3}{2}e<\frac{1}{2}-c$ and $e<\frac{a+h+i}{2}+c-g$. Adding these two, we get
$e<\frac{1+a+h+i}{5}-\frac{2g}{5}$. Now $\phi_{37}<0$ for $e<\frac{1+a+h+i}{5}-\frac{2g}{5}$, $\frac{2.54}{9}<g<0.53127$,
$1.48a-i<h\leq a$ and $0.48a<i\leq a$ and $0.3<a<0.35$. Hence we have $A<1.3434$.\vspace{2mm}

Repeating the cycle we get\vspace{2mm}

\noindent\textbf{Claim(vii)} $E<0.812$ and  $B\leq1.1838$, $g<0.456b+0.48(h+i)$, $C\leq1.053$, $B\leq1.1818$ and $A\leq1.3393$\vspace{2mm}

Repeating the cycle once again we get\vspace{2mm}

\noindent\textbf{Claim(viii)} $g<0.436b+0.48(h+i)$, $d>0.054$ and  $B\leq1.174$ (proving here $\phi_{3}<0$
for $f<b+\frac{h+i}{2}-0.054-\frac{g}{2}$, $g<0.436b+0.48(h+i)$) \vspace{1mm}

\noindent \textbf{Claim(ix)} $A\leq1.3315$, $g<0.412b+0.47(h+i)$, $d>0.072$\\
\noindent \textbf{Final contradiction}\vspace{2mm}

Now $\frac{A^3}{BCD}>\frac{1.32^3}{1.174\times1.053\times 0.928}>2$, and  $\phi_{37}<0$ working as in Claim (x).This gives a contradiction.\vspace{1mm}

\noindent {\center \textbf{Case VI: 1.35 $\leq A < $ 1.38}}\vspace{1mm}

We work in a similar way as in Case V. First we can assume $B>1.1542$ and $C\leq 1.134$. \vspace{1mm}


\noindent  \textbf{Subcase I:~~~$ h+i<1.48a$}\vspace{1mm}

Here we get that $E<0.8312$, $B<1.2354$ ($\phi_{29}<0$ for $B\geq 1.2354$ by Note II), $g<0.536b+0.474(h+i)$ and $d>0.085$ \vspace{1mm}


Finally  we find that $B^4AFGHI>2$ for $B>1.1542$ and $\phi_{3}<0$ for $f<b+\frac{h+i}{2}-0.015-\frac{g}{2}$, $g<0.536b+0.474(h+i)$.
This gives a contradiction. \vspace{1mm}

\noindent  \textbf{Subcase II:~~~$ h+i\geq 1.48a$}\vspace{1mm}

\noindent\textbf{(i)} $E<0.915$ \\
\noindent\textbf{(ii)} $D<0.9499$  and $B<1.2539$ ($\phi_{29}<0$ for $B\geq 1.2539$ by Note II)\\
\noindent\textbf{(iii)} $g<0.75b+0.5(h+i)$ and $C\leq1.0899$\\
\noindent\textbf{(iv)}  $B\leq1.2374$ (proving here $\phi_{3}<0$)and $A\leq 1.3684$ \\
\noindent\textbf{(v)} $g<0.674b+0.5(h+i)$ and $C<1.0761$\\
\noindent\textbf{(vi)} $E<0.838$ and  $B\leq1.1966$ (proving here $\phi_{3}<0$)\\
\noindent \textbf{Final contradiction :}  $\frac{A^3}{BCD}>2$ and  $\phi_{37}<0$. \vspace{1mm}

\noindent {\center \textbf{Case VII: 1.38 $\leq A < $ 1.41}}\vspace{1mm}

We work in a similar way as in Case V. Firstly  proving $\phi_6^*<0$ we can assume $B>1.161$ and then working as in Claim (i)
we can take  $C\leq 1.128$ \vspace{2mm}

\noindent \textbf{Subcase I:~~~$ h+i<1.48a$}\vspace{1mm}

We get that $E<0.828$, $B<1.2316$, $g<0.534b+0.48(h+i)$ and $d>0.066$. \vspace{1mm}


Finally  we find that $B^4AFGHI>2$  for $B>1.161$ and $\phi_{3}<0$ for $f<b+\frac{h+i}{2}-0.034-\frac{g}{2}$,
$g<0.534b+0.48(h+i)$. This gives a contradiction. \vspace{1mm}

\noindent  \textbf{Subcase II:~~~$ h+i\geq 1.48a$}\vspace{1mm}

We get that $E<0.901$, $D<0.9353$ and $B\leq1.2388$ (by showing that $\phi_{3}<0$ for $f<b+\frac{h+i}{2}-0.0647-\frac{g}{2}$,
$0<g<b+\frac{h+i}{2}-0.0647-\frac{0.099}{2})$.\vspace{1mm}

Finally we find that $\frac{A^3}{BCD}>2$ but  $\phi_{37}<0$. This gives a contradiction.\vspace{2mm}

\noindent {\center \textbf{Case VIII: 1.41 $\leq A < $ 1.43}}\vspace{1mm}

Working as in Claim (i) we can take  $C\leq 1.1217$ \vspace{2mm}

Suppose $B<1.2494$ then we have $\frac{A^3}{BCD}>\frac{1.41^3}{1.2494\times1.1217}>2$. If $G>\frac{6.46}{9}$, we find that
 $\phi_{37}<0$ for  $G>\frac{6.46}{9}$, $E\geq0.46873A$,  $1<H,I\leq A$ and $1<A\leq 2$. Let now $G\leq\frac{6.46}{9}$, i.e.
 $\frac{2.54}{9}<g<0.53127$. If $h+i<a$, then $\phi_{37}<0$ for $0<e<\frac{a+h+i}{2}+c-g<\frac{a+h+i}{2}+0.1217-g$,
  $\frac{2.54}{9}<g<0.53127$, $0<h+i\leq a$. Now take $h+i>a$. As in Claim (x) of Case V, we also have $e<\frac{1+a+h+i}{5}-\frac{2g}{5}$.
   Now $\phi_{37}<0$ for $e<\frac{1+a+h+i}{5}-\frac{2g}{5}$, $\frac{2.54}{9}<g<0.53127$, $a<h+i\leq 2a$ and $0.41<a<0.43$.
   Hence we can take  $B\geq1.2494$.\vspace{2mm}

Suppose $E\geq0.851$, i.e. $e<0.149$. Let first $h+i<1.35a$
We find that $\phi_{19}<0$ for $0<g<\frac{a+h+i}{2}+c-e$, $0<e<0.149$, $0<h\leq 1.35a-i$, $0<i\leq a$, $0<c<0.1217$ and and $0.41<a<0.43$.
When  $h+i\geq 1.35a$, $\phi_{35}<0$  for $G\geq\frac{2}{3}E$, $E\geq0.851$,
 $1.35a<h+i<2a$, $0.2494<b<0.322$ and $0.41<a<0.43$. Hence we must have $E<0.851$. \vspace{2mm}

Finally  $\frac{B^3}{CDE}>\frac{1.2494^3}{1.1217\times0.851}>2$ and  $\phi_{29}<0$. This gives a contradiction.\vspace{1mm}

\noindent {\center \textbf{Case IX: 1.43 $\leq A < $ 1.465}}\vspace{1mm}

Working as above we first find that $C\leq 1.1177$. Then $\frac{B^3}{CDE}>2$  for $B>1.3076$ and  $\phi_{29}<0$. So we can take $B\leq 1.3076$.
But now $\frac{A^3}{BCD}>2$ and  $\phi_{37}<0$. This gives a contradiction.\vspace{2mm}

{\noindent \bf Proposition 44.} Case (15) i.e. $ A>1,~ B>1, ~C> 1, ~D>1,~ E\leq 1, ~F\leq 1,~ G \leq 1, ~H \leq 1, ~I \leq 1$
 does not arise.\vspace{2mm}

{\noindent \bf Proof.} Here  $ c\leq \frac{1}{2},~ d\leq \frac{1}{ 3}$.\vspace{2mm}
Using the weak inequalities $(2,1,2,2,2), ~(2,1,2,1,2,1)$, $(2,1,2,2,1,1),~
(2,2,2,2,1)$, ~$(1,2,2,2,1,1)$, ~$(1,2,2,2,2)$, $(1,2,2,1,2,1)$, $(1,2,2,1,1,1,1)$ and $(2,2,2,1,1,1)$   we get
\begin{equation}2b+c-2e-2g-2i>0,\vspace{-3mm}\end{equation}
\begin{equation}2b+c-2e-f-2h-i>0,,\vspace{-2mm}\end{equation}
\begin{equation}2b+c-2e-2g-h-i>0,\vspace{-2mm}\end{equation}
\begin{equation}2b+2d-2f-2h-i>0,\vspace{-2mm}\end{equation}
\begin{equation}a+2c-2e-2g-h-i>0\vspace{-2mm}\end{equation}
\begin{equation}a+2c-2e-2g-2i>0\vspace{-2mm}\end{equation}
\begin{equation}a+2c-2e-f-2h-i>0\vspace{-2mm}\end{equation}
\begin{equation}a+2c-2e-f-g-h-i>0\vspace{-2mm}\end{equation}
\begin{equation}2b+2d-2f-g-h-i>0\end{equation}

\noindent{\bf Claim (i)} $D < 1.26$

Suppose $D\geq 1.26$, then $\frac{D^3}{EFG}>2$. Now $\phi_1<0$ for $H>0.46873D$ and $1.26\leq D\leq A<2.1327$. \vspace{2mm}

\noindent{\bf Claim (ii)} $B< 1.7132$

Suppose $B\geq 1.7132$, then $\frac{B^3}{CDE}>2$. When $\frac{F^3}{GHI}>2$ we find $\phi_2<0$ for $F>0.46873B$ and $1.7132\leq B\leq A<2.1327$.
 When $\frac{F^3}{GHI}\leq 2$ we find $\phi_3<0$ for $GHI>\frac{F^3}{2}$ and $1.7132\leq B\leq A<2.1327$, $\frac{2}{3}<F\leq 1$. \vspace{2mm}

\noindent{\bf Claim (iii)} $A>1.155$

Suppose $A\leq 1.155$. If $g+h+i<1.6d+0.546a$, we find that $\phi_{25}<0$ for $1<D\leq A.$ If $g+h+i\geq 1.6d+0.546a$, which together
with (6.16) gives $f<b+0.2d-0.273a$, we find $\phi_{22}^{(1)}<0$ for $h<b+d-f-\frac{i}{2}, 0<i<b+d-f, 0<b\leq a$,
$0<d\leq a.$ \vspace{3mm}

We divide the range of $A$ into many subintervals and  arrive at a contradiction in each.\vspace{3mm}

\noindent \textbf{Case I: 1.9 $\leq A <$2.1327.}

  Here $\frac{A^3}{BCD}>2$ and $E\geq 0.46873A>0.8905$. We find that $\phi_4<0$ for $G\geq 0.72$, $I\geq \frac{2}{3}G$ and $E>0.46873A$
  or for $I\geq 0.75$, $G\geq \frac{2}{3}E$ and $E>0.46873A$. Also $\phi_5<0$ for $F\geq 0.86$, $H\geq \frac{2}{3}F$ and $E>0.46873A$ or
  for $H\geq 0.705$, $F\geq \frac{3}{4}E$ and $E>0.46873A$. So we can assume  $G<0.72,~I<0.75,~F<0.86,~H<0.705$. Now $\frac{D^3}{EFG}>2$
   and $\phi_9<0$ for $D\geq 1.2$. Finally  $\phi_4^*<0$ for $1\geq x\geq \frac{1}{1.7132 \times 1.5 \times 1.2A}$ and $1.9 \leq A <2.1327$.
 \vspace{2mm}

\noindent \textbf{Case II: 1.865 $\leq A < $ 1.9}\vspace{2mm}

 We work as in Case I
  to get $G<0.731,~I<0.85,~F<0.87,~H<0.725$ and $D< 1.214$.
Let first $\frac{2}{3}A\leq C <1.455$. Then $\phi_4^*<0$ for $1\geq x\geq \frac{1}{1.7132 \times 1.455 \times 1.214A}$ and $1.865 \leq A <1.9$. \\
When $1.455 \leq C \leq 1.5$, we find that  $\frac{C^3}{DEF}>2$. We note that $\phi_{11}$ is a decreasing function  of $G$.
We use $G>0.46873C$ and $I\geq \frac{2}{3}G$  to get
$$\phi_{11} <\phi_{11}^{(1)}=4A-\frac{2A^2}{B}+4C+\frac{14}{3}0.46873C-2(0.46873C)^2(\frac{2}{3}C^5AB)^{1/2}-9.$$
but $\phi_{11}^{(1)}$ which is a decreasing function of  $C$ is at most $0$ for  $1 \leq B <1.7132$ and $1.865 \leq A <1.9$.\\

\noindent \textbf{Case III: 1.755 $\leq A < $ 1.865}\\

\noindent \textbf{Subcase I: $1.674 \leq B <  1.7132, 1<C \leq 1.5$}\vspace{2mm}

 Here $\frac{B^3}{CDE}>2$. We find that $\phi_6<0$ for $H\geq 0.585$, $\phi_7<0$ for $G\geq 0.719$, $I\geq \frac{2}{3}G$, $F>0.46873B$.
 So we can take $H< 0.585$ and $G< 0.719$. Then $F^2>GH$ and  $\phi_8<0$ for $I\geq 0.578$ and $F>0.46873B$. So we can take $I< 0.578$.
 Now when $\frac{F^3}{GHI}>2$ we find $\phi_2<0$. When  $\frac{F^3}{GHI}\leq 2$ we have $I>\frac{F^3}{2GH}>0.5742$ and $G> \frac{F^3}{2HI}>0.714$
 and find that $\phi_7<0$. \vspace{2mm}

\noindent \textbf{Subcase II: $1.56 \leq B <  1.674, ~\frac{2}{3}A<C \leq 1.335$}\vspace{3mm}

Here $\frac{B^3}{CDE}>2$. We find that $\phi_6<0$ for $H\geq 0.648$ and $F>0.46873B$ or for $F>0.866$ and $H\geq \frac{2}{3}F$;
$\phi_7<0$ for $G\geq 0.769$, $I\geq \frac{2}{3}G$; and  $\phi_8<0$ for $I\geq 0.834$ and $F>0.46873B$. So we can assume that
$H< 0.648$, $F<0.866$, $G< 0.769$ and $I< 0.834$. Then  $\frac{D^3}{EFG}>2$ and $\phi_{23}<0$ for $D> 1.152$. So can take $D\leq  1.152$.

Suppose $F\geq E(>0.46873 A >0.8226)$. Then  $\phi_8<0$ for $I\geq 0.555$. So can take  $I< 0.555$. But then $\frac{F^3}{GHI}>2$ and
$\phi_{12}<0$ for
$C\leq 1.32$ or $B\geq 1.6$. If $C> 1.32$ and $B< 1.6$ then $\phi_{23}<0$. Therefore we can assume $F<E$.

Suppose $I\geq E( >0.8226)$. Then  $\phi_7<0$ for $G\geq 0.62$  and   $\phi_{8}<0$ for
$F\geq 0.735$ So can take  $G< 0.62$ and $F< 0.735$. Then $\phi_{23}<0$ Therefore we can assume $I<E$.

Now $E$ is greater than each of $F,G,H,I$, so $(4,5^*)$ holds as $\frac{A^3}{BCD}>2$. But then $\phi_4^*<0$ for
$1\geq x\geq \frac{1}{1.674 \times 1.335 \times 1.152A}$ and $1.755 \leq A <1.865$.  \vspace{3mm}

\noindent \textbf{Subcase III: $1.56 \leq B <  1.65, ~1.335<C \leq 1.5$}\vspace{3mm}

Let $E< \lambda$ and $G<\mu$, $I<\nu$, ~ $\lambda, \mu ,\nu $ be some rational numbers. The weak inequality (1,2,2,2,2) gives
$A+2C+2\lambda+2G+2I>9$ which further gives $G> \frac{9-A}{2}-\lambda-C-I$. We also have $I>\frac{9-A}{2}-\lambda-C-G> \frac{9-A}{2}-\lambda-C-\mu$. \\
We use
$$G>\left\{\begin{array}{lll}0.46873C & {\rm if }& I\geq \frac{9-A}{2}-\lambda-1.46873C\\\frac{9-A}{2}-\lambda-C-I
& {\rm if }& \frac{9-A}{2}-C-\lambda-\mu<I<\frac{9-A}{2}-\lambda-1.46873C\end{array}\right. $$
In the first case we replace $G$ by $0.46873C$ and $I$ by $\frac{9-A}{2}-\lambda-1.46873C$ (whenever $\phi_{11}$ is decreasing function of $I$)
 to get $\phi_{11}\leq$ $$ \phi_{11}^{(2)}=\frac{7}{2}A-\frac{2A^2}{B}+3(1.46873)C-\lambda-2(0.46873)^{\frac{3}{2}}C^4 (AB)^{\frac{3}{2}}
 (\frac{9-A}{2}-\lambda-1.46873C)^{\frac{1}{2}}-\frac{9}{2}$$

In the second case we replace $G$ by $\frac{9-A}{2}-\lambda-C-I$ to get $\phi_{11}\leq
\phi_{11}^{(3)}$ and then consider it as a function of $I$, where $\frac{9-A}{2}-C-\lambda-\mu<I<\frac{9-A}{2}-\lambda-1.46873C$.\vspace{2mm}

Here $\frac{B^3}{CDE}>2$. We find that $\phi_6<0$  for $F>0.866$. For $D<1.26$ and $F<0.866$ we get $\frac{C^3}{DEF}>2$. But then $\phi_{11}<0$
for $I\geq \frac{2}{3} G$,  $G\geq 0.726$ and $C>1.335$. When $G<0.726=\mu$,
 we get $\phi_{11}^{(2)}<0$ and $\phi_{11}^{(3)}<0$ on taking $E\leq 1=\lambda$. \vspace{3mm}

\noindent \textbf{Subcase IV: $1.65 \leq B <  1.674, ~1.335<C \leq 1.39$}

As before, using $\phi_6<0$ and $\phi_7<0$, we get that $F<0.845, H<0.6$ and $G<0.73=\mu$. Now $\frac{D^3}{EFG}>2$ and
$\phi_{12}<0$ for $D\geq 1.142$; $\frac{E^3}{FGH}>2$ and $\phi_{10}<0$ for $E\geq 0.987$. So can take $D< 1.142$ and $E< 0.987=\lambda.$ Now
$\phi_{11}^{(2)}<0$ and $\phi_{11}^{(3)}<0$.
 \vspace{3mm}

 \noindent \textbf{Subcase V: $1.65 \leq B <  1.674, ~1.39<C \leq 1.484$}\vspace{2mm}

Here, as above, using  $\phi_7<0$, we get that  $G<0.73=\mu$. Now $\phi_{11}^{(2)}<0$ and $\phi_{11}^{(3)}<0$ on taking $E\leq 1=\lambda$.
\vspace{3mm}

 \noindent \textbf{Subcase VI: $1.65 \leq B <  1.674, ~1.484<C \leq 1.5$}\vspace{2mm}

Here we find  $\phi_7<0$ and get that  $G<0.73=\mu$, $\phi_{11}<0$ and get that $E<0.99=\lambda$. Now $\phi_{11}^{(2)}<0$ and
$\phi_{11}^{(3)}<0$. \vspace{3mm}

 \noindent \textbf{Subcase VII: $\frac{3}{4}A \leq B <  1.56, ~1.365<C \leq 1.5$}\vspace{2mm}

Here $\frac{C^3}{DEF}>2$ and find that $\phi_{11}<0$ for $I\geq \frac{2}{3}$ and $G\geq 0.68$. Take now $G< 0.68=\mu$. Then
$\phi_{11}^{(2)}<0$ and $\phi_{11}^{(3)}<0$ on taking $E\leq 1=\lambda$. \vspace{2mm}

 \noindent \textbf{Subcase VIII: $\frac{3}{4}A \leq B <  1.495, ~\frac{2}{3}A \leq C \leq 1.365$}\vspace{2mm}

Here $\frac{A^3}{BCD}>2$ and $E\geq 0.46873A>0.82$. We find that $\phi_4<0$ for $G\geq 0.775$, $I\geq \frac{2}{3}G$ and
$E>0.46873A$. Also $\phi_5<0$ for $F\geq 0.892$, $H\geq \frac{2}{3}F$ and $E>0.46873A$ or for $H\geq 0.84$, $F\geq \frac{3}{4}E$
and $E>0.46873A$. So we can assume  $G<0.775,~F<0.892,~H<0.84$. Now $\frac{D^3}{EFG}>2$  and $\phi_9<0$ for $D\geq 1.165$.
Therefore $C\geq \frac{2}{3}A>1.17>D$ and so $(2,7^*)$ holds. But  $\phi_3^*<0$. \vspace{3mm}

 \noindent \textbf{Subcase IX: $1.495 \leq B <  1.56, ~\frac{2}{3}A \leq C \leq 1.365$}\vspace{2mm}

Proceeding as in Subcase VIII, first we obtain that $\frac{A^3}{BCD}>2$, $G<0.775,~F<0.892,~H<0.84$ and $D< 1.165$. Now $\frac{B^3}{CDE}>2$.
  Working as in Subcase II, using $\phi_6$; we can take that $H< 0.697$ and $F<0.885$.

Suppose $F\geq E(>0.46873 A >0.82)$. Then  $\phi_8<0$ for $I\geq 0.585$. So can take  $I< 0.585$. But then $\frac{F^3}{GHI}>2$
and $\phi_{12}<0$ for
$C\leq 1.365$ and  $B\geq 1.495$. Therefore we can assume $F<E$.

Suppose $I\geq E( >0.82)$. Then  $\phi_7<0$ for $G\geq 0.652$  and   $\phi_{8}<0$ for
$F\geq 0.743$ So can take  $G< 0.652$ and $F< 0.743$.  But then $\frac{D^3}{EFG}>2$  and $\phi_{23}<0$ for $D\geq 1$.
Therefore we can assume $I<E$.

Now $E$ is greater than each of $F,G,H,I$, so $(4,5^*)$ holds. But then $\phi_4^*<0$ for
$1\geq x\geq \frac{1}{1.56 \times 1.365 \times 1.165A}$ and $1.755 \leq A <1.865$.  \vspace{3mm}

\noindent \textbf{Case IV: 1.54 $\leq B\leq A < $ 1.755}\vspace{2mm}

\noindent \textbf{Subcase I: $1.66 \leq B\leq A$}

Working as in Subcase I of Case III,  we can take $H< 0.59$ and $G< 0.725$, $I< 0.593$. Now when $\frac{F^3}{GHI}>2$ we find
$\phi_2<0$. When  $\frac{F^3}{GHI}\leq 2$ we have $I>\frac{F^3}{2GH}>0.55$. If $G>0.68$ we  find that $\phi_7<0$.
If $G \leq 0.68$, we have $I>\frac{F^3}{2GH}>0.586$ and $G>\frac{F^3}{2HI}>0.6732$. Again we  find that $\phi_7<0$. \vspace{2mm}

\noindent \textbf{Subcase II: $1.54 \leq B <1.66,~ 1.375<C \leq 1.5$}\vspace{2mm}

Here $\frac{C^3}{DEF}>2$ and find that $\phi_{11}<0$ for $A \geq B$, $I\geq \frac{2}{3}$ and $G\geq 0.74$.
Take now $G< 0.74=\mu$. Then  $\phi_{11}^{(2)}<0$ and $\phi_{11}^{(3)}<0$ on taking $E\leq 1=\lambda$. \vspace{2mm}

\noindent \textbf{Subcase III: $1.54 \leq B <1.66,~ 1<C \leq 1.375$}\vspace{2mm}

If $C\leq 1.2422$, we find that $\phi_2^*<0$, so can take $C>1.2422$. Here $\frac{B^3}{CDE}>2$. We find that $\phi_6<0$ for
$H\geq 0.654$ and $F>0.46873B$ or for $F>0.8677$ and $H\geq \frac{2}{3}F$ ;  $\phi_7<0$ for $G\geq 0.7763$, $~
I\geq \frac{2}{3}G$; and  $\phi_8<0$ for $I\geq 0.8565$ and $F>0.46873B$. So we can assume that $H< 0.654$, $F<0.8677$,
$G< 0.7763=\mu$ and $I< 0.8565$. Then $\frac{E^3}{FGH}>2$ for $E\geq 0.9604$ and $\phi_{23}<0$. Therefore can take $E<0.9604$.
Now $\frac{D^3}{EFG}>2$ for $D> 1.1045$ and $\phi_{12}<0$ for $D> 1.18$ and $C<1.375$ or for $D>1.1045$ and $C<1.31$. So can
take $D\leq  1.18$ if $C\geq 1.31$ or $D \leq 1.1045$ if $C<1.31$.  Now $\frac{C^3}{DEF}>2$ for $C>1.2422$. \vspace{2mm}

If $I \leq 0.633=\nu$ we find that $\phi_{11}<0$ with $\lambda=1, ~\mu=0.7763$.
When $0.633<I<0.721=\nu$, then $\phi_7<0$ for $G\geq 0.7$ and $\phi_8<0$ for $F\geq 0.78$ or for $B\geq 1.636$. Therefore we can
assume that $G<0.7$, $F<0.78$ and $B<0.636$. Then $\frac{E^3}{FGH}>2$ for $E \geq 0.953$ as $H <0.654$. But then $\phi_{23}<0$,
so we can take $E<0.953$. But now we find that $\phi_{11}<0$ with $\lambda=0.953$ and at end points of $I$.\vspace{2mm}

When $0.721<I<0.8565=\nu$, then $\phi_7<0$ for $G\geq 0.66$ and $\phi_8<0$ for $F\geq 0.752$ or for $B\geq 1.592$. Therefore we can
assume that $G<0.66$, $F<0.752$ and $B<1.592$. Then $\frac{E^3}{FGH}>2$ for $E \geq 0.866$ as $H <0.654$. But then $\phi_{23}<0$,
so we can take $E<0.866$. But now we find that $\phi_{11}^{(2)}<0$ with $\lambda=0.866$ and at end points of $I$.\vspace{2mm}

\noindent \textbf{Case V:   1.67 $\leq A < $ 1.755, 1 $\leq B < $1.54,}\vspace{2mm}

\noindent \textbf{Subcase I: $1 \leq B < 1.448$}

Suppose first $\frac{C^3}{DEF}>2$, then $\phi_{11}<0$ for $A>1.67$, $I\geq \frac{2}{3}G$ and $G>0.72$. When $G\leq 0.72=\mu$,
we find that  $\phi_{11}^{(2)}<0$ and $\phi_{11}^{(3)}<0$ on taking $E\leq 1=\lambda$. Therefore we can take $\frac{C^3}{DEF}\leq 2$.
 This gives $C\leq 1.361$. If $B\leq 1.395$, we find $\phi_1^*<0.$ Similarly if $C\leq 1.15$, we find $\phi_2^*<0.$ Therefore we can
 take $B>1.395$ and $C>1.15$. This is true for all $B$, $1\leq B<1.54$.\\

Suppose $e+g<0.45b+0.25c$. Using (6.8) we get $i<b+\frac{c}{2}-e-g$. Then we find that $\phi_{23}<0$ for $A>1.67$.
Therefore
we can take $e+g\geq 0.45b+0.25c$ which together with (6.8) gives $i<0.55b+0.25c$. Now $D^4ABCHI>2$ for
$h<b+\frac{c}{2}-\frac{i}{2}$ (from (6.9)) and $D>1.145$, but then $\phi_9<0$. Hence we can take $D\leq 1.145$
 which is less than $C$. So $(2,7^*)$ holds, but $\phi_3^*<0$ for $A>1.67$ and $B<1.448$.

\noindent \textbf{Subcase II: $1.448 \leq B < 1.54$}

As in case VI, we can take $\frac{C^3}{DEF}\leq 2$ and $C\leq 1.361$. If $C\leq 1.184$, we find $\phi_2^*<0.$ Therefore we can also
take $C>1.184$.\\

Suppose  $\frac{B^3}{CDE}>2$, then we
 find that $\phi_6<0$ for $H\geq 0.728$ and $F>0.46873B$ or for $F>0.897$ and $H\geq \frac{2}{3}F $;  $\phi_7<0$ for $G\geq 0.831$, $~
I\geq \frac{2}{3}G$. So we can assume that $H< 0.728$, $F<0.897$ and $G< 0.813$. Now $\frac{C^3}{DEF}\leq 2$ implies $C< 1.313$.
Now $\frac{D^3}{EFG}>2$, for $D >1.15$.  Also $\phi_9<0$  for $D >1.15$, $B\leq 1.5$ and $\phi_{12}<0$  for $D >1.15$, $B> 1.5$.
Therefore we can take $D\leq 1.15$. But then $\frac{C^3}{DEF}\leq 2$ implies $C< 1.2731$. This gives $\frac{A^3}{BCD}>2$. Suppose $E\geq 0.897$.
 If  $I>E$, we get $\phi_4<0$ and if $I \leq E$ we find $\phi_4^*<0$. Therefore we can take $E< 0.897$. Now $D\leq 1.15$, $E< 0.897$, $F<0.897$
  and $\frac{C^3}{DEF}\leq 2$ implies $C< 1.2278$. This further restricts $D$ to $D\leq 1.1018$ using $\phi_{12}<0$. Now using $C>1.184$ and
  $\frac{C^3}{DEF}\leq 2$ we find that $E\geq \frac{C^3}{2DF}>0.839$. Similarly $F>0.839$. But then we
 find that $\phi_6<0$ for $H\geq 0.62$ and $F>0.839$;  $\phi_7<0$ for $G\geq 0.775$,
 $~I\geq \frac{2}{3}G$ and $F>0.839$; $\phi_8<0$ for $I\geq 0.54$ and $F>0.839$.  Therefore we can assume $H< 0.62$, $G<0.775$ and $I< 0.54$.
 Then $\frac{F^3}{GHI}>2$ and $\phi_{23}<0$ for $F>0.839$. This gives a contradiction.\vspace{2mm}

Hence we must have $\frac{B^3}{CDE}\leq 2$. This gives $B\leq 1.509$ using $C\leq 1.361$ and $D \leq 1.26$. As $E\geq \frac{B^3}{2CD}$,
$D\geq \frac{B^3}{2CE}$, $C\geq \frac{B^3}{2DE}$, we get $C>1.2047$, $D >1.1153$ and $E >0.8852$. If $D \geq 1.215$ and $\frac{D^3}{EFG}>2$,
we find that $\phi_9<0$. If $D \geq 1.215$ and $\frac{D^3}{EFG}\leq 2$, we get $E\geq \frac{D^3}{2}>0.8968$ and similarly $F>0.8968$.
Then $F^4ABCDE>2$ but
$\phi_{23}<0$. Hence we must have $D\leq 1.215$. Similarly if $C \geq 1.31$, and $\frac{C^3}{DEF}\leq 2$ we find that $E$ and $F$ are
$\geq \frac{C^3}{2D}>0.9251$.  But then $F^4ABCDE>2$ and $\phi_{23}<0$. Hence we can take $C<1.31$.\vspace{2mm}

Thus we are left now with $\frac{B^3}{CDE}\leq 2$, $D\leq 1.215$, $C<1.31$. Together with $B>1.448$, we get $B<1.4711,~D>1.1588,~C>1.2494$
and $E>0.9537$. Finally as in Subcase IV, we can take $e+g\geq 0.403b+0.25c$ which  gives $i<0.597b+0.25c$. Now $D^4ABCHI>2$ for
$h<b+\frac{c}{2}-\frac{i}{2}$  and $D>1.1588$, but then $\phi_9<0$.\vspace{2mm}

\noindent \textbf{Case VI: 1.58 $\leq A < $ 1.67, ~1 $\leq B < $1.54}\vspace{2mm}

\noindent \textbf{Subcase I: $1.361 \leq C < 1.5$}\vspace{2mm}

Here $\frac{C^3}{DEF}>2$, and $\phi_{11}<0$ for $A>1.58$, $I\geq \frac{2}{3}G$ and $G>0.74$. When $G\leq 0.74=\mu$,
we find that  $\phi_{11}^{(2)}<0$ and $\phi_{11}^{(3)}<0$ on taking $E\leq 1=\lambda$. \vspace{2mm}

\noindent \textbf{Subcase II: $1 <B \leq 1.398,
1 \leq C < 1.361$}\vspace{2mm}

If $B\leq 1.347$, we find $\phi_1^*<0.$ Similarly if $C\leq 1.138$, we find $\phi_2^*<0.$ Therefore we can take $B>1.347$ and $C>1.138$. \\
Suppose $e+g<0.499b+0.25c$. Using  $i<b+\frac{c}{2}-e-g$ from (6.8), we find that $\phi_{23}<0$ for $A>1.58$.
 Therefore
we can take $e+g\geq 0.499b+0.25c$ which together with (6.10) gives $h+i<1.002b+0.5c$. Now $D^4ABCHI>2$ for
$h+i<1.002b+0.5c$ and $D>1.125$, but then $\phi_9<0$. Hence we can take $D\leq 1.125$ which is less than $C$. So $(2,7^*)$ holds, but
 $\phi_3^*<0$ for $A>1.58$ and $B<1.398$.\vspace{2mm}

\noindent \textbf{Subcase III: $1.398 <B \leq 1.54,
1 \leq C < 1.361, ~\frac{B^3}{CDE}>2, ~~\frac{C^3}{DEF}>2$}\vspace{2mm}

Firstly $\phi_2^*<0$, if $C\leq 1.171$. Therefore we can take  $C>1.171$. \\
Now $\phi_{11}<0$ for $I\geq \frac{2}{3}G$ and $G\geq 0.85$. Let therefore $G<0.85=\mu$. If $I\leq 0.753=\nu$ then $\phi_{11}<0$
with $E\leq 1= \lambda$. So let $I\geq 0.753$. Now $\phi_{11}<0$  for $G>0.46873C$ and $C\geq 1.293$. Therefore can take $C< 1.293$\vspace{2mm}

Suppose $F>0.0.74$. As $\frac{B^3}{CDE}>2$,  we
 find that $\phi_6<0$ for $H\geq 0.68$  and $\phi_7<0$ for $G\geq 0.665$. When $H< 0.68$ and  $G< 0.665$ then $F^2>GH$ and
 $\phi_8<0$ for $F>0.0.74$ and $I\geq 0.753$. Therefore we can take $F\leq 0.74$.\vspace{2mm}

Suppose $E>0.9248$.  We
 find that $\phi_6<0$ for $H\geq 0.76$  and $\phi_7<0$ for $G\geq 0.703$. When $H< 0.76$ and  $G< 0.703$ then $\frac{E^3}{FGH}>2$ and
 $\phi_{23}<0$ for $E>0.9248$ and $I\geq 0.753$. Therefore we can take $E\leq 0.9248$.\vspace{2mm}

 Finally $\phi_{23}<0$ for $g<\frac{a}{2}+c-e-\frac{h}{2}-\frac{i}{2}< \frac{a}{2}+c-0.0752-\frac{h}{2}-\frac{i}{2}$,
 $0<h<\frac{a}{2}+c-e-\frac{i}{2}-\frac{f}{2}<\frac{a}{2}+c-0.0752-\frac{i}{2}-\frac{0.26}{2}$, $0<i<0.247$.\vspace{2mm}

 \noindent \textbf{Subcase IV: $1.398 <B \leq 1.54,
1 \leq C < 1.361, ~\frac{B^3}{CDE}\leq 2, ~~\frac{C^3}{DEF}>2$}\vspace{2mm}

As in Subcase III, we get $C>1.171$.  $\frac{B^3}{CDE}\leq 2$ implies $B<1.509$. Now $\phi_{11}<0$ for $I\geq \frac{2}{3}G$ and
$G\geq 0.83$. Let therefore $G<0.83=\mu$. If $I\leq 0.777=\nu$ then $\phi_{11}<0$ with $E\leq 1= \lambda$. So let $I\geq 0.777$.
 Now $\phi_{11}<0$  for $G>0.46873C$ and $C\geq 1.25$. Therefore can take $C< 1.25$. But then $B<1.4659$. We repeat the cycle and
 find that $I\geq 0.952$, then $\phi_{11}<0$  for $G>0.46873C$ and $C\geq 1.171$.\vspace{2mm}

\noindent \textbf{Subcase V: $1.398 <B \leq 1.54,
1 \leq C < 1.361, ~\frac{B^3}{CDE}> 2, ~~\frac{C^3}{DEF}\leq 2$}\vspace{2mm}

As in Subcase III, we get $C>1.171$. Suppose first $F>0.821$.  We
 find that $\phi_6<0$ for $H\geq 0.636$  and $\phi_7<0$ for $G\geq 0.797$. When $H< 0.636$ and  $G< 0.797$ then $F^2>GH$ and
 $\phi_{8}<0$ for  $I\geq 0.545$. Therefore we can take $I\leq 0.545$. But then  $\frac{F^3}{GHI}>2$ and
 $\phi_{2}<0$. Therefore we can take $F\leq 0.821$. This implies $C< 1.2743$.\vspace{2mm}

Secondly $\phi_7<0$ for $G\geq 0.86$, $F\geq \frac{2}{3}$ and  $~
I\geq \frac{2}{3}G$. So we can take $G<0.86$. Now $\frac{D^3}{EFG}>2$ for $D>1.122$, but $\phi_{12}<0$ for $D>1.165$.  If $D\leq 1.165$,
we have $C\leq 1.2414$ and then $\phi_{12}<0$ for $1.122<D\leq1.165$. Hence we can take $D\leq 1.122$ which gives $C \leq 1.226$.\vspace{2mm}

Now we have $E>\frac{C^3}{2DF}>0.8715$ and $F>\frac{C^3}{2DE}>0.7155$ for $C\geq 1.171$ and $D\leq 1.122$. For  $F>0.7155$  we
 find that $\phi_6<0$ for $H\geq 0.7$  and $\phi_7<0$ for $G\geq 0.835$. Therefore $H< 0.7$  and $G< 0.835$. Also $F\leq 0.821$.
 Then for $E \geq 0.9864$ we find  $\frac{E^3}{FGH}> 2$  but $\phi_{23}<0$. Hence we can take $E<0.9864$. Now $\frac{D^3}{EFG}>2$
 for $D>1.1059$, but $\phi_{12}<0$. This gives $D<1.1059$, $C< \sqrt[3]{2DEF}<1.2145$. Further if  $B\geq 1.465$, we find $\phi_2^*<0.$
  Therefore can take $B<1.465$.\vspace{2mm}

 Suppose $I\geq E( >0.8715)$. We find that  $\phi_7<0$ for $G\geq 0.64$. Therefore can take $G< 0.64$. Then $F^2>GH$ as $F>0.7155$
 and $H<0.7$ but $\phi_8<0$. Therefore we must have $I<E$. Now $E$ is larger than each of $F,~G,~H,~I$ and $\frac{A^3}{BCD}> 2$
  for $A>1.58,~B<1.465,C<1.2145,~D<1.1059$. Therefore $(4,5^*)$ holds but $\phi_4^*<0$. This gives a contradiction.\vspace{2mm}

 \noindent \textbf{Subcase V: $1.398 <B \leq 1.54,
1 \leq C < 1.361, ~\frac{B^3}{CDE}\leq 2, ~~\frac{C^3}{DEF}\leq 2$}\vspace{2mm}

As in Subcase III, we get $C>1.171$. Then we get $B<1.509$ using $\frac{B^3}{CDE}\leq 2$. \\
Suppose first $D\geq 1.217$. If $\frac{D^3}{EFG}> 2$ also, we get a contradiction as $\phi_9<0$. When $\frac{D^3}{EFG}\leq 2$ we get
 $E,~F \geq \frac{D^3}{2}>0.9$. If $F>0.92$ and $E>0.9$ we find $F^4ABCDE>2$ and $\phi_{20}<0$ for  $B<1.509$. If $F\leq 0.92$,
 $\frac{D^3}{EFG}\leq 2$ gives $E>0.979$. Also we get $D<1.2254,~C<1.3113,$ and hence $~B<1.4757$. Now again $F^4ABCDE>2$ for $E>0.979$,
  $F>0.9$ and $\phi_{20}<0$ for  $B<1.4757$.
This gives a contradiction. Hence we can take $D<1.217$ and therefore $C<1.3452,$  $B<1.485$.\vspace{2mm}

Similarly if we take $C \geq 1.297$, we find $D>1.0909,~E>0.8963,~F>0.8963$. If $F>0.92$ and $E>0.8963$ we find $F^4ABCDE>2$ and
$\phi_{20}<0$ for  $B<1.485$. If $\leq 0.92$ we get $D>1.1866$ and $E>0.975$ using $\frac{C^3}{DEF}\leq 2$. But then $F^4ABCDE>2$ for
$E>0.975$, $F>0.8963$ and $\phi_{20}<0$ for  $B<1.485$.
This gives a contradiction. Hence we can take $C<1.297$ and therefore $B<1.467$, $D>1.0532, E>0.8654$.\vspace{2mm}

Suppose $e+g<0.3b+0.2c$. Using  $i<b+\frac{c}{2}-e-g$, from (6.8), we find that $\phi_{19}<0$ for $A>1.58$.
Therefore
we can take $e+g\geq 0.3b+0.2c$ which together with (6.8)) gives $i<0.7b+0.3c$.
Now $D^4ABCHI>2$ for $h<b+\frac{c}{2}-\frac{i}{2}$ (from (6.9)) and $D>1.1776$, but then $\phi_9<0$. Hence we
can take $D\leq 1.1776.$ This implies $B<1.451$. \vspace{2mm}

If $F\geq 0.883$, we find $\phi_{22}<0$ using  $h<b+d-f-\frac{i}{2}$ (from (6.11)), $0<i<0.7b+0.3c$, $1.0532<D<1.1776$. Therefore we must have $F<0.883$
which gives $C<1.2765, B<1.4433$. We repeat the cycle by proving $e+g\geq 0.385b+0.25c$, $i<0.615b+0.25c$, $D<1.1267$, $C<1.2578$ and
$B<1.4152$. Also $B>1.398$ gives $D>\frac{B^3}{2C}>1.0861$. Further we get $F<0.785$ by proving $\phi_{22}<0$. But then we must have
$C<1.2094, B<1.398$.\vspace{2mm}

\noindent \textbf{Case VII: 1.5 $\leq A < $ 1.58, ~1 $\leq B < $1.54}\vspace{2mm}

\noindent \textbf{Subcase I: $1.361 \leq C < 1.5$}\vspace{2mm}

Here $\frac{C^3}{DEF}>2$, and $\phi_{11}<0$ for $A>1.5$, $I\geq \frac{2}{3}G$ and $G>0.76$. When $G\leq 0.76=\mu$,
we find that  $\phi_{11}^{(2)}<0$ and $\phi_{11}^{(3)}<0$ on taking $E\leq 1=\lambda$. \vspace{2mm}

\noindent \textbf{Subcase II: $1 <B \leq 1.352,
1 \leq C < 1.361$}\vspace{2mm}

We work as in Subcase II of Case VI. If $B\leq 1.3$, we find $\phi_1^*<0.$ Similarly if $C\leq 1.128$, we find $\phi_2^*<0.$ Therefore
we can take $B>1.3$ and $C>1.128$. \\
Suppose $e+g<0.527b+0.27c$. Using $i<b+\frac{c}{2}-e-g$ we find that $\phi_{19}<0$ for $a>0.5$. Therefore
we can take $e+g\geq 0.527b+0.27c$ which  gives $h+i<0.946b+0.46c$. Now $D^4ABCHI>2$ for  $D>1.1$, but then $\phi_9<0$. Hence we can take
$D\leq 1.1$ which is less than $C$. So $(2,7^*)$ holds, but $\phi_3^*<0$ for $A>1.5$ and $B<1.352$.\vspace{2mm}

\noindent \textbf{Subcase III: $1.352 <B \leq 1.4,~1 \leq C < 1.361,  ~~\frac{C^3}{DEF}>2$}\vspace{2mm}

Firstly $\phi_2^*<0$, if $C\leq 1.162$. Therefore we can take  $C>1.162$. \\
Now $\phi_{11}<0$ for $I\geq \frac{2}{3}G$ and $G\geq 0.75$. Let therefore $G<0.75=\mu$. If $I\leq 0.87=\nu$ then $\phi_{11}<0$ with $E\leq 1= \lambda$. So let $I\geq 0.87$. Now $\phi_{11}<0$  for $G>0.46873C$ and $C\geq 1.162$. \vspace{2mm}

\noindent \textbf{Subcase IV: $1.352 <B \leq 1.4, ~1 \leq C < 1.361,  ~~\frac{C^3}{DEF}\leq 2$}\vspace{2mm}

Firstly as in Subcase III, we can take  $C>1.162$. Suppose $e+g<0.4b+0.209c$. Using  $i<b+\frac{c}{2}-e-g$
we find that $\phi_{23}<0$ for $A>1.5$. Therefore
we can take $e+g\geq 0.4b+0.209c$ which together with (6.8) gives $i<0.6b+0.291c$.
Now $D^4ABCHI>2$ for $h<b+\frac{c}{2}-\frac{i}{2}$ ((6.9)) and $D>1.159$, but then $\phi_9<0$. Hence we can take
$D\leq 1.159.$ This implies $C<1.3235$ as $\frac{C^3}{DEF}\leq 2$. \vspace{2mm}

 We find $\phi_{23}<0$ using  $H>1-b-d+f+\frac{i}{2}, 0<i<0.6b+0.291c$, when  $F\geq 0.87$, $1<D\leq 1.1$ or when $F\geq 0.839$, $1.1<D<1.159$.
 Therefore we must have $F<0.87$ when $1<D\leq 1.1$ and $F< 0.839$ when  $1.1<D<1.159$. which gives $C<1.2483$. We repeat the cycle by
 proving $e+g\geq 0.43b+0.25c$, $i<0.57b+0.25c$, $D<1.119$. Further we get $F<0.81$ by proving $\phi_{23}<0$. But then we must have $C<1.2194$.\\

 Further if $E\geq 0.896$, we find $\phi_{23}<0$ using  $g<b+\frac{c}{2}-e-\frac{h+i}{2}$, (from (6.10)), $0<h<b+\frac{c}{2}-e-\frac{f+i}{2}$, from (6.9) and $
 0<i<0.57b+0.25c$.
 Therefore we can take $E<0.896$ which gives $C<1.1755$. Finally for $1.162<C<1.1755$, $D>\frac{C^3}{2EF}>1.08$, we find $D^4ABCHI>2$ and
 $\phi_9<0$ using $h<b+\frac{c}{2}-e-\frac{f+i}{2}$.\vspace{2mm}

\noindent \textbf{Subcase V: $1.4 <B \leq 1.54,
1 \leq C < 1.361, ~\frac{B^3}{CDE}>2, ~~\frac{C^3}{DEF}>2$}\vspace{2mm}

Firstly $\phi_2^*<0$, if $C\leq 1.193$. Therefore we can take  $C>1.193$. \\
Now $\phi_{11}<0$ for $I\geq \frac{2}{3}G$ and $G\geq 0.865$. Let therefore $G<0.865=\mu$. If $I\leq 0.769=\nu$ then $\phi_{11}<0$ with
$E\leq 1= \lambda$. So let $I\geq 0.769$. \vspace{2mm}

Suppose $F>0.73$. As $\frac{B^3}{CDE}>2$,  we
 find that $\phi_6<0$ for $H\geq 0.675$  and $\phi_7<0$ for $G\geq 0.642$. When $H< 0.675$ and  $G< 0.642$ then $F^2>GH$ and
 $\phi_8<0$ for $F>0.73$ and $I\geq 0.769$. Therefore we can take $F\leq 0.73$.\vspace{2mm}

Suppose $E>0.9014$.  We
 find that $\phi_6<0$ for $H\geq 0.742$  and $\phi_7<0$ for $G\geq 0.676$. When $H< 0.742$ and  $G< 0.676$ then $\frac{E^3}{FGH}>2$ and
 $\phi_{23}<0$ for $E>0.9014$ and $I\geq 0.769$. Therefore we can take $E\leq 0.9014$.\vspace{2mm}

 Finally $\phi_{23}<0$ for $g<\frac{a}{2}+c-e-\frac{h}{2}-\frac{i}{2}< \frac{a}{2}+c-0.0986-\frac{h}{2}-\frac{i}{2}$,
 $0<h<\frac{a}{2}+c-e-\frac{i}{2}-\frac{f}{2}<\frac{a}{2}+c-0.0986-\frac{i}{2}-\frac{0.27}{2}$, $0<i<0.231$.\vspace{2mm}

\noindent \textbf{Subcase VI: $1.4 <B \leq 1.54,
1 \leq C < 1.361, ~\frac{B^3}{CDE}\leq 2, ~~\frac{C^3}{DEF}>2$}\vspace{2mm}

As in Subcase V, we get $C>1.193$.  $\frac{B^3}{CDE}\leq 2$ implies $B<1.509$. Now $\phi_{11}<0$ for $I\geq \frac{2}{3}G$ and $G\geq 0.86$.
Let therefore $G<0.86=\mu$. If $I\leq 0.788=\nu$ then $\phi_{11}<0$ with $E\leq 1= \lambda$. So let $I\geq 0.788$. Now $\phi_{11}<0$
 for $G>0.46873C$ and $C\geq 1.298$. Therefore can take $C< 1.298$. But then $B<1.4845$. We repeat the cycle and find that $I\geq 0.87$,
  then $\phi_{11}<0$  for $G>0.46873C$ and $C\geq 1.252$. Therefore can take $C< 1.252$. But then $B<1.4667$. We repeat the cycle once again
   and find that $I\geq 0.936$, then $\phi_{11}<0$  for $G>0.46873C$ and $C\geq 1.214$. Therefore can take $C< 1.214$. But then $B<1.4517$.
   Finally we find that $I\geq 0.99$, then $\phi_{11}<0$  for $G>0.46873C$ and $C\geq 1.193$.\\

\noindent \textbf{Subcase VII: $1.4 <B \leq 1.54,
1 \leq C < 1.361, ~\frac{B^3}{CDE}>2, ~~\frac{C^3}{DEF}\leq 2$}\vspace{2mm}

As in Subcase V, we get $C>1.193$. Suppose first $F>0.815$.  We
 find that $\phi_6<0$ for $H\geq 0.63$  and $\phi_7<0$ for $G\geq 0.79$. When $H< 0.63$ and  $G< 0.79$ then $F^2>GH$ and
 $\phi_{8}<0$ for  $I\geq 0.531$. Therefore we can take $I\leq 0.531$. But then  $\frac{F^3}{GHI}>2$ and
 $\phi_{2}<0$. Therefore we can take $F\leq 0.815$. This implies $C< 1.2712$.\vspace{2mm}

Secondly $\phi_7<0$ for $G\geq 0.853$, $F\geq \frac{2}{3}$ and  $~I\geq \frac{2}{3}G$. So we can take $G<0.853$. Now $\frac{D^3}{EFG}>2$ for $D>1.1162$, but $\phi_{12}<0$.  If $D\leq 1.1162$, we have
$C\leq 1.2208$. Also $F\geq \frac{C^3}{2D}>0.7605$ as $C>1.193$.\vspace{2mm}

For  $F>0.7605$,  we find that $\phi_6<0$ for $H\geq 0.655$  and $\phi_7<0$ for $G\geq 0.81$. Therefore can take $H< 0.655$  and $G< 0.81$. Also $F\leq 0.815$.
 Then for $E \geq 0.9528$ we find  $\frac{E^3}{FGH}> 2$  but $\phi_{18}<0$. Hence we can take $E<0.9528$. Now $D\geq \frac{C^3}{2EF}>1.09$ as
 $C>1.193$. Therefore $\frac{D^3}{EFG}>2$ for $D>1.09$, but $\phi_{13}<0$. This gives a contradiction.\vspace{2mm}

\noindent \textbf{Subcase VIII: $1.4 <B \leq 1.54,
1 \leq C < 1.361, ~\frac{B^3}{CDE}\leq 2, ~~\frac{C^3}{DEF}\leq 2$}\vspace{2mm}

Here $C>1.193$ and $B<1.509$.
Suppose first $C\geq 1.31$.  Then $\frac{C^3}{DEF}\leq 2$ gives  $D \geq 1.124,~E>0.892,~F>0.892$. If $F>0.928$  we find $F^4ABCDE>2$
and $\phi_{14}<0$ for  $B<1.509$. If $0.892<F\leq 0.928$, we get $C<1.3274$,  $D>1.2112,~E>0.9613$. Now again $F^4ABCDE>2$.
But $\phi_{14}<0$ for  $F>0.899$. For $F<0.899$, we get $C<1.3134$ and $\phi_{14}<0$.
This gives a contradiction. Hence we can take  $C<1.31,$ and so  $B<1.489$.\vspace{2mm}

Suppose $e+g<0.223a+0.3c$. Using (6.13) we get $I>1-c-\frac{a}{2}+e+g$. Then we find that $\phi_{19}<0$
for $1.5<A<1.58$, $1.4 <B \leq 1.489,~ 1.193 \leq C < 1.31$. Therefore
we can take $e+g\geq 0.223a+0.3c$ which together with (6.12) gives $i<0.277a+0.7c$.
Now $D^4ABCHI>2$ for $h<c+\frac{a}{2}-\frac{i}{2}$ (from (6,14)) and $D>1.1715$, but then $\phi_{16}<0$.
Hence we can take $D\leq 1.1715.$ This implies $B<1.4533$ and $E>0.894$. \vspace{2mm}

Suppose  $C\geq 1.286$.  Then $\frac{C^3}{DEF}\leq 2$ gives  $D \geq 1.0633,~E>0.9077,~F>0.9077$. If $F>0.9356$  we find $F^4ABCDE>2$
and $\phi_{14}<0$ for  $B<1.4533$. If $F\leq 0.9356$, we get   $D>1.136,~E>0.9697$. Now again $F^4ABCDE>2$ for  $F>0.9077$ and $\phi_{14}<0$.
This gives a contradiction. Hence we can take  $C<1.286,$ and so  $B<1.4444, ~D>1.0668,~E>0.9106$.\vspace{2mm}

Suppose $e+g<0.25a+0.365c$. Working as above, we find that $\phi_{19}<0$ for $1.5<A<1.58$, $1.4 <B \leq 1.4444, ~1.193 \leq C < 1.286$.
 Therefore
we can take $e+g\geq 0.25a+0.365c$ which  gives $i<0.25a+0.635c$.
Now $D^4ABCHI>2$ for $h<c+\frac{a}{2}-\frac{i}{2}$  and $D>1.155$, but then $\phi_{16}<0$. We note here that $\phi_{16}<0$ even for $D>1$.
Hence we can take $D\leq 1.155.$ This implies $B<1.4376$.\vspace{2mm}

Suppose $e+g<0.28b+0.198c$. Using  $i<b+\frac{c}{2}-e-g$ we find that $\phi_{19}<0$
for $1.5<A<1.58$, $1.4 <B \leq 1.4376, ~1.193 \leq C < 1.286$. Therefore
we can take $e+g\geq 0.28b+0.198c$ which together with (6.8) gives $i<0.72b+0.302c$. Further we find that
$\phi_{22}<0$ for $0.0668<d<0.155,~0<f<0.118$, $0<i<0.72b+0.302c$, $0.5<a<0.58$, $0.4<b<0.4376$ and $0.193<c<0.286$.
Therefore we can take $f\geq 0.118$ i. e. $F\leq 0.882$. This together with $\frac{B^3}{CDE}\leq 2, ~~\frac{C^3}{DEF}\leq 2$
implies $C<1.2678$ and $B<1.43072$.\vspace{2mm}

Repeating the cycle we find that $\phi_{19}<0$ for $e+g<0.315b+0.213c$, $A>1.5$, $1.4 <B \leq 1.43072, 1.193 \leq C < 1.2678$. Therefore
we can take $e+g\geq 0.315b+0.213c$ which together with (6.8) gives $i<0.685b+0.287c$. Now $D^4ABCHI>2$
for $h<b+d-f+\frac{i}{2}<b+d-0.118+\frac{i}{2}$, $0<i<0.685b+0.287c$, $A>1.5$, $d>0.1335$, $0.193<c<0.2678$ and $0.4<b<0.43072$.
 Already $\phi_{16}<0$ even for $D>1$. Therefore we can take $D\leq 1.1335$ which gives $C<1.26$, $B<1.4189$ and $D>1.0888$.\vspace{2mm}

We again repeat the cycle to find that  $\phi_{22}<0$ for $0.0888<d<0.1335,~0<f<0.165$, $0<i<0.685b+0.287c$, $0.5<a<0.58$,
$0.4<b<0.4189$ and $0.193<c<0.26$. Therefore we can take $f\geq 0.165$ i. e. $F\leq 0.835$. But now we use
$h<b+d-f+\frac{i}{2}<b+d-0.165+\frac{i}{2}$ to find that
 $D^4ABCHI>2$ for $d>0.103$. So we can take $D\leq 1.103$. But then $\frac{C^3}{DEF}\leq 2$ implies $C<1.2258$ and
 $\frac{B^3}{CDE}\leq 2$ implies $B<1.394$. This gives a contradiction as $B>1.4$ in this subcase.\vspace{2mm}

\noindent \textbf{Case VIII: 1.155 $\leq A < $ 1.5, ~1 $\leq B \leq A $}\vspace{2mm}

\noindent \textbf{Subcase I: $\frac{C^3}{DEF}> 2 $}\vspace{2mm}

 Suppose first $1.27 \leq C < 1.5$.
If $G\geq 0.79$, we find that $\phi_{12}<0$ for $I\geq \frac{2}{3}G$. If $G<0.79=\mu$, we work as in Subcase (iii) of Case III to find that
the corresponding $\phi_{12}^{(2)}<0$, $\phi_{12}^{(3)}<0$ with $\lambda=1$, for   $1<B\leq A$, $1<A<1.5$ and $1.27<C\leq 1.5$.\vspace{2mm}

Let now $1 \leq C < 1.27$.
 Suppose $e+g<0.27a+0.305c$. Using (6.13) we get $I>1-c-\frac{a}{2}+e+g$. Then  we find that $\phi_{19}<0$ for $1<A<1.5$,
 $1 <B \leq A, 1 \leq C < 1.27$. Therefore
we can take $e+g\geq 0.27a+0.305c$ which  gives $i<0.73a+0.195c$. Now $\phi_{21}<0$ for $g<\frac{a}{2}+c-\frac{h+i}{2},~0<h<\frac{a}{2}+c-\frac{i}{2},$
$0<i<0.73a+0.195c$, $0<b\leq a, 0<a<0.5$ and $0<c<0.27$.\vspace{2mm}

\noindent \textbf{Subcase II: $\frac{C^3}{DEF}\leq 2,~1.35<B\leq A $ }\vspace{2mm}

We note that $\frac{C^3}{DEF}\leq  2$ implies $C< 1.361$.
Working as in  subcase V of Case VII, we get a contradiction if we have $\frac{B^3}{CDE}> 2 $. Therefore we can assume $\frac{B^3}{CDE}\leq 2 $.
Also $\phi_2^*<0$, if $C\leq 1.179$. Therefore we can take  $C>1.179$. \vspace{2mm}

Suppose first $C\geq 1.33$.  Then $\frac{C^3}{DEF}\leq 2$ gives  $D \geq 1.1763,~E>0.933,~F>0.933$. Then  we find $F^4ABCDE>2$ and
$\phi_{15}<0$ for  $1.35<B\leq A<1.5$. This gives a contradiction. Hence we can take  $C<1.33$.\vspace{2mm}

 Suppose $e+g<0.18a+0.32c$. Using  $I>1-c-\frac{a}{2}+e+g$ we find that $\phi_{19}<0$ for
 $1<A<1.5$, $1.35 <B \leq A, 1.179 \leq C < 1.361$. Therefore
we can take $e+g\geq 0.18a+0.32c$ which  gives $i<0.32a+0.68c$. Now $D^4ABCHI>2$ for $h<c+\frac{a}{2}-\frac{i}{2}$  and $D>1.185$, but
then $\phi_{16}<0$. Hence we can take $D\leq 1.185.$ This implies $C<1.308$ and so  $B<1.4577$.\vspace{2mm}

We repeat the cycle to get $e+g\geq 0.24a+0.32c$, $i<0.26a+0.68c$, $D\leq 1.17.$\\If $F\geq 0.952$ and $1<D<1.05$ or if $F\geq 0.927$
and $1.05<D<1.17$, we find $\phi_{22}<0$ using  $H>1-b-d+f+\frac{i}{2}, 0<i<0.26a+0.68c$. Therefore we must have $F< 0.952$ when
$1<D<1.05$ and $F< 0.927$ when $1.05<D<1.17$ which gives $C<1.2864$. This further gives $B<1.4439$.\vspace{2mm}

We repeat the cycle again to get $e+g\geq 0.279a+0.32c$, $i<0.221a+0.68c$, $D\leq 1.1463$ using $h<c+\frac{a}{2}-\frac{i}{2}-\frac{f}{2}$,
 $F< 0.907$ when $1<D<1.07$ and $F< 0.882$ when $1.07<D<1.1463$ which gives $C<1.2563$ and $B<1.4228$. If $e+g<0.317a+0.34c$,
 we find that $\phi_{19}<0$ for $1.35<A<1.5$, $1.35 <B \leq \min(A,1.4228), 1.179 \leq C < 1.2864$. Therefore
we can take $e+g\geq 0.317a+0.34c$ which  gives $i<0.183a+0.66c$.\vspace{2mm}

If $E\geq 0.881$, we find $\phi_{23}<0$. So we can take $E<0.881$ which gives $C<1.213$ and $B<1.4$. We repeat the cycle once again
to get $e+g\geq 0.317a+0.475c$, $i<0.183a+0.525c$, $h+i<0.366a+0.105c$, $E<0.835$, $D<1.05$ which gives $C<1.179$. This gives a contradiction. \vspace{2mm}

\noindent \textbf{Subcase III: $\frac{C^3}{DEF}\leq 2,~1<B\leq 1.35,~1.37< A< 1.5, ~1<C<1.23 $ }\vspace{2mm}

Firstly if $B\leq 1.222$, we find $\phi_1^*<0.$ Similarly if $C\leq 1.093$, we find $\phi_2^*<0.$ Therefore we can take $B>1.222$ and $C>1.093$.\\
 Suppose $e+g<0.312a+0.5c$.  Then  we find that $\phi_{19}<0$ for $1.37<A<1.5$, $1.222 <B \leq 1.35, 1.093 \leq C < 1.23$. Therefore
we can take $e+g\geq 0.312a+0.5c$ which  gives $i<0.188a+0.5c$. Now $D^4ABCHI>2$ for $h<c+\frac{a}{2}-\frac{i}{2}$  and $D>1.155$,
 but then $\phi_{16}<0$. Hence we can take $D\leq 1.155.$\vspace{2mm}

If $F\geq 0.904$  we find $\phi_{22}<0$ using  $H>1-b-d+f+\frac{i}{2}, 0<i<0.188a+0.5c$. Therefore we must have $F< 0.904$.
Here $C^2>1.093^2>DE$. If $f+g+h+i\leq 1.5c+0.416a$, we find $\phi_{24}<0$ for $1.222<B\leq 1.35$, $1.37< A< 1.5, ~1.093<C<1.23 $.
Therefore we can take $f+g+h+i> 1.5c+0.416a$
which gives (using (6.15)) $e<0.292a+0.25c$. Now $\phi_{23}<0$ for $ A< 1.5, ~C<1.23 $. This gives a contradiction.\vspace{2mm}

\noindent \textbf{Subcase IV: $\frac{C^3}{DEF}\leq 2,~1<B\leq 1.35,~1.37< A< 1.5, ~1.23<C<1.361 $ }\vspace{2mm}

Firstly if $B\leq 1.222$, we find $\phi_1^*<0.$  Therefore we can take $B>1.222$.\\
 Suppose $e+g<0.34b+0.227c$.  Then  we find that $\phi_{19}<0$ for $1.37<A<1.5$, $1.222 <B \leq 1.35, 1.23 \leq C < 1.361$. Therefore
we can take $e+g\geq 0.34b+0.227c$ which  gives $i<0.66b+0.273c$. Now $D^4ABCHI>2$ for  $D>1.1692$ using $h<b+\frac{c}{2}-\frac{i}{2}$,
 but then $\phi_{17}<0$. Hence we can take $D\leq 1.1692.$\vspace{2mm}

Suppose  $C\geq 1.31$.  Then $\frac{C^3}{DEF}\leq 2$ gives  $D \geq 1.124,~E>0.96,~F>0.96$.  Now  $F^4ABCDE>2$. But then $\phi_{15}<0$.
This gives a contradiction. Hence we can take  $C<1.31$. We repeat the cycle  to get $e+g\geq 0.353b+0.25c$, $i<0.647b+0.25c$.
If $F\geq 0.885$ and $1<D<1.11$ or if $F\geq 0.849$ and $1.11<D<1.1692$, we find $\phi_{22}<0$ using  $H>1-b-d+f+\frac{i}{2}, 0<i<0.647b+0.25c$.
Therefore we must have $F< 0.885$ when $1<D<1.11$ and $F< 0.849$ when $1.11<D<1.1692$ which gives $C<1.2569$. We repeat the cycle once again
to get
$e+g\geq 0.37b+0.26c$, $i<0.63b+0.24c$, $D<1.137$ (using $h<b+\frac{c}{2}-\frac{f}{2}-\frac{i}{2}$), $F<0.851$, $F<0.824$ when
 $1.1<D<1.1692$, $C<1.233$. Finally if $F\geq 0.84$ and $1<D<1.137$ or if $F\geq 0.815$ and $1.1<D<1.137$, we find $\phi_{22}<0$ using
  $H>1-b-d+f+\frac{i}{2}, 0<i<0.63b+0.24c$. Therefore we must have $F< 0.84$ when $1<D<1.1$ and $F< 0.815$ when $1.1<D<1.137$ which gives
   $C<1.23$, a contradiction.\vspace{2mm}

\noindent \textbf{Subcase V: $\frac{C^3}{DEF}\leq 2,~1<B\leq 1.35,~1.155< A< 1.37, ~B\leq C $ }\\

Firstly if $B\leq 1.075$, we find $\phi_1^*<0.$  Therefore we can take $B>1.075$.\\
Suppose $e+g<0.3b+0.192c$.  Then  we find that $\phi_{19}<0$ for $I>1-b-\frac{c}{2}+e+g$, $\max\{1.155,C\}<A<1.37$,
$1.075 <B \leq \min\{1.35,C\}$. Therefore
we can take $e+g\geq 0.3b+0.192c$ which  gives $i<0.7b+0.308c$. Now $D^4ABCHI>2$ for  $D>1.182$ using $h<b+\frac{c}{2}-\frac{i}{2}$,
but then $\phi_{16}<0$. Hence we can take $D\leq 1.182.$ This gives $C<1.3322$.\\
If $F\geq 0.935$ for $1<D<1.09$ or if $F\geq 0.893$ for $1.09<D<1.182$, we find $\phi_{22}<0$ using  $H>1-b-d+f+\frac{i}{2}, 0<i<0.7b+0.308c$.
Therefore we must have $F< 0.935$ when $1<D<1.09$ and $F< 0.893$ when $1.09<D<1.182$ which gives $C<1.2829$.
We repeat the cycle once again to get
$e+g\geq 0.466b+0.25c$, $i<0.534b+0.25c$, $h+i<1.0684b+0.5c$. \vspace{2mm}

If $E\geq 0.912$ for $F<0.935$ or if $E\geq 0.9$ for $F<0.9893$, we find $\phi_{23}<0$ using $g<b+\frac{c}{2}-e-\frac{h+i}{2}$,
$0<h <b+\frac{c}{2}-e-\frac{f+i}{2}$,  $0<i<0.534b+0.25c$. Therefore we must have $E< 0.912$ when $F<0.935$ and $E< 0.9$ when $F<0.9893$
which gives $C<\sqrt[3]{2DEF}<1.2386$. Now $D^4ABCHI>2$ for  $D>1.145$ using
$h<b+\frac{c}{2}-\frac{i}{2}-\frac{f}{2}-e< b+\frac{c}{2}-\frac{i}{2}-\frac{0.065}{2}-0.088$, but then $\phi_{16}<0$.
Hence we can take $D\leq 1.145.$ \\
If $g+h+i<1.6d+0.59a$, we find that $\phi_{25}<0$ for $1<D\leq 1.145.$ Therefore we can take $g+h+i\geq 1.6d+0.59a$ which together
with (6.16) gives $f<b+0.2d-0.295a$. If $B<1.177$, we find that
$\phi_{22}<0$. Therefore can take that $B\geq 1.177$ and hence $A\geq C\geq B\geq 1.177$. Suppose $e+g<0.54b+0.3c$.  Then  we find
that $\phi_{19}<0$ for $1.177<A<1.37$, $1.177 <B \leq \min\{1.35,C\}$. Therefore
we can take $e+g\geq 0.54b+0.3c$ which  gives $i<0.46b+0.2c$, $h+i<0.92b+0.4c$.\vspace{2mm}

If $F\geq 0.864$ for $1<D<1.1$ or if $F\geq 0.84$ for $1.1<D<1.145$, we find $\phi_{22}<0$ using  $H>1-b-d+f+\frac{i}{2}, 0<i<0.46b+0.2c$.
 Therefore we must have $F< 0.864$ when $1<D<1.1$ and $F< 0.84$ when $1.1<D<1.145$ which gives $C<1.192$. Finally $E\geq \frac{C^3}{2DF}>0.847$
  and $\phi_{23}<0$ for   $h+i<0.92b+0.4c$, $1.177\leq B \leq C<1.192$. This gives a contradiction.\\

\noindent \textbf{Subcase VI: $\frac{C^3}{DEF}\leq 2,~1<B\leq 1.35,~1.155< A< 1.37, ~B> C $ }\vspace{2mm}

As in Subcase IV, we get $B>1.075$.  Suppose $e+g<0.2a+0.328c$. Using (6.13) we get $I>1-c-\frac{a}{2}+e+g$.
Then  we find that $\phi_{19}<0$ for $1.155<A<1.37$, $1<C <B \leq \min\{A,1.35\}$. Therefore
we can take $e+g\geq 0.2a+0.328c$ which  gives $i<0.3a+0.672c$. Now $D^4ABCHI>2$ for $h<c+\frac{a}{2}-\frac{i}{2}$  and $D>1.178$,
 but then $\phi_{16}<0$. Hence we can take $D\leq 1.178.$ This implies $C<1.331$.\vspace{2mm}

If $F\geq 0.95$ for $1<D<1.088$ or if $F\geq 0.91$ for $1.088<D<1.178$, we find $\phi_{22}<0$ using
$H>1-\frac{a+b}{2}-d+f+\frac{i}{2}, 0<i<0.3a+0.672c$. Therefore we must have $F< 0.95$ when $1<D<1.088$ and $F< 0.91$ when
$1.088<D<1.178$ which gives $C<1.2895$. Further if $B\leq 1.092$, we find $\phi_1^*<0.$  Therefore we can take $B>1.092$. If
$e+g<0.25a+0.432c$,  we find that $\phi_{19}<0$ for $1 <C \leq \min\{B,1.2895\}$,  $1 <B \leq \min\{A,1.35\}$. Therefore
we can take $e+g\geq 0.25a+0.432c$ which  gives $i<0.25a+0.568c$. Now $D^4ABCHI>2$ for  $D>1.154$ using
$h<c+\frac{a}{2}-\frac{i}{2}-\frac{f}{2}< c+\frac{a}{2}-\frac{i}{2}-\frac{0.05}{2}$, but then $\phi_{16}<0$. Hence we can take $D\leq 1.154.$\vspace{2mm}

If $g+h+i<1.55d+0.6a$, we find that $\phi_{25}<0$ for $1<D\leq 1.154.$ Therefore we can take $g+h+i\geq 1.55d+0.6a$, which together
with (6.16) gives $f<b+0.225d-0.3a$. If $B<1.16$, we find that
$\phi_{22}<0$. Therefore can take that $B\geq 1.16$ and then $ C\geq  1.078$, as $\phi_2^*<0$, if $C\leq 1.078$.\vspace{2mm}

If $F\geq 0.9$ for $1<D<1.095$ or if $F\geq 0.866$ for $1.095<D<1.154$, we find $\phi_{22}<0$ using  $H>1-b-d+f+\frac{i}{2}$. Therefore
we must have $F< 0.9$ when $1<D<1.095$ and $F< 0.866$ when $1.095<D<1.154$ which gives $C<1.2597$.
We repeat the cycle once again to get
$e+g\geq 0.25a+0.5c$, $i<0.25a+0.5c$, $h+i<0.5a+c$. \vspace{2mm}

If $E\geq 0.876$, we find $\phi_{23}<0$ using  $g<c+\frac{a}{2}-e-\frac{h+i}{2}$, $0<h <c+\frac{a}{2}-e-\frac{f+i}{2}$,
$0<i<0.25a+0.5c$ and $F<0.9$. Therefore we must have $E< 0.876$  which gives $C<\sqrt[3]{2DEF}<1.2053$. If $f+g+h+i<1.5c+0.448a$
we find that $C^2>DE$ and $\phi_{24}<0$. Therefore we can take $f+g+h+i\geq 1.5c+0.448a$ which together with (6.15) gives $e<0.25c+0.276a$. Now $\phi_{23}<0$ for $1.078<C<1.2053$. This gives a contradiction.\\

\noindent \textbf{Proposition 45.} Case (14) i.e. $A>1$, $B>1$, $C>1$, $D>1$, $E\leq 1$, $F\leq1$, $G\leq1$, $H\leq1$, $I>1$ does not arise.\vspace{2mm}\\
{\noindent \bf Proof.}
Here $B\leq 2$, $C\leq 1.5$, $D\leq\frac{4}{3}$.
Using the weak inequalities
(1,2,2,1,2,1), ~$(1,2,2,2,1,1)$, ~$(2,2,2,2,1)$, ~(2,1,2,1,2,1), ~(2,1,2,2,1,1), ~$(2,2,2,1,1,1)$ ~and ~(1,2,1,2,1,1,1)  we get
\begin{equation}a+2c-2e-f-2h+i>0,\vspace{-2mm}\end{equation}
\begin{equation}a+2c-2e-2g-h+i>0\vspace{-2mm}\end{equation}
\begin{equation}2b+2d-2f-2h+i>0,\vspace{-2mm}\end{equation}
\begin{equation}2b+c-2e-f-2h+i>0\vspace{-2mm}\end{equation}
\begin{equation}2b+2d-2f-g-h+i>0\vspace{-2mm}\end{equation}
\begin{equation}a+2c-d-2f-g+h+i>0\end{equation}

\noindent\textbf{Claim(i) $D\leq1.155$ and $B\leq1.7325$}\vspace{2mm}

Suppose $D>1.155$. If $B\geq 1.5$, we have $D^4ABCHI>D^4(1.5)^2(\frac{4}{9}D)>2$. Let now $B<1.5$. Using (6.20) we have $h<b+\frac{c+i}{2}$. Also $H>\frac{4}{9}D$. So we have
$$h<\left\{\begin{array}{lll}b+\frac{c+i}{2} & {\rm if }& 0<c+i\leq \frac{10}{9}-\frac{8}{9}d-2b\\\frac{5}{9}-\frac{4}{9}d
& {\rm if }& c+i>\frac{10}{9}-\frac{8}{9}d-2b. \end{array}\right. $$

We find that $D^4ABCHI>(1+d)^4(1+a)(1+b)(1+c+i)(1-h)>2$ for $h,~c+i$ as above and $0<b< \min\{a,0.5\}$, $0.155<d\leq a\leq 1.1327$.
 Also  $\phi_{38}<0$ for $H>0.46873D$, $1\leq I\leq A$ and $1.155<D\leq A<2.1327$.  This gives a contradiction. Now $B\leq \frac{3}{2}D<1.7325$. \vspace{2mm}

\noindent\textbf{Claim(ii) $C\leq1.322$ and $A \leq 1.983$}\vspace{2mm}

Suppose $C>1.322$, then $\frac{C^3}{DEF}>2$. But $\phi_{11}<0$ for $G>0.46873C$, $1\leq I\leq A$, $A\geq B$ and $B<1.7325$. This gives a contradiction. Now $A\leq \frac{3}{2}C<1.983$.\vspace{2mm}

\noindent\textbf{Claim(iii) $B\leq1.533$}\vspace{2mm}

Suppose $B>1.533$, then $\frac{B^3}{CDE}>2$. Now if $F^2>GH$, then $\phi_{8}<0$ for $F\geq 0.46873B$, $1<I\leq A$ and $1.533<B\leq A<1.983$. If $F^2\leq GH$, then $\phi_{3}<0$ using $GH>F^2$,  $F\geq0.46873B$, $1<I\leq A$ and $1.44<B\leq A<1.983$. Hence $B\leq1.533$.\vspace{2mm}

\noindent\textbf{Claim(iv) $A\leq 1.685$}\vspace{2mm}

Suppose $A>1.685$, then $\frac{A^3}{BCD}>2$. Now if $\frac{E^3}{FGH}>2$, then $\phi_{26}<0$ for  $1<I\leq A$, $0.46873A<E\leq 1$ and $1.685<A<1.983$.\\

 Let now $\frac{E^3}{FGH}\leq 2$. If $G\geq 0.6082$ then $\phi_{4}<0$ using $G>0.6082$,  $I>1$, $0.46873A<E\leq 1$ and $1.685<A<1.983$. So we must have $G<0.6082$. This gives $H>\frac{E^3}{2F\times 0.6082}$. If $F>0.691$, we find that $\phi_{27}<0$ for $H>\frac{E^3}{2F\times 0.6082}$, $1<I\leq A$, $0.691<F<1$, $0.46873A<E\leq 1$ and $1.685<A<1.983$. So we must have $F\leq 0.691$. Now $E^2>(0.46873A)^2>FG$ and so (4,3,1,1) holds. But $\phi_{28}<0$ for $H>\frac{E^3}{2\times 0.691\times 0.6082}$, $1<I\leq A$, $0.46873A<E\leq 1$ and $1.685<A<1.983$. This gives a contradiction. Hence we must have $A\leq 1.685$. \vspace{2mm}

\noindent\textbf{Claim(v) $B\leq 1.498$}\vspace{2mm}

Using $A\leq 1.685$ instead of $A<1.983$ and working as in Claim (iii) we can take $B\leq 1.498$.

\noindent\textbf{Claim(vi) $B\leq 1.451$ if $A<1.6$}\vspace{2mm}

Suppose $B>1.451$, then $\frac{B^3}{CDE}>\frac{(1.451)^3}{1.322\times 1.155}>2$.  If $F>0.74$, we see that $\phi_{29}<0$ for $H>\frac{2}{3}F$, $1<I\leq A<1.6$. Therefore we can take $F\leq 0.74$. If $G>0.62$, we find that $\phi_{7}<0$ for $F>0.46873B$, $1<I\leq A<1.6$. Therefore we can take $G<0.62$. Now  $\frac{D^3}{EFG}>\frac{1}{ 0.74\times0.62}>2$.
If $D>1.045$ or if $C<1.23$ we find that  $\phi_{16}<0$  for $B>1.451$, $H>0.46873D$, $1<I\leq A$, $1<A<1.6$. Therefore we can take $D\leq 1.045$ and $C\geq 1.23$.  Now $\frac{C^3}{DEF}>\frac{(1.23)^3}{1.045\times 0.74}>2$. If $A>1.53$ or if $C>1.26$, we find that  $\phi_{11}<0$  for $G>0.46873C$, $1<I\leq A$. Therefore we can take $A\leq 1.53$ and $C\leq 1.26$. But then $\phi_{16}<0$. This gives a contradiction. \vspace{2mm}

\textbf{We divide the range of $A$ into several
 subintervals and  arrive at a contradiction in each.}\vspace{2mm}

\noindent \textbf{Case I: 1.6 $\leq A < $ 1.685}\vspace{2mm}

\noindent If $B\leq 1.26$, we find that $\phi_7^*<0$, so can take $B>1.26$. \vspace{2mm}

\noindent \textbf{Subcase I: $1.26 \leq B <  1.34$}

Here $\frac{A^3}{BCD}>2$. We find that $\phi_4<0$ for $G\geq 0.636$, $1<I\leq A$, $0.46873A<E\leq 1$ and $1.6<A<1.685$. Also $\phi_{27}<0$ for $H>\frac{2}{3}F$, $F\geq 0.782$, $1<I\leq A$, $0.46873A<E\leq 1$ and $1.6<A<1.685$. So we can take $G< 0.636$ and $F< 0.782$. Now $\frac{D^3}{EFG}> \frac{D^3}{0.0.782\times0.636}>2 $ for $D>1$ but $\phi_{17}<0$ for $H>0.46873D$, $1<I\leq A$, $A>1.6$, $1<C<1.322$ and $1<B<1.34$. \vspace{2mm}

\noindent \textbf{Subcase II: $1.34 \leq B <  1.498, ~1<C<1.183, ~1.045<D\leq 1.155$}\vspace{2mm}

Here again $\frac{A^3}{BCD}>2$. As in Subcase I, we get that $G< 0.636$ and $F< 0.782$ and $\frac{D^3}{EFG}> \frac{D^3}{0.0.782\times0.636}>2 $ for $D>1$.   But $\phi_{16}<0$ for $B>1.34$, $H>0.46873D$, $1<I\leq A$, $1.6<A<1.685$, $C<1.183$ and $1.045<D<1.155$. \vspace{2mm}

\noindent \textbf{Subcase III: $1.34 \leq B <  1.498, ~1<C<1.183, ~1<D\leq 1.045$}\vspace{2mm}

Here again $\frac{A^3}{BCD}>2$. As in Subcase I, we get that $G< 0.636$ and $F< 0.782$ and  so $E^2>FG$. But then $\phi_{28}<0$ for $H>0.46873$, $E>0.865$, $1<I\leq A$  and $1.6<A<1.685$. So we can take $E\leq 0.865$. Further $\frac{D^3}{EFG}>2 $ for $D>1$. If $B>1.415$, $\phi_{16}<0$ for $B>1.415$, $H>0.46873D$, $1<I\leq A$, $A>1.6$, $C<1.183$ and $1<D<1.045$. Also if $C<1.12$ and $B>1.34$ we get
$\phi_{16}<0$. Therefore we can take $B\leq 1.415$ and $C\geq 1.12$.\vspace{2mm}

 Now  $\frac{C^3}{DEF}>\frac{(1.2)^3}{1.045\times 0.865\times 0.782}>2$. We see that $\phi_{11}<0$ for $G>0.46873C$, $1<I\leq A$, if $A>1.65,~B<1.415,~C>1.12$ or if $A>1.6,~B<1.39,~C>1.12$ or if $A>1.6,~B<1.415,~C>1.155$. Therefore we can assume that $A<1.65$, $B>1.39$ and  $C<1.155$. Now $\phi_{16}<0$. This gives a contradiction. \vspace{2mm}

 \noindent \textbf{Subcase IV: $1.34 \leq B <  1.395, ~1.183<C<1.322$}\vspace{2mm}

If $E>0.945$ then $E^4ABCDI>(0.945)^4(1.6)(1.34)(1.183)>2$. But then $\phi_8^*<0$  for $1<I\leq A$, $1.6\times1.34\times1.183<x=ABCD<1.685\times1.498\times1.322\times1.155$, $0.945<E<1$. Therefore we can take $E\leq 0.945$. If $D>1.11$ then
$D^4ABCHI>(1.11)^4(1.6)(1.34)(1.183)(0.46873\times 1.11)>2$, but then $\phi_{16}<0$.  So we can take $D\leq 1.11$. Now $\frac{A^3}{BCD}>2$. As in Subcase I and III, we get that $G< 0.636$ and $F< 0.782$ and  $E\leq 0.865$.

Now  $\frac{C^3}{DEF}>\frac{(1.183)^3}{1.11\times 0.865\times 0.782}>2$. We see that $\phi_{11}<0$ for $G>0.46873C$, $1<I\leq A$, $1.6<A<1.685,~1.34<B<1.395,~1.183C<1.322$. This gives a contradiction. \vspace{2mm}

 \noindent \textbf{Subcase V: $1.395 \leq B <  1.498, ~1.183<C<1.322$}\vspace{2mm}

As in Subcase IV, we can take $E<0.945$ and $D<1.11$.  If $C>1.281$, we have  $\frac{C^3}{DEF}>\frac{(1.281)^3}{1.11\times 0.945}>2$ and $\phi_{11}<0$. Therefore we can take $C\leq 1.281$. Now $\frac{B^3}{CDE}>\frac{(1.395)^3}{1.281\times 1.11\times 0.945}>2$.  If $F>0.76$, we see that $\phi_{29}<0$ for $H>\frac{2}{3}F$, $1<I\leq A$. Therefore we can take $F<0.76$. If $G>0.65$, we find that $\phi_{7}<0$ for $F>0.46873B$, $1<I\leq A$. Therefore we can take $G<0.65$. Now $\frac{C^3}{DEF}>\frac{(1.183)^3}{1.11\times 0.945\times0.76}>2$ and $\frac{D^3}{EFG}>\frac{1}{0.945\times 0.76\times0.65}>2$.
If $C>1.21$ or if $B<1.45$ we find that  $\phi_{11}<0$. Therefore we can take $C\leq 1.21$ and $B\geq 1.45$. But then $\phi_{16}<0$. This gives a contradiction. \vspace{2mm}

\noindent \textbf{Case II: 1.5 $\leq A < $ 1.6}\vspace{2mm}

\noindent If $B\leq 1.22$, we find that $\phi_7^*<0$, so can take $B>1.22$. \vspace{2mm}

\noindent \textbf{Subcase I: $1.22 \leq B <  1.27$}\vspace{2mm}

Working as in Claim (i), we find that $D^4ABCHI>2$ for $h<\min\{b+\frac{c+i}{2}, \frac{5}{9}-\frac{4}{9}d\}$, $b<0.27$ and $D>1.098$. Also  $\phi_{38}<0$ for $H>0.46873D$, $1\leq I\leq A$ and $1<D\leq A<1.6$.  This gives a contradiction. So we can take $D\leq 1.098$.\vspace{2mm}

If $C>D$ then $(2,6^*,1)$ holds, but $\phi_5^*<0$ for $A>1.5$ and $B<1.27$. Therefore we can take $C\leq D\leq 1.098$. But then  $\frac{A^3}{BCD}>2$. As in Subcase I of Case I, we can take  that $G< 0.665$ and $F< 0.82$. Also we find that $\phi_4<0$ for $E\geq 0.89$, $G\geq \frac{2}{3}E$, $1<I\leq A$ and $1.5\leq A<1.6$. So we can take $E<0.89$. Now $\frac{D^3}{EFG}>2$ for $D>1$ and $\phi_{17}<0$ for $H>0.46873D$, $1\leq I\leq A$ and $1<C\leq D\leq 1.098$.  This gives a contradiction. \vspace{2mm}

\noindent \textbf{Subcase II: $1.27 \leq B <  1.35, ~1.22<C\leq 1.322$}\vspace{2mm}

Here $C>D$. If $C>I$, then $(2,7^*)$ holds but $\phi_3^*<0$ for $A>1.5$ and $B<1.35$. So we can take $I\geq C$. If $C>1.28$, we find that  $E^4ABCDI>2$ for $E>0.895$; and if $1.22<C\leq 1.28$, $E^4ABCDI>2$ for $E>0.917$. Also we find that $\phi_{10}<0$ for $1<I\leq A$, $A>1.5$ and $C,~E$ satisfying the respective bounds.  Therefore we can take $E\leq 0.895$ if $C>1.28$ and $E\leq 0.917$ if $1.22<C\leq 1.28$. Further if  $\frac{C^3}{DEF}>2$ we find that  $\phi_{11}<0$ for $1.22<C\leq 1.322$. Also $\frac{C^3}{DEF}>2$ for $C>1.28$ and $E\leq 0.895$. Therefore we can take $\frac{C^3}{DEF}\leq 2$ and are left with the case  $1.22<C\leq 1.28$ and so $E\leq 0.917$. Now $D^4ABCHI>2$, $\phi_{17}<0$ for $H>0.46873D$ and $D>1.086$; So we can take $D<1.086$. But then $F>\frac{C^3}{2DE}> \frac{1.22^3}{2\times 1.086\times 0.917}>0.911$, $G>\frac{3}{4}F>0.683$, $H> \frac{2}{3}F>0.607$ and $\frac{C^3}{DEF}=C^4ABGHI>2$. This gives a contradiction. \vspace{2mm}

\noindent \textbf{Subcase III: $1.27 \leq B <  1.35, ~1.125<C\leq 1.22,~\frac{A^3}{BCD}\leq 2 $}\vspace{2mm}

Working as in Claim (i), we find that $D^4ABCHI>2$ for $h<\min\{b+\frac{c+i}{2}, \frac{5}{9}-\frac{4}{9}d\}$, $b<0.35$ and $D>1.1104$. Also  $\phi_{17}<0$ for $H>0.46873D$, $1\leq I\leq A$ and $1.1104<D\leq 1.155$.  So we can take $D\leq 1.1104$.
Here $C>D$ and working as in Subcase II, we can take $I\geq C$. Again $E^4ABCDI>2$ for $E>0.955$ and $C>1.125$ but $\phi_{10}<0$ for $1<I\leq A$, $A>1.5$. Therefore we can take $E \leq 0.955$ i.e. $e>0.045$. If $F>0.885$, we see that $\phi_{22}<0$ for $H> \frac{2}{3}F$, $B>1.27$, $C<I<A$, $1.125<C<1.22$ and $1.5<A<1.6$. Therefore we can take $F<0.885$ i.e. $f>0.115$.  Now using $h<b+\frac{c+i}{2}-e-\frac{f}{2}$ i.e. $h<\min\{b+\frac{c+i}{2}-0.1075, \frac{5}{9}-\frac{4}{9}d\}$, we find that $D^4ABCHI>2$ and $\phi_{17}<0$ for $H>0.46873D$,  $D>1.071$, $C\leq I\leq A$. So we can take $D<1.071$.\vspace{2mm}

 Further $\frac{A^3}{BCD}\leq 2 $ implies $A<1.523$, $D>1.024$, $C>1.167$, $B>1.291$. With bounds on $A,~B,C,~D$ improved, we repeat the cycle to get
$E \leq 0.93$ i.e. $e>0.07$, $F<0.873$ i.e. $f>0.127$ and $D<1.06$. If $C>1.199$, $\frac{C^3}{DEF}>2$ and $\phi_{11}<0$. Therefore can take $C\leq 1.199$. This gives, using $\frac{A^3}{BCD}\leq 2 $, that $A<1.509$, $D>1.042$, $C>1.179$, $B>1.327$. Repeating the cycle once again we get $E \leq 0.912$ i.e. $e>0.088$, $F<0.865$ i.e. $f>0.135$ and $D<1.053$. But then $A\leq (2BCD)^{1/3}<1.5$. This gives a contradiction.\vspace{2mm}

\noindent \textbf{Subcase IV: $1.27 \leq B <  1.35, ~1.125<C\leq 1.22,~\frac{A^3}{BCD}>2 $}\vspace{2mm}

We start as in Subcase I of Case I. As $\frac{A^3}{BCD}>2$, we find that $\phi_4<0$ for $G\geq 0.665$, $1<I\leq A$, $0.46873A<E\leq 1$  or for $E\geq 0.89$, $G>\frac{2}{3}E$, $1<I\leq A$ and $1.5<A<1.6$.  Also $\phi_{27}<0$ for $H>\frac{2}{3}F$, $F\geq 0.82$, $1<I\leq A$, $0.46873A<E\leq 1$ and $1.5<A<1.6$. So we can take $G< 0.665$, $E<0.89$ and $F< 0.82$. Now $\frac{D^3}{EFG}> 2 $ for $D>1$ but $\phi_{17}<0$ for $H>0.46873D$, $1<I\leq A$, $A>1.5$, $1<C\leq 1.22$ and $1.27<B<1.35$.  \vspace{2mm}

\noindent \textbf{Subcase V: $1.27 \leq B <  1.35, ~1<C\leq 1.125 $}\vspace{2mm}

Working as in Subcase III, we can take $D\leq 1.1104$. Then $\frac{A^3}{BCD}>2$ for $A>1.5,~B<1.35,~C<1.125$.  We get a contradiction working as in Subcase IV.\vspace{2mm}

\noindent \textbf{Subcase VI: $1.35 \leq B <  1.4, ~\frac{B^3}{CDE}>2  $}\vspace{2mm}

We find that $\phi_{29}<0$ if $F>0.77$ and using $H>\frac{2}{3}F$ or if $H>0.69$ and using $F\geq \frac{2}{3}$ and for  $1<I\leq A$, $1.5<A<1.6,~1.35 \leq B <  1.4$ . Therefore we can take $F\leq 0.77$ and $H\leq 0.69$. Further if $G>0.65$, we find that $\phi_{7}<0$ for $F\geq \frac{2}{3}$, $1<I\leq A$, $1.5<A<1.6,~1.35 \leq B <  1.4$. Therefore we can take $G<0.65$. Now $\frac{E^3}{FGH}>2$ and $\phi_8^*<0$  for $1<I\leq A$, $1.5\times1.35<x=ABCD<1.6\times1.4\times1.322\times1.155$ and for $0.92<E<1$. Therefore we can take $E\leq 0.92$. \vspace{2mm}

Now  $\frac{C^3}{DEF}>2$ and  $\phi_{11}<0$ for $C>1.186$, $G>0.46873C$,  $1\leq I\leq A$, $1.5<A\leq 1.6$ and $1.35<B<1.4$. Therefore we can take $C\leq 1.186$. Further $\frac{D^3}{EFG}>2$ for $D>1, E\leq 0.92, ~F\leq 0.77~,G\leq 0.65$.  If $D>1.025$ or if $C<1.165$ or if $A<1.55$ or if $B>1.38$, we find that  $\phi_{16}<0$  for  $H>0.46873D$, $1<I\leq A$. Therefore we can take $D\leq 1.025$, $C\geq 1.165$, $A\geq 1.55$ and $B\leq 1.38$.  But then  $\frac{C^3}{DEF}>2$ and  $\phi_{11}<0$  for $C>1.165$, $G>0.46873C$,  $1\leq I\leq A$. This gives a contradiction. \vspace{2mm}

\noindent \textbf{Subcase VII: $1.35 \leq B <  1.4, ~\frac{B^3}{CDE}\leq 2  $}\vspace{2mm}

If $F>0.905$, we see that $\phi_{22}<0$ for $H> \frac{2}{3}F$, $B>1.35$, $D>1$, $1<I\leq A$, $1<C<1.322$ and $1.5<A<1.6$. Therefore we can take $F<0.905$ i.e. $f>0.095$. If $C\leq 1.05$, we find that $\phi_9^*<0$, so can take $C>1.05$. Now if $E>0.985$, $E^4ABCDI>E^4(1.5)(1.35)(1.05)>2$ and $\phi_8^*<0$  for $1<I\leq A$. Therefore we can take $E\leq 0.985$ i. e. $e>0.015$. Now using $h<b+\frac{c+i}{2}-e-\frac{f}{2}$ i.e. $h<\min\{b+\frac{c+i}{2}-0.0625, \frac{5}{9}-\frac{4}{9}d\}$, we find that $D^4ABCHI>2$ and $\phi_{16}<0$ for $H>0.46873D$,  $D>1.095$, $1\leq I\leq A$. So we can take $D<1.095$. Further if $C>1.25$, then $\frac{C^3}{DEF}>2$ and $\phi_{11}<0$ for $G>0.46873C$, $1\leq I\leq A$. Therefore we can take $C\leq 1.25$. With better bound on $C$, we find $\phi_{22}<0$ for $F\geq 0.895$. Thus can take $F<0.895$ i.e. $f>0.105$. \vspace{2mm}

As $\frac{B^3}{CDE}\leq 2$, we have $B<\sqrt[3]{2CDE}<1.392$, $C\geq \frac{B^3}{2DE}>1.14056$.  We repeat the cycle to get $E<0.9647$ i.e. $e>0.03
53$, $D<1.0862$, $C\leq 1.234$. Then $B<1.373$, $C\geq 1.174$, Also $D\geq \frac{B^3}{2CE}>1.033$. Now if $E>0.95$, $E^4ABCDI>E^4(1.5)(1.35)(1.174)(1.033)>2$ and $\phi_8^*<0$  for $1<I\leq A$. Therefore we can take $E\leq 0.95$ i. e. $e>0.05$.  $D<1.079$ and $C\leq 1.225$. We repeat the cycle once again to get $B<1.3593$, $C\geq 1.2$, $D\geq \frac{B^3}{2CE}>1.057$. $E<0.9394$ i.e. $e>0.0606$, $D<1.073$ and $C\leq 1.218$. But then  $B<\sqrt[3]{2CDE}<1.3491<1.35$, a contradiction.\vspace{2mm}

\noindent \textbf{Subcase VIII: $1.4 \leq B <  1.451 $}\vspace{2mm}

If $C\leq 1.085$, we find that $\phi_9^*<0$, so can take $C>1.085$. If $F>0.886$, we see that $\phi_{22}<0$ for $H> \frac{2}{3}F$, $B>1.4$, $D>1$, $1<I\leq A$, $1.085<C<1.322$ and $1.5<A<1.6$. Therefore we can take $F<0.886$ i.e. $f>0.114$.  Now if $E>0.968$, $E^4ABCDI>E^4(1.5)(1.4)(1.085)>2$ and $\phi_8^*<0$  for $1<I\leq A$. Therefore we can take $E\leq 0.968$ i. e. $e>0.032$. Now using $h<b+\frac{c+i}{2}-e-\frac{f}{2}$ i.e. $h<\min\{b+\frac{c+i}{2}-0.089, \frac{5}{9}-\frac{4}{9}d\}$, we find that $D^4ABCHI>2$ and $\phi_{16}<0$ for $H>0.46873D$,  $D>1.096$, $1\leq I\leq A$. So we can take $D<1.096$. Further if $C>1.235$, then $\frac{C^3}{DEF}>2$ and $\phi_{11}<0$ for $G>0.46873C$, $1\leq I\leq A$. Therefore we can take $C\leq 1.235$.\vspace{2mm}

Now we see that  $\frac{B^3}{CDE}>\frac{1.4^3}{1.235\times 1.096\times 0.968}>2$. Working as in Subcase VI, we can take $F\leq 0.77$, $H\leq 0.69$ and $G<0.65$. Then $\frac{E^3}{FGH}>2$ for $E>0.93$ and $\phi_8^*<0$  for $1<I\leq A$, $(1.5)(1.4)( 1.085)<x=ABCD<(1.6)(1.451)(1.235)(1.1096)$ and for $0.93<E<1$. Therefore we can take $E\leq 0.93$. Now if $C>1.22$, then $\frac{C^3}{DEF}>2$ and $\phi_{11}<0$ for $G>0.46873C$, $1\leq I\leq A$. Therefore we can take $C\leq 1.22$. But then $\frac{D^3}{EFG}>\frac{1}{0.93\times 0.77\times 0.65}>2$ for $D>1$, whereas $\phi_{16}<0$ for either $D>1.103$ or for $D>1$ and $C<1.2$ or for $D>1$ and $A<1.55$ or for $D>1$ and $B>1.425$. Therefore we can take $D\leq 1.03$, $C\geq 1.2$, $A\geq 1.55$ and $B\leq 1.425$. But then $\frac{C^3}{DEF}>2$ and $\phi_{11}<0$. This gives a contradiction.\vspace{2mm}

\noindent \textbf{Case III: 1.4 $\leq A < $ 1.5}\vspace{2mm}

\noindent If $B\leq 1.18$, we find that $\phi_7^*<0$, so can take $B>1.18$. \vspace{2mm}

\noindent \textbf{Subcase I: $1.18 \leq B <  1.29, ~ 1.23\leq C<1.322$}\vspace{2mm}

Here $C>D$. If $C>I$, then $(2,7^*)$ holds but $\phi_3^*<0$ for $A>1.4$ and $B<1.29$. So we can take $I\geq C$. As $C\geq 1.23$, we find that  $E^4ABCDI>2$ for $E>0.946$; and   $\phi_{10}<0$.  Therefore we can take $E\leq 0.946$ i.e. $e>0.054$. If $F>0.925$, we see that $\phi_{22}<0$ for $H> \frac{2}{3}F$, $B>1.18$, $C<I<A$, $1.23<C<1.322$ and $1.4<A<1.5$. Therefore we can take $F<0.925$ i.e. $f>0.075$.  Further if  $\frac{C^3}{DEF}>2$ we find that  $\phi_{11}<0$. Therefore we can take $\frac{C^3}{DEF}\leq 2$. This gives $C<1.265$, $D>1.063$. Again we find that  $E^4ABCDI>2$ for $E>0.932$; and   $\phi_{10}<0$.  Therefore we can take $E\leq 0.932$ i.e. $e>0.068$. Now using $h<b+\frac{c+i}{2}-e-\frac{f}{2}$ i.e. $h<\min\{b+\frac{c+i}{2}-0.1075, \frac{5}{9}-\frac{4}{9}d\}$, we find that $D^4ABCHI>2$ and $\phi_{17}<0$ for $H>0.46873D$,  $D>1.08$, $C\leq I\leq A$. So we can take $D<1.08$. But then $\phi_{22}<0 $ for $F>0.9$; so we can take $F\leq 0.9$. Now  $C<(2DEF)^{\frac{1}{3}}<1.23$. This gives a contradiction. \vspace{2mm}

\noindent \textbf{Subcase II: $1.18 \leq B <  1.29, ~ 1<C<1.23, ~\frac{A^3}{BCD}>2 $}\vspace{2mm}

We find that $\phi_4<0$ for $G\geq 0.65$, $E\geq\frac{3}{4}D>\frac{3}{4}$ or for $G\geq \frac{2}{3}E$ and $E>0.92$. Therefore we can take $G\leq 0.65$ and $E\leq 0.92$.  Also $\phi_{27}<0$ for $H>\frac{2}{3}F$, $F\geq 0.8$, $1<I\leq A$ and $1.4<A<1.5$. So we can take  $F< 0.8$. Now $\frac{D^3}{EFG}>\frac{1}{0.92\times 0.8\times 0.65}>2$ for $D>1$, whereas $\phi_{17}<0$ for $1<I\leq A$, $A>1.4$, $1<C<1.23$, $0.46873D<H<1$, $1<D<1.155$ and $1.18\leq B<1.29$. This gives a contradiction.\vspace{2mm}

\noindent \textbf{Subcase III: $1.18 \leq B <  1.29, ~ 1.14<C<1.23, ~\frac{A^3}{BCD}\leq 2,~i<0.85a $}\vspace{2mm}

Here if  $\frac{C^3}{DEF}>2$ we find that  $\phi_{11}<0$ for $0<i\leq a$, $1<C<1.23$. Therefore we can take $\frac{C^3}{DEF}\leq 2$. Further if  $\frac{D^3}{EFG}>2$ we find that  $\phi_{17}<0$ for $0<i\leq a$. Therefore we can take $\frac{D^3}{EFG}\leq 2$. Also using  $h<\min\{b+\frac{c+i}{2}, \frac{5}{9}-\frac{4}{9}d\}$, we find that $D^4ABCHI>2$  for   $D>1.1181$. So we can take $D<1.1181$. This gives $C>D$. Now $\phi_5^*<0$ for $A>1.4$, $B<1.243$ and $i<0.85a$. Therefore we can take $B\geq 1.243$. If $C>I$, then $(2,7^*)$ holds but $\phi_3^*<0$ for $A>1.4$ and $B<1.29$. So we can take $I\geq C$. As $C\geq 1.14, B \geq 1.243$, we find that  $E^4ABCDI>2$ for $E>0.9698$; and   $\phi_{10}<0$.  Therefore we can take $E\leq 0.9698$ i.e. $e>0.0302$.\vspace{2mm}

If $F>0.91$, we see that $\phi_{22}<0$ for $H> \frac{2}{3}F$, $B>1.243$, $C<I\leq A$, $1.14<C<1.23$ and $1.4<A<1.5$. Therefore we can take $F<0.91$.
Further making use of $h<b+d-f+\frac{i}{2}$, from (), in place of $H> \frac{2}{3}F$ and noting that $\phi_{22}^{(1)}$ is an increasing function of $C$ here, we find that $\phi_{22}^{(1)}<0$ for $c<0.23$, $0.09<f\leq 0.23$, $0<d<0.1181$ and $0<i<0.85a$. Therefore we can take $f>0.23$. Now using $h<b+\frac{c+i}{2}-e-\frac{f}{2}$ i.e. $h<\min\{b+\frac{c+i}{2}-0.1452, \frac{5}{9}-\frac{4}{9}d\}$, we find that $D^4ABCHI>2$ for $D>1.069$. So we can take $D<1.069$. Further $\frac{C^3}{DEF}\leq 2$ gives $C<1.1688$. Repeating the cycle we get  $f>0.258$, $D<1.064$, $C<1.1527$ and $E>\frac{C^3}{2DF}>0.938$. Now $\phi_{19}^{(2)}$ is a decreasing function of $G,E,C$ and of $A$. For $E>0.938$, $G>0.72$, $C>1.14$, $A>1.4$, $1<I\leq A$, we find that $\phi_{19}^{(2)}<0$ for $1.243<B<1.29$ and $1<D<1.064$. Therefore we must have $G \leq 0.72$. But now $D^3<2EFG$ implies $D<1.012$ and $C^3<2DEF$ implies $C<1.136$. This gives a contradiction.\vspace{2mm}

\noindent \textbf{Subcase IV: $1.18 \leq B <  1.29, ~ 1.14<C<1.23, ~\frac{A^3}{BCD}\leq 2,~i\geq 0.85a $}\vspace{2mm}

As in Subcase III, we can take $\frac{C^3}{DEF}\leq 2$, $\frac{D^3}{EFG}\leq 2$, $D<1.1181$. This gives $C>D$. Now $\phi_5^*<0$ for $A>1.4$, $B<1.232$ and $i<a$. Therefore we can take $B\geq 1.232$. Now we find that  $E^4ABCDI>2$ for $I>1+0.85a>1+0.85\times 0.4$ and $E>0.934$; but  $\phi_{10}<0$.  Therefore we can take $E\leq 0.934$ i.e. $e>0.066$.

Again working as in Subcase III  we can take firstly $F<0.89$
and then  $f>0.215$. Now using $h<b+\frac{c+i}{2}-e-\frac{f}{2}$ i.e. $h<\min\{b+\frac{c+i}{2}-0.1452, \frac{5}{9}-\frac{4}{9}d\}$, we find that $D^4ABCHI>2$ for $D>1.059$. So we can take $D<1.059$. Further $\frac{C^3}{DEF}\leq 2$ gives $C<1.1581$. Repeating the cycle we get  $f>0.25$
and $E>\frac{C^3}{2DF}>0.891$. For $E>0.831$, $G>0.71$,$C>1.14$, $A>1.4$, $1<I\leq A$, we find that $\phi_{19}^{(2)}<0$ for $1.232\leq B<1.29$ and $1<D<1.059$. Therefore we must have $G \leq 0.71$. But now $D^3<2EFG$ implies $D<1$. This gives a contradiction.\vspace{2mm}

\noindent \textbf{Subcase V: $1.18 \leq B <  1.29, ~ 1<C<1.14, ~\frac{A^3}{BCD}\leq 2,~i< 0.85a $}\vspace{2mm}

As in Subcase III, we can take $\frac{C^3}{DEF}\leq 2$, $\frac{D^3}{EFG}\leq 2$, $D<1.1181$.  Now $\phi_7^*<0$ for $C<1.14$, $i<0.85a$, $1.4<A<1.5$, and $B<1.222$. Therefore we can take $B\geq 1.222$. If $F>0.95$, we see that $\phi_{22}<0$ for $H> \frac{2}{3}F$, $B>1.222$, $0<i<0.85a$, $1<C<1.14$ and $1.4<A<1.5$. Therefore we can take $F<0.95$. Further making use of $h<b+d-f+\frac{i}{2}$ we find that $\phi_{22}^{(1)}<0$ for $c<0.14$, $0.05<f\leq 0.24$, $0<d<0.1181$ and $0<i<0.85a$. Therefore we can take $f>0.24$. Now using $h<b+\frac{c+i}{2}-e-\frac{f}{2}$ i.e. $h<\min\{b+\frac{c+i}{2}-0.1452, \frac{5}{9}-\frac{4}{9}d\}$, we find that $D^4ABCHI>2$ for $D>1.076$. So we can take $D<1.076$. Moreover $\frac{A^3}{BCD}\leq 2$ gives $A<1.469$.

If $e+g\leq 0.71c+0.3(a+i)$, we find that $\phi_{19}^{(1)}<0$ for $0<i<0.85a$, $0<c<0.14$, $0<b<0.29$ and $0.4<a<0.469$. So we can take $e+g> 0.71c+0.3(a+i)$ which together with () gives $h<0.58c+0.4(a+i)$. Now $\phi_{23}<0$ for $b<0.29$, $0<g<c+\frac{a+i}{2}-e-\frac{h}{2}, ~0<h<0.58c+0.4(a+i)$,$0<e<0.1,~0<i<0.85a$, $0.4<a<0.469$ and $0<c<0.14$. Therefore we can take $e>0.1$. Now using $h<c+\frac{a+i}{2}-e-\frac{f}{2}<c+\frac{a+i}{2}-0.1-\frac{0.24}{2}$ we find that $\frac{D^3}{EFG}=D^4ABCHI>2$ for $0<i<0.85a$, $0.4<a<0.469$, $0<c<0.14$ and $d>0.04$. Therefore we can take $D\leq 1.04$. Then $\frac{C^3}{DEF}\leq 2$ implies $C<1.125$ and $\frac{A^3}{BCD}\leq 2$ gives $A< 1.446$. We repeat the cycle to get $e+g> 0.73c+0.311(a+i)$, $h<0.54c+0.378(a+i)$, $e>0.142$, $C<1.107$, $D\leq 1.026$, $C<1.102$ and $A< 1.429$.

Further making use of $h<b+d-f+\frac{i}{2}$ we find that $\phi_{22}^{(1)}<0$ for $c<0.102$, $0.24<f\leq 0.254$, $0<d<0.026$ and $0<i<0.85a$. Therefore we can take $f>0.254$, i.e. $F<0.746$ which gives $C<1.096$ and $A<1.427$. We repeat the cycle once again to get $e+g> 0.74c+0.33(a+i)$, $h<0.52c+0.34(a+i)$, $e>0.165$, $C<1.086$. This gives  $B>\frac{A^3}{2CD}>1.2313$. Now $\frac{B^3}{CDE}> 2$ and $\phi_7<0$ for $F>\frac{2}{3}$, $G>0.6$, $0<i<0.85a$, $0.4<a<0.427$ and $B<1.29$. Therefore we can take $G<0.6$. But then $D<(2EFG)^{\frac{1}{3}}<1$.  This gives a contradiction.

\noindent \textbf{Subcase VI: $1.18 \leq B <  1.29, ~ 1<C<1.11, ~\frac{A^3}{BCD}\leq 2,~i\geq 0.85a $}\vspace{2mm}

As in Subcase III, we can take $\frac{C^3}{DEF}\leq 2$, $\frac{D^3}{EFG}\leq 2$, $D<1.1181$. Also $\frac{A^3}{BCD}\leq 2$ gives $A< 1.474$. Now $\phi_7^*<0$ for $A>1.4$, $C<1.11$, $B<1.213$ and $i<a$. Therefore we can take $B\geq 1.213$. Now we find that  $E^4ABCDI>2$ for $I>1+0.85a>1+0.85\times 0.4$ and $E>0.9683$; but  $\phi_{10}<0$.  Therefore we can take $E\leq 0.9683$ i.e. $e>0.0317$.

Again working as in Subcase III first we can take $F<0.891$.
and then  $f>0.258$. Now using $h<\frac{a+i}{2}+c-e-\frac{f}{2}<\frac{a+i}{2}+c-0.1607$, we find that $D^4ABCHI>2$ for $D>1.044$. So we can take $D<1.044$.
Further $\frac{A^3}{BCD}\leq 2$ gives $A<1.441$. If $e+g\leq 0.71c+0.305(a+i)$, we find that $\phi_{19}^{(1)}<0$ for $0<i<a$, $0<c<0.11$, $0<b<0.29$ and $0.4<a<0.441$. So we can take $e+g> 0.71c+0.305(a+i)$ which together with () gives $h<0.58c+0.39(a+i)$. Now $\phi_{23}<0$ for $b<0.29$, $0<g<c+\frac{a+i}{2}-e-\frac{h}{2}, ~0<h<0.58c+0.39(a+i)$,$0<e<0.095,~0<i<a$, $0.4<a<0.441$ and $0<c<0.11$. Therefore we can take $e>0.095$. Again using $h<c+\frac{a+i}{2}-e-\frac{f}{2}<c+\frac{a+i}{2}-0.095-\frac{0.258}{2}$ we find that $\frac{D^3}{EFG}=D^4ABCHI>2$ for  $0.4<a<0.441$, $0<c<0.11$ and $d>0.016$. Therefore we can take $D\leq 1.016$.

If $\frac{B^3}{CDE}\leq 2$ we get $B<1.269$ and so $A<1.4199$.
Repeating the cycle we get  $E<0.885$, $D<1.009$, $C<1.099$, $B<1.236$ and so $A<1.4$.
 This gives a contradiction.

 If $\frac{B^3}{CDE}> 2$, we find that $\phi_7<0$ for $G>0.68$, $F>\frac{2}{3}$, $0.85a<i\leq a$. Therefore we can take $G<0.68$. But then $D< (2EFG)^{\frac{1}{3}}<1$. This gives a contradiction.\vspace{2mm}

\noindent \textbf{Subcase VII: $1.18 \leq B <  1.29, ~ 1.11<C<1.14, ~\frac{A^3}{BCD}\leq 2,~i\geq 0.85a $}\vspace{2mm}

As in Subcase III, we can take $\frac{C^3}{DEF}\leq 2$, $\frac{D^3}{EFG}\leq 2$, $D<1.1181$. Also $\frac{A^3}{BCD}\leq 2$ gives $A< 1.488$. Now $\phi_7^*<0$ for $A>1.4$, $1.11<C<1.14$, $B<1.21$ and $i<a$. Therefore we can take $B\geq 1.21$. Now we find that  $E^4ABCDI>2$ for $I>1+0.85a>1+0.85\times 0.4$ and $E>0.9439$; but  $\phi_{10}<0$.  Therefore we can take $E\leq 0.9439$ i.e. $e>0.0561$.

Again working as in Subcase III first we can take $F<0.891$.
and then  $f>0.249$. Now using $h<\frac{a+i}{2}+c-e-\frac{f}{2}$, we find that $D^4ABCHI>2$ for $D>1.043$. So we can take $D\leq 1.043$.
Further  $\frac{A^3}{BCD}\leq 2$ gives $A<1.4531$ and $G>\frac{D^3}{2EF}$ gives $G>0.705$
If $E>0.91$ we find $\phi_{19}^{(2)}<0$ for $G>0.705$,$C>1.11$, $1<I<A$, $A>1.4$. Therefore we can take $E\leq 0.91$. As $\frac{C^3}{DEF}\leq 2$, we get $C<1.126$, and then $A<1.4471$. Repeating the cycle we get  $D<1.023$, $G>0.73$, $E<0.86$. But now $C^3<2DEF$ implies $C<1.1$. This gives a contradiction.\vspace{2mm}

\noindent \textbf{Subcase VIII: $1.29 \leq B <  1.35, 1<C<1.322, ~i\leq 0.6a $}\vspace{2mm}

\noindent {Claim (i)} $C>1.092$. If $C\leq 1.092$, we find that $\phi_9^*<0$ for $0<i\leq 0.6a$.

\noindent {Claim (ii)} $F<0.935$. If $F>0.935$, we see that $\phi_{22}<0$ for $H> \frac{2}{3}F$, $B>1.29$, $1<I\leq A$, $D>1$, $1<C\leq 1.322$ and $1.4<A<1.5$.

\noindent{Claim (iii)} $D\leq 1.1165$ using $h<\min\{b+\frac{c+i}{2}-\frac{f}{2}, \frac{5}{9}-\frac{4}{9}d\}$, and proving  $D^4ABCHI>2$ and $\phi_{17}<0$.

\noindent {Claim (iv)} $C\leq 1.2781$  proving  $\frac{C^3}{DEF}>2$ and $\phi_{11}<0$.

\noindent {Claim (v)} $f>0.25$. Using $h<b+d-f+\frac{i}{2}$, in place of $H> \frac{2}{3}F$  we find that $\phi_{22}^{(1)}<0$ for $c<0.2781$, $0.065<f\leq 0.25$, $0<d<0.1165$ and $0<i<0.6a$.

\noindent {Claim (vi)} $C\leq 1.188$  proving  $\frac{C^3}{DEF}>2$ and  $\phi_{11}<0$.

\noindent {Claim (vii)} $D\leq 1.062$ using $h<c+\frac{a+i}{2}-\frac{f}{2}$ and proving  $D^4ABCHI>2$ and $\phi_{17}<0$. Further $C<1.167$.

\noindent {Claim (viii)} $g+h\geq 1.5d+0.526(a+i)$  proving  $\phi_{25}<0$ for  $0<i\leq 0.6a$, $0<d\leq 0.062$, $0.4<a<0.5$. This gives $f<b+d+\frac{i}{2}-\frac{1.5d+0.526(a+i)}{2}$.

\noindent {Claim (ix)} $i>0.4a$ and $b>0.32$. If $i<0.4a$ and $b<0.35$ or $i<0.6a$ but $b<0.32$ we find that $\phi_{22}^{(1)}<0$ for $c<0.167$, $0.065<f\leq 0.25$, $0<d<0.062$

\noindent {Claim (x)} $f>0.28$, $D<1.04$, $C<1.145$. (by repeating claims (v),(vi) and (vii)).

\noindent {Claim (xi)}$E<0.962$ for otherwise $E^4ABCDI>2$ using $i>0.4a$, $b>0.32$ and $\phi_{10}<0$.

\noindent {Claim (xii)} $G<0.69$. If $G\geq 0.69$, $\phi_7<0$ for $F\geq \frac{2}{3}, ~0.4a<i<a$

\noindent {Final contradiction :} $\frac{D^3}{EFG}>2$ for $D>1$ but $\phi_{17}<0$. \vspace{2mm}

\noindent \textbf{Subcase IX: $1.29 \leq B <  1.35, ~i> 0.6a,~1.17<C<1.322 $}\vspace{2mm}

\noindent {Claim (i)} $F<0.89$. If $F>0.89$, we see that $\phi_{22}<0$ for $H> \frac{2}{3}F$, $B>1.29$, $0.6a<i\leq a$, $D>1$, $1.17<C\leq 1.322$ and $1.4<A<1.5$.

\noindent {Claim (ii)}$E<0.9348$ for otherwise $E^4ABCDI>2$ using $i>0.6a$, $c>0.17$, $b>0.29$ and $\phi_{10}<0$.

\noindent {Claim (iii)} $D\leq 1.083$ using $h<c+\frac{a+i}{2}-e-\frac{f}{2}$ and proving  $D^4ABCHI>2$ and $\phi_{17}<0$.

\noindent {Claim (iv)} $C\leq 1.217$  proving  $\frac{C^3}{DEF}>2$ and  $\phi_{11}<0$.

\noindent {Claim (v)} $f>0.205$. Using $h<b+d-f+\frac{i}{2}$, we find that $\phi_{22}^{(1)}<0$ for $c<0.217$, $0.11<f\leq 0.205$, $0<d<0.083$ and $0.6a<i<a$.

\noindent {Claim (vi)} $D\leq 1.065$ using $h<c+\frac{a+i}{2}-e-\frac{f}{2}$ and proving  $D^4ABCHI>2$ and $\phi_{17}<0$.

\noindent {Final contradiction :}  $C<(2DEF)^{\frac{1}{3}}<1.17$. This gives a contradiction. \vspace{2mm}

\noindent \textbf{Subcase X: $1.29 \leq B <  1.35, ~i> 0.6a,~1<C\leq 1.17 $}\vspace{2mm}

\noindent {Claim (i)} $C>1.053$. If $C\leq 1.053$, we find that $\phi_9^*<0$ for $0.6a<i\leq a$.

\noindent {Claim (ii)} $\frac{A^3}{BCD}\leq 2$. For if $\frac{A^3}{BCD}> 2$, working as in Subcase II of Case III we can take  $G\leq 0.65$ and $E\leq 0.92$  and  $F< 0.8$. Now $\frac{D^3}{EFG}>\frac{1}{0.92\times 0.8\times 0.65}>2$ for $D>1$, whereas $\phi_{16}<0$ for $B>1.29$. $1<I\leq A$, $1<C<1.17$, $H>0.46873D$, and $1<D<1.155$ . This gives a contradiction.

\noindent {Claim (iii)} Working as in Subcase VIII, we can take $F<0.935$,   $D\leq 1.1165$.

\noindent {Claim (iv)} $f>0.185$. Using $h<b+d-f+\frac{i}{2}$,   we find that $\phi_{22}^{(1)}<0$ for $c<0.17$, $0.065<f\leq 0.185$, $0<d<0.1165$ and $0<i<a$.


\noindent {Claim (v)}$E<0.9597$ for otherwise $E^4ABCDI>2$ using $i>0.6a$, $c>0.053$, $b>0.29$ and $\phi_{10}<0$.

\noindent {Claim (vi)} $D\leq 1.063$ using $h<c+\frac{a+i}{2}-e-\frac{f}{2}$ and proving  $D^4ABCHI>2$ and $\phi_{17}<0$.

\noindent {Claim (vii)} $f>0.242$. Using $h<b+d-f+\frac{i}{2}$,  we find that $\phi_{22}^{(1)}<0$ for $c<0.17$, $0.185<f\leq 0.242$, $0<d<0.063$ and $0<i<a$.

\noindent {Claim (viii)} $D\leq 1.047$ using $h<c+\frac{a+i}{2}-e-\frac{f}{2}$ and proving  $D^4ABCHI>2$ and $\phi_{17}<0$.

\noindent {Claim (ix)} $C\leq 1.151$  proving  $\frac{C^3}{DEF}>2$ and  $\phi_{32}^{(1)}<0$.

\noindent {Claim (x)} $\frac{B^3}{CDE}\leq 2$. For if $\frac{B^3}{CDE}> 2$, working as in Subcase VI of Case II we can take  $G\leq 0.66$   and  $F< 0.76$. Now $\frac{D^3}{EFG}>\frac{1}{0.9597\times 0.76\times 0.66}>2$ for $D>1$, whereas $\phi_{16}<0$ for $B>1.29$. $1<I\leq A$, $1<C<1.151$, $H>0.46873D$, and $1<D<1.047$ .

\noindent {Claim (xi)}$B<1.3226$,  $A<1.4718$, using Claims (x) and  (ii).

\noindent {Claim (xii)} $f>0.262$, $D<1.032$, $C<1.135$, $B<1.311$ and $A<1.4536$. With $A,B$ and $D$ reduced, we repeat the arguments in claims (vii)-(ix), (xi) and obtain these.

\noindent {Claim (xiii)} $C>1.073$. If $C\leq 1.073$, we find that $\phi_9^*<0$ for $0.6a<i\leq a$.

\noindent {Claim (xiv)}$E<0.9552$, $D<1.0245$ by repeating arguments in claims (v) and (vi).

\noindent {Claim (xv)} $ i>0.78a$. If $i\leq 0.78a$ and $g+h\leq 1.5d+0.517(a+i)$ we find that   $\phi_{25}<0$ for  $0<i\leq 0.78a$, $0<d\leq 0.0245$, $0.4<a<0.4536$. If  $g+h>1.5d+0.517(a+i)$, we get  $f<b+d+\frac{i}{2}-\frac{1.5d+0.517(a+i)}{2}$. But then $\phi_{22}^{(1)}<0$ for $h<b+d+\frac{i}{2}-f$, $0<c<0.135$, $0<i\leq 0.78a$, $0<d\leq 0.0245$, $0.29<b<0.311$, $0.4<a<0.4536$.

\noindent {Claim (xvi)}$E<0.942$, $D<1.019$, $f>0.271$, $D<1.017$, $C<1.118
$, $B<1.29$, by repeating arguments in claims (v)-(ix) and using claim (x).
 This contradicts that $B>1.29$.
 \vspace{2mm}

\noindent \textbf{Subcase XI: $1.35 \leq B <  1.451, 1<C<1.322, ~\frac{B^3}{CDE}> 2 $}\vspace{2mm}

We find that $\phi_{29}<0$ if $F>0.784$ and using $H>\frac{2}{3}F$ or if $H>0.64$ and using $F\geq \frac{2}{3}$ and for  $1<I\leq A$, $1.4<A<1.5,~1.35 \leq B <  1.451$ . Therefore we can take $F\leq 0.784$ and $H\leq 0.64$. Further if $G>0.66$, we find that $\phi_{7}<0$ for $F\geq \frac{2}{3}$, $1<I\leq A$, $1.4<A<1.5,~1.35 \leq B <  1.451$. Therefore we can take $G<0.66$. Now $\frac{E^3}{FGH}>2$ and $\phi_8^*<0$  for $1<I\leq A$, $1.4\times1.35<x=ABCD<1.5\times1.451\times1.322\times1.155$ and for $0.92<E<1$. Therefore we can take $E\leq 0.92$. \vspace{2mm}

Now  $\frac{C^3}{DEF}>2$ for $C>1.2$ and  $\phi_{11}<0$ for $C>1.26$, $G>0.46873C$,  $1\leq I\leq A$, $1.4<A\leq 1.5$ and $1.35<B<1.451$. Therefore we can take $C\leq 1.26$. Further $\frac{D^3}{EFG}>2$ for $D>1, E\leq 0.92, ~F\leq 0.784~,G\leq 0.66$.  If $D>1.05$ or if $C<1.21$ or if $B>1.425$, we find that  $\phi_{16}<0$  for  $H>0.46873D$, $1<I\leq A$. Therefore we can take $D\leq 1.05$, $C\geq 1.21$ and $B\leq 1.425$.  Further  $\phi_{11}<0$ for $C>1.24$ and $1.35<B\leq 1.425$.  So we can take $C<1.24$. If $A<1.444$ or $B>1.394$, we find that $\phi_{16}<0$. So we can take $A\geq 1.444$ and $B<1.394$. But then   $\phi_{11}<0$  for  $A\geq 1.444$, $C>1.21$ and $1.35<B\leq 1.394$. This gives a contradiction. \vspace{2mm}

\noindent \textbf{Subcase XII: $1.35 \leq B <  1.451, 1<C<1.322, ~\frac{B^3}{CDE}\leq 2, ~i\leq 0.45a $}\vspace{2mm}

\noindent {Claim (i)} $C>1.144$. If $C\leq 1.144$, we find that $\phi_9^*<0$ for $0<i\leq 0.45a$.

\noindent {Claim (ii)} $F<0.915$. If $F>0.915$, we see that $\phi_{22}<0$ for $H> \frac{2}{3}F$, $B>1.35$, $1<I\leq A$, $D>1$, $1<C\leq 1.322$ and $1.4<A<1.5$.

\noindent {Claim (iii)}$E<0.9807$ for otherwise $E^4ABCDI>2$ using  $C>1.144$, $B>1.35$, $A>1.4$ and $\phi_{8}^*<0$..

\noindent{Claim (iv)} $D\leq 1.1275$ using $h<\min\{b+\frac{c+i}{2}-e-\frac{f}{2}, \frac{5}{9}-\frac{4}{9}d\}$, and proving  $D^4ABCHI>2$ and $\phi_{16}<0$ for $D>1.1275$.

\noindent {Claim (v)} $C\leq 1.2649$  proving  $\frac{C^3}{DEF}>2$ and $\phi_{11}<0$. Hence $B<(2CDE)^{\frac{1}{3}}<1.4091$.

\noindent {Claim (vi)} $f>0.24$. Using $h<b+d-f+\frac{i}{2}$, we find that $\phi_{22}^{(1)}<0$ for $c<0.2649$, $0.085<f\leq 0.24$, $0<d<0.1275$,  $0<i<0.45a$ and $b<\min\{a,0.451\}$.

\noindent {Claim (vii)} $D\leq 1.058$ using $h<c+\frac{a+i}{2}-e-\frac{f}{2}$ and proving  $D^4ABCHI>2$ and $\phi_{16}<0$ for $0<i<0.45a$. Hence $B<1.3795$ and $C>\frac{B^3}{2DE}>1.185$.

\noindent {Final contradiction :} $\frac{C^3}{DEF}>2$ for $C>1.185$ but $\phi_{11}<0$. \vspace{2mm}

\noindent \textbf{Subcase XIII: $1.35 \leq B <  1.451, 1<C<1.322, ~\frac{B^3}{CDE}\leq 2, ~i> 0.45a $}\vspace{2mm}

\noindent {Claim (i)} $C>1.093$. If $C\leq 1.093$, we find that $\phi_9^*<0$ for $0<i\leq a$.

\noindent {Claim (ii)} $F<0.88$. If $F>0.88$, we see that $\phi_{22}<0$ for $H> \frac{2}{3}F$, $B>1.35$, $1+0.45a<I\leq A$, $D>1$, $1<C\leq 1.322$ and $1.4<A<1.5$.

\noindent {Claim (iii)}$E<0.9518$ for otherwise $E^4ABCDI>2$ using  $C>1.093$, $B>1.35$, $A>1.4$, $i>0.45a$ and $\phi_{8}^*<0$.

\noindent{Claim (iv)} $D\leq 1.1075$ using $h<\min\{b+\frac{c+i}{2}-e-\frac{f}{2}, \frac{5}{9}-\frac{4}{9}d\}$, and proving  $D^4ABCHI>2$ and $\phi_{16}<0$. Hence $B<1.408$.

\noindent {Claim (v)} $C\leq 1.235$  proving  $\frac{C^3}{DEF}>2$ and $\phi_{11}<0$. Hence $B<1.376$. $C>\frac{B^3}{2DE}$ gives $C>1.167$. $D>\frac{B^3}{2CE}$ gives $D>1.046$

\noindent {Claim (vi)}$E<0.926$ for otherwise $E^4ABCDI>2$ using  $C>1.167$, $B>1.35$, $A>1.4$, $D>1.046$, $i>0.45a$ and $\phi_{10}<0$.

\noindent {Claim (vii)} $D\leq 1.0855$ using $h<c+\frac{a+i}{2}-e-\frac{f}{2}$ and proving  $D^4ABCHI>2$ and $\phi_{16}<0$. Further $B<1.3541$.

\noindent {Claim (viii)} $C\leq 1.2095$  proving  $\frac{C^3}{DEF}>2$ and  $\phi_{11}<0$. But then $B<(2CDE)^{\frac{1}{3}}<1.35$. This gives a contradiction.
\vspace{2mm}

\noindent \textbf{Case IV: 1.3 $\leq A < $ 1.4}\vspace{2mm}

\noindent \textbf{Subcase I: $1.3<B\leq A,~  \frac{B^3}{CDE}>2$}\vspace{2mm}

Working as in Subcase VI of Case II we can take $F<0.77$, $H<0.61$ and $G<0.62$. Now $\frac{E^3}{FGH}>2$ and $\phi_8^*<0$  for $0.91<E<1$ and for $1<I\leq A$, $1.3\times1.3<x=ABCD<1.4\times1.4\times1.322\times1.155$. Therefore we can take $E\leq 0.91$. \vspace{2mm}

Now  $\frac{C^3}{DEF}>2$ for $C>1.22$ and  $\phi_{32}^{(1)}<0$ for $B\leq A$, $G>0.46873C$, $h<\frac{a+i}{2}+c-e-\frac{f}{2}< \frac{a+i}{2}+c-0.205$, $1<I\leq A$ and $1.22<C<1.322$. Therefore we can take $C\leq 1.22$. Finally $\frac{D^3}{EFG}>2$ for $D>1, E\leq 0.91, ~F\leq 0.77~,G\leq 0.62$; but  $\phi_{16}<0$  for  $0.46873D<H<0.61$, $1<I\leq A$, $B>1.3$ and $C\leq 1.22$. This gives a contradiction. \vspace{2mm}

\noindent \textbf{Subcase II: $1.3<B\leq A, \frac{B^3}{CDE}\leq 2,~ i<0.6a$}\vspace{2mm}

\noindent {Claim (i)} $F<0.941$. If $F>0.941$, we see that $\phi_{22}<0$ for $H> \frac{2}{3}F$, $B>1.3$, $1<I\leq A$, $D>1$, $1<C\leq \min\{A,1.322\}$ and $1.3<A<1.4$.

\noindent{Claim (ii)} $D\leq 1.142$ using $h<\min\{b+\frac{c+i}{2}-\frac{f}{2}, \frac{5}{9}-\frac{4}{9}d\}$, and proving  $D^4ABCHI>2$ and $\phi_{16}<0$.

\noindent {Claim (iii)} $C\leq 1.291$  proving  $\frac{C^3}{DEF}>2$ and $\phi_{11}<0$.

\noindent {Claim (iv)} $f>0.214$. Using $h<b+d-f+\frac{i}{2}$, in place of $H> \frac{2}{3}F$  we find that $\phi_{22}^{(1)}<0$ for $c<0.291$, $0.059<f\leq 0.214$, $0<d<0.142$ and $0<i<0.6a$.

\noindent {Claim (v)} $C\leq 1.216$  proving  $\frac{C^3}{DEF}>2$ and $\phi_{32}^{(1)}<0$ using $h<\frac{a+i}{2}+c-\frac{f}{2}< \frac{a+i}{2}+c-0.107$.

\noindent {Claim (vi)} $D\leq 1.072$ using $h<c+\frac{a+i}{2}-\frac{f}{2}$ and proving  $D^4ABCHI>2$ and $\phi_{16}<0$. Further $C<1.191$.

\noindent {Claim (vii)} $B<1.367$ using $\frac{B^3}{CDE}\leq 2$.

\noindent {Claim (viii)} $e+g\geq 0.705c+0.296(a+i)$  proving  $\phi_{19}^{(1)}<0$ for  $0<i\leq 0.6a$, $0<c\leq 0.191$, $0.3<b\leq \min\{a,0.367\}$. This gives $h<0.59c+0.408(a+i)$.

\noindent {Claim (ix)} $e\geq 0.132$ i.e. $E<0.868$ by proving   $\phi_{23}<0$ for $b<0.367$, $i<0.6a$. Then $C<1.136$

%
%
%

\noindent {Final contradiction :} $B<(2CDE)^{\frac{1}{3}}<1.3$ \vspace{2mm}

\noindent \textbf{Subcase III: $1.3<B\leq A,~ i\geq 0.6a,~ \frac{B^3}{CDE}\leq 2$}\vspace{2mm}

\noindent {Claim (i)} $C>1.1$, for if $C\leq 1.1$, we find that $\phi_9^*<0$.

\noindent {Claim (ii)} $F<0.941$, $D\leq 1.142$ , $C\leq 1.291$ as in Subcase II.

\noindent {Claim (iii)} $E<0.9772$ if $C>1.1$ and $E<0.9562$ if $C>1.2$ for otherwise we have $E^4ABCDI>2$ and $\phi_8^*<0$.

\noindent {Claim (iv)}  $F<0.905$ if $1.2<C<1.291$ and $F<0.895$ if $C<1.2$ for otherwise we have $\phi_{22}<0$ for $H> \frac{2}{3}F$, $B>1.3$, $D>1$, $1+0.6a<I\leq A$ and $1.3<A<1.4$.

\noindent {Claim (v)} If $C>1.2$, we must have $C<1.255$ for otherwise $\frac{C^3}{DEF}>2$ and $\phi_{11}<0$.

\noindent {Claim (vi)}$D<1.095$ if $1.2<C<1.255$ and $D<1.085$ if $C<1.2$. Here we use $e+\frac{f}{2}> 0.0913 $ if $1.2<C<1.255$ and $e+\frac{f}{2}>0.0753$ if $C<1.2$. Using $h<\frac{a+i}{2}+c-e-\frac{f}{2}$, we get $D^4ABCHI>2$ for $D\geq 1.097$ if $1.2<C<1.255$ and $D\geq 1.085$ if $C<1.2$ and $\phi_{16}<0$.

\noindent {Claim (vii)}$B< 1.366$ if $C<1.2$.

\noindent {Claim (viii)} If $C>1.2$, we must have $C<1.2376$ for otherwise $\frac{C^3}{DEF}>2$ and $\phi_{11}<0$. Hence $B<1.374$.

\noindent {Claim (ix)} $f>0.2$ if  $1.2<C<1.2376$. Using $h<b+d-f+\frac{i}{2}$,  we find that $\phi_{22}^{(1)}<0$ for $0.2<c<0.2376$, $0.095<f\leq 0.2$, $0<d<0.095$ and $0<i<a$. Similarly $f>0.238$ if $C<1.2$.

\noindent {Claim (x)} If $C>1.2$, we must have $C<1.22$; so $B<1.33$. If $C<1.2$, we must have $D<1.0496$, $C<1.161$ and so $B<1.336$.

\noindent {Claim (xi)} If $C>1.2$, we have $\frac{C^3}{DEF}>2$ and $\phi_{32}^{(1)}<0$, giving a contradiction.

\noindent {Claim (xii)} $e+g\geq 0.7c+0.28(a+i)$  proving  $\phi_{19}^{(1)}<0$. This gives $h<0.6c+0.44(a+i)$.

\noindent {Claim (xiii)} $e\geq 0.08$ by proving   $\phi_{23}<0$. Then $C<1.1375$. But then $B<(2CDE)^{\frac{1}{3}}<1.3$, a contradiction. \vspace{2mm}

\noindent \textbf{Case IV: 1 $\leq A < $ 1.19}\vspace{2mm}

Let first $A<1.164$. If $g+h<1.5d+0.466(a+i)$, we find that $\phi_{25}<0$ for $0<i<a$ and $0<d<\min \{a,0.155\}$. When $g+h\geq 1.5d+0.466(a+i)$, we get from (6.21), that $f<b+d+\frac{i}{2}-\frac{1}{2}(1.5d+0.466(a+i))$. Then $\phi_{22}^{(2)}<0$ for $h<b+d+\frac{i}{2}-f$, $0<f<b+d+\frac{i}{2}-\frac{1}{2}(1.5d+0.466(a+i))$, $0<d<\min \{a,0.155\}$, $0<i<a$ and $0<b\leq a$.\vspace{2mm}

Let now $1.164\leq A<1.19$. If $g+h<1.5d+0.47(a+i)$, we find that $\phi_{25}<0$ for $0<i\leq a$ and $0<d<0.155$. When $g+h\geq 1.5d+0.47(a+i)$, we get that  $\phi_{22}^{(2)}<0$ for $h<b+d+\frac{i}{2}-f$, $0<f<b+d+\frac{i}{2}-\frac{1}{2}(1.5d+0.47(a+i))$, $0<d<0.155$, $0<i<a$ and $0<b\leq 0.164$. Therefore we can take  $b> 0.164$. If $e+g<0.71c+0.325(a+i)$, we find that $\phi_{19}^{(1)}<0$ for  $0<i\leq a$, $0<c\leq a$, $0.164<b\leq a$,and $0.164<a<0.19$. When $e+g\geq 0.71c+0.325(a+i)$
we have $h<0.58c+0.35(a+i)$ from $(1,2,2,2,1,1)_w$. But then $D^4ABCHI> 2$ for $B>1.164$, $D>1.12$, $0<i\leq a$ and $0<c\leq a$ . Also
$\phi_{38}<0$ for $H>0.46872D$, $1<I\leq A$, $1.12<D<1.155$, $1.164\leq A<1.19$. Therefore we can take $D\leq 1.12$. With $d$ now improved, we arrive at a contradiction by proving $\phi_{25}<0$ when $g+h<1.5d+0.528(a+i)$ and $\phi_{22}^{(2)}<0$ when $g+h\geq 1.5d+0.528(a+i)$. \vspace{2mm}

\noindent \textbf{Case V: 1.19 $\leq A < $ 1.25}\vspace{2mm}

Working as in Claim (i), we find that $D^4ABCHI>2$ for $h<\min\{b+\frac{c+i}{2}, \frac{5}{9}-\frac{4}{9}d\}$, $b<a$ and $D>1.1492$. Also  $\phi_{38}<0$ for $H>0.46873D$, $1\leq I\leq A$ and $1<D\leq A<1.25$.  This gives a contradiction. So we can take $D\leq 1.1492$.

\noindent \textbf{Subcase I: $i<0.8a$}\vspace{2mm}

\noindent {Claim (i)} $B\geq 1.19$.
Suppose $B<1.19$. If $g+h<1.5d+0.509(a+i)$, we find that $\phi_{25}<0$. If  $g+h\geq 1.5d+0.509(a+i)$, we get that  $\phi_{22}^{(2)}<0$.\vspace{2mm}

\noindent {Claim (ii)} $f\geq 0.164$.
Suppose $f<0.164$. Then $\phi_{22}^{(1)}<0$

\noindent {Claim (iii)} $D\leq 1.109$ using $h<c+\frac{a+i}{2}-\frac{f}{2}$ and proving  $D^4ABCHI>2$ and $\phi_{17}<0$.

\noindent{Claim (iv)} $C\leq 1.229$  proving  $\frac{C^3}{DEF}>2$ and $\phi_{11}<0$.

\noindent {Claim (v)} $e+g\geq 0.69c+0.29(a+i)$  proving  $\phi_{19}^{(1)}<0$. This gives $h<0.62c+0.42(a+i)$.

\noindent {Claim (vi)} $e\geq 0.09$  proving   $\phi_{23}<0$.

\noindent{Claim (vii)} $D\leq 1.067$ using $h<c+\frac{a+i}{2}-e-\frac{f}{2}$ and proving  $D^4ABCHI>2$ and $\phi_{17}<0$.

\noindent {Claim (viii)} $C\leq 1.176$  proving  $\frac{C^3}{DEF}>2$ and $\phi_{32}<0$.

\noindent {Final contradiction :} If $g+h<1.5d+0.58(a+i)$, we find that $\phi_{25}<0$. If  $g+h\geq 1.5d+0.58(a+i)$, we get that  $\phi_{22}^{(1)}<0$.\vspace{2mm}

\noindent  \textbf{Subcase II: $i\geq 0.8a$}\vspace{2mm}

\noindent {Claim (i)} $B\geq 1.174$.
Suppose $B<1.174$. If $g+h<1.5d+0.47(a+i)$, we find that $\phi_{25}<0$. If  $g+h\geq 1.5d+0.47(a+i)$, we get that  $\phi_{22}^{(2)}<0$.\vspace{2mm}

\noindent {Claim (ii)} $f\geq 0.185$
by proving  $\phi_{22}^{(1)}<0$.

\noindent {Claim (iii)} $D\leq 1.103$ using $h<c+\frac{a+i}{2}-\frac{f}{2}$ and proving  $D^4ABCHI>2$ and $\phi_{38}<0$.

\noindent {Claim (iv)} $C\leq 1.216$  proving  $\frac{C^3}{DEF}>2$ and $\phi_{11}<0$.

\noindent {Claim (v)} $e+g\geq 0.69c+0.29(a+i)$  proving  $\phi_{19}^{(1)}<0$. This gives $h<0.62c+0.42(a+i)$.

\noindent {Claim (vi)} $e\geq 0.08$  proving   $\phi_{23}<0$.

\noindent {Claim (vii)} $D\leq 1.068$ using $h<c+\frac{a+i}{2}-e-\frac{f}{2}$ and proving  $D^4ABCHI>2$ and $\phi_{17}<0$.

\noindent {Claim (viii)} $C\leq 1.17$  proving  $\frac{C^3}{DEF}>2$ and $\phi_{32}<0$.

\noindent {Claim (ix)} $e+g\geq 0.72c+0.315(a+i)$  proving  $\phi_{19}^{(1)}<0$. This gives $h<0.56c+0.37(a+i)$.

\noindent {Claim (x)} $e\geq 0.13$  proving   $\phi_{23}<0$.

\noindent {Claim (xi)} $D\leq 1.05$ using $h<c+\frac{a+i}{2}-e-\frac{f}{2}$ and proving  $D^4ABCHI>2$ and $\phi_{17}<0$.

\noindent {Claim (xii)} $C\leq 1.142$  proving  $\frac{C^3}{DEF}>2$ and $\phi_{32}<0$.

\noindent {Final contradiction :} If $g+h<1.5d+0.548(a+i)$, we find that $\phi_{25}<0$. If  $g+h\geq 1.5d+0.548(a+i)$, we get that  $\phi_{22}^{(1)}<0$.\vspace{2mm}

\noindent  \textbf{Case VI: 1.25 $\leq A < $ 1.3}\vspace{2mm}

Working as in Claim (i), we find that $D^4ABCHI>2$ for $h<\min\{b+\frac{c+i}{2}, \frac{5}{9}-\frac{4}{9}d\}$, $b<a$ and $D>1.143$. Also  $\phi_{38}<0$ for $H>0.46873D$, $1\leq I\leq A$ and $1<D\leq A<1.3$.  This gives a contradiction. So we can take $D\leq 1.143$.

\noindent \textbf{\bf Subcase I: $i<0.7a$}\vspace{2mm}

\noindent {Claim (i)} $B\geq 1.218$.
Suppose $B<1.218$. If $g+h<1.5d+0.509(a+i)$, we find that $\phi_{25}<0$. If  $g+h\geq 1.5d+0.509(a+i)$, we get that  $\phi_{22}^{(2)}<0$.\vspace{2mm}

\noindent {Claim (ii)} $f\geq 0.182$
proving  $\phi_{22}^{(1)}<0$

\noindent {Claim (iii)} $C\leq 1.232$  proving  $\frac{C^3}{DEF}>2$ and $\phi_{11}<0$.

\noindent{Claim (iv)} $D\leq 1.098$ using $h<c+\frac{a+i}{2}-\frac{f}{2}$ and proving  $D^4ABCHI>2$ and $\phi_{17}<0$.

\noindent {Claim (v)} $C\leq 1.216$  proving  $\frac{C^3}{DEF}>2$ and $\phi_{11}<0$.

\noindent {Claim (vi)} $e+g\geq 0.7c+0.3(a+i)$  proving  $\phi_{19}^{(1)}<0$. This gives $h<0.6c+0.4(a+i)$.

\noindent {Claim (vii)} $e\geq 0.12$  proving   $\phi_{23}<0$.

\noindent {Claim (viii)} $D\leq 1.057$ using $h<c+\frac{a+i}{2}-e-\frac{f}{2}$ and proving  $D^4ABCHI>2$ and $\phi_{17}<0$.

\noindent {Claim (ix)} $C\leq 1.151$  proving  $\frac{C^3}{DEF}>2$ and $\phi_{32}<0$.

\noindent {Final contradiction :} If $g+h<1.5d+0.585(a+i)$, we find that $\phi_{25}<0$. If  $g+h\geq 1.5d+0.585(a+i)$, we get that  $\phi_{22}^{(1)}<0$.\vspace{2mm}

\noindent  \textbf{Subcase II: $i\geq 0.7a$} \vspace{2mm}

\noindent{Claim (i)} $B\geq 1.185$ proving $\phi_{25}<0$ for $g+h<1.5d+0.465(a+i)$, and $\phi_{22}^{(2)}<0$ for  $g+h\geq 1.5d+0.465(a+i)$.

\noindent {Claim (ii)} $f\geq 0.2$
by proving  $\phi_{22}^{(1)}<0$.

\noindent {Claim (iii)} $C\leq 1.223$  proving  $\frac{C^3}{DEF}>2$ and $\phi_{11}<0$.

\noindent {Claim (iv)} $D\leq 1.093$ using $h<c+\frac{a+i}{2}-\frac{f}{2}$ and proving  $D^4ABCHI>2$ and $\phi_{38}<0$.

\noindent {Claim (v)} $C\leq 1.205$  proving  $\frac{C^3}{DEF}>2$ and $\phi_{11}<0$.

\noindent {Claim (vi)} $e+g\geq 0.7c+0.285(a+i)$  proving  $\phi_{19}^{(1)}<0$. This gives $h<0.6c+0.43(a+i)$.

\noindent {Claim (vii)} $e\geq 0.07$  proving   $\phi_{23}<0$.

\noindent {Claim (viii)} $C\leq 1.177$  proving  $\frac{C^3}{DEF}>2$ and $\phi_{32}<0$.

\noindent {Claim (ix)} $D\leq 1.06$ using $h<c+\frac{a+i}{2}-e-\frac{f}{2}$ and proving  $D^4ABCHI>2$ and $\phi_{17}<0$.

\noindent {Claim (x)} $B\geq 1.235$.
proving $\phi_{25}<0$ for $g+h<1.5d+0.528(a+i)$, and $\phi_{22}^{(2)}<0$ for  $g+h\geq 1.5d+0.528(a+i)$.\vspace{2mm}

\noindent {Claim (xi)} $f\geq 0.23$
by proving  $\phi_{22}^{(1)}<0$.

\noindent {Claim (xii)} $C\leq 1.494$  proving  $\frac{C^3}{DEF}>2$ and $\phi_{32}<0$.

\noindent {Claim (xiii)} $D\leq 1.038$ and   $C<1.142$.

\noindent {Claim (xiv)} $e+g\geq 0.72c+0.315(a+i)$ by proving  $\phi_{19}^{(1)}<0$. This gives $h<0.56c+0.37(a+i)$.

\noindent {Claim (xv)} $e\geq 0.132$ and   $C<1.116$.

\noindent {Claim (xvi)} $B\geq 1.265$.

\noindent {Final contradiction :} For $B\geq 1.265$, $1<I\leq A$,  $1<C<1.116$, $1.25<A\leq 1.3$ we find that $\phi_{9}^*<0$. \vspace{2mm}

\noindent \textbf{Proposition 46.} Case (6) i.e. $A>1$, $B>1$, $C>1$, $D>1$, $E>1$, $F\leq1$, $G\leq1$, $H\leq1$, $I\leq 1$ does not arise.\vspace{2mm}\\
{\noindent \bf Proof.} Here $B\leq A<2.1326324$, ~$C\leq 2$, $D\leq 1.5$, $E\leq\frac{4}{3}$.\\
Using the weak inequalities
$(1,2,2,2,1,1)$, ~$(2,2,2,2,1)$, ~$(1,2,2,2,2)$, ~$(2,2,2,1,2)$ and ~$(2,1,1,2,1,2)$   we get
\begin{equation}a+2c+2e-2g-h-i>0\vspace{-2mm}\end{equation}
\begin{equation}2b+2d-2f-2h-i>0,\vspace{-2mm}\end{equation}
\begin{equation}a+2c+2e-2g-2i>0\vspace{-2mm}\end{equation}
\begin{equation}2b+2d-2f-g-2i>0\vspace{-2mm}\end{equation}
\begin{equation}2b+c+d-2f-g-2i>0\end{equation}

\noindent \textbf{Claim(i) $E\leq1.155$}\vspace{2mm}

Suppose $E>1.155$. If $B\geq 1.5$, we have $E^4ABCDI>E^4(1.5)^2(\frac{4}{9}E)>2$. Let now $B<1.5$. Using $(6.27)$ We have $i<b+\frac{c+d}{2}$. Also $I>\frac{4}{9}E$. So we have
$$i<\left\{\begin{array}{lll}b+\frac{c+d}{2} & {\rm if }& 0<c+d\leq \frac{10}{9}-\frac{8}{9}e-2b\\\frac{5}{9}-\frac{4}{9}e
& {\rm if }& c+d>\frac{10}{9}-\frac{8}{9}e-2b. \end{array}\right. $$

We find that $E^4ABCDI>(1+e)^4(1+a)(1+b)(1+c+d)(1-i)>2$ for $i,~c+d$ as above and $0<b< \min\{a,0.5\}$, $0.155<e\leq a$.
 Also  $\phi_8^*=4(x)^{1/4}+4E-\frac{1}{2}E^5xI+I-9>0$, where $x=ABCD$, has maximum at $x=(\frac{2}{E^5I})^{\frac{4}{3}}$. So $\phi_8^*<4E+3(\frac{2}{E^5I})^{\frac{1}{3}}+I<0$ for $E>1.08$ and $0.46873\leq I\leq 1$. This gives a contradiction. \vspace{2mm}

 \noindent \textbf{Claim(ii) $D\leq1.322$}\vspace{2mm}

 Suppose $D>1.322$, then $\frac{D^3}{EFG}>\frac{1.322^3}{1.155}>2$. Now $\phi_9<0$ for $H>0.46873D$, $1<C<A$ and for $A>B$.\vspace{2mm}

  \noindent \textbf{Claim(iii) $C\leq1.498$}\vspace{2mm}

  Suppose $C>1.498$, then $\frac{C^3}{DEF}>2$. If $G^2>HI$, then $\phi_{43}<0$ using $G>0.46873C$. If $G^2<HI$, then $\phi_{32}<0$ using $HI>G^2$,  $G\geq0.46873C$ and $B\leq A<2.14$. Hence $C\leq1.498$.\vspace{2mm}

  \noindent \textbf{Claim(iv) $B\leq1.71$}\vspace{2mm}

  Suppose $B>1.71$, then $\frac{B^3}{CDE}>\frac{1.71^3}{1.498\times1.322\times1.155}>2$. If $\frac{F^3}{GHI}>2$, then $\phi_2<0$. Let now $\frac{F^3}{GHI}\leq2$.  If $H \geq 0.58$ then $\phi_{6}<0$ using $F>0.46873B$. So we must have $H<0.58$. Now $I>\frac{F^3}{2GH}>\frac{F^3}{2G(0.58)}$. Using this lower bound of $I$,
  $G>\frac{2}{3}$ and $1<A<2.14$ we find $\phi_7<0$. Therefore we must have  $B\leq1.71$. \vspace{2mm}

  \noindent \textbf{Claim(v) $\frac{A^3}{BCD}<2$ and $A\leq1.893$}\vspace{2mm}

  Suppose $\frac{A^3}{BCD}>2$. If $1<A<1.57$, then $\phi_4^*<0$ using that $\phi_4^*$ has its maximum at $x=\{\frac{2}{A^5}\}^{\frac{5}{4}}$. If $\frac{E^3}{FGH}>2$, then $\phi_{26}<0$ for $0.46873<I<1$, $1<E<1.155$ and $1<A<2.14$. Now consider $A\geq1.57$ and $\frac{E^3}{FGH}<2$. If $G>0.76$ we find that  $\phi_4<0$. So we must have $G<0.76$. So $H>\frac{E^3}{2(0.76)F}>0.657$. Now $\phi_5<0$ for $H>0.657$, $F>\frac{3}{4}$, $1<E<1.155$ and $1<A<2.14$. Hence  $\frac{A^3}{BCD}<2$. This gives  $A<(2BCD)^{1/3}<(2\times1.71\times1.498\times1.322)^{1/3}<1.893.$\vspace{2mm}

   \noindent \textbf{Claim(vi) $B<1.678$}\vspace{2mm}

    Working as in Claim(iv), we first get that $\frac{F^3}{GHI}\leq2$, $H<0.585$. If $I>0.56$, we find that $\phi_7<0$ for $G>\frac{2}{3}$, $F>0.46873B$ and $1<A<2.14$. Therefore we can take $I\leq 0.56$. Now $G>\frac{F^3}{2HI}>\frac{F^3}{2(0.585)(0.56)}>0.74$ for $F>0.46873B$, $B\geq  1.678$  and $\phi_7<0$ for $I>\frac{2}{3}G$, $1.678\leq B\leq A<2.14$.  Therefore we must have  $B< 1.678$.\vspace{2mm}

   \textbf{We divide the range of $A$ into several
 subintervals and  arrive at a contradiction in each.}\vspace{2mm}

 \noindent \textbf{Case I: $1.71\leq A\leq1.893$}\vspace{2mm}

Suppose first that $\frac{E^3}{FGH}<2$. As $\frac{A^3}{BCD}<2$, we get $A^4E^4I<4$, which gives $A<(\frac{4}{0.46873})^{\frac{1}{4}}<1.71$, a contradiction. So we must have $\frac{E^3}{FGH}\geq2$.\vspace{2mm}

\noindent \textbf{Subcase I:}  $B\leq1.55$, $C>1.33$ \vspace{2mm}

As $C>{\mbox{each of}}~\{D,E,F,G,H,I\}$, the inequality $(2,7^*)$ holds and $\phi_3^*<0 $ for $B\leq 1.47$. So we can take $B>1.47$.  For $F>0.908$, $F^4ABCDE>(0.908)^4(1.71)(1.47)(1.33)>2$ and $\phi_{12}^*<0$. So we can take $F<0.908$. For $E>1.08$, $E^4ABCDI>0.46873E^5ABC>2$ and $\phi_8^*<0$. So $E<1.08$. Now for $D>1.252$, $\frac{D^3}{EFG}>2$ and $\phi_9<0$. So $D<1.252$. Now  $\phi_{10}<0$ for $1<E<1.08$ and $B<1.55$.\vspace{2mm}

\noindent \textbf{Subcase II:} $B\leq1.55$, $C\leq1.33$.\vspace{2mm}

As $\frac{A^3}{BCD}<2$, we get  $A<1.7599$, $B>1.421$, $C>1.22$, $D>1.212$. As in Subcase (i)  we have $F<0.908$, $E<1.08$ and $D<1.252$. Using $I>\frac{2}{3}G$, we find that $\phi_{19}^{(2)}<0$ for $G>0.82$. So  we have $G<0.82$. For $D>1.212$, we have $\frac{D^3}{EFG}>2$ and $\phi_9<0$. This gives a contradiction. \vspace{2mm}

\noindent \textbf{Subcase III:} $1.55<B\leq1.59$, $C\leq1.33$.\vspace{2mm}

 Here $\frac{A^3}{BCD}<2$ gives $A<1.774$, $C>1.189$, $D>1.182$. As in Subcase (i) we have $F<0.908$, $E<1.08$ and $D<1.252$. Again $\frac{A^3}{BCD}<2$ gives  $C>1.255$. Now using $I>\frac{2}{3}G$, we find that $\phi_{19}^{(2)}<0$ for $G>0.82$. So now we have $G<0.82$. For $D>1.182$, we have $\frac{D^3}{EFG}>2$ and $\phi_{13}<0$. This gives a contradiction. \vspace{2mm}

\noindent \textbf{Subcase IV:} $1.55<B\leq1.59$, $C>1.33$.\vspace{2mm}

As above, we have $F<0.908$, $E<1.08$ and $D<1.252$ and $G<0.8$. Now for $D>1.208$, $\frac{D^3}{EFG}>2$ and $\phi_9<0$. Therefore $D<1.208$.  For these bounds on $A,~B,~C,~D$ and $E$ we find that $\phi_{10}<0$.\vspace{2mm}

\noindent \textbf{Subcase V:} $1.59<B\leq1.678$, $C>1.43$.\vspace{2mm}

Here $F<0.908$, $E<1.08$ and $D<1.252$ and $G<0.79$. So for $D>1.239$, $\frac{D^3}{EFG}>2$ and $\phi_9<0$. Therefore $D<1.239.$ Now $\phi_{10}$ gives $D>1.21$, $C<1.465$, $A<1.775$. For these bounds on $A,~B,~C,~D$ and $E,~F,~G$ we find $\frac{D^3}{EFG}>2$ and $\phi_{13}<0$.\vspace{2mm}

\noindent \textbf{Subcase VI:} $1.59<B\leq1.678$, $C\leq1.43$.\vspace{2mm}

Here $\frac{A^3}{BCD}<2$ gives $C>1.127$, $D>1.04$. Also we have $F<0.908$, $E<1.08$ and $D<1.252$. Further $A<(2BCD)^{1/3}<1.818$ and $\frac{B^3}{CDE}>2$. We find $\phi_7,0$ for $G\geq 0.77$. Therefore we can take $G<0.77$.  For $D>1.198$, $\frac{D^3}{EFG}>2$ and $\phi_{13}<0$. So we have $D<1.198$. Then $A<(2BCD)^{1/3}<1.792$. Now we get $C>1.35$ by proving $\phi_{18}<0$. Further $\phi_{10}<0$ if either $B\leq 1.615$ or $C\geq 1.416$ or $D\leq 1.144$. So we can take $B>1.615$, $C<1.416$, $D>1.144$. For $D>1.178$, $\frac{D^3}{EFG}>2$ and $\phi_{13}<0$. So we have $D<1.178$. Then $A<(2BCD)^{1/3}<1.776$. We repeat the cycle to get that $C>1.374$, $B>1.655$ and  $C<1.394$. But then $\phi_{18}<0$. This gives a contradiction.\vspace{2mm}

 \noindent \textbf{Case II: $1.6<A\leq1.71$}\vspace{2mm}

\noindent \textbf{Subcase I:}  $B\leq1.35$\vspace{2mm}

Suppose $C<D$, then we must have $D>1.23$ as for $C<D<1.23$ we get $A<(2BCD)^{1/3}<1.6$. For $D>1.23$, we have $\phi_1^*<0$ as $B<1.35$. Therefore we must have $C>D$. Now suppose $C<E$. Using $i<\min\{b+\frac{c+d}{2},\frac{5}{9}-\frac{4}{9}e\}$, we get $E^4ABCDI>2$ for $E>1.093$ and  $\phi_8^*<0$.
  So $C\leq E\leq 1.093$ but then $A<(2BCD)^{1/3}<1.6$, a contradiction. Therefore we must have $C>E$. Now  $C>{\mbox{each of}}\{D,E,F,G,H,I\}$, therefore $(2,7^*)$ holds but  then $\phi_3^*<0$. This gives a contradiction. \vspace{2mm}

 \noindent \textbf{Subcase II:}  $1.35<B\leq1.4$ \vspace{2mm}

 Here also we will prove that $C>D$ and $C>E$. First suppose $C\leq D$, then we must have $D>1.2$, else $A$ will be $<1.6$. For $D>1.2$, $(1,2,6^*)$ holds but  $\phi_2^*<0$ for  $C<1.13$. When $C>1.13$, and also $E>1.08$, $E^4ABCDI>(0.46873)E^5(1.6)(1.35)(1.13)(1.2)>2$ but $\phi_8^*<0$. Also for $F>0.908$, $F^4ABCDE>2$ but $\phi_{12}^*<0$. Therefore we must have  $E<1.08$ and $F<0.908$. Also $G<0.86$ otherwise $\phi_{19}^{(2)}<0$ using $I>\frac{2}{3}G$. Now $\frac{D^3}{EFG}>2$ and $\phi_9<0$ for $D>1.2$. This gives a contradiction. So we must have $C>D$.\\
 Now suppose $C\leq E$. Using $i<\min\{b+\frac{c+d}{2},\frac{5}{9}-\frac{4}{9}e\}$, we get $E^4ABCDI>2$ for $E>1.105$ and $\phi_8^*<0$. So we get $E\leq1.105$. Hence $C\leq 1.093$ but then  $A<(2BCD)^{1/3}<1.6$, a contradiction. Therefore we must have $C>E$. Now  $C>{\mbox{each of}}\{D,E,F,G,H,I\}$, therefore $(2,7^*)$ holds but  then $\phi_3^*<0$. This gives a contradiction. \vspace{2mm}

 \noindent \textbf{Subcase III:}   $B>1.4$, $C>1.35$ \vspace{2mm}

%
%
Working as in subcase II we can take $E<1.08$, $F<0.908$. If $D>1.252$, then we find that $\frac{D^3}{EFG}>2$ and $\phi_9<0$. So $D<1.252$.
  Then $\frac{C^3}{DEF}>2$. Now $\phi_{11}<0$ if  $G>0.75$ and $\phi_{42}<0$ if $H>0.65$. Therefore we can take $G\leq0.75$ and $H\leq0.65$. Then $\frac{E^3}{FGH}>2$ for $E>1$. So $\phi_{10}$, $\phi_{18}$ holds. We see that $\phi_{10}<0$ if $1.4<B\leq1.5$.

  If $1.5<B\leq1.55$, we see that $\phi_{10}<0$ if  $ C\geq1.42$ or  $D\leq1.205$. So we can take $C<1.42$, $D>1.205$. For $D>1.205$, we have $\frac{D^3}{EFG}>2$. But $\phi_{13}<0$ if  $D>1.22, ~C<1.42$ or $D>1.205,~C<1.39$. So we can take $D<1.22$ and $C>1.39$. Now for these bounds we see that $\phi_{10}<0$. This gives a contradiction.

  If $1.55<B<1.6$, $1.35<C<1.4$, then $\phi_{10}<0$ for $D<1.165$, so we can take $D>1.165$.  So $\frac{D^3}{EFG}>2$. But $\phi_{13}<0$ if  $D>1.188$ or $C<1.37$ So we can take  $D<1.188, ~C>1.37$. Subsequently $\phi_{10}<0$ if $C>1.38$ and $D<1.18$. So we can take $D>1.18,~C<1.38$. Finally $\phi_{13}<0$. This gives a contradiction.

  Similarly we deal with  the ranges~ $1.55<B<1.6$, $1.4<C<1.465$; $1.55<B<1.6$, $1.465<C<1.498$; $1.6<B<1.678$, $1.35<C<1.4$; $1.6<B<1.678$, $1.4<C<1.45$; $1.6<B<1.678$, $1.45<C<1.498$ separately and using the inequalities $\phi_{10}, ~\phi_{13}$ and $\phi_{18}$ get a contradiction in each.\vspace{2mm}

 \noindent \textbf{Subcase IV:}   $B>1.4$, $C\leq1.35$, $\frac{E^3}{FGH}>2$ \vspace{2mm}

 Working as in Subcase III we get contradiction here. \vspace{2mm}

  \noindent \textbf{Subcase V:}   $B>1.4$, $C\leq1.35$, $\frac{E^3}{FGH}\leq 2$, $\frac{B^3}{CDE}>2$ \vspace{2mm}

  Here $FGH>\frac{1}{2}$. Now $\phi_6<0$ for $F\geq0.912$ or $H\geq0.676$ and $\phi_7<0$ for $G\geq0.822$. So we can take $F<0.912$, $H<0.676$ and $G<0.822$. Using $FGH>\frac{1}{2}$ we get $F>0.899$ and $H>0.666$. Now $\phi_6<0$.\vspace{2mm}

  \noindent \textbf{Subcase VI:}   $B>1.4$, $1.25<C\leq1.35$, $\frac{E^3}{FGH}\leq 2$, $\frac{B^3}{CDE}\leq 2$, $D>1.14$ \vspace{2mm}

  Here we can take $F<0.908$ and $FGH>\frac{1}{2}$. We must have $I>0.523$ as for $I\leq0.523$ we have $FGH<0.908(\frac{3}{2}I)(\frac{4}{3}I)<0.5.$ Now $\frac{A^3}{BCD}<2$ and $\frac{E^3}{FGH}<2$ implies $A^4E^4I<4,$ i.e. $E<1.04$. $\frac{B^3}{CDE}<2$ and $\frac{E^3}{FGH}<2$ implies $AB^4E^3I<4$, i.e. $B<1.478$. Using $\phi_{19}$ we have $G<0.83$. It gives $D<1.162$ otherwise $\frac{D^3}{EFG}>2$ and $\phi_{13}<0$. Also $F<0.889$ using $\phi_{15}.$ It further gives $G<0.82$ and $D<1.1493$. If $E>1.01$, we get $\frac{E^3}{FGH}=E^4ABCDI>2$ using $i<b+d-f-\frac{g}{2}$. So we must have $E<1.01$ and finally for $D>1.14$, $\frac{D^3}{EFG}>2$ and $\phi_{13}<0$. \vspace{2mm}

  Working similarly, we get a contradiction in the following subcases.\vspace{2mm}

   \noindent \textbf{Subcase VII:}   $B>1.4$, $1.25<C\leq1.35$, $\frac{E^3}{FGH}\leq 2$, $\frac{B^3}{CDE}\leq 2$, $D\leq 1.14$ \vspace{2mm}

    \noindent \textbf{Subcase VIII:}   $B>1.4$, $C\leq1.25$, $\frac{E^3}{FGH}\leq 2$, $\frac{B^3}{CDE}\leq 2$, $D>1.14$ \vspace{2mm}

     \noindent \textbf{Subcase IX:}   $B>1.4$, $C\leq1.25$, $\frac{E^3}{FGH}\leq 2$, $\frac{B^3}{CDE}\leq 2$, $D\leq 1.14$ \vspace{2mm}

%
%

\noindent \textbf{Case III: 1 $< A \leq $ 1.16}\vspace{2mm}

If $h+i<1.5e+0.465(a+d)$, we find that $\phi_{39}<0~$ for $~0<e<\min\{a,0.155\}$ and $0<d<a$. If $h+i\geq1.5e+0.465(a+d)$,  we get from  $(6.23)$ that $g<\frac{a}{2}+c+e-\frac{h+i}{2}<\frac{a}{2}+c+e-\frac{1.5e+0.465(a+d)}{2}$. Then $\phi_{19}^{(3)}<0~$ for $i<\frac{a}{2}+c+e-g$, (from(6.25)),
~$0<g<\frac{a}{2}+c+e-\frac{1.5e+0.465(a+d)}{2}$, ~$0<e<\min\{a,0.155\}$, ~$d\leq a$ and ~$0<c\leq a$.\vspace{2mm}

\noindent \textbf{Case IV: 1.16 $< A \leq $ 1.2}\vspace{2mm}

If $h+i<1.5e+0.47(a+d)$, we find that $\phi_{39}<0$ for $0<e<0.155$ and $0<d<a$. If $h+i\geq1.5e+0.47(a+d)$, then we get from $(6.23)$ that $g<\frac{a}{2}+c+e-\frac{1.5e+0.47(a+d)}{2}$. Then $\phi_{19}^{(3)}<0$ for $i<\frac{a}{2}+c+e-g$, ~$0<g<\frac{a}{2}+c+e-\frac{1.5e+0.47(a+d)}{2}$, ~$0<e<0.155$, ~$d\leq a$ and ~$0<c\leq 0.165$. So we can take $c>0.165$. Also $a\geq c>0.165$. Now if $f+h\leq 0.66(b+d)$, then $\phi_{40}^{(1)}<0$ for $0<d\leq a$, ~$0<b\leq a$ and $0.165<a<0.2$. So we have $f+h> 0.66(b+d)$. Now $(6.24)$ gives $i<2(b+d)-2(f+h)<0.68(b+d)$. Using $i<0.68(b+d)$, $a\geq c>0.165$, $0<b<0.2$, $0<d<0.2$, we get $E^4ABCDI>2$ for $E\geq1.102$ and $\phi_{8}^*<0$. So $E<1.102$. Now $\phi_{39}<0$  if $h+i<1.5e+0.546(a+d)$. Therefore  $h+i\geq 1.5e+0.546(a+d)$ which gives $g<\frac{a}{2}+c+e-\frac{1.5e+0.546(a+d)}{2}$. Then $\phi_{19}^{(3)}<0$ for $0<c\leq a$.\vspace{2mm}

\noindent \textbf{Case V: 1.2 $< A \leq $ 1.25}\vspace{2mm}

Working as in Claim(i) on page 113, we get $E^4ABCDI>2$ for $i<\min\{\frac{5}{9}-\frac{4}{9}e,b+\frac{c+d}{2}\}$, $b\leq a$, $E\geq 1.1488$ and $\phi_{8}^*<0$. So we can take $E<1.1488$. Now we have following claims:\vspace{2mm}\\
Claim(i)  $C\geq1.13$ and $D>1.045$. Suppose $c<0.13$. If $h+i<1.5e+0.475(a+d)$, then $\phi_{39}<0$. If $h+i\geq1.5e+0.475(a+d)$, then $\phi_{19}^{(3)}<0$. So $C>1.13$. Now $\phi_{10}^*<0$ for $D\leq1.045$. So we can take $D>1.045$.\\
Claim(ii) $g>0.18$ and $g>0.21$ if $e>0.095$. For $0<e<0.1488$, $0<g<0.18$, $\phi_{19}^{(3)}<0$ and if $0.095<e<0.1488$, then $\phi_{19}^{(3)}<0$ for $0<g<0.21$.\\
Claim(iii) $D<1.2199$. Suppose $D\geq1.2199$. If $1<E<1.095$, then $\frac{D^3}{EFG}>2$, using $G<1-g<0.82$. If $1.095\leq E<1.1488$, then $\frac{D^3}{EFG}>2$, using $G<1-g<0.79$. Also $\phi_9<0$ for $D>1$, giving thereby a contradiction.\\
Claim(iv) $f+h>0.63(b+d)$, proving $\phi_{40}^{(2)}<0$. Now $(6.24)$ gives $i<0.74(b+d)$.\\
Claim(v) $E<1.105$. Using $i<0.74(b+d)$ we get $E^4ABCDI>2$ for $E\geq1.105$, then $\phi_{8}^*<0$.\\
Claim(vi) $f>0.095$ and $D<1.1793$. For $f\leq 0.095$ we get $\phi_{22}^{(1)}<0$. So $F<0.905$. Now for $D>1.1793$, $\frac{D^3}{EFG}>2$ and $\phi_9<0$.\\
Claim(vii) $E<1.073$. Using $(6.26)$ we get $i<b+d-0.095-\frac{0.18}{2}$, then $E^4ABCDI>2$ for $E>1.073$. Also $\phi_{10}<0$ for $E\geq 1.073$.\\
Claim(viii) $C>1.17$, $D>1.073$ and $B>1.14$.\\
Suppose $C\leq1.17$. If $h+i<1.5e+0.6(a+d)$, then $\phi_{39}<0$. If $h+i\geq1.5e+0.6(a+d)$, then $\phi_{19}^{(3)}<0$. So $C>1.17$. Now using $\phi_{10}^*$  we can take $D>1.073$. So $D>E$. Now $\phi_1^*<0$ for $B<1.14$ and $A>\max\{C,1.2\}$. Hence we can take $B>1.14$.\\
Claim(ix) $f+h>0.663(b+d)$ and $f>0.135$. Using $\phi_{40}^{(2)}$ we get $f+h>0.663(b+d)$, then $(6.24)$ gives $i<0.674(b+d)$. Now $\phi_{22}^{(1)}<0$ for $f\leq 0.135$. Hence $F<1-0.135=0.865$, Now for $D>1.151$, $\frac{D^3}{EFG}>2$ and $\phi_9<0$. So we can take $D\leq 1.151$.\\
Claim(x) $E<1.0425$ proving $E^4ABCDI>2$ for $i<b+d-0.135-\frac{0.18}{2}$.\\
Using $E<1.0422$, $F<0.865$, $G<0.82$, we get $\frac{D^3}{EFG}>2$, then $\phi_9$ gives contradiction.\\
Claim(xi) $g>0.208$ proving $\phi_{19}^{(2)}<0$. So $G\leq 0.792$, which further gives $D<1.127$.\\
Claim(xii) $C>1.205$, $D>1.097$ and $B>1.185$.\\
Claim(xiii) $f+h>0.705(b+d)$ and $f>0.152$.\\
Claim(xiv) $E<1.024$, $D<1.113$.\\
Claim(xv) $C<1.235$, proving $\phi_{10}^*<0$ for $C\geq 1.235$.\\
Claim(xvi) $C>1.212$, $D>1.103$. If $h+i<1.82e+0.652(a+d)$, then $\phi_{39}<0$. If $h+i\geq1.82e+0.652(a+d)$, then $\phi_{19}^{(3)}<0$ for $C\leq 1.212$. So $C>1.212$. Now using $\phi_{10}^*$  we can take $D>1.103$.\\
Final: Using $\phi_{19}^{(2)}$ we get $g>0.23$, i.e. $G\leq
0.77$. Now for $D>1.103$, $\frac{D^3}{EFG}>2$ and $\phi_9<0$.\vspace{2mm}

\noindent \textbf{Case VI: 1.25 $ <A \leq $ 1.3}\vspace{2mm}\\
\textbf{Subcase I:} $D>1.15$\vspace{2mm}

Working as in Case V, we get $E^4ABCDI>2$ for $i<min\{\frac{5}{9}-\frac{4}{9}e,b+\frac{c+d}{2}\}$, $b\leq a$ and $E\geq 1.143$. So we can take $E<1.143$ using $\phi_{8}^*$. If $h+i<1.5e+0.468(a+d)$, then $\phi_{39}<0$. If $h+i\geq1.5e+0.468(a+d)$, then $\phi_{19}^{(3)}<0$ for $C<1.183$. So we can take $C\geq1.183$. Now $g>0.195$ proving $\phi_{19}^{(3)}<0$. So $G<1-0.195=0.805$. Now for $D>1.226$, $\frac{D^3}{EFG}>2$ and $\phi_9<0$. So $D\leq1.226$. Using $\phi_{40}^{(2)}$ we get $f+h>0.602(b+d)$ and further $i<0.796(b+d)$, which gives $E<1.092$ and hence $D<1.207$. Now $f+h>0.617(b+d)$, i.e. $i<0.766(b+d)$. Now $\phi_{22}^{(1)}<0$ for $f<0.062$. Hence $F<1-0.062<0.938$. It gives $D<1.182$. Now proving $E^4ABCDI>2$ and $\phi_{10}<0$ we get $E<1.06$. Further we get $g>0.218$ using $\phi_{19}^{(2)}$. It gives $D<1.159$. $\phi_{22}^{(1)}$ gives $f>0.09$, i.e. $F<0.91$. Now for $D>1.15$, $\frac{D^3}{EFG}>2$ and $\phi_9<0$. \vspace{2mm}

\noindent \textbf{Subcase II:} $1<D\leq1.15$, $C>1.21$\vspace{2mm}

If $1<D\leq1.08$ and $C>1.21$ we find that  $\phi_{10}^*<0$. So we can take $D>1.08$. As in Subcase I, we have $E<1.143$.\\
Claim(i) $g>0.215$ otherwise $\phi_{19}^{(2)}<0$\\
Claim(ii) $E<1.07$. For $E>1.07$, $E^4ABCDI>2$ for $i<b+d-\frac{0.215}{2}$ and $\phi_{10}<0$.\\
Claim(iii) $f+h>0.648(b+d)$, $i<0.704(b+d)$, otherwise $\phi_{40}^{(2)}<0$.\\
Claim(iv) $f>0.115$ otherwise $\phi_{22}^{(1)}<0$.\\
Claim(v) $E<1.029$ and $D<1.127$\\
Claim(vi) $g>0.237$ and $D<1.116$.\\
Claim(vii) $C<1.262$ otherwise $\phi_{10}^*<0.$\\
Claim(viii) $f+h>0.72(b+d)$, otherwise $\phi_{40}^{(2)}<0$.\\
Claim(ix) $f>0.155$ otherwise $\phi_{22}^{(1)}<0$.\\
Claim(x) $E<1.014$ and $D<1.094$\\
Claim(xi)$C<1.228$ and $B>1.27$. Since $\phi_{10}^*<0$ for $C\geq1.228$ or $B\leq1.27$. \\
Claim(xii) $g>0.25$ otherwise $\phi_{19}^{(2)}<0$\\
Claim(xiii) $f+h>0.757(b+d)$, otherwise $\phi_{40}^{(2)}<0$.\\
Claim(xiv) $f>0.17$ otherwise $\phi_{22}^{(1)}<0$.\\
Finally $\frac{D^3}{EFG}>2$ for $D>1.08$ and $\phi_9<0$. This gives a contradiction. \vspace{2mm}

\noindent \textbf{Subcase III:} $1<D\leq1.08$, $C\leq1.21$ \vspace{2mm}

Here $E<1.143$. If $h+i<1.5e+0.61(a+d)$, then $\phi_{39}<0$. If $h+i\geq1.5e+0.61(a+d)$, then $\phi_{19}^{(3)}<0$ for $C<1.161$. So we can take $C\geq1.161$.
Now $D>1.045$ and $B>1.22$ otherwise $\phi_{10}^*<0$. Using this we get $g>0.205$ by proving $\phi_{19}^{(2)}<0$. For $E>1.077$, $E^4ABCDI>2$ for $i<b+d-\frac{0.205}{2}$ and $\phi_{10}<0$. Next we get $f+h>0.75(b+d)$ and then $f>0.164$. It gives $E<1.023$. Now if $h+i<1.79e+0.68(a+d)$, then $\phi_{39}<0$.
If $h+i\geq1.79e+0.68(a+d)$, then $\phi_{19}^{(3)}<0$ for $c<0.21$.\vspace{2mm}

\noindent \textbf{Subcase IV:} $1.08<D\leq1.15$, $C\leq1.21$\vspace{2mm}

Here also we get contradiction by working similar to as in Subcase (III).\vspace{2mm}\\
\noindent \textbf{Case VII: 1.3 $< A \leq$ 1.35}\vspace{2mm}\\
\textbf{Subcase I:} $B\leq1.3$, $1.2<C\leq1.3$\vspace{2mm}

Here $ D>1.08$ for otherwise $\phi_{10}^*<0$. We get $E^4ABCDI>2$ for $i<min\{\frac{5}{9}-\frac{4}{9}e,b+\frac{c+d}{2}\}$, $b\leq a$ and $E\geq 1.1404$. Then $\phi_{8}^*<0$. So we can take $E<1.1404$. \\
Claim(i) $g>0.205$ otherwise $\phi_{19}^{(2)}<0$\\
Claim(ii) $D<1.2195$ as for $D>1.2195$, $\frac{D^3}{EFG}>2$ and $\phi_9<0$\\
Claim(iii) $E<1.1$ otherwise $E^4ABCDI>2$ using $i<b+d-\frac{g}{2}$ and $\phi_{8}^*<0$. Then $D<1.205$\\
Claim(iv) $g>0.226$ otherwise  $\phi_{19}^{(2)}$ and then $D<1.195$\\
Claim(v)  $E<1.075$ and $D<1.185$\\
Claim(vi) $f+h>0.598(b+d)$, otherwise  $\phi_{40}<0$.\\
Claim(vii) $f>0.056$ otherwise $\phi_{22}<0$.\\
Claim(viii) $D<1.1625$ as for $D>1.1625$, $\frac{D^3}{EFG}>2$ and $\phi_9<0$\\
Claim(ix) $E<1.047$ and $D<1.1533$\\
Claim(x) $f+h>0.66(b+d)$ and $i<0.68(b+d)$.\\
Claim(xi) $f>0.128$, $F<0.872$ and then $D<1.123$\\
Claim(xii) $E<1.019$ and $D<1.113$.\\
Claim(xiii) $C<1.255$ and $B>1.2$ as  $\phi_{10}^*<0$ for $C\geq1.255$ or $B\leq1.2$. \\
Claim(xiv) $f+h>0.73(b+d)$ and $i<0.54(b+d)$.\\
Claim(xv) $f>0.16$, $F<0.84$ and then $D<1.099$\\
Claim(xvi) $E<1.007$ and $D<1.095$.\\
Claim(xvii) $C<1.23$ and $B>1.25$ as  $\phi_{10}^*<0$ for $C\geq1.23$ or $B\leq1.25$. \\
Claim(xviii) $f+h>0.755(b+d)$ and $i<0.49(b+d)$.\\
Claim(xix) $f>0.182$, $F<0.818$ and then $D<1.085$\\
Finally we have $E^4ABCDI>2$ for $E>1$ and $i<b+d-f\frac{g}{2}$ and $\phi_{10}<0$. Hence we get a contradiction.\vspace{2mm}

\noindent \textbf{Subcase II:} $B\leq1.3$, $1.3<C\leq A$\\
Here $ D>1.14$ for otherwise $\phi_{10}^*<0$. As in Subcase I,  we can take $E<1.1404$.\\
Claim(i) $g>0.17$ otherwise $\phi_{19}^{(2)}<0$\\
Claim(ii) $D<1.25$ as for $D>1.25$, $\frac{D^3}{EFG}>2$ and $\phi_9<0$\\
Claim(iii) $E<1.089$ using $i<b+d-\frac{g}{2}$  and then $D<1.205$\\
Claim(iv) $g>0.208$ otherwise $\phi_{19}^{(2)}<0$ and then $D<1.1993$\\
Claim(v)  $E<1.059$ using $i<b+d-\frac{g}{2}$ and then $D<1.189$\\
Claim(vi) $f+h>0.645(b+d)$, otherwise $\phi_{40}<0$.\\
Claim(vii) $f>0.123$ otherwise $\phi_{22}<0$. Hence $\frac{D^3}{EFG}>2$ and $\phi_9<0$.\vspace{2mm}

\noindent \textbf{Subcase III:} $B\leq1.3$, $C\leq1.2$, $D>1.13$\\
Here also $E<1.1404$\\
Claim(i) $g>0.185$ otherwise $\phi_{19}^{(2)}<0$\\
Claim(ii) $D<1.2301$ as for $D>1.2301$, $\frac{D^3}{EFG}>2$ and $\phi_9<0$\\
Claim(iii) $h+i>1.5e+0.43(a+d)$ otherwise $\phi_{39}<0$.\\
Claim(iv) $C>1.13$ otherwise $\phi_{19}<0$.\\
Claim(v) $g>0.205$ otherwise $\phi_{19}^{(2)}<0$ and then $D<1.2195$\\
Claim(vi)  $E<1.102$ using $i<b+d-\frac{g}{2}$ and then $D<1.2056$\\
Claim(vii) $f+h>0.63(b+d)$, otherwise $\phi_{40}<0$.\\
Claim(viii) $f>0.075$ otherwise $\phi_{22}<0$ and $D<1.175$.\\
Claim(ix)  $E<1.062$ using $i<b+d-f-\frac{g}{2}$ and then $D<1.161$\\
Claim(x) $h+i>1.5e+0.62(a+d)$ otherwise $\phi_{39}<0$.\\
Claim(xi) $C>1.19$ otherwise $\phi_{19}<0$.\\
Claim(xii) $f+h>0.64(b+d)$, otherwise $\phi_{40}<0$.\\
Claim(xiii) $f>0.11$ otherwise $\phi_{22}<0$ and $D<1.146$.\\
Claim(xiv)  $E<1.035$ using $i<b+d-f-\frac{g}{2}$ and then $D<1.136$\\
Claim(xv) $h+i>1.65e+0.639(a+d)$ otherwise $\phi_{39}<0$.\\
Finally $\phi_{19}<0$.\vspace{2mm}

\noindent
\textbf{Subcase IV:} $B\leq1.3$, $C\leq1.2$, $D\leq1.13$\\
Proof is similar.\\
\noindent \textbf{Subcase V:} $B>1.3$, $C>1.3$\vspace{2mm}

Here $D>1.122$ otherwise $\phi_{10}^*<0$. Also $E<1.1404$. For $F>0.9491$, $F^4ABCDE>2$ and $\phi_{12}^*<0$. So $F\leq0.9491$. $g>0.135$ otherwise $\phi_{19}^{(3)}<0$. So $G<0.865$, which gives $D<1.24$ otherwise $\frac{D^3}{EFG}>2$ and $\phi_9<0$. It gives $E<1.07$. It further gives $g>0.205$. If $\frac{E^3}{FGH}>2$, then $\phi_{10}<0$ for $E>1$. So we can have $\frac{E^3}{FGH}\leq2$, i.e $H\geq\frac{1}{2FG}>0.66$. Now if $\frac{D^3}{EFG}>2$, then $\phi_9<0$. Hence $\frac{D^3}{EFG}\leq2.$ Now using $i<b+d-f-\frac{g}{2}$, we find $\frac{E^3}{FGH}=E^4ABCDI>2$ for $E\geq1.05$. So we can take $E<1.05$ and so $D<1.166$. Repeating the cycle we get $D<1.1497$. Now $C<1.313$ otherwise $\phi_{10}^*<0$. We get $f+h>0.655(b+d)$ using $\phi_{40}<0$. Now $\phi_{22}^{(3)}<0$ for $0<f\leq 0.12$. So we can take $f>0.12$. Finally we have $\frac{D^3}{EFG}>2$ for $D>1.122$ and $\phi_9<0$. This gives contradiction.\vspace{2mm}

\noindent \textbf{Subcase VI:} $B>1.3$, $1.2<C\leq1.3$, $D>1.14$\vspace{2mm}

Here $E<1.1404$. For $F>0.9666$, $F^4ABCDE>2$ and $\phi_{12}^*<0$. So $F<0.9666$. Now $\frac{D^3}{EFG}>2$ for $D>1.302$ and $\phi_9<0$, so we get $D<1.302$. Further we get $g>0.2$ otherwise  $\phi_{19}^{(3)}<0$. With this we get $D<1.208$ and $E<1.082$. With these improved values we get $g>0.21$, $E<1.075$ and $D<1.1795$. Now $\phi_{40}$ gives $f+h>0.635(b+d)$, i.e. $i<0.73(b+d)$. Now $\phi_{22}$ gives $f>0.07$, i.e. $F<0.93$. Now $D<1.165$ and $E<1.055$. Repeating this cycle we get a contradiction.\vspace{2mm}

\noindent \textbf{Subcase VII:} $B>1.3$, $1.2<C\leq1.3$, $D\leq1.14$\vspace{2mm}

Here also $E<1.1404$. Also $D>1.053$ otherwise $\phi_{10}^*<0$. Now for $F>0.9838$, $F^4ABCDE>2$, so proving $\phi_{12}^*<0$ we get $F<0.9838$. We have \\
Claim(i) $g>0.214$ otherwise $\phi_{19}^{(3)}<0$.\\
Claim(ii) $\frac{E^3}{FGH}\leq 2$ for otherwise $\phi_{10}<0$ for $E>1$ gives a contradiction. So we can take $FGH>\frac{1}{2}.$ Also  $E<1.069$ otherwise $E^4ABCDI>2$ using  $i<b+d-f-\frac{g}{2}$. \\
Claim(iii) $f+h>0.652(b+d)$ otherwise $\phi_{40}<0$, i.e. $i<0.696(b+d)$.\\
Claim(iv) $f>0.118$ otherwise $\phi_{22}^{(3)}<0$. It gives $E<1.032$.
$FGH>\frac{1}{2}$ gives $H>0.71$.\\
Claim(v) $\frac{D^3}{EFG}\leq 2$ for otherwise $\phi_{9}<0$ for $D>1$ gives a contradiction. So  $D<(2EFG)^{\frac{1}{3}}<1.127$.\\
Claim(vi) $f+h>0.668(b+d)$ and $i<0.664(b+d)$.\\
Claim(vii) $D<1.111$ and $E<1.02$.\\
Claim(viii) $\frac{C^3}{DEF}<2$ and $C<1.249$.\\
Suppose $\frac{C^3}{DEF}\geq2$, then $\phi_{42}<0$ if $H\geq 0.595$ and $C>1.2$. For $H< 0.595$  we find $\frac{E^3}{FGH}>2$ for $E>1$. Then $C<1.249$ for $D<1.111$, $E<1.02$ ~{\mbox and}~ $F<0.855$.\\
Claim(ix) $f+h>0.698(b+d)$ and $i<0.604(b+d)$.\\
Claim(x) $f>0.18$ otherwise $\phi_{22}^{(3)}<0$ and so $\frac{C^3}{DEF}<2$ implies $C<1.22$.\\
Claim(xi) $f+h>0.735(b+d)$ and $i<0.53(b+d)$.\\
Claim(xii) $f>0.21$. It gives $D<1.08$.\\
Finally $\frac{C^3}{DEF}<2$ implies $C<1.2$, a contradiction.\vspace{2mm}

\noindent \textbf{Subcase VIII:} $B>1.3$, $C\leq1.2$.\vspace{2mm}

\noindent Proof is similar.

\noindent \textbf{Case VIII: 1.35 $< A \leq$ 1.4}\vspace{2mm}

\noindent \textbf{Subcase I:} $B\leq1.26$\\
Proving $\phi_3^*<0$ we get $C<D$ or $C<E$. First suppose $C<D$, then we have\\
Claim(i) $g>0.17$  otherwise $\phi_{19}<0$.\\
Claim(ii) $D<1.245$ as for $D\geq1.245$, $\frac{D^3}{EFG}>2$ and $\phi_9<0$. So $C<D<1.245$.\\
Claim(iii) $f+h>0.665(b+d)$  otherwise $\phi_{40}<0$. It gives $i<0.67(b+d).$\\
Claim(iv) $E<1.104$ using $i<0.67(b+d)$ and proving $\phi_8^*<0$.\\
Claim(v) $C<D<1.224$ as for $D>1.224$, $\frac{D^3}{EFG}>2$ and $\phi_9<0$.\\
Claim(vi) $f>0.155$ otherwise $\phi_{22}<0$. This gives $C<D<1.157$.\\
Claim(vii) $E<1.055$, this gives $C<D<1.1396$.\\
Claim(viii) $h+i>1.5e+0.63(a+d)$ otherwise  $\phi_{39}<0$.\\
Now $\phi_{19}<0$. This gives a contradiction.\\ Working similarly we get a contradiction when $C<E$.\vspace{2mm}

\noindent \textbf{Subcase II:} $1.26<B\leq1.35$, $1<C\leq1.144$\\
If $D>E$, then $\phi_2^*<0$. So $D<E<1.155$. It gives $f+h>0.7(b+d)$, i.e. $i<0.6(b+d)$. Using this we get $E^4ABCDI>2$ for $E>1.09$. So $E<1.09$ by proving $\phi_8^*<0$. For
$D<E<1.09$ we get $h+i>1.5e+0.645(a+d)$. For this we get $\phi_{19}<0$.\vspace{2mm}

\noindent \textbf{Subcase III:} $1.26<B\leq1.35$, $1.144<C\leq1.26$\\
Claim(i) $g>0.205$ otherwise $\phi_{19}^{(3)}<0$.\\
Claim(ii) $D<1.225$ as for $D>1.225$, $\frac{D^3}{EFG}>2$ and $\phi_9<0$.\\
Claim(iii) $f+h>0.615(b+d)$\\
Claim(iv) $E<1.088$ using $i<0.77(b+d)$. It further gives $D<1.2005$.\\
Claim(v) $f+h>0.647(b+d)$ and $i<0.706(b+d)$.\\
Claim(vi) $f>0.12$ and $D<1.151$.\\
Claim(vii) $E<1.036$ using $i<b+d-f-\frac{g}{2}$. Then $D<1.1318$.\\
Claim(viii) $f+h>0.688(b+d)$ and $i<0.624(b+d)$.\\
Claim(ix) $f>0.17$ and $D<1.1099$.\\
Claim(x) $E<1.017$\\
Claim(xi) $g>0.212$, $G<0.788$ and $D<1.0999$\\
Claim(xii) $C\leq1.21$. For if $C>1.21$, then for $E>1$, $E^4ABCDI>2$ and $\phi_{10}<0$. Now $f+h>0.735(b+d)$, $f>0.2$ and $D<1.087$. If $b<0.32$, then also $E^4ABCDI>2$. So $B>1.32$. Then $\phi_{22}<0$ for $f<0.21$. So $f>0.21$ then also $E^4ABCDI>2$ for $E>1$, but $\phi_{10}<0$. This gives a contradiction.\vspace{1mm}\\
\noindent \textbf{Subcase IV:} $1.26<B\leq1.35$, $1.26<C\leq A$, $\frac{C^3}{DEF}>2$.\\
$\phi_{11}$ and $\phi_{42}$ gives $G<0.78$ and $H<0.6$. It gives $D<1.25$ as for $D>1.25$, $\frac{D^3}{EFG}>2$ and $\phi_9<0$. Also for $E>1$, $\frac{E^3}{FGH}>2$ and
$\phi_{10}<0$. \vspace{2mm}

\noindent \textbf{Subcase V:} $1.26<B\leq1.35$, $1.26<C\leq1.35$, $\frac{C^3}{DEF}\leq 2$.\vspace{2mm}

Here $D>1.095$ otherwise $\phi_{10}^*<0$. For $F>0.9609$, we get $F^4ABCDE>2$ and $\phi_8^*<0$. For $E>1.127$, $E^4ABCDI>2$ using $I>0.46873E$. and $\phi_{12}^*<0$. Therefore  $F<0.9609$ and $E<1.127$.
For $D>1.294$, $\frac{D^3}{EFG}>2$ and $\phi_9<0$. So $D<1.294$. Now we have following claims:\\
Claim(i) $G<0.91$ using $\phi_{19}$\\
Claim(ii) $D<1.254$\\
Claim(iii) $f+h>0.59(b+d)$ and $i<0.82(b+d)$\\
Claim(iv) $E<1.084$\\
Claim(v) $D<1.238$ as for $D>1.238$, we get $\frac{D^3}{EFG}>2$ and $\phi_9<0$\\
Claim(vi) $g>0.195$ and $G<0.805$ using $\phi_{19}$. then we get $D<1.189$.\\
Claim(vii) $f+h>0.62(b+d)$, i.e. $i<0.76(b+d)$\\
Claim(viii) $E<1.056$ and $D<1.178$\\
Now $\frac{C^3}{DEF}<2$ implies $C<1.338$. This gives $f+h>0.635(b+d)$, i.e. $i<0.73(b+d)$.\\
Claim(ix) $f>0.09$ and $D<1.157$\\
$\frac{C^3}{DEF}<2$ implies $C<1.306$.\\
Claim(x) $E<1.025$ and $D<1.146$. Then $C<1.289$.\\
Claim(xi) $f+h>0.675(b+d)$, i.e. $i<0.65(b+d)$\\
Claim(xii) $f>0.17$, $F<0.83$\\
$\frac{C^3}{DEF}<2$ implies $C<1.26$. This gives a contradiction.\vspace{2mm}

\noindent \textbf{Subcase VI:} $1.26<B\leq1.35$, $1.35<C\leq A$, $\frac{C^3}{DEF}\leq 2$.\\
$\phi_{10}^*$ implies $D>1.15$. For $F>0.933$, we get $F^4ABCDE>2$ and for $E>1.101$, $E^4ABCDI>2$ using $I>0.46873E$. $\phi_8^*$ and $\phi_{12}^*$ implies $F<0.933$ and $E<1.101$.
For $D>1.272$, $\frac{D^3}{EFG}>2$ and $\phi_9<0$. So $D<1.272$. Now $\phi_{19}$ gives $G<0.88$. Now $D<1.219$. $\frac{C^3}{DEF}<2$ implies $C<1.359$, $D>1.198$. Now we
get for $F>0.924$, $F^4ABCDE>2$ and for $E>1.092$, $E^4ABCDI>2$ using $I>0.46873E$. $\phi_8^*$ and $\phi_{12}^*$ implies $F<0.924$ and $E<1.092$. $\frac{C^3}{DEF}<2$ implies $C<1.35$.\vspace{2mm}

\noindent \textbf{Subcase VII:} $B>1.35$, $C\leq1.2$\\
If $D>E$, $(2,1,6^*)$ holds but  $\phi_1^*<0$.  Therefore $D\leq E$.\\
Claim(i) $g>0.165$ using $\phi_{19}^{(1)}$ and $D\leq E<1.155$.\\
Claim(ii) $f+h>0.65(b+d)$, i.e. $i<0.7(b+d)$ using $\phi_{40}^{(3)}$. \\
Claim(iii) $f>0.18$, i.e $F<0.82$ using $\phi_{40}^{(2)}$\\
Claim(iv) $E<1.068$ using $i<b+d-f-\frac{g}{2}$.\\
Claim(v) $h+i>1.75d+0.64(a+e)$ using $\phi_{39}$.\\
Claim(vi) $C>1.16$ using $\phi_{19}$.\\
Repeating the cycle with this improved $C$ we get $g>0.215$, $D\leq E<1.015$. It gives $h+i>1.95d+0.7(a+e)$ otherwise $\phi_{39}<0$. Then $\phi_{19}<0$ if $C\leq 1.175$. So $C>1.175$.
Finally we get $\phi_2^{*}<0$.\vspace{2mm}

\noindent \textbf{Subcase VIII:} $B>1.35$, $1.2<C\leq1.25$\vspace{2mm}

\noindent Claim(i) $g>0.2$  otherwise  $\phi_{19}^{(1)}<0$.\\Claim(ii) $F<0.978$, for if $F>0.978$, then $F^4ABCDE>2$ and  $\phi_{12}^*<0$.\\
 Claim(iii) $D<1.219$ as $\frac{D^3}{EFG}>\frac{1.219^3}{1.155\times 0.8\times 0.9798}>2$
and $\phi_{13}<0$, so $D<1.219.$ \\Claim(iv) $E<1.097$. Using  $i<b+d-f-\frac{g}{2}$ we get $E^4ABCDI>2$ for $E>1.097$ and $\phi_8^*<0$ so $E<1.097$ which further gives $D<1.198$.\\ Claim(v) $F<0.955$ and $H<0.76$ using  $\phi_{44}$.\\Claim(vi) $g>0.212$\\Claim(vii) $E<1.045$ and $D<1.162$\\ Claim(viii) $f+h>0.655(b+d)$ otherwise $\phi_{40}<0$, i.e. $i<0.69(b+d)$.\\
Claim(ix) $f>0.165$, $F<0.835$\\
Now $\frac{E^3}{FGH}>2$ and $\phi_{10}<0$.\vspace{2mm}

\noindent \textbf{Subcase IX:} $B>1.35$, $1.25<C\leq A$\vspace{2mm}

\noindent Proof is similar.\vspace{2mm}

\noindent \textbf{Case IX: 1.4 $ <A \leq $ 1.5}\vspace{2mm}\\
\noindent \textbf{Subcase I:} $1<B\leq1.29$\\
If $C>\mbox{each of}~\{D,E,F,G,H,I\} $, then $\phi_3^*<0$, which gives a contradiction. So we have either $C<D$ or $C<E$. Suppose first that $C<D$, i.e. $C<D<1.322$. We have following claims:\\
Claim(i) $E<1.118$\\
Claim(ii) $D<1.308$ as for $D\geq1.308$, $\frac{D^3}{EFG}>2$ and $\phi_9<0$.\\
Claim(iii) $g>0.16$ and $D<1.234$. For $g<0.16$, $\phi_{19}^{(2)}<0$. So $G<1-0.16=0.84$, which gives $D<1.234$.\\
Claim(iv) $f+h>0.66(b+d)$, by proving $\phi_{40}^{(2)}<0$.\\
Claim(v) $f>0.16$ by proving  $\phi_{22}^{(1)}<0$ for $f<0.16$.\\
Claim(vi) $D<1.165$, as for $D>1.165$, $\frac{D^3}{EFG}>2$ and $\phi_9<0$.\\
Claim(vii) $E<1.05$, otherwise proving $E^4ABCDI>2$ and $\phi_{41}<0$.\\
It further gives $D<1.141$.\\
Claim(viii) $h+i>1.5e+0.6(a+d)$ otherwise proving $\phi_{39}<0$. \\
Claim(ix) $C>1.1$, as $\phi_{19}^{(3)}<0$ for $C\leq 1.1$.\\
Claim(x) $f+h>0.735(b+d)$, otherwise proving $\phi_{40}^{(2)}<0$.\\
Claim(xi) $f>0.19$ by proving  $\phi_{22}^{(1)}<0$ for $f<0.19$.\\
Claim(xii) $D<1.127$, as for $D>1.127$, $\frac{D^3}{EFG}>2$ and $\phi_9<0$.\\
Now $\phi_{39}<0$. This gives a contradiction.\\
Proceeding similarly we get contradiction if $C<E$.\vspace{2mm}

\noindent \textbf{Subcase II:} $1.29<B\leq1.4$, $C\leq1.2$\\
Here we have $g>0.17$ using  $\phi_{19}^{(1)}$. Using $i<\frac{a}{2}+c+e-\frac{0.17}{2}$ and $\phi_8^*$ we get $E<1.08$. Now $D<1.215$ as for $D>1.215$, $\frac{D^3}{EFG}>2$ and $\phi_{13}<0$. Proving $\phi_{40}^{(2)}<0$ we get $f+h>0.652(b+d)$, i.e. $i<0.696(b+d)$. Using  $\phi_{22}^{(1)}$ we get $f>0.19$, i.e. $F<0.81$. Now proving $\frac{D^3}{EFG}>2$ and $\phi_{13}<0$ we get $D<1.133$. Now $\phi_{11}^*$ gives $C>1.12$. We further get $g>0.205$ using $\phi_{19}^{(1)}$ and then $E<1.015$, $D<1.094$. By proving $f+h>0.73(b+d)$, we get $f>0.215$, i.e. $F<0.785$ and $D<1.083$. $E^4ABCDI>2$ for $B<1.38$ or $C>1.13$. Also $\phi_{18}<0$ for $E>1$. So we get $B>1.38$ and $C<1.13$. Now $\phi_{11}^*<0$, a contradiction.\vspace{2mm}

\noindent \textbf{Subcase III:} $1.29<B\leq1.4$, $1.2<C\leq1.25$\\
Working  as above we get a contradiction in this subcase also.\vspace{1mm}\\
\noindent \textbf{Subcase IV:} $1.29<B\leq1.4$, $1.25<C\leq1.498$, $\frac{E^3}{FGH}>2$.\\
Here $\phi_{10}$ gives $D>1.21$. Now $F^4ABCDE>2$ for $F>0.926$ and $\phi_{12}^*<0$. So $F<0.926$. For $E>1.094$, $E^4ABCDI>2$ and $\phi_8^*<0$, so $E<1.094$. Now for $D>1.266$, $\frac{D^3}{EFG}>2$ and $\phi_9<0$, so $D<1.266$. Further $\phi_{10}$ implies $C<1.355$ and then $G<0.905$ using $\phi_{19}^2$. It further gives $D<1.23$ and $C<1.285$ using $\phi_9$ and $\phi_{10}$ respectively. For these values $\phi_{18}$ gives that $B<1.368$. Now $\phi_{10}<0$.\vspace{1mm}\\
\noindent  \textbf{Subcase V:} $1.29<B\leq1.4$, $1.25<C\leq1.498$, $~\frac{E^3}{FGH}\leq 2$,  $~\frac{C^3}{DEF}>2$.\vspace{1mm}\\
 Here $\phi_{11}$ gives $G<0.8$ and $\phi_{42}$ gives $H<0.62$. So $\frac{E^3}{FGH}>2$ for $E>1$, a contradiction.
 \vspace{2mm}

\noindent \textbf{ Subcase VI:} $1.29<B\leq1.4$, $1.25<C\leq1.3$, $~\frac{E^3}{FGH}\leq 2$,  $~\frac{C^3}{DEF}\leq 2$\\

 Using $\phi_{10}^*$ we get $D>1.068$. Now $F^4ABCDE>2$ for $F>0.955$ and $\phi_{12}^*<0$, so $F<0.955$. For $E>1.121$, we have $E^4ABCDI>2$ and $\phi_8^*<0$, so $E<1.121$. Now for $D>1.289$, $\frac{D^3}{EFG}>2$ and $\phi_9<0$, so $D<1.289$. Now $\phi_{19}^{(1)}$ implies $g>0.173$, i.e $G<0.827$. For $D>1.2099$, $\frac{D^3}{EFG}>2$. For $B<1.35$, $\phi_9<0$ and for $B>1.35$, $\phi_{13}<0$. So we have $D<1.2099$.\\
 For $E>1.067$ and $i<b+d-f-\frac{g}{2}$, we get $\frac{E^3}{FGH}=E^4ABCDI>2$, a contradiction. So We have $E<1.067$. It further gives $g>0.195$, $D<1.1795$ and $E<1.053$. Then $f+h>0.63(b+d)$ otherwise $\phi_{40}^{(3)}<0$. Therefore  $i<0.74(b+d)$. Now $f>0.063$ using $\phi_{22}^{3}$, i.e. $F<0.937$. Now if $H<0.65$, then $\frac{E^3}{FGH}>2$, so we must have $H>0.65$. For $D>1.169$, $\frac{D^3}{EFG}>2$ and $\phi_9<0$. So $D<1.169$. Now repeating this cycle we get $F<0.91$, $E<1.028$, $D<1.147$, which gives $C<(2DEF)^{\frac{1}{3}}<1.2899$. Then $f+h>0.65(b+d)$, $f>0.14$, $E<1.01$, $D<1.112$ and hence $C<1.25$.\vspace{2mm}

\noindent \textbf{Subcase VII:} $1.29<B\leq1.4$, $1.3<C\leq1.498$, $~\frac{E^3}{FGH}\leq 2$,  $~\frac{C^3}{DEF}\leq 2$\\
Here we have the following:\\
Claim(i) $D>1.1$ otherwise   $\phi_{10}^*<0$.\\
Claim(ii) $E<1.106$ otherwise  $\phi_8^*<0.$\\
Claim(iii) $F<0.939$ otherwise  $\phi_{12}^*<0.$
Then $C<(2DEF)^{\frac{1}{3}}<1.401$\\
Claim(iv) $D<1.276$ otherwise  $\phi_9<0.$\\
Claim(v) $G<0.89$ otherwise  $\phi_{19}^{2}<0$ using $I>\frac{2}{3}G$.\\
Claim(vi) $E<1.0875$ otherwise  $\phi_8^*<0$.
Now $C<1.377$.\\
 Claim(vii) $D<1.225$ otherwise  $\phi_9<0$. Then
$C<1.358$.\\
Claim(viii) $g>0.226$ otherwise  $\phi_{19}^{3}<0$.\\
Now $FGH>\frac{1}{2}$ implies $H>0.687$.\\
Claim(ix) $D<1.165$  otherwise  $\phi_9<0$ using $H>0.687$. \\
Claim(x) $E<1.026$\\
Claim(xi) $D<1.143$\\
Then $E<1.0215$ and hence $C<(2DEF)^{\frac{1}{3}}<1.3$, a contradiction.\vspace{2mm}

\noindent\textbf{Subcase VIII:} $1.4<B\leq1.5$, $1<C\leq1.185$\vspace{2mm}

If $D>\max\{E,F,G,H,I\}$, then $\phi_2^*<0$. So we can take $D<E<1.155$. Now $\phi_{11}^*<0$.\vspace{2mm}

\noindent \textbf{Subcase IX:} $1.4<B<1.5$, $1.185<C<1.35$, $\frac{E^3}{FGH}>2$\\
Claim(i) $C>1.28$ otherwise  $\phi_{18}<0$.\\
Claim(ii) $D>1.17$ otherwise  $\phi_{10}.$\\
Claim(iii) $E<1.08$ as for $E>1.08$, $E^4ABCDI>(1.08)^4\times1.4^2\times1.28\times1.17\times0.46873>2$ and $\phi_{8}^*<0$\\
Claim(iv) $F<0.909$ as for $F>0.909$, $F^4ABCDE>2$ and $\phi_{12}^*<0$\\
Claim(v) $D<1.253$ as for $D>1.253$, $\frac{D^3}{EFG}>2$ and $\phi_{9}<0$
Claim(vi) $C>1.297$ otherwise  $\phi_{18}<0$.\\
Claim(vii) $D>1.178$ otherwise  $\phi_{18}<0$.\\
Claim(viii) $G<0.865$ otherwise  $\phi_{19}^2<0.$\\
Claim(ix) $D<1.205$ otherwise  $\phi_{13}<0$.\\
Claim(x) $B>1.45$ otherwise  $\phi_{18}<0.$\\
Claim(xi) $E<1.0599$ as for $E>1.0599$, $E^4ABCDI>2$ and $\phi_{18}<0$\\
Claim(xii) $G<0.848$ otherwise  $\phi_{19}^2<0.$\\
Now $\frac{D^3}{EFG}>2$ and $\phi_{13}<0$\vspace{2mm}

\noindent\textbf{Subcase X:} $1.4<B\leq1.5$, $1.35<C\leq1.498$, $\frac{E^3}{FGH}>2$\vspace{2mm}

Here $E<1.08$, $F<0.908$, $D<1.253$.
 If $D\leq 1.21$ or $C\geq 1.42$ we find that  $\phi_{18}<0$; so  we have $D>1.21$ and $C<1.42$ .
 If $G\geq 0.85$ we find $\phi_{19}^2<0$. So we must have $G<0.85$.
 Now $\frac{D^3}{EFG}>2$ for $D>1.21$. By proving $\phi_{9}<0$  we have $D<1.235$.
 Next we find $\phi_{10}<0$ for $C\geq1.4$. So we have $C<1.4$. Now
 $\phi_{13}<0$ for $B\geq1.455$ or $C<1.358$. So we have $B<1.455$ and $C>1.358$. Finally  $\phi_{10}<0$. This gives a contradiction.\vspace{2mm}

\noindent\textbf{Subcase XI:} $1.4<B\leq1.5$, $\frac{E^3}{FGH}\leq 2$,  $\frac{B^3}{CDE}>2$\vspace{2mm}

If $H\geq 0.66$ or $F\geq 0.91$ we find $\phi_6<0$. So  $H<0.66$ and $F<0.91$.
If $G\geq 0.81$, we see that $\phi_7<0$. Therefore  $G<0.81$. Now $\frac{E^3}{FGH}>2$ for $E>1$, a contradiction.\vspace{2mm}

\noindent\textbf{Subcase XII:} $1.4<B\leq1.5$, $\frac{E^3}{FGH}\leq 2$,  $\frac{B^3}{CDE}\leq 2$, $~1.185<C\leq1.25$, $D>1.16$\vspace{2mm}

\noindent Here we have \\
Claim(i) $E<1.097$ as for $E>1.097$, $\frac{E^3}{FGH}=E^4ABCDI>(1.097)^4\times 1.4^2\times1.185\times1.16\times0.46873>2$.\\
Claim(ii) $F<0.929$ as for $F>0.929$, $F^4ABCDE>2$ and $\phi_{12}^*<0$\\
Claim(iii) $g>0.2$ otherwise  $\phi_{19}^3<0$.\\
Claim(iv) $D<1.178$ as for $D>1.178$, $\frac{D^3}{EFG}>2$ and $\phi_{13}<0$.\\
Claim(v) $B<(2CDE)^{1/3}$ implies $B<1.479$.\\
Now using $i<b+d-f-\frac{g}{2}$, we get $\frac{E^3}{FGH}=E^4ABCDI>2$ for $E>1.063$. So $E<1.063$. It gives $D<1.165$. Now $g>0.21$, i.e. $G<0.79$. Now $\frac{D^3}{EFG}>2$ for $D>1.16$ and $\phi_{13}<0$.\vspace{2mm}

\noindent\textbf{Subcase XIII:} $1.4<B\leq1.5$, $\frac{E^3}{FGH}\leq 2$,  $\frac{B^3}{CDE}\leq 2$, $1.185<C\leq1.25$, $1.1<D\leq1.16$\\
Here we have\\
 Claim(i) $g>0.21$ and $F>0.9407$.\\
 Claim(ii) $E<1.065$. As for $i<b+d-f-\frac{g}{2}$ and $E>1.065$, we find that $\frac{E^3}{FGH}=E^4ABCDI>2$. Now $B<(2CDE)^{1/3}$ implies $B<1.457$, which further gives $E<1.058$. \\
 Claim(iii) $f+h>0.637(b+d)$ otherwise $\phi_{40}^{(3)}<0$, i.e. $i<0.726(b+d)$.\\
  Claim(iv) $f>0.105$ otherwise $\phi_{22}^{(3)}<0$.
  Now we get $E<1.038$.\\
  Claim(v) $D<1.137$ as for $D\geq1.137$, we find $\frac{D^3}{EFG}>2$ and $\phi_{13}<0$.
  Now $B<1.435$.\\
   Claim(vi) $f+h>0.658(b+d)$ so $i<0.684(b+d)$.\\
  Claim(vii) $f>0.18$ otherwise  $\phi_{22}^{(3)}<0$.\\
  Now for $i<b+d-f-\frac{g}{2}$ and $E\geq1.01$, we find that $\frac{E^3}{FGH}=E^4ABCDI>2$. So $E<1.01$.\\
  Finally for $D>1.1$, we find $\frac{D^3}{EFG}>2$ and $\phi_{13}<0$. This gives a contradiction.\vspace{2mm}

\noindent\textbf{Subcase XIV:} $1.4<B\leq1.5$, $\frac{E^3}{FGH}\leq 2$,  $\frac{B^3}{CDE}\leq 2$, $1.185<C\leq1.25$, $1.0<D\leq1.1$\vspace{1mm}\\
Here $g>0.21$, $E<1.1294$, $F<0.9634$. It gives $B<1.4591$ and further $E<1.055$. Again  using $B<(2CDE)^{1/3}$ we have $B<1.427$. Now using $\phi_{40}^{(3)}$ we get $f+h>0.67(b+d)$, i.e. $i<0.66(b+d)$. $\phi_{22}^{(3)}$ gives $f>0.185$. Now using $i<b+d-f-\frac{g}{2}$, we get $E^4ABCDI>2$ for $E>1$, a contradiction. \vspace{2mm}

\noindent\textbf{Subcase XV:} $1.4<B\leq1.5$, $\frac{E^3}{FGH}\leq 2$,  $\frac{B^3}{CDE}\leq 2$,~$C>1.25$, $\frac{C^3}{DEF}>2$\\
Here $\phi_{11}$ gives $G<0.8$ and $\phi_{42}$ gives $H<0.67$. Also for $F>0.951$, $F^4ABCDI>0.951^4\times1.4^2\times1.25>2$. So $F<0.951$ using $\phi_{12}^*$. $E^3<2FGH$ implies $G>\frac{1}{2FH}>0.784$. Now $\phi_{42}$ gives $H<0.62$. Now we get $FGH<\frac{1}{2}$, we get contradiction as in Subcase III.\vspace{2mm}

\noindent \textbf{Subcase XVI:} $1.4<B\leq1.5$, $C>1.25$,$~ \frac{E^3}{FGH}\leq 2$,  $\frac{B^3}{CDE}\leq 2$, $\frac{C^3}{DEF}\leq 2$
\\
Working as above we get contradiction here.\vspace{2mm}

\noindent \textbf{Case X: $1.5<A\leq1.6$}\vspace{2mm}

\noindent\textbf{ Subcase I:}  $B\leq1.35$\vspace{2mm}

For $i<\min\{b+\frac{c+d}{2},\frac{5}{9}-\frac{4}{9}e\}$ and for $E>1.113$ we get $E^4ABCDI>2$. Using $\phi_8^*$ we get $E<1.113$. Now for $D>1.307$, $\frac{D^3}{EFG}>\frac{1.307^3}{1.113}>2$ and $\phi_9<0$. So we have $D<1.307$.
If $C>{\mbox{each of}}\{D,E,F,G,H,I\}$, then $\phi_3^*<0$. So we must have $C<D$ or $C<E$. If $D<1.115$, then $C<1.115$ and so $A<(2BCD)^{1/3}<(2\times1.35\times1.115\times1.115)^{1/3}$. So $D\geq1.115>E$. So $(2,1,6^*)$ holds. But $\phi_1^*$ implies $B>1.31$. Now if $C>1.25$, then $D>1.25$. For $E>1.069$, we have $E^4ABCDI>(0.46873)(1.069)^5(1.5)(1.31)(1.25)^2>2$ and for $F>0.9$ we have $F^4ABCDE>2$. Using $\phi_{41}$ we get $E<1.069$ and using $\phi_{20}$ we get $F<0.9$. So $\frac{D^3}{EFG}>2$ and using $\phi_9$ we get contradiction. So we have $C<1.25$. Now $\phi_1^*$ implies $B>1.32$ and $A<1.555$, then  $\phi_2^*$ implies $C>1.145$. Now we have $F<0.941$ as $F^4ABCDE>2$ and $\phi_{12}^*<0$, it further gives $D<1.2795$. Using $i<\frac{a}{2}+c+e-g$ and $\phi_{19}$ we get $g>0.255$ and further we get $D<1.16$. Now $\phi_1^*$ implies $B>1.332$ and $A<1.529$. Now $\phi_1^*<0$. Now suppose $C<E$. If $D<E$, then $A<(2BCD)^{1/3}<1.5$. If $D>E$, then $\phi_1^*$ implies $B>1.31$ and for this $\phi_2^*<0$.\vspace{2mm}

\noindent\textbf{Subcase II:} $1.35<B\leq1.45$, $C>1.35$\vspace{2mm}

$\phi_{10}^*$ implies $D>1.125$. Then using $\phi_{12}^*$ we have $F<0.908$ and using $\phi_8^*$ we get $E<1.08$. Now for $D>1.25$, $\frac{D^3}{EFG}>2$ and $\phi_9<0$. We have $\frac{C^3}{DEF}>2$. $\phi_{11}$ gives $G<0.74$ and $\phi_{42}$ gives $H<0.63$. Now for $E>1$, $\frac{E^3}{FGH}>2$ and $\phi_{10}<0$.\vspace{2mm}

\noindent\textbf{Subcase III:} $1.35<B\leq1.45$, $C\leq1.35$, $\frac{E^3}{FGH}>2$\vspace{2mm}

We see that  $\phi_{18}<0$ for $1<C<1.225$. Now consider $1.225<C<1.35$. Here
$\phi_{10}<0$ for $D<1.19$.So $D>1.19$ Then we will have $F<0.908$ and $E<1.08$. Now $\phi_{19}^{(2)}$ gives $G<0.862$.  For $D>1.21$, $\frac{D^3}{EFG}>2$ and $\phi_9<0$, so we have $D<1.21$. Now $\phi_{10}$ implies $B>1.43$, $C<1.3$ and $E<1.065$. Finally $\phi_{13}<0$.\vspace{2mm}

\noindent\textbf{Subcase IV:} $1.35<B\leq1.45$, $C\leq1.35$, $\frac{E^3}{FGH}\leq 2$, $\frac{B^3}{CDE}>2$\vspace{2mm}

Using $\phi_6$ and $\phi_7$, we have $F<0.93$, $G<0.84$, $H<0.683$. Here
 $FGH>\frac{1}{2}$, which gives $F>0.871$ and $H>0.64$. Now $\phi_6<0$.\vspace{2mm}

\noindent\textbf{Subcase V:} $1.35<B\leq1.45$, $C\leq1.155$\\
For $D>E$, $\phi_2^*<0$. Do $D<E$. For $E<1.06$ we get $A<1.5$. So $E>1.06$. Now using $i<\frac{a}{2}+c+e-g$, $\phi_{19}$ gives $g>0.2$. $\phi_{44}$ gives $H<0.84$. Now we get contradiction as $FGH<\frac{1}{2}$.\vspace{2mm}

\noindent\textbf{Subcase VI:} $1.35<B\leq1.45$, $1.155<C\leq1.25$, $D\leq1.14$\\
Here $E<1.128$ using $I>0.46873E$ and $\phi_8^*$. Using $\phi_{12}^*$, we have $F<0.9617$. Using $\phi_{19}$ we have $g>0.17$ and further $E<1.07$. Now $\phi_{42}^{(2)}$ implies $f+h>0.645(b+d)$, i.e. $i<0.71(b+d)$. $\phi_{22}$ implies $f>0.15$.Using $FGH>\frac{1}{2}$and working like this we get contradiction.\vspace{2mm}

\noindent\textbf{Subcase VII:} $1.35<B\leq1.45$, $1.155<C\leq1.25$, $D>1.14$.\vspace{2mm}

\noindent
Claim(i) $F<0.931$ otherwise $\phi_{12}^*<0$ \\
Claim(ii) $E<1.099$ otherwise  $\phi_{8}^*<0$ \\
Claim(iii) $g>0.185$ otherwise  $\phi_{19}<0$ \\
Claim(iv) $D<1.186$ otherwise  $\phi_{9}<0$ \\
Claim(v) $E<1.064$ and $D<1.174$\\
Claim(vi) $f+h>0.635(b+d)$ otherwise  $\phi_{42}^{(2)}<0$ \\
Claim(vii) $f>0.11$ otherwise  $\phi_{22}^{(2)}<0$ \\
Claim(viii) $D<1.156$ otherwise  $\phi_{9}<0$ \\
Claim(ix) $E<1.038$ and $F<0.86$\\
Claim(x) $D<1.14$ otherwise  $\phi_{9}<0$.\\
Working similarly we get contradiction in following subcases also.\vspace{2mm}

\noindent\textbf{Subcase VIII:} $1.35<B\leq1.45$, $1.25<C\leq1.25$, $D>1.14$\vspace{2mm}

\noindent\textbf{Subcase IX:} $1.35<B\leq1.45$, $1.25<C\leq1.25$, $D>1.14$ \vspace{2mm}

\noindent\textbf{Subcase X:} $B>1.45$, $\frac{E^3}{FGH}>2$ \vspace{2mm}

\noindent Claim(i) $C<1.285$. Here $\phi_{18}<0$.\\
Claim(ii) $1.285<C<1.325$, $B>1.5$.
Here $E<1.08$, $F<0.912$ and then $D<1.254$ using $\phi_9$.
then $\phi_{18}<0$.\\
Claim(iii) $1.285<C<1.325$, $B<1.5$.
$\phi_{10}^*$ gives $D>1.05$.
$\phi_{45}$ gives $G<0.855$ (using $I>\frac{2}{3}G$).
Then $E<1.08$, $F<0.909$ and $D<1.1885$. Now we get contradiction using $\phi_{18}$ and $\phi_{10}$.\\
Claim(iv) $1.325<C<1.4$, $B>1.5$\\
Claim(v) $1.325<C<1.4$, $B<1.5$\\
Claim(vi) $1.4<C<1.498$\vspace{2mm}

\noindent \textbf{Subcase XI:} $B>1.45$, $\frac{E^3}{FGH}\leq 2$,  $\frac{B^3}{CDE}>2$ \vspace{2mm}

If $F\geq 0.89$ or $H\geq 0.65$, one finds $\phi_6<0$ and if $H\geq 0.65$, $\phi_7<0$. So we must have  $F<0.89$, $H<0.65$ and $G<0.81$. Now $FGH>\frac{1}{2}$, a contradiction.\vspace{2mm}

\noindent \textbf{Subcase XII:} $B>1.51$, $\frac{E^3}{FGH}\leq 2$,  $\frac{B^3}{CDE}\leq 2$, $C>1.25$ \vspace{2mm}

If $\frac{B^3}{CDE}>2$ or $\frac{E^3}{FGH}>2$, then we get contradiction from Subcases X and XI. So $\frac{B^3}{CDE}<2$ and $\frac{E^3}{FGH}<2$. If $I<0.52$, then $FGH<\frac{1}{2}$, so $I>0.52$. Also $\frac{A^3}{BCD}<2$. Therefore $AB^4E^3I<4$, i.e. $B<(\frac{4}{AI})^{\frac{1}{4}}<1.51.$\vspace{2mm}

\noindent \textbf{Subcase XIII:} $B>1.51$, $\frac{E^3}{FGH}\leq 2$,  $\frac{B^3}{CDE}\leq 2$, $C<1.25$ \vspace{2mm}

If $D>E$, then $\phi_2^*<0$. If $D<E<1.155$, then $B<(2CDE)^{\frac{1}{3}}<1.51$\vspace{2mm}

\noindent \textbf{Subcase XIV:} $1.45<B<1.51$, $\frac{E^3}{FGH}\leq 2$,  $\frac{B^3}{CDE}\leq 2$, $C<1.22$ \vspace{2mm}

As above, $D<E$. $\frac{B^3}{CDE}<2$ implies $D>1.08$. For $I>0.46873E$, and for $E>1.109$, $E^4ABCDI>2$, so $E<1.109$ by proving $\phi_8^*<0$. So $B<(2CEE)^{\frac{1}{3}}<1.45$.\vspace{2mm}

\noindent \textbf{Subcase XV:} $1.45<B<1.51$,$~\frac{E^3}{FGH}\leq 2$,  $\frac{B^3}{CDE}\leq 2$, $1.22<C<1.25$ \vspace{2mm}

 For $I>0.46873E$, and for $E>1.0997$, $E^4ABCDI>2$, so $E<1.0997$ as $\phi_8^*<0$. $\frac{B^3}{CDE}<2$ implies $D>1.108$. Then $E<1.078$ and $F<0.908$. Now we get $G<0.83$ by proving $\phi_{45}<0$. Now $D<1.176$, $E<1.045$ and $B<1.45$, a contradiction.\vspace{2mm}

\noindent Similarly we get contradiction when $1.45<B<1.51$, $\frac{E^3}{FGH}\leq 2$,  $\frac{B^3}{CDE}\leq 2$, $1.25<C<1.498$. \vspace{2mm}

\noindent The Theorem, i.e. Conjecture III for $n=9$, follows from Propositions 1-46.\vspace{2mm}

\noindent{\bf Acknowledgements:} The authors are grateful to Professors R. P. Bambah, R. J. Hans-Gill and Ranjeet Sehmi for various discussions throughout the preparation of this paper.\vspace{2mm}\\
\textbf{\large{References}}
\begin{enumerate}
\item{R. P. Bambah, V. C. Dumir and R. J. Hans-Gill,  Non-homogeneous problems: Conjectures of Minkowski and Watson, Number Theory, Trends in Mathematics, Birkhauser Verlag, Basel, (2000) 15-41.}
\item{R.P. Bambah, A.C. Woods, {\it On a theorem of Dyson}, J. Number Theory 6 (1974), 422-433.}
\item{R.P. Bambah, A.C. Woods, {\it Minkowski's conjecture for n = 5. A theorem of Skubenko}, J. Number Theory 12 (1980), 27-48.}
\item{B.J. Birch, H.P.F. Swinnerton-Dyer, {\it On the inhomogeneous minimum of the product of n linear forms}, Mathematika 3 (1956),
25-39.}
\item{H. F. Blichfeldt, {\it The minimum values of positive quadratic forms in six, seven and eight variables}, Math. Z. 39 (1934), 1-15.}
\item{H. Cohn and N. Elkies,  {\it New upper bounds on sphere packings}, I. Ann. of Math. (2) 157 (2003), no. 2, 689-714.}
\item{J. H. Conway and N. J. A. Sloane, Sphere packings, Lattices and groups, Springer-verlag, 2nd edition
.}
\item{P. Gruber, {\it Convex and discrete geometry}, Springer Grundlehren  Series, Vol.336 (2007).}
\item{P. Gruber and C. G. Lekkerkerker, Geometry of Numbers, Second Edition, North Holland (1987).}
\item{R. J. Hans-Gill, Madhu Raka , Ranjeet Sehmi and Sucheta, {\it A unified simple proof of Woods' conjecture for $n \leq 6$}, J. Number Theory,129 (2009) 1000-1010.}
\item{R. J. Hans-Gill, Madhu Raka and Ranjeet Sehmi, {\it On conjectures of Minkowski and Woods for $n= 7$}, J. Number Theory, 129 (2009), 1011-1033.}
\item{R. J. Hans-Gill, Madhu Raka and Ranjeet Sehmi, {\it Estimates On Conjectures of Minkowski and Woods}, Indian Jl. Pure Appl. Math.,41(4) (2010), 595-606.}
\item{R.J. Hans-Gill, Madhu Raka, Ranjeet Sehmi, {\it On Conjectures of Minkowski and  Woods for $n=8$}, Acta Arithmetica, 147(4) (2011), 337-385.}
\item{R. J. Hans-Gill, Madhu Raka and Ranjeet Sehmi, {\it Estimates On Conjectures of Minkowski and Woods II}, Indian Jl. Pure Appl. Math.,42(5) (2011), 307-333.}
\item{Leetika Kathuria, R. J. Hans-Gill and Madhu Raka,{\it On a question of Uri Shapira and Barak Weiss} to appear in Indian Jl. Pure Appl. Math.}
\item{A. Korkine, G. Zolotareff, {\it Sur les formes quadratiques}, Math. Ann. 6 (1873), 366-389; ~{\it Sur les formes quadratiques positives}, Math. Ann. 11 (1877), 242-292.}
\item{C.T. McMullen, {\it Minkowski's conjecture, well rounded lattices and topological dimension}, J. Amer. Math. Soc. 18 (2005),
711-734.}
\item{H. Minkowski, {\it $\ddot{U}$ber die Annäherung an eine reelle Grösse durch rationale Zahlen}, Math. Ann. 54 (1901), 91-124.}
\item{R.A. Pendavingh and S.H.M. Van Zwam, {\it New Korkine-Zolotarev inequalities}, SIAM J. Optim. 18 (2007), no. 1, 364-378.}
\item{Uri Shapira and Barak Weiss, {\it On the stable lattices and the diagonal group}, arXiv:1309.4025v1 [math.DS] 16 September 2013.}
\item{A.C. Woods, {\it The densest double lattice packing of four spheres}, Mathematika 12 (1965) 138-142.}
\item{A.C. Woods, {\it Lattice coverings of five space by spheres}, Mathematika 12 (1965) 143-150.}
\item{A.C. Woods, {\it Covering six space with spheres}, J. Number Theory 4 (1972) 157-180.}
\end{enumerate}

\end{document}